\documentclass[10pt]{amsart}
\usepackage{palatino}
\usepackage[letterpaper, left=1in, right= 1in]{geometry} 
\usepackage{FancyMath}

\title[Koszul dual algebras part 1]{Koszul dual $\A_{\infty}$-algebras from star-shaped diagrams -- part 1}
\author[Isabella Khan]{Isabella Khan}
\date{\today}

\begin{document}

\maketitle

\begin{abstract}
	By slicing the Heegaard diagram for a given $3$-manifold in a particular way, it is possible to construct $\A_{\infty}$-bimodules, the tensor product of which retrieves the Heegaard Floer homology of the original 3-manifold. The first step in this is to construct algebras corresponding to the individual slices. Here, we use the graphical calculus for $\A_{\infty}$-structures introduced in~\cite{DiagBible} to construct Koszul dual weighted $\A_{\infty}$-algebras $\A$ and $\B$, and dualizing bimodules for a particular star-shaped class of slice. The duality result is then proved in the sequel. 
\end{abstract}\blfootnote{The author was partially supported by a NSF Graduate Research Fellowship, and by the Simons Collaboration on New Structures in Low Dimensional Topology.}

\tableofcontents

\section{Introduction} 
	
	Heegaard Floer homology is a three-manifold invariant which can be used to retrieve classical data about the underlying manifold (see e.g.~\cite{Heeg1} and~\cite{Heeg2}). It is constructed using a Heegaard diagram for the given three-manifold, as well as pseudo-holomorphic disks associated to this Heegaard diagram. One advantage of the Heegaard Floer package is that it is often more easily computable than other existing three-manifold invariants, while also retrieving much of the same data. 

	The aim in bordered Heegaard Floer constructions, such as e.g.~\cite{Bim2014},~\cite{BordBook},~\cite{Torus}, is to retrieve this data by associating algebraic objects to subsections of the original Heegaard diagram rather than the whole. More specifically, we first cut up the Heegaard diagram into two or more pieces and associate an algebraic object to each piece. Then we recombine these objects (usually by taking a tensor product) to get back the original Heegaard Floer homology. Since the algebraic objects corresponding to the cut-up pieces of the Heegaard diagram can be more easily computable than the whole Heegaard Floer complex, this often speeds up computations for the whole three-manifold.
 
	Knot Floer homology first defined by Ozsv\'ath and Szabo in~\cite{KnotFloer}, and, separately, by Rasmussen in~\cite{Ras}, is a knot invariant closely related to Heegaard Floer homology, which also admits a bordered construction. Bordered knot Floer homology allows us to retrieve the knot Floer homology corresponding to a given knot $K$ from modules and bimodules associated to different parts of a braid presentation of $K$ -- see e.g.~\cite{Kauff}, ~\cite{AlgMat},~\cite{Khov}. More specifically, we start with a knot diagram for $K$ and take a tubular neighborhood, labelled so as to make it a Heegaard diagram. For instance, in the case of the the left-handed trefoil, this looks like
	\begin{center}
		\includegraphics[width = 12cm]{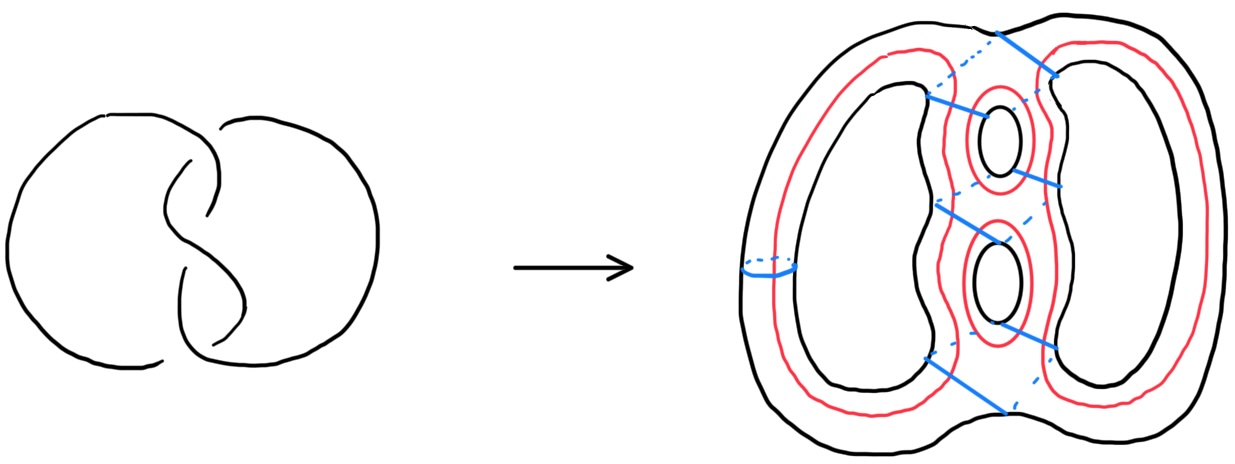}
	\end{center}
	Notice that the blue $\beta$-circles encode the crossings in this diagram. Then we can slice the diagram horizontally:
	\begin{center}
		\includegraphics[width = 12cm]{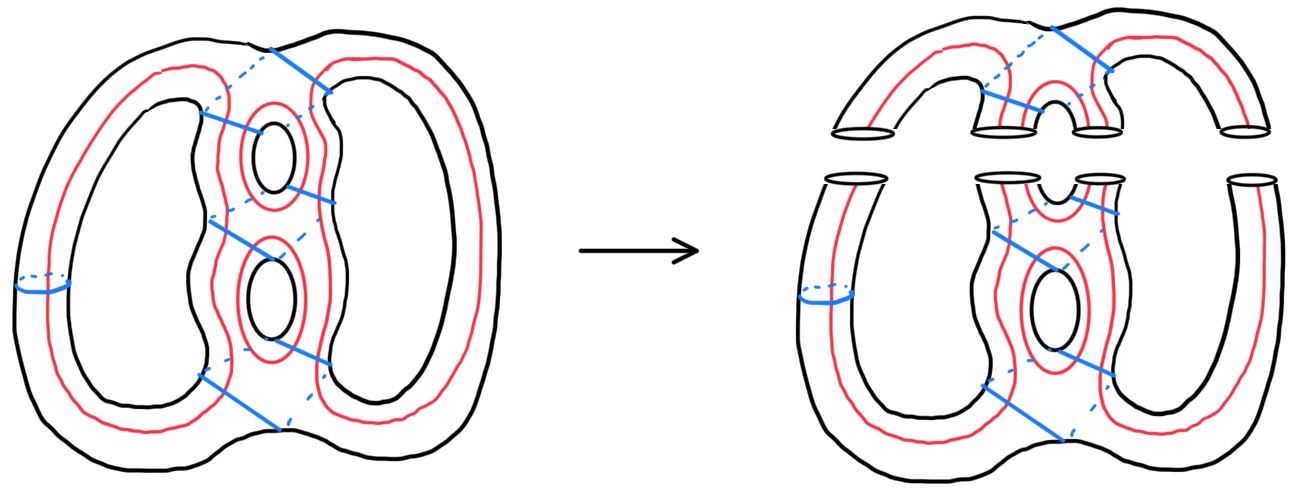}
	\end{center}
	and in this case, the top portion of the diagram becomes a sphere with four punctures, labelled with three red $\alpha$-arcs and a $\beta$-circle, that is
	\begin{center}
		\includegraphics[width = 12cm]{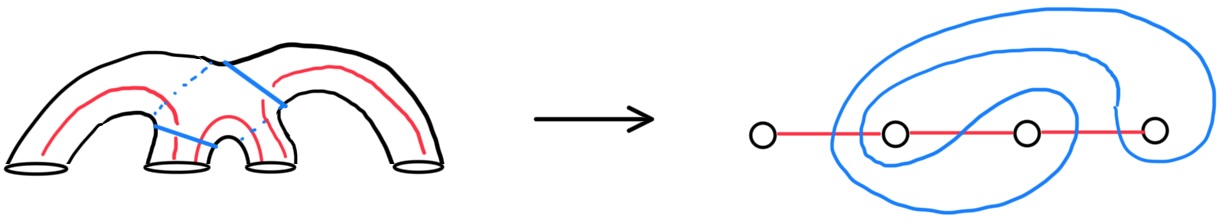}
	\end{center}
	The first step in the bordered knot Floer construction is to associate an $\A_{\infty}$-algebra $\A$ to the horizontal slice. The $\beta$-circle does not have any bearing on $\A$, so we are really just associating an algebra to the planar graph
	\begin{figure}[H]
		\includegraphics[width = 5cm]{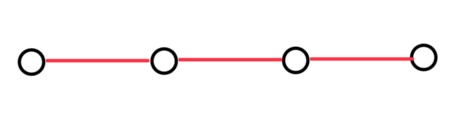}

		\caption{}\label{tref4'}
	\end{figure}
	\hspace{-0.45cm}See~\cite{KnotFloer},~\cite{Kauff},~\cite{AlgMat} for this construction. One way to further understand this algebra is to construct a Koszul dual algebra $\B$ for $\A$. In~\cite{Pong1}, the authors do just this, associating an algebra $\B$ to the blue portion of the following figure:
	\begin{figure}[H]
		\includegraphics[width =7cm]{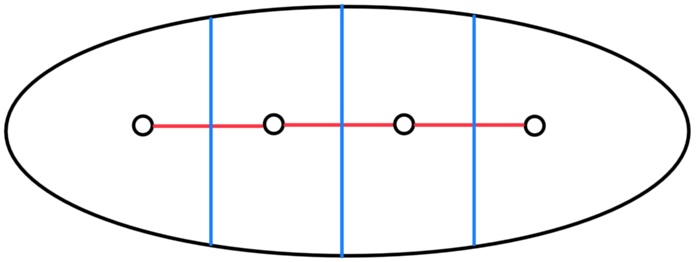}
		\caption{}\label{pongy}
	\end{figure}
	\hspace{-0.5cm}just as $\A$ is associated to the red portion. They then prove that $\A$ and $\B$ are Koszul dual.

	It is also of interest to consider more complicated planar graphs, which may arise in horizontal slices corresponding to other tangles. In this paper, we consider one such class of planar graph -- namely star shaped graphs as in
	\begin{figure}[H]
		\includegraphics[width = 6cm]{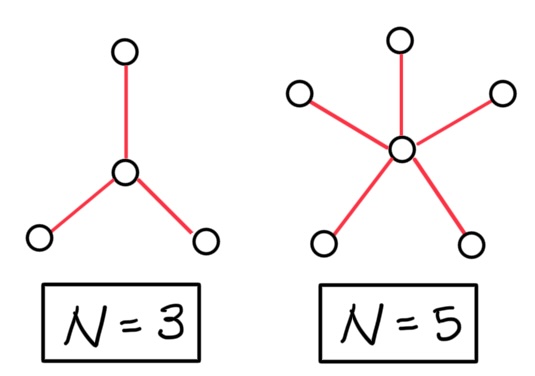}
		\caption{}\label{star71}
	\end{figure}
	\hspace{-0.5cm}where $N$ denotes the number of $1$-valent vertices, and is usually assumed greater than two. (It is important to remember that these ``vertices'' are in fact punctures in the bordered Heegaard diagram, and correspond to Reeb chords and Reeb orbits picked up by the pseudo-holomorphic disks used in the construction.) 

	This paper is part one of a pair of papers which together prove the following duality relation. 

    \begin{theorem}\label{duality}
        Let $\A$ and $\B$ be the weighted $\A_{\infty}$-algebras defined in Sections~\ref{alph} and~\ref{bb1}, respectively. Then there exists a DD-bimodule $\:^{\A} X^{\B}$, constructed in Section~\ref{DDbim}, and a weighted AA-bimodule $\:_{\B} Y_{\A}$ constructed in Section~\ref{aabim}, such that 
        \begin{equation}\label{diag4}
            \:^{\A} X^{\B} \boxtimes \:_{\B} Y_{\A} \cong \:^{\A}\id_{\A}
        \end{equation}
        and
        \begin{equation}\label{diag5}
            \:^{\B} X^{\A} \boxtimes \:_{\A} Y_{\B} \cong \:^{\B}\id_{\B}
        \end{equation}
        where $\:^{\A} \id_{\A}$ and $\:^{\B} \id_{\B}$ are the identity DA-bimodules over $\A$ and $\B$ respectively (as defined in Section~\ref{kos}), and $\boxtimes$ is a notion of tensor product defined in the sequel. 
    \end{theorem}
    
    In this paper, we start with a graph as in Figure~\ref{star71}, and construct a corresponding weighted $\A_{\infty}$-algebra $\A$. Then, in order to better understand these higher operations, we consider diagrams as in
	\begin{figure}[H]
		\includegraphics[width = 8cm]{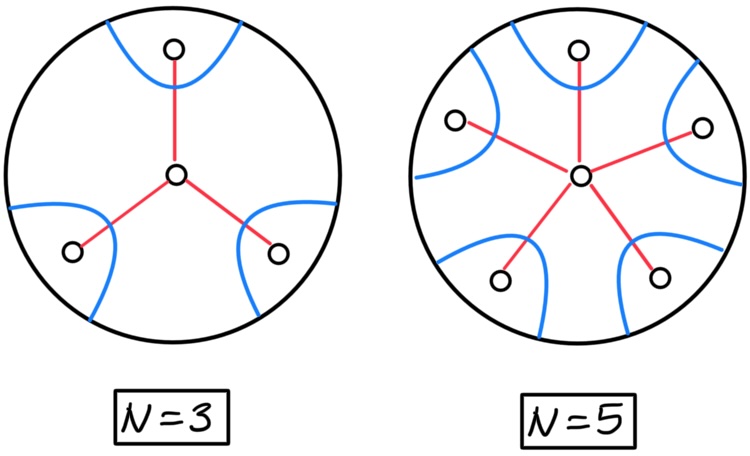}
		\caption{}\label{stardiag}
	\end{figure}
	\noindent and associate to the blue portion of the diagram a second weighted $\A_{\infty}$-algebra $\B$. Finally, we construct in Sections~\ref{aabim} and~\ref{DDbim} a pair of dualizing bimodules, which we use in the sequel to verify the duality relations in Theorem~\ref{duality}.

	This Koszul duality result can be placed in a larger framework having to do with various folklore duality conjectures. It is expected that to any surface with a handle decomposition one can associate an appropriate weighted algebra. Then look at a dual decomposition of the surface. It is then expected that the algebra associated with this dual decomposition will be Koszul dual to one another. For examples of existing results in the literature, see for instance the work of Ozsv\'ath and Szab\'o in~\cite{Pong1}, Zarev's work on the bordered sutured construction,~\cite{Zar}, or the work of Roberts~\cite{Rob} or Manion~\cite{ManAlg} on bordered algebras. The duality result of Theorem~\ref{duality} can be interpreted as a special case of this general principle.
	
	From an algebraic standpoint, this paper exhibits a number of new phenomena. Primarily, it gives combinatorially computable examples of $\A_{\infty}$-structures of various types -- operationally bounded, non-operationally bounded, and weighted -- all with a clearly discernible geometric foundation. 
	
	These more complicated $\A_{\infty}$-structures arise when we allow disks in a bordered Heegaard diagram to cross punctures in the partial Heegaard diagram, and this paper also gives explicitly computed examples of what happens on the algebraic structures when we do this. (For analogous computations in the bordered Heegaard Floer case, see e.g.~\cite{Torus}.) In particular, these orbits are what force us to include higher multiplications and weighted operations on both algebras. One side (the $\B$-side) is always operationally bounded, but on the $\A$-side, we have both weighted operations and non-operationally bounded operations. To keep track of all of the algebra structures, we therefore have to develop new notation. See Sections~\ref{alph} and~\ref{bb1} for more details.
	
	In a similar vein, this paper also gives an exposition of (and concrete examples for) the system of graphical calculus by which $\A_{\infty}$-operations and relations can be represented using planar trees. While this is discussed in full in~\cite{DiagBible}, a reader encountering this notation for the first time could also view Section~\ref{definitions} as a very brief introduction to this way of managing the arithmetic of $\A_{\infty}$-structures.
	
	Also worth noting are the possible connections between the algebras constructed here, and the fully wrapped Fukaya category of $\mathbb{C}P^1$ -- see e.g.~\cite{AurFuk1} and~\cite{AurFuk2}. In the related case of the Pong algebra considered by Ozsv\'ath and Szab\'o, this connection is explored in~\cite{Pong2}, and is also discussed in~\cite{graph}. It is expected that the algebras constructed here may likewise be represented using the fully wrapped Fukaya category; this will be explored further in a forthcoming paper.
	
	The structure of this paper is as follows. In Section~\ref{definitions}, we define the notions of $\A_{\infty}$-algebras, maps, and bimodules used in the rest of the paper. 
    In Sections~\ref{alph} and~\ref{bb1}, we construct the pair of dual $\A_{\infty}$-algebras $\A$ and $\B$ for this class of diagram and verify that these satisfy the $\A_{\infty}$-relations defined in Section~\ref{definitions}. In Sections~\ref{aabim} and~\ref{DDbim}, we construct the bimodules needed to prove the duality relation. 
	
	\begin{ackno}
		\emph{I would like to thank Peter Ozsv\'ath and Zolt\'an Szab\'o for many interesting discussions and a great deal of very helpful advice in the preparation of this paper. I would also like to thank Robert Lipshitz and Andy Manion for many helpful comments on early drafts of this paper.}
	\end{ackno}
	
	\section{Definitions}\label{definitions}
	
	The goal of this section is to set up the machinery and graphical calculus that will be used to construct the algebraic objects in the following sections. It should be emphasized again that this section is meant rather to introduce the graphical calculus for $\A_{\infty}$-algebras used in this paper, than to serve as a comprehensive introduction to $\A_{\infty}$-structures. For such an introduction, see for instance Chapters 2, 6, and 7 of~\cite{BordBook}. 
	
	While much of the exposition in this section follows~\cite{BordBook},~\cite{DiagBible},~\cite{Kauff},~\cite{AlgMat}, etc., some of the notation and definitions are specific to this paper. The reader less familiar with the $\A_{\infty}$-structures discussed here should therefore concentrate mainly on the arithmetic of trees used to calculate $\A_{\infty}$-operations and relations: it is this which will be used in the body of the paper. 
	
	\subsection{$\A_{\infty}$ algebras, and operations as trees}
		
		In this subsection we will introduce both unweighted and weighted $\A_{\infty}$-algebras. We start with the definition of an unweighted $\A_{\infty}$-algebra. 
		
		Start with a ring $R$ of characteristic two, usually $\F[V_0,\ldots, \: V_{N +1}]$, where $\F = \F_2$. An \emph{unweighted $\A_{\infty}$ algebra} is a $R$-module $\A$ equipped with multiplication maps
		\[
			\mu_n : \A^{\otimes n} \to \A
		\]
		where $n \in \Z_{\geq 0}$ and the tensor products are over $R$. In particular, the unweighted $\mu_0$ is a map $R \to \A$ and is determined by its action on $1 \in R$, so it can be viewed as  an element of $\A$. In more familiar language, $\mu_0$ is curvature (if it is nonzero), $\mu_1$ is the differential, and $\mu_2$ is the usual multiplication. The rest are just generalized operations. We require that composition of the $\mu_n$'s satisfies the relation
		\begin{equation}\label{defs1}
			\sum_{\begin{matrix}\substack {1 \leq r \leq n \\ 1 \leq j \leq (n - r)} \end{matrix}} \mu_{n - r + 1} (a_1, \ldots, a_{j}, \mu_r(a_{j + 1}, \ldots, a_{j + r}), a_{j + r + 1}, \ldots, a_{n}) = 0
		\end{equation}
		The \emph{$\A_{\infty}$ relation for $\A$ with $n$ inputs} is defined by~\eqref{defs1}. 
		
		We usually work with the case where $\mu_0$, the unweighted curvature, is zero. To illustrate some of the basic $\A_{\infty}$-relations, we assume for the rest of the section that $\mu_0 = 0$. In this case, for small $n$, the $\A_{\infty}$ relations are, in more familiar language
		\begin{itemize}			
			\item $n = 1$: $\mu_1 \circ \mu_1 \equiv 0$ i.e. $\de^2 = 0$;
			
			\item $n = 2$: $(\mu_1 \circ \mu_2)(a, b) = \mu_2(\mu_1( a), b) + \mu_2(a, \mu_1 (b))$ is the Leibniz rule, i.e
			\[
				\de (ab) = (\de a) b + a (\de b)
			\]
		\end{itemize}
				
		A \emph{corolla} is defined to be a planar tree with a single internal  vertex, $n$ input edges and one output edge, as in.
		\begin{figure}[H]
			\includegraphics[width = 4.5cm]{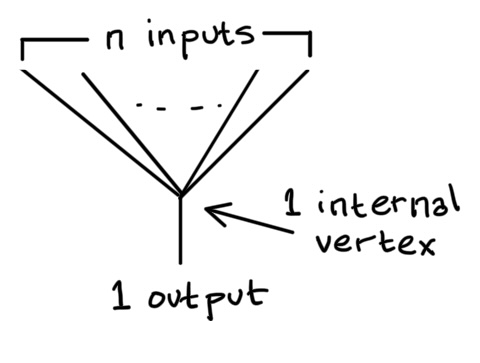}
		\caption{}\label{cor1'}
		\end{figure}
		\hspace{-0.45cm}In this paper, we will most commonly present $\A_{\infty}$-operations and relations in terms of trees. Operatins, in particular, are written in terms of corollas: the corolla from Figure~\ref{cor1'} is defined to represent $\mu_n$.
		
		The next step is to represent the $\A_{\infty}$-relations in terms of trees. We state that any tree determines an $\A_{\infty}$-relation as the sum of all the ways to add a single internal edge to $T$ by replacing an internal vertex with an edge. For instance, in the following cases:
		\begin{itemize}
			\item $\de^2 = 0$: 
			\begin{center}
				 \includegraphics[width = 4cm]{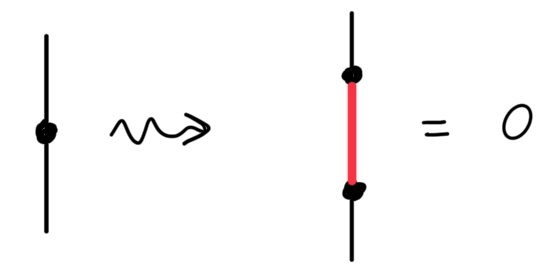}
			\end{center}
			(The red edge represents the new, added edge.)
			
			\item Leibniz rule:
			\begin{center}
				\includegraphics[width = 8cm]{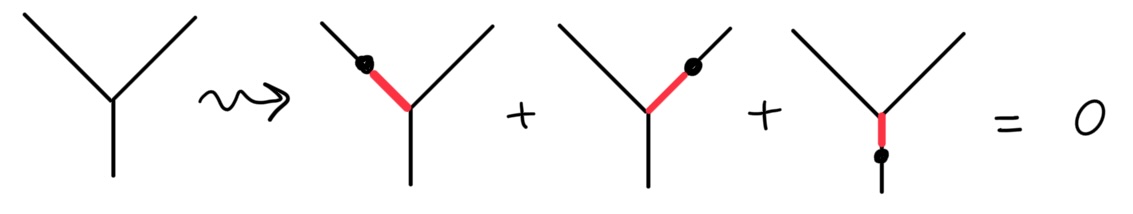}
			\end{center}
			
			\item Associativity fails, but holds up to homology:
			\begin{center}
				\includegraphics[width = 10cm]{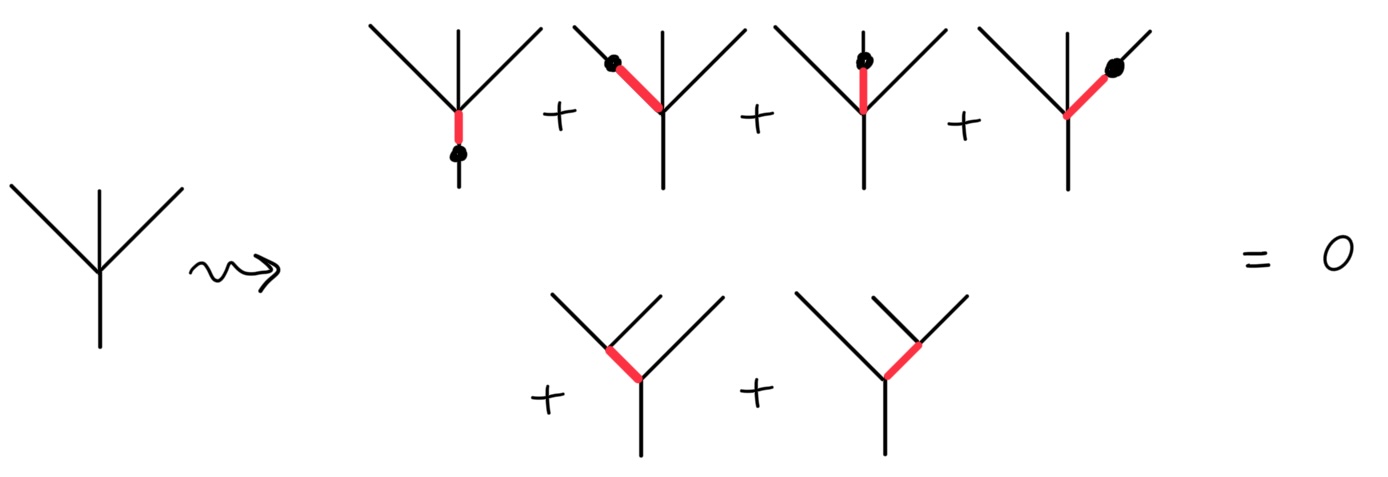}
			\end{center}
			Again, normal associativity would say that the last two terms sum to zero. If we are dealing with cycles $(\mu_1 (a_i) = 0$ for each $i$) and setting boundaries to zero, i.e. working on homology, associativity still holds.
		\end{itemize}
		Occasionally, as we did above, we will not label the input elements at the top of a tree and just write out the trees. 
		
		Next, we define a \emph{weighted $\A_{\infty}$-algebra}: a weighted $\A_{\infty}$ algebra is an $R$-module $\A$, equipped with \emph{operations}, i.e. maps 
		\[
			\mu_n^{\w}: \A^{\otimes n} \to \A.
		\]
		where $\w$ is a weight vector in a finite-dimensional vector-space, usually over $\F_2$. It should be emphasized that while the weight in a weighted operation does have geometric significance, it is, at this point, a purely algebraic construct, with properties to be described in the following paragraphs.
		
		In the case of weighted algebras, we still write operations as trees, except now weighted operations will be labelled with weight at the vertex, i.e.
		\begin{center}
			\includegraphics[width = 8cm]{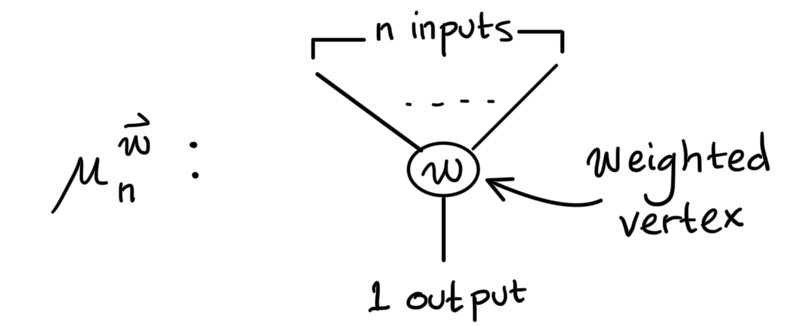}
		\end{center}
		 We will most commonly encounter weighted $\mu_0$s -- namely, $\mu_0^{\e_i} = U_i$, for $0 \leq i \leq N + 1$ -- which will correspond to picking up an orbit in the bimodule. 
		
		Again, \emph{the $\A_{\infty}$-relation for a given (weighted) tree} is the sum over the number of ways to push out an edge. This is equivalent to a weighted version of~\eqref{defs1}, but we will only use this graphical definition in this paper. There are four types of terms that can appear in the $\A_{\infty}$ relation for a given tree.
		 First, recall the notions of a pull, a split, a push, or a differential. In terms of trees, these look like
	\begin{itemize}
		\item Pull:
		\begin{figure}[H]
			\includegraphics[width = 6cm]{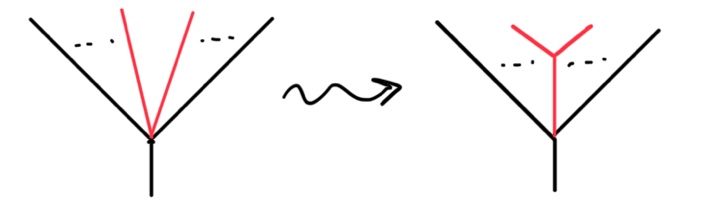}
			
		\end{figure}
		
		\item Split: (Where we are assuming the internal operation is a higher multiplication)
		\begin{figure}[H]
			\includegraphics[width= 12cm]{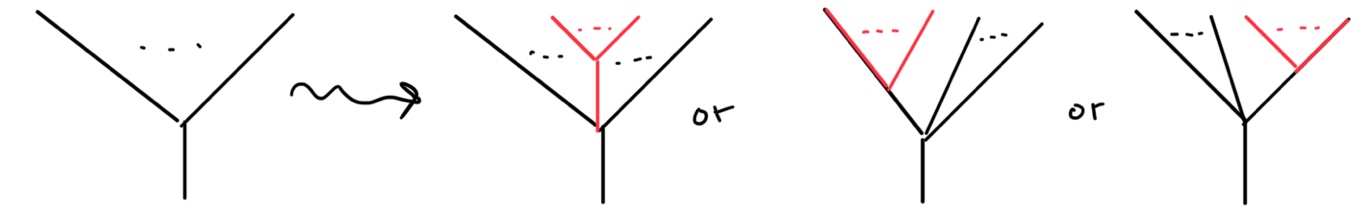}
			
		\end{figure}
		
		\item Push:
		\begin{figure}[H]
			\includegraphics[width = 6cm]{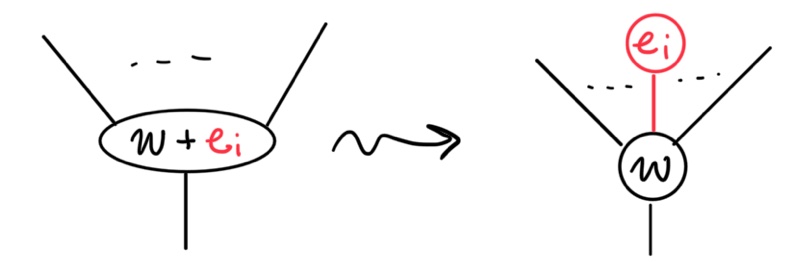}
		\end{figure}
		
		\item Differential:
		\begin{figure}[H]
			\includegraphics[width = 9cm]{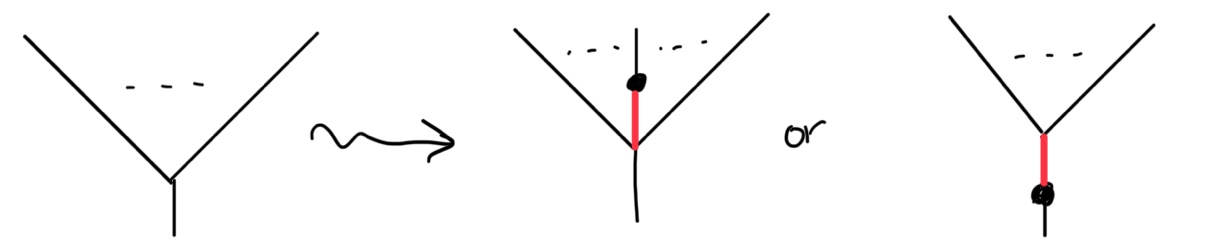}
		\end{figure}
	\end{itemize}
		
		\subsection{AA bimodules}\label{aabimdef}
		
		The basic data for any type of $\A_{\infty}$-bimodule is a $R$-module $X$ and two $\A_{\infty}$-algebras, $\A$, and $\B$, over $R$. An \emph{$AA$-bimodule} $\:_{\B} Y_{\A}$ is an $R$-module $Y$ equipped with maps
				\[
					m_{j|\x|i} : {\B}^{\otimes j} \otimes Y \otimes \A^{\otimes i} \to Y, 
				\]
				written (in terms of trees) as
				\begin{figure}[H]
					\includegraphics[width = 4cm]{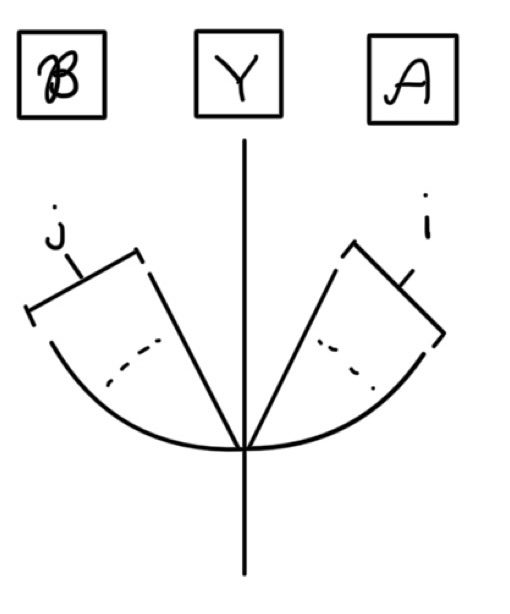}
					
					\caption{}\label{bimdefpic}
				\end{figure}
				\hspace{-0.45cm}that is, with $i$ inputs on the $\A$-side, and $j$ inputs on the $\B$-side. As in the case of $\A_{\infty}$ algebras, the \emph{$\A_{\infty}$ relation for a tree in an $AA$-bimodule} (such as the one above) says that all ways to add a single edge to this tree sum to zero. In more traditional notation, for a fixed $i, j \in \N$, $\x \in Y$, and $a_1, \ldots, a_i \in \A$, $b_1, \ldots, b_j \in \B$, the corresponding $\A_{\infty}$ relation is
				\begin{multline}\label{defs1'}
					\sum_{i'\leq i''} m_{j |\x | i - (i'' - i')}(b_j, \ldots b_1, \x, a_1, \ldots, \mu_{i'' - i' + 1}(a_{i'}, \ldots, a_{i''}), \ldots, a_i) \\+ \sum_{j'\leq j''} m_{j - (j'' - j') |\x | i}(b_j, \ldots, \mu_{j'' - j' + 1}(b_{j''}, \ldots, b_{j'}), \ldots, b_1, \x, a_1, \ldots, a_i)  = 0
				\end{multline}
				We can also have weighted $AA$-bimodules, just as we can have weighted $\A_{\infty}$-algebras. In this case, the $\A_{\infty}$-relation for a given tree $T$ is still the sum over all the ways to add an edge to $T$. Again, there are four ways an edge can be added: namely, a pull, a push, a split, or a differential. In the case of $AA$-bimodules, these look a little different than they did in the case of algebras, so we give examples here:
			\begin{description}
				\item[Pull] Brings together a subset of the elements on either side, as in
				\begin{center}
					\includegraphics[width = 8cm]{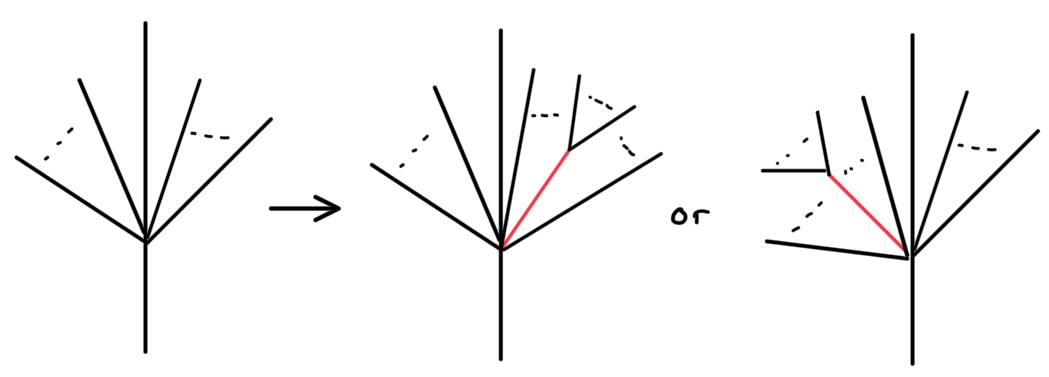}
				\end{center}
				
				\item[Split] Decomposes the original tree into a composition of two bimodule operations, as in:
				\begin{center}
					\includegraphics[width = 5cm]{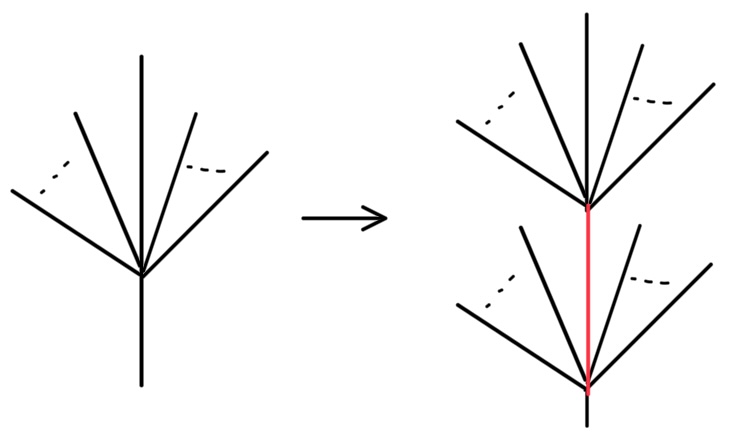}
				\end{center}
				
				\item[Push] Decomposes a weighted tree by pushing out a weighted operation on one side or another, as in:
				\begin{center}
					\includegraphics[width = 10cm]{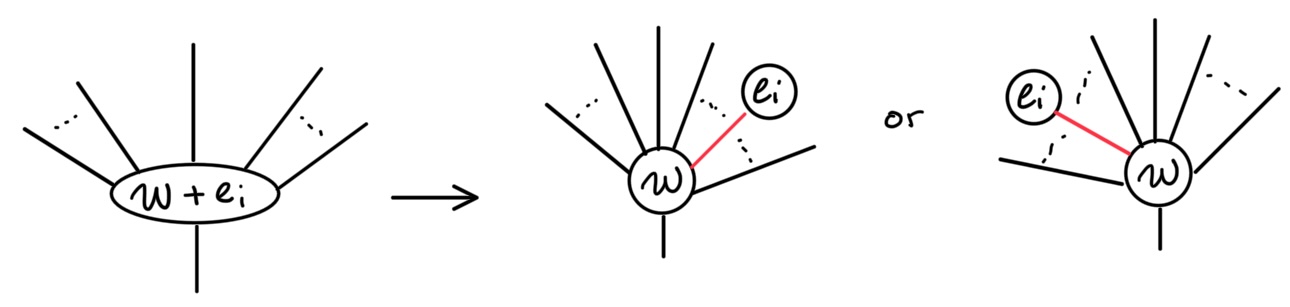}
				\end{center}
				
				\item[Differential] Takes a differential of one of the input elements on either side, as in:
				\begin{center}
					\includegraphics[width = 8cm]{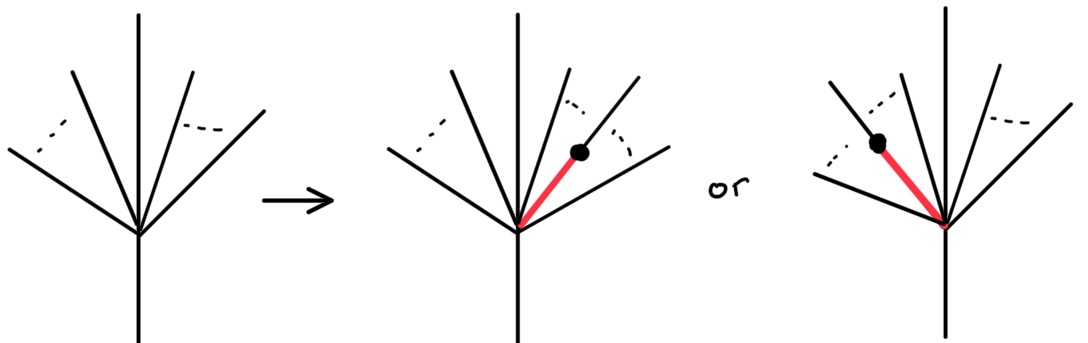}
				\end{center}
				
			\end{description}
		
			\subsection{DD Bimodules}\label{ddorigin}
			
				 A \emph{DD bimodule} $^{\A}X^{\B}$ over $\A$ and $\B$ is an $R$-module with a single basic operation:
		\[
			\delta^1: X \to \A \otimes \B \otimes X.
		\]
		We apply $\delta^1$ repeatedly to get
		\[
			\delta^n = \underbrace{\delta^1 \circ \cdots \circ \delta^1}_{n \text{ times}} : X \to (\A \otimes \B)^{\otimes n} \otimes X.
		\]
		Writing $\I$ for the identity on $X$, the $\A_{\infty}$ relation for $X$ is
		\begin{equation}\label{dd2}
			\sum_{n, \w} (\mu_n^{\w} \otimes \I) \circ \delta^n = 0
		\end{equation}
		where the sum is over all $n \geq 0$ and finite linear combinations $\w$, and the $\mu_n^{\w}$ denote the (weighted) operations on $\A \otimes \B$, as defined below, in Section~\ref{tens}.
		
		\begin{remark}\label{dd135} \emph{For readers familiar with the notion of type D modules, it is worth noting that a $DD$ bimodule over $\A$ and $\B$ is actually just a type-$D$ module over $\A \otimes \B$.}
		\end{remark}

        \subsection{Koszul Duality}\label{kos}

        This subsection adapts the relevant definition from Section 8 of~\cite{Morph} to the current setting. The aim here is to define what it means for two $\A_{\infty}$ algebras to be Koszul dual. We give a more complete definition in the sequel, but this section explains in brief the terms which appear in Theorem~\ref{duality}.
		
		Start with weighted $\A_{\infty}$ algebras $\A$ and $\B$ over a single ground ring $R$ and weight space $\Lambda$, as usual. 

        Define the \emph{identy bimodule} $\:^{\A} \id_{\A}$ to be an $(R,R)$-bimodule with a single non-zero operation $\delta_2^1: \id \otimes \A \to \A \otimes \id$ given by $\delta_2^1(\x \otimes a) = a \otimes \x$. 

        We will define the precise type of box product which appears in Theorem~\ref{duality} in Section 2.5 of the sequel. In brief, this notion of product, $\boxtimes$, gives a means by which to take a tensor product of a DD-bimodule over weighted algebras and a weighted AA-bimodule over weighted algebras, and obtain a different type of weighted $\A_{\infty}$-bimodule called a DA-bimodule.

        With this framework, we say that a DD-bimodule $\:^{\A} X^{\B}$ is \emph{quasi-invertible} if and only if there exists an AA-bimodule $\:_{\B} Y_{\A}$ such that
        \begin{equation}\label{aux1}
            \:^{\A} X^{\B} \boxtimes \:_{\B} Y_{\A} \simeq \:^{\A}\id_{\A}
        \end{equation}
        and
        \begin{equation}\label{aux2}
            \:^{\B} X^{\A} \boxtimes \:_{\A} Y_{\B} \simeq \:^{\B}\id_{\B}.
        \end{equation}
        
        Now, define a \emph{Koszul dualizing bimodule} to be a DD bimodule $^{\A}X^{\B}$ which is quasi-invertible. We say that $\A$ and $\B$ are \emph{Koszul dual} if they admit a Koszul dualizing bimodule $X$.

        Notice first that~\cite{Morph} deals only with the case of bimodules over dg-algebras with ground ring $\F_2$. This makes the algebra significantly simpler. In our case, the ground ring is somewhat larger, and the algebras involved are not only bona-fide $\A_{\infty}$-algebras with non-zero higher operations, but also weighted $\A_{\infty}$-algebras. Proving a duality theorem such as Theorem~\ref{duality} in this case requires a certain amount of new algebraic machinery, which is defined in Section 2 of the sequel. 

        In the sequel, we also verify a stronger result, namely that the bimodules on the right hand sides of~\eqref{aux1} and~\eqref{aux2} are actually isomorphic, as weighted DA-bimodules, to $\:^{\A}\id_{\A}$ and $\:^{\B}\id_{\B}$, respectively, rather than just chain homotopy equivalent as in~\eqref{aux1} and~\eqref{aux2}.

        \section{Diagonals and tensor products}\label{diag}
	
	The goal of this section is to define \emph{weighted diagonals}, especially in the context of our algebras $\A$ and $\B$, and discuss the role they play in describing operations on $\A_{\infty}$-tensor products, both of algebras (e.g. $\A \otimes \B$) and bimodules (e.g. $X \boxtimes Y$). This will be relevant to the verification of the $\A_{\infty}$-relations for the $DD$-bimodule $\:^{\A} X^{\B}$ in Section~\ref{DDbim}, since this module is in fact just a type-$D$ module over $\A \otimes \B$.

	\subsection{Diagonals}\label{diagsfirst}
	
	The main source for this is Section 6 of~\cite{DiagBible}. In this section, we summarize the definitions and results most relevant the constructions that follow, modified to fit this situation. The main definition of this section is a \emph{weighted diagonal}, given in~\eqref{diagstartshere}. The other definitions that lead up to this point in the section are intended entirely to prepare the reader to understand the stipulations of this most important definition.
	
	The basic set-up is as follows. Throughout, a \emph{weight vector} $\w$ is defined to be a finite $\Z_{\geq 0}$-linear combination of the vectors $\{\e_i\}_{i = 0}^{N + 1}$, viewed an element of the $\Z_{\geq 0}$-linear space $V$ generated by these basis vectors. A \emph{weighted tree} $T$ is a planar tree where each internal vertex $v$ of $T$ is labelled with a particular weight vector $\w(v) \in V$. The \emph{total weight} $\w(T)$ of $T$ is the sum of all the $\w(v)$, over all the internal vertices of $T$. Define $\wt(T)$ to denote the magnitude of the total weight $\w(T)$ of $T$. Define the \emph{formal dimension} $\dim T$ of a tree $T$ with $n$ inputs, total weight $\w(T)$, and $|v|$ internal vertices as
	\begin{equation}\label{diag9}
		\dim T = n + 2\wt(T) - |v| -1
	\end{equation}
	In this section, we will consider only a subset of the total set of trees, namely \emph{stably weighted trees}. To define \emph{stably weighted}, we start by looking at a given vertex $v$ in a tree $T$, with corresponding weight $\w(v)$. This vertex $(v, \w(v))$ is defined to be \emph{stable} when $v$ is $3$-valent, $|\w(v)| > 0$, or both. The tree $T$ is \emph{stably weighted} if all its vertices are stable. So, for instance, considering
	\begin{center}
		 \includegraphics[width = 9cm]{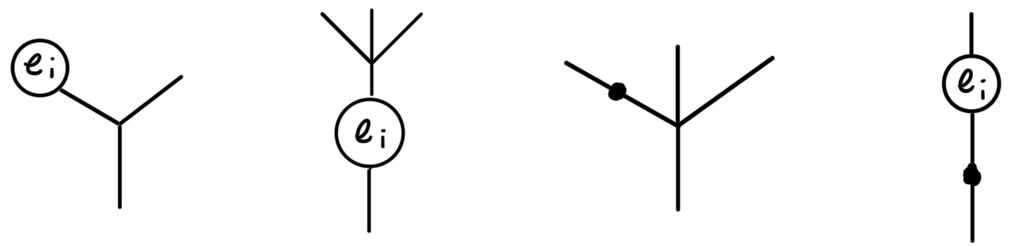}
	\end{center}
	the left two trees are stably weighted, and the right two are not.
	
	Define $X^{n, \w}_k$ to be the set of stably weighted trees with $n$ inputs, total weight $\w$, and formal dimension $k$. Define a \emph{differential} on $X^{n, \w}_k$ as follows. Let $\del T$ be sum of all the ways to replace a vertex of $T$ with an edge in such a way that the result is still a stably weighted tree. For instance:
	\begin{center}
		\includegraphics[width = 6cm]{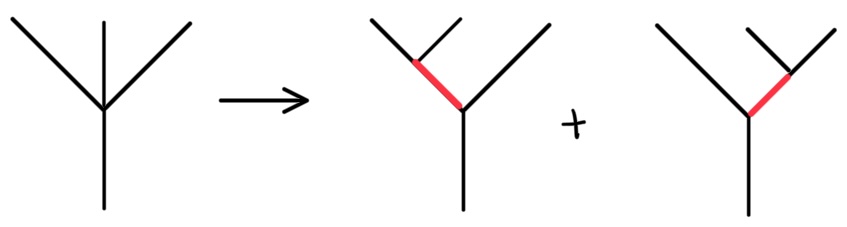}
	\end{center}
	or
	\begin{center}
		\includegraphics[width = 6cm]{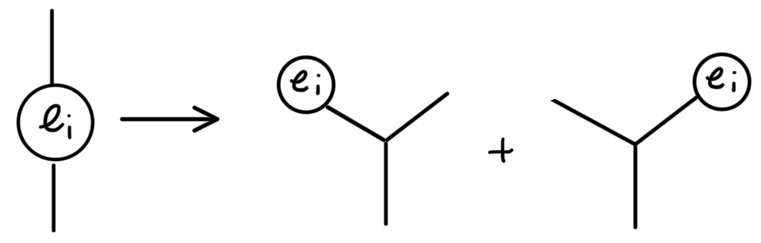}
	\end{center}
	Notice that when we look at the ways to push out an edge from the 3-input corolla, above, we do not get
	\begin{center}
		 \includegraphics[width = 2cm]{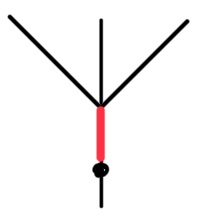}
	\end{center}
	as a summand of the result, because this tree is not stably weighted (it has a 2-valent vertex). $\del$ is clearly a differential, i.e. $\del^2 = 0$, and drops dimension by 1 in each case. 
	
	Common conventions are as follows: we usually want our trees to be in dimension 0, and we usually suppress the dimension and speak of $X_*^{n, \w} = \bigoplus_{k \in \Z} X_{k}^{n, \w}$ as a chain complex. The \emph{input leaves} of a tree $T$ are defined to be the edges adjacent to the 1-valent input vertices of $T$. The \emph{output leaf} of $T$ is the edge adjacent to the 1-valent output vertex of $T$. 
	
	 We also use special symbols for the \emph{corollas} -- trees with one internal vertex, either unweighted or weighted, and any number of inputs. We let $\Psi_n^{\w}$ denote corolla with $n$ inputs and weight $\w$ placed at its single vertex. We permit $n = 0$, and in particular, the popsicles
	\begin{center}
		\includegraphics[width = 3cm]{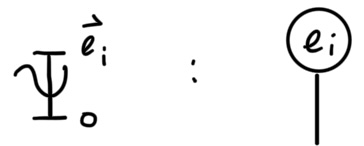}
	\end{center}
	as well as $\Psi_1^0$, which will give back the differential of an element.
	
	We also want to compose trees, e.g.
	\begin{center}
		\includegraphics[width = 7cm]{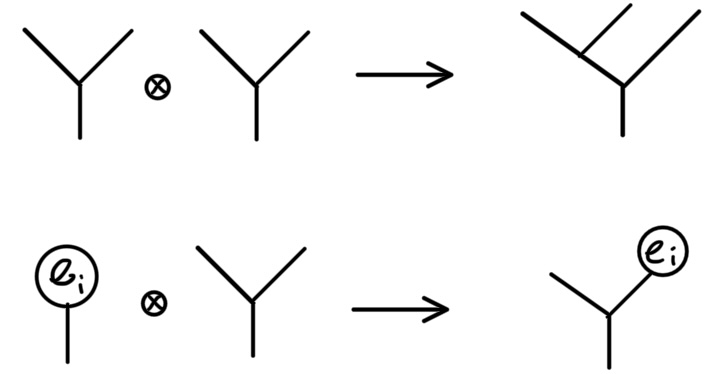}
	\end{center}
	We do this explicitly with \emph{stacking maps} $\phi_{i,j,n;\vv,\w}$. To define these, we need a little more notation, namely the \emph{gluing map} $\circ_i$; for each $S \in X^{m, \vv}_*, \: T \in X^{n, \w}_*$, $T \circ_i S$ is defined as the result when we attach the output leaf of $S$ to the $i$-th input of $T$. Here, we are not considering the actual labelling of the inputs / outputs, but only the trees themselves. For instance, the two stacked trees above are:
	\begin{center}
		\includegraphics[width = 6cm]{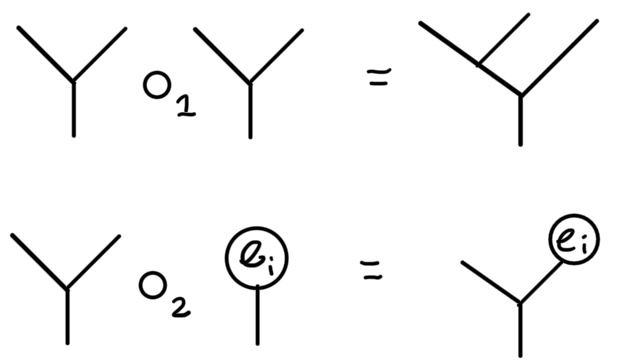}
	\end{center}
	With this notation, we further define the \emph{stacking map}
	\begin{equation}\label{diag10}
		\phi_{i,j,n;\vv,\w} : X_*^{(j - i + 1),\vv} \otimes X_*^{(n + i - j), \w} \to X_*^{n, \vv + \w}
	\end{equation}
	for $1 \leq i \leq n$ (so we can choose any branch to stack into) and $i-1 \leq j \leq n$, as
	\begin{equation}\label{diag11}
		\phi_{i,j,n;\vv,\w}(S \otimes T) = T \circ_i S 
	\end{equation}
	So the examples given above can be written in these terms as
	\begin{center}
		\includegraphics[width = 8cm]{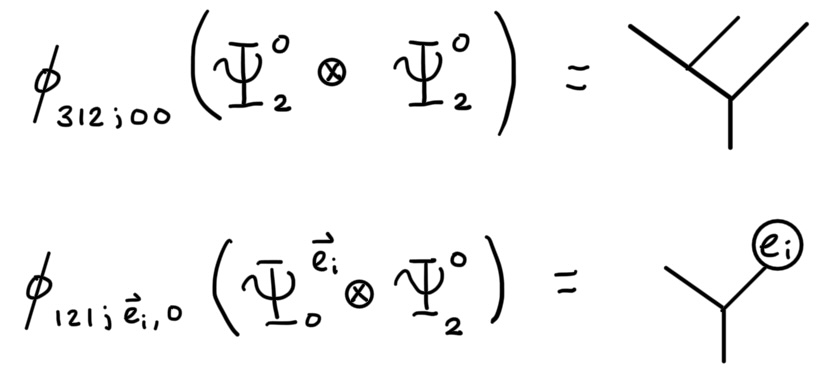}
	\end{center}
	When we define weighted diagonals, we will work over our usual ground ring $R = \F_2[V_0, \ldots, V_{N + 1}]$. We define \emph{weight for tensor products of trees} as 
	\begin{equation}\label{diag6}
		\begin{cases}
			\wt_1(V_0^{a_0} \cdots V_{N + 1}^{a_{N + 1}} \cdot S \otimes T) = a_0 + \wt (S) \\
			\wt_2(V_0^{a_0} \cdots V_{N + 1}^{a_{N + 1}} \cdot S \otimes T) = \sum_{i = 1}^{N +1} a_i + \wt (T)
		\end{cases}
	\end{equation}
	where again $\wt (S), \wt (T)$ just denote the magnitudes of the total weights of $S, T$, respectively.
	
	The last piece of the set-up is the \emph{generalized set of trees} we will use in our definition of the diagonal. We start with $X_*^{n, \w}$ and add
	\begin{enumerate}[label = (\roman*)]
		\item The \emph{stump}, $\top$, with no vertices and no inputs. Making the (formal) stipulation that $\top$ actually has $-1$ internal vertices,~\eqref{diag9} gives $\dim \top = 0$. We write $\tilde{\mathcal{T}}_{0,0}$ for this stump. 
		
		\item The \emph{shoot}, $\downarrow$, with one input and no internal vertices, so that $\dim \downarrow = 0$. We write $\tilde{\mathcal{T}}_{1,0}$ for this shoot.
	\end{enumerate}
	We then define 
	\begin{equation}\label{diag12}
		\tilde{X}_*^{n, \w} = \begin{cases} X^{n, \w}_* & n > 0 \\ R(\tilde{\mathcal{T}}_{0,0}) & n = 0, \w = 0; \\ R(\tilde{\mathcal{T}}_{1,0}) & n = 1, \w = 0; \end{cases}
	\end{equation}
	We extend the differential to the complex $\tilde{X}_*^{n, \w}$ by letting it be the same as usual when $n > 0$, and vanishing on $\tilde{X}^{0,0}_*, \tilde{X}^{1,0}_*$. To extend the gluing / composition map $\circ_i$, we make the following definitions:
	\begin{enumerate}[label = (\roman*)]
		\item For $\top$, and any suitable $T$, we set $T \circ_i \top = 0$ unless the vertex adjacent to the $i$-th input leaf of $T$ is 3-valent -- i.e. near this vertex, $T$ looks like the two-pronged corolla $\Psi_2^0$ -- in which case, $T \circ_i \top$ is the tree obtained from $T$ by removing the $i$-th input and the adjacent leaf. For instance:
		\begin{center}
			\includegraphics[width = 7cm]{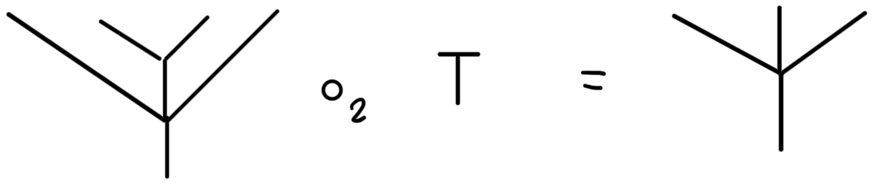}
		\end{center}
		
		\item Any composition with $\downarrow$ is the identity $\tilde{X}^{n, \w} \to \tilde{X}^{n, \w}$. 
	\end{enumerate}
	
	We are now ready to define diagonals. A \emph{weighted diagonal} is a family of maps
	\begin{equation}\label{diagstartshere}
		\Gamma^{n, \w}_*: X^{n, \w}_* \to \bigoplus_{\tiny |\w - \w_1- \w_2| \geq 0 } \tilde{X}^{n, \w_1} \otimes \tilde{X}^{n, \w_2} \otimes R
	\end{equation}
	satisfying conditions~\ref{presdim} through~\ref{basecase}, where $n \in \Z_{\geq 0}$ and $\w$ ranges over all finite sums of the basic weight vectors $\{\e_i\}_{i = 0}^{N +1}$. The rules we require $\Gamma$ to satisfy are as follows:
	\begin{enumerate}[label = \textbf{WD\arabic*.}]
		\item\label{presdim} (Preserves dimension) For $n, \w$, and $T \in X^{n, \w}_*$,
		\begin{equation}\label{diag3}
			\dim \Gamma(T) = \dim T
		\end{equation}
		with the convention that $\dim (S \otimes T) = \dim S + \dim T$ for any tensor product of trees;
		
		\item\label{weighthom} (Preserves weight) For any tree $T\in X^{n, \w}$, we have
		\begin{equation}\label{diag5}
			\wt_1(\Gamma^{n, \w}(T)) = \wt_2(\Gamma^{n, \w}(T)) = \wt(T) = |\w|.
		\end{equation}
		where the weights on each side are defined as in~\eqref{diag6}.
		
		\item\label{stack} (Stacking) The requirement is
		\begin{equation}\label{diag7}
			\Gamma^{n, \vv + \w} \circ \phi_{i,j,n;\vv,\w} = \sum_{{\tiny \begin{matrix} |\vv - \vv_1 - \vv_2| \geq 0 \\ |\w - \w_1 - \w_2| \geq 0\end{matrix}}} (\phi_{i,j,n;\vv_1,\w_1} \otimes \phi_{i,j,n;\vv_2,\w_2})\circ(\Gamma^{(j - i + 1), \vv} \otimes \Gamma^{(n + i - j + 1), \w})
		\end{equation}
		
		\item\label{basecase} (Non-degeneracy / Base-case requirements)
		\begin{enumerate}[label = (\alph*)]
			\item (Base case for unweighted) By~\ref{weighthom}, $\Gamma^{2,0}$ is a map from $X^{2,0} \to X^{2,0} \otimes X^{2,0}$. But $X^{2,0}$ is canonically isomorphic to $R$, with generator the unique two pronged corolla. As already noted, we then require that $\Gamma^{2,0}$ is the canonical isomorphism:\vspace{3pt}
			\begin{center}
				\includegraphics[width = 8cm]{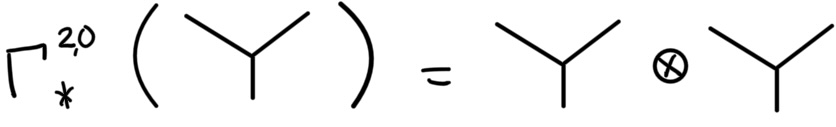}
			\end{center}
			
			\item (Base case for weighted) We require that $\Gamma^{0,\e_i}$ is the predetermined seed for $\e_i$,
			\begin{center}
				\includegraphics[width = 10cm]{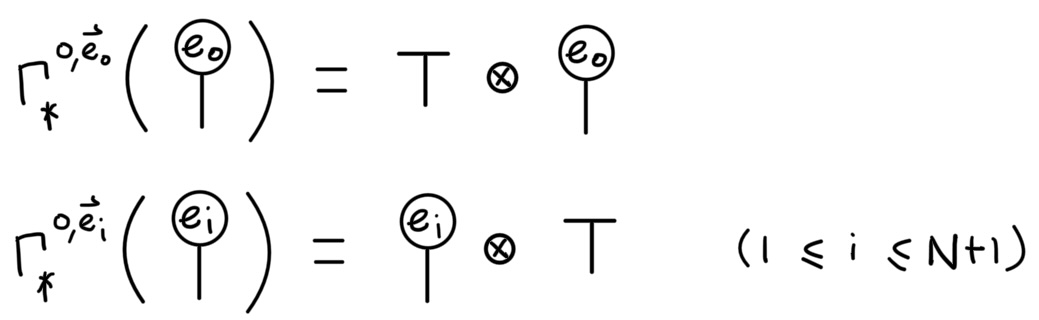}
			\end{center}
			
			\item (Ruling out bad cases) At most one of each pair of trees in the diagonal is a stump or a shoot, that is, the image of $\Gamma^{n, \w}$ is in
			\begin{equation}\label{diag8}
				\bigoplus_{\tiny{ |\w - \w_1 - \w_2| \geq 0}} ( \tilde{X}^{n, \w_1}_* \otimes X^{n, \w_2}_* \otimes R + X^{n, \w_1}_*\otimes \tilde{X}^{n, \w_2}_* \otimes R)
			\end{equation}
		\end{enumerate}
	\end{enumerate}
	
	In practice, diagonals satisfying~\ref{presdim}-~\ref{basecase} are constructed using the following recipe, rather than found in nature. Starting from base-cases, or \emph{seeds}, which we choose manually at the beginning, a diagonal is built by induction on the number of inputs and total weight. The seed for the unweighted part of the diagonal (i.e. the $\Gamma^{n, 0}_*$, which have only unweighted trees in their domains) is by definition always 
	\begin{center}
		\includegraphics[width = 7cm]{diag1}
	\end{center}
	From this comes the unweighted part of the diagonal, via an acyclic models type construction: for an example and further details of this part of the construction, see the discussion, surrounding~\eqref{indstep}, below. 
	
	To construct the weighted part of the diagonal, i.e. the maps $\Gamma^{n, \w}$ with $|\w| > 0$, we have to manually choose our seed trees. We choose a set of seed trees $T_i$, each of which is required to be a linear combination, with coefficients either 0 or 1, of the following elements:
	\begin{equation}\label{diag1}
		V_{i} \cdot \Psi_0^{\e_i} \otimes \top, \quad V_{0} \cdot \top \otimes \Psi_0^{\e_i} 
	\end{equation}
	where $i$ ranges over $1 \ldots, N + 1$. We then define $\Gamma^{0, \e_i}(\Psi_0^{\e_i}) = T_i$ for each $\e_i$, $0 \leq i \leq N + 1$. To construct the rest of the diagonal, we again apply an inductive acyclic models type construction; that is, choose $n \geq 0$ and a weight $\w$, and a tree $T \in X^{n, \w}_*$. Assume we have already defined $\Gamma^{m, \vv}_*$ for the cases
	\begin{itemize}
		\item  $m \leq n$, $|\vv| < |\w|$ and
		\item $m < n$, $|\vv| \leq |\w|$
	\end{itemize}
	This means that we know what $\Gamma^{*, *}_* \del T$ is, and one can show that it is a cycle. Then, since the associahedron is acyclic, we can guarantee that $\Gamma^{*,*}_* \del T = \del T'$ for some 
	\begin{equation}\label{indstep}
		T' \in \bigoplus_{\tiny{ |\w - \w_1 - \w_2| \geq 0}} \tilde{X}^{n, \w_1}_* \otimes\tilde{X}^{n, \w_2}_* \otimes R.
	\end{equation}
	We then choose such a $T'$, and define $\Gamma^{n, \w} T = T'$.

	In the remainder of this subsection, we give further details about the inductive step by which we construct the diagonal from a the base-case. We work with the unweighted case for this example. The base case is given from~\ref{basecase}, so the next step is to define $\Gamma^{3,0} \Psi_3^0$, that is, when we ask
	\begin{center}
		\includegraphics[width = 4cm]{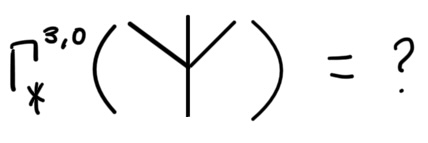}
	\end{center}
	we have
	\begin{center}
		\includegraphics[width = 7cm]{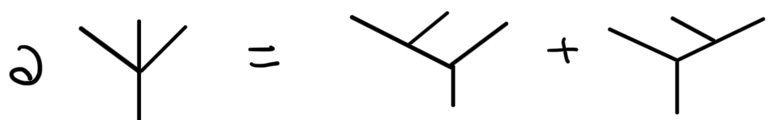}
	\end{center}
	The stacking rules give us
	\begin{center}
		\includegraphics[width = 8cm]{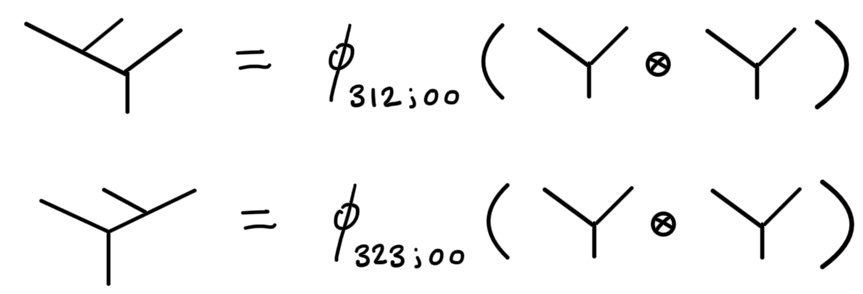}
	\end{center}
	so that
	\begin{center}
		\includegraphics[width = 14cm]{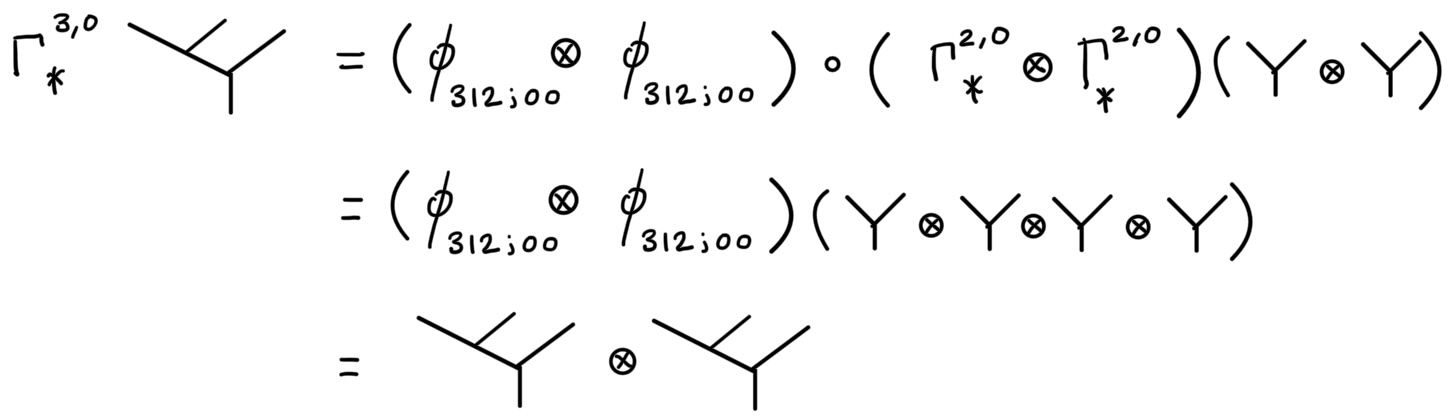}
	\end{center}
	and likewise,
	\begin{center}
		\includegraphics[width = 8cm]{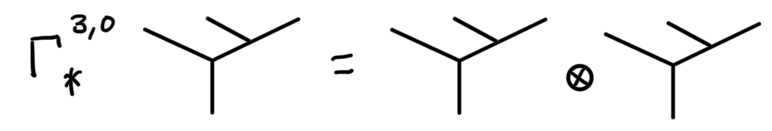}
	\end{center}
	Incidentally, this holds (by induction) for all compositions of $\mu_2$s, so e.g.
	\begin{center}
		\includegraphics[width = 10cm]{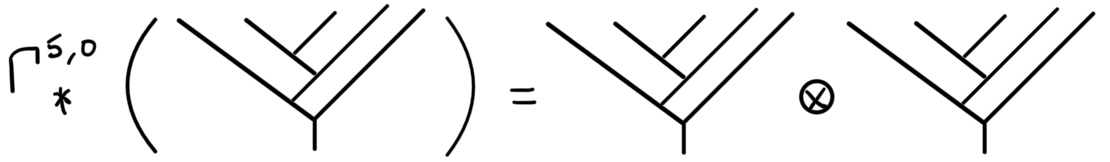}
	\end{center}
	and so on, for any number of inputs, provided that all operations involved are $\mu_2$s. But, using this for the case $n = 3$, we get
	\begin{center}
		\includegraphics[width = 12cm]{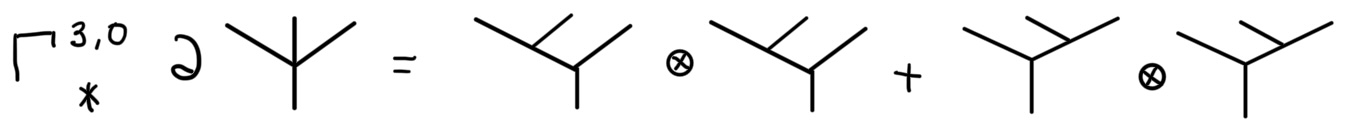}
	\end{center}
	And we have 
	\begin{center}
		\includegraphics[width = 14cm]{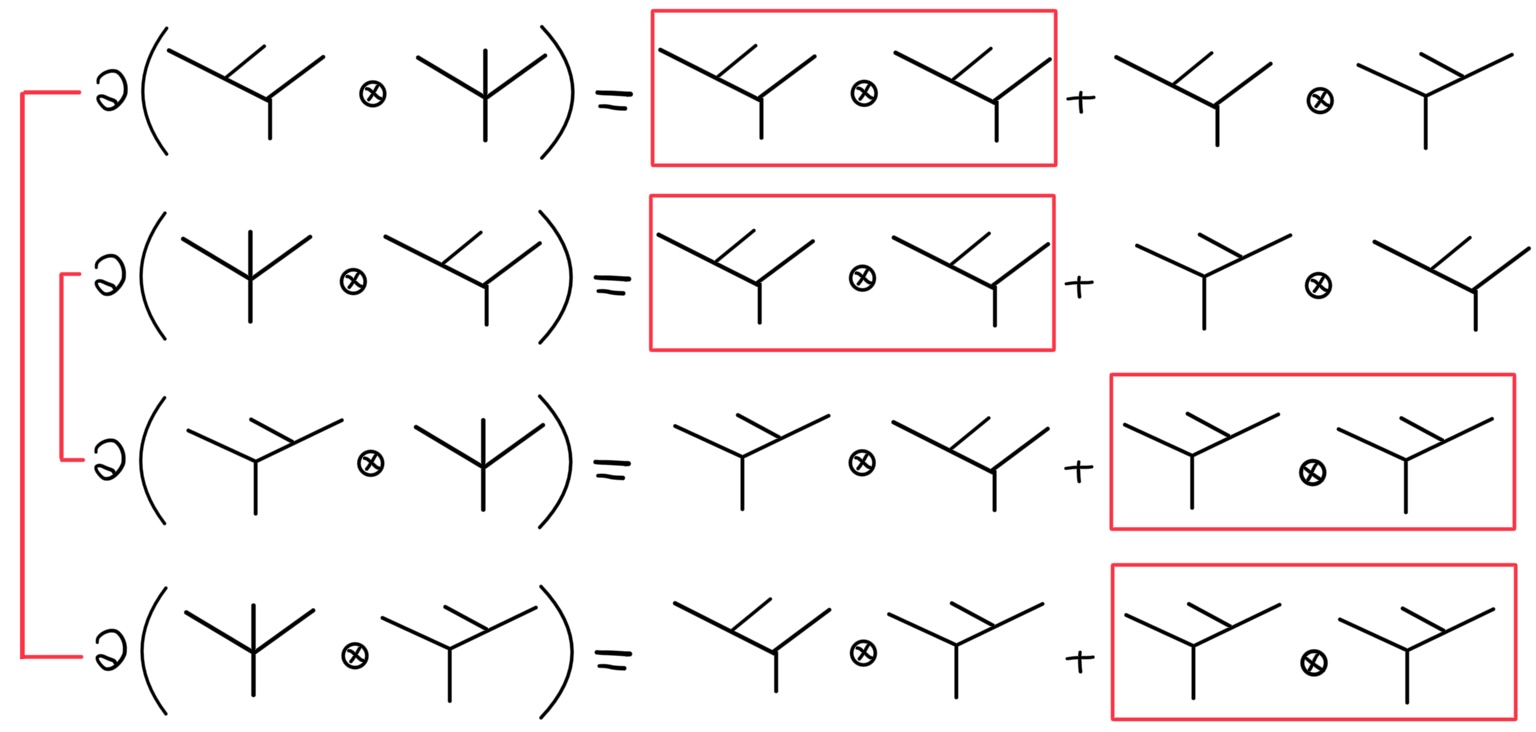}
	\end{center}
	that is, we have two choices that give the correct differential, namely:
	\begin{center}
		\includegraphics[width= 14cm]{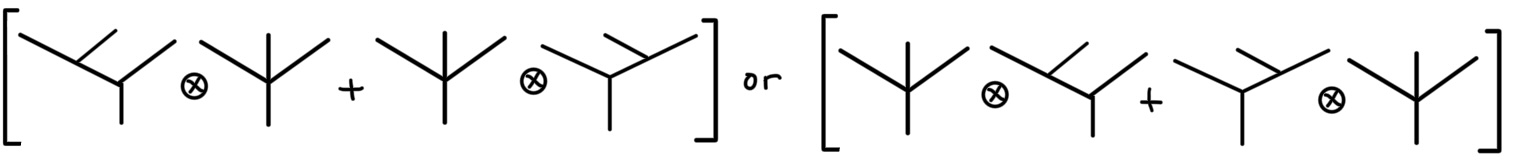}
	\end{center}
	We choose:
	\begin{center}
		\includegraphics[width = 10cm]{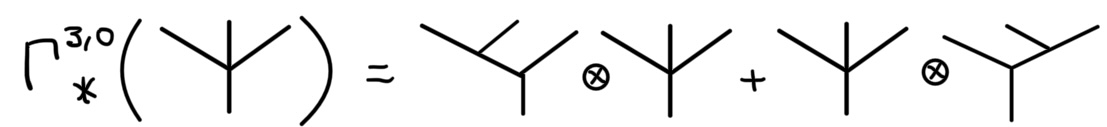}
	\end{center}
	This is what is called the \emph{right-moving} choice, and in general, we will always choose the right-moving options. For general trees, the notion of being \emph{right-moving} is defined as follows. The basic right-moving pairs of trees are
	\begin{center}
		\includegraphics[width = 10cm]{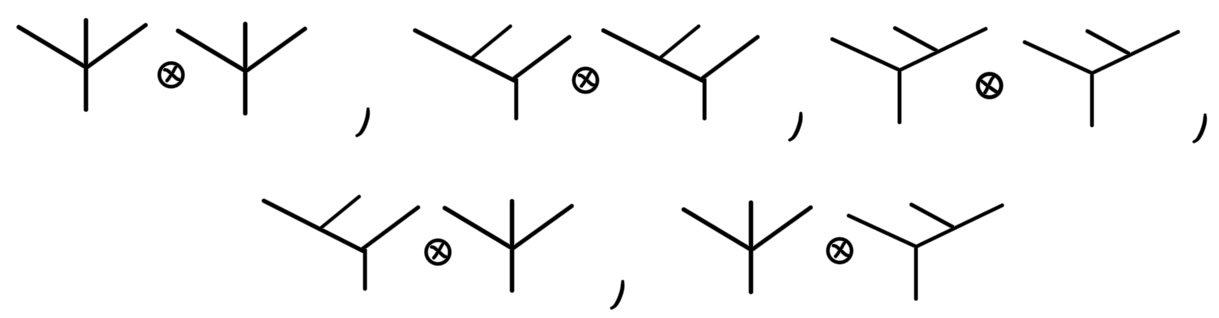}
	\end{center}
	For more general trees, we need to use the 	notion of a \emph{profile tree} corresponding to a subset of the inputs. For $T \in X_*^{n, \w}$ and $I := \{i_1 < \cdots < i_k\}$ with $i_1 \geq 1, \: i_k \leq n$, the profile tree corresponding to $I$, written $T(I)$ is the tree that results when we delete the inputs whose indices are not in $I$, and then delete any 2-valent vertices that result. For instance, with $I = (1, 3, 4)$, and a particular choice of $T$, this looks like:
	\begin{center}
		\includegraphics[width = 7cm]{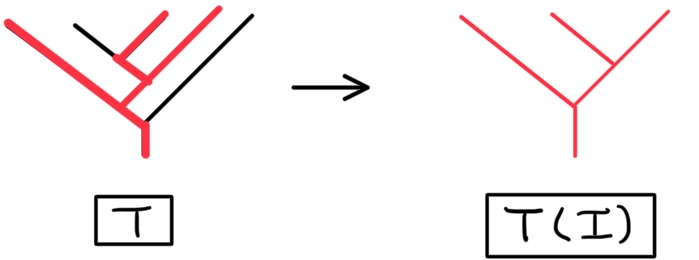}
	\end{center}
	The profile tree for a product $S \otimes T$ of trees (both with $n$ inputs) corresponding to the set $I$ is written $(S \otimes T) (I)$, and is defined to be is $S(I) \otimes T(I)$. A product of trees $S \otimes T$ is defined to be \emph{right-moving} if every profile tree $(S \otimes T)(I)$ corresponding to a 3-element subset $I$ is one of the basic right-moving pairs. So, for instance,
	\begin{center}
		\includegraphics[width = 6cm]{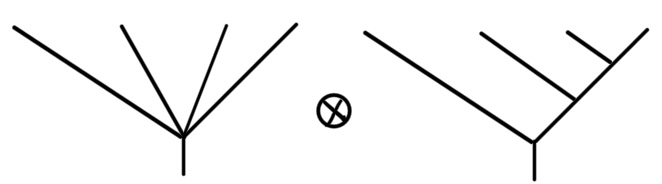}
	\end{center}
	is right moving, while
	\begin{center}
		\includegraphics[width = 6cm]{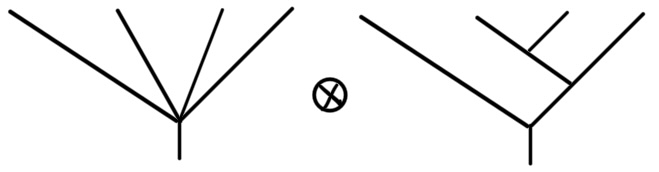}
	\end{center}
	is not. Note that ``right-moving or not'' does not depend on the weight, or its distribution in the tree. So when determining whether a tree is or is not right-moving, we can ignore weight. It is a result of Masuda-Thomas-Tonks-Vallette,~\cite{MasAssoc}, cited in~\cite{DiagBible}, that we can find a right-moving diagonal in the unweighted case, which is all we will need to prove Theorem~\ref{duality}. 
	
	\subsection{$\A_{\infty}$-structures for tensor products of algebras}\label{tens}
	
	So far, all discussion of diagonals has been entirely on the level of trees, with no reference to the actual inputs (which in our case, will be algebra elements). The notion of a weighted diagonal allows us to define the $\A_{\infty}$-structure for a tensor product of algebras $\A \otimes \B$, in the following way. 
	
	First, it is clear from the construction that a weighted diagonal $\Gamma$ is a chain map on the complex of trees. Next, we label the inputs with algebra elements. The data for an operation on $\A \otimes \B$ is as follows 
	\begin{itemize}
		\item $\Gamma^{n, \w}_* T$, for $T \in X^{n, \w}_*$, $n \in \Z_{\geq 0}$, and $\w$ ranging over finite linear combinations of the basic weight vectors;
		
		\item A collection $\{a_i \otimes b_i\}_{i = 1}^n$ of pairs in $\A \otimes \B$.
	\end{itemize}
	The corresponding operation is now the sum over all products of trees which make up $\Gamma^{n, \w}_* T$, with the inputs of the left tree (of each pair) labelled with $\{a_i\}_{i = 1}^n$, and the inputs of the right tree labelled with $\{b_i\}_{i = 1}^n$. We use the usual maps
	\begin{align*}
		\mu&: X^{n, \w}_* \to \mathrm{Mor}(\A^{\otimes n}, \A) \\
		\nu&: X^{n, \w}_* \to \mathrm{Mor}(\B^{\otimes n}, \B)
	\end{align*}
	to go from trees to valid operations, on the $\A$-side and the $\B$-side respectively. We further stipulate that $\mu(\top) = V_0$, and $\nu(\top) = \sum_{i =1}^N V_i$. For instance, 
	\begin{center} 
	\end{center}
	or
	\begin{center}
		\includegraphics[width = 11cm]{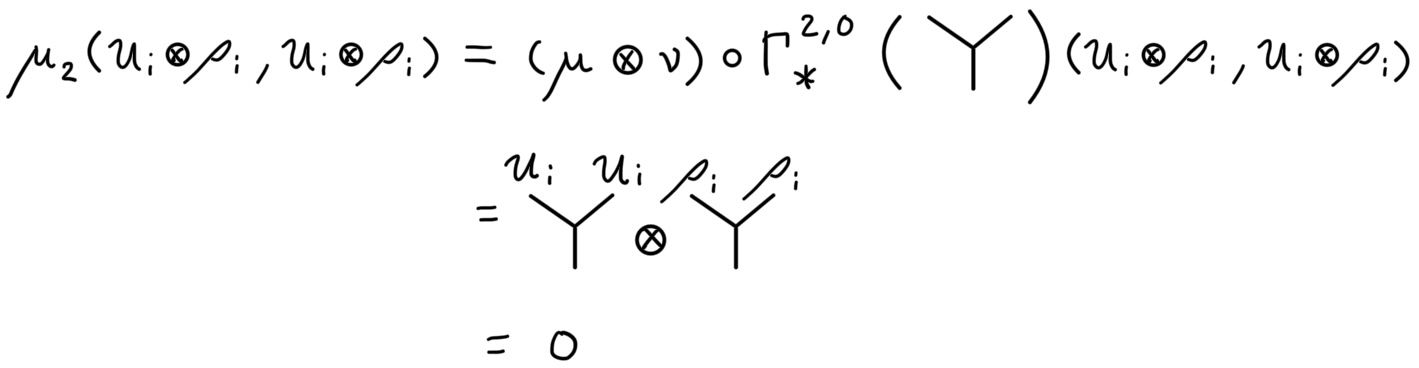}
	\end{center}
	That the operations on $\A \otimes \B$ determined in this way satisfy $\A_{\infty}$ relations follows directly from the fact that $\Gamma$ is a chain map, and from the original $\A_{\infty}$ relations on $\A$ and $\B$.

\section{The $\alpha$-bordered algebra $\A$}\label{alph}

		The object of this section is to define a weighted $\A_{\infty}$ algebra $\A$, for star diagrams of the form
		
	\begin{figure}[h!]
		\includegraphics[width = 5cm]{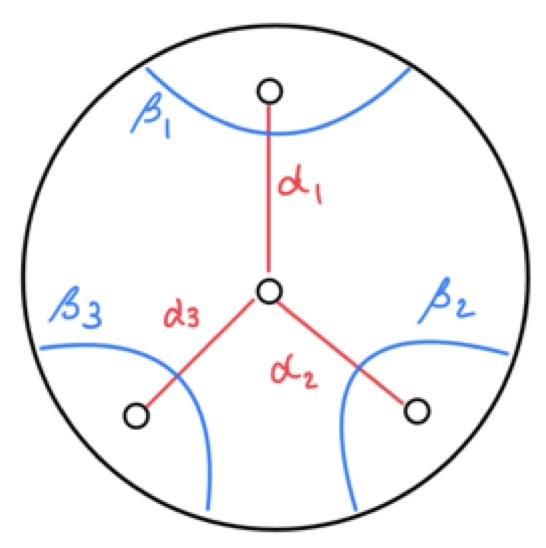}
		
		\caption{}
	\end{figure}
	That is, we are given, for $N \geq 3$, a diagram consisting of $N + 2$ boundary components (the black circles), $N$ red spokes (the $\alpha$-arcs) and $N$ blue arcs going out to the outer boundary component (the $\beta$-arcs). The algebra $\A$ is called ``$\alpha$-bordered'' because the generators, some of which are pictured below, in Figure~\ref{fig10'}, for $N = 3$, correspond to Reeb chords on the boundary circles, which are adjacent to $\alpha$-arcs:
	
	\begin{figure}[h!]
		\includegraphics[width = 5cm]{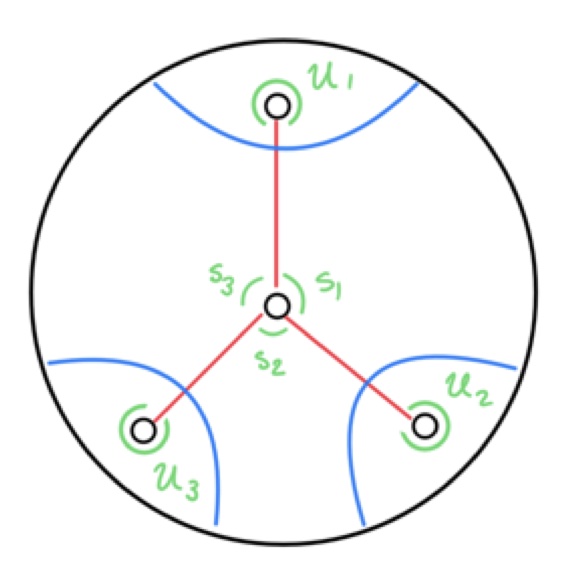}
		\caption{}\label{fig10'}
	\end{figure}
	There is a corresponding $\beta$-bordered algebra, $\B$, which will be defined in Section~\ref{bb1}.
	
	\subsection{Algebra elements and multiplication}
		Fix $N \geq 3$, and let $\F$ be a field (again, usually of characteristic $2$). The weighted $\A_{\infty}$ algebra $\A$ is, discounting $\A_{\infty}$ operations, an algebra (in the usual sense) over the ring $R = \F[V_0, \ldots, V_{N + 1}]$. Its generators are $\{U_1, \ldots, U_N\}$ and $\{s_{i (i + 1)}\}_{i = 1}^N$, where the indices for the $s_{ij}$ are counted $\emm N$. We write 
		\[
			s_{ij} = s_{i (i + 1)} s_{(i + 1)(i + 2)} \cdots s_{(j - 1) j},
		\]
		for each $0 < j - i < N$ (counted in $\N$, not $\emm N$). More generally, we stipulate that
		\begin{align*}
			U_i U_j &= \delta_{ij}, \\
			U_i s_{jk} &= s_{ij} U_k = 0 \quad \text{for all } i, j, k \\
			s_{ij} s_{k \ell} &= \begin{cases} s_{i \ell} & \text{ j = k} \\ 0 & \text{ otherwise}\end{cases} 
		\end{align*}
		and we define
		\[
			U_{N + 1} = \sum_{i = 1}^{N} s_{ii}
		\]
		so
		\begin{align*}
			U_i U_{N + 1} &= 0 \quad \text{ for each } 1 \leq i \leq N \\
			s_{ij} U_{N + 1} &= s_{ij} s_{jj} \\
			U_{N + 1} s_{ij} &= s_{ii} s_{ij} \\
		\end{align*}
		Throughout the following exposition, we say that two elements are \emph{multipliable} if they have non-zero product, and \emph{cannot be multiplied} or are \emph{non-multipliable} if they have product zero. 
		
		Graphically, each $U_i$ with $1 \leq i \leq N$ is the long chord around the $i$-th boundary circle (which is at the end of a spoke), while $s_{ij}$ with $0 < j - i < N$ is a short chord around the central boundary circle, for instance \vspace{-0.25cm}
		
		\begin{figure}[h!]
			\includegraphics[width = 5cm]{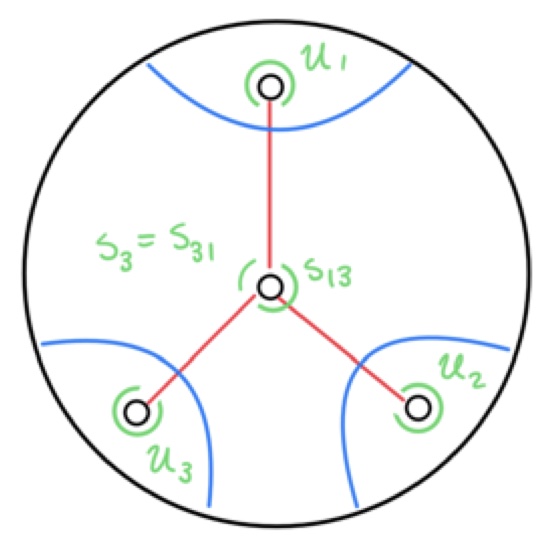}
			
			\caption{}
		\end{figure}
		This gives a geometric reason why the $s_{ij}$ and $U_i$ (except $U_{N + 1}$) cannot be multiplied. Multiplication corresponds to concatenation of adjacent chords, and even the nearest two nearest $U_i$ to a given $s_{jk}$ are each separated from it by an $\alpha$-arc. In the pseudo-holomorphic picture of things, to say that $U_i$ and $s_{jk}$ are multipliable is to say that this $\alpha$-arc can be compressed to a point, which is not allowed. This is also the underlying reason why it makes sense that $U_i$ and $U_j$ cannot be multiplied. 
		
		Each $s_{ii}$ denotes the long chord around the central boundary circle beginning and ending at the central spoke, as in the picture at right.
		
		\begin{wrapfigure}{r}{4cm}\vspace{-0.5cm}
			\includegraphics[width = 4cm]{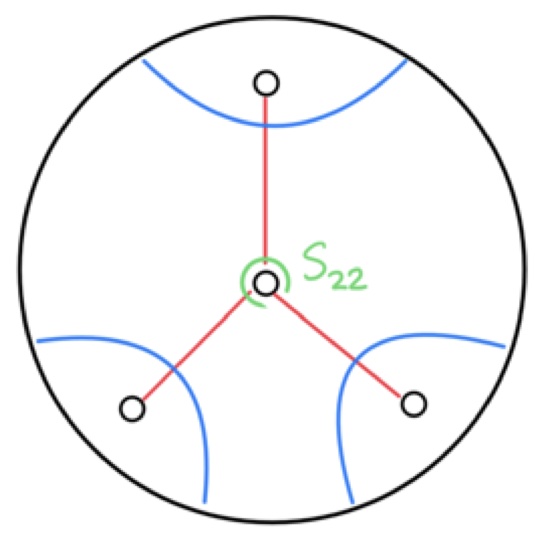}
			
			\caption{}\vspace{-1.7cm}
		\end{wrapfigure}
		This then makes $U_{N + 1}$ in some sense ``the'' long chord around the central boundary circle. 
		
		In keeping with the view of our algebra elements as Reeb chords, we define \emph{initial and final idempotents} for each element. We say that there are $N$ idempotents in total, corresponding to the $N$ spokes in the star picture above. We define each $U_i$ ($1 \leq i \leq N$) to both start and end in the $i$-th idempotent, and define $s_{ij}$ to start in the $i$-th idempotent and end in the $j$-th idempotent. It follows that $U_{N + 1}$ has summands in each idempotent state $1 \leq i \leq N$. In view of these definitions, it is clear from the above multiplication rules that two elements $a_1, a_2 \in \A$ are multipliable (i.e. $a_1 a_2 \neq 0$) \emph{only if} the final idempotent of $a_1$ is the same as the initial idempotent of $a_2$. This is, however, not always sufficient, as we see from the $U_i s_{ij} = 0$ (for $1 \leq i \leq N)$. 
	
	\subsection{Higher operations}
		
		We have now defined the multiplication on $\A$ (i.e., in terms of $\A_{\infty}$ operations, the unweighted $\mu_2$). The object of this subsection is to define the other (weighted and unweighted) operations which will make $\A$ into a bona-fide $\A_{\infty}$ algebra. Recall that operations correspond to certain trees whose inputs are labelled with algebra elements, and which may or may not contain weight. The $\A_{\infty}$ relations are in one-to-one correspondence with the set of labelled trees; each tree $T$ gives rise to an $\A_{\infty}$ relation, as the number of ways to add a single internal edge to $T$ by replacing one internal vertex with a new edge. We now need to define the other operations on $\A$ so that for each $T$, the sum over all such ways to add an edge vanishes. 
		
		Operations (both weighted and unweighted) will correspond to so-called \emph{allowable labelled graphs}, analogous to the construction of the weighted bordered torus algebra in~\cite{Torus}. First, a \emph{labelled graph} is a planar graph in the disk, with each sector around each vertex labelled with an element of $\A$. 
		
		For fixed $N$, the \emph{basic allowable (unweighted) planar graphs} are a set of $2N$-valent graphs, each with one internal vertex, and one marked  vertex on the boundary called the \emph{root}. Each such graph corresponds to a corolla with $2N$ inputs and one output. For example, in the case $N = 3$,
		
		\begin{figure}[H] 
			\includegraphics[width = 4cm]{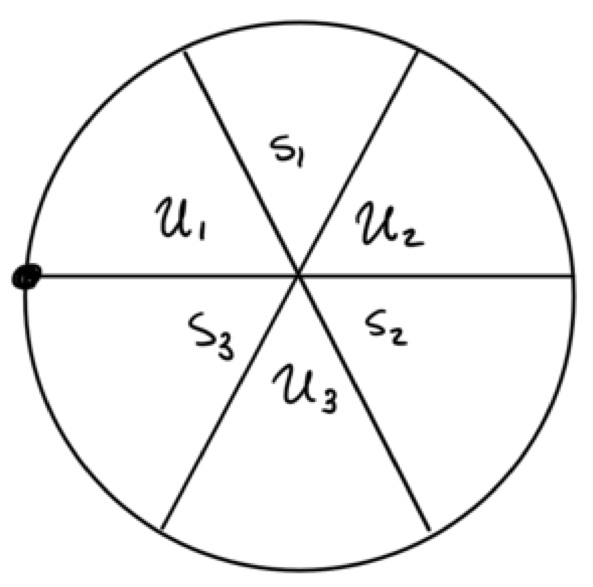}\hspace{0.5cm}\includegraphics[width = 6cm]{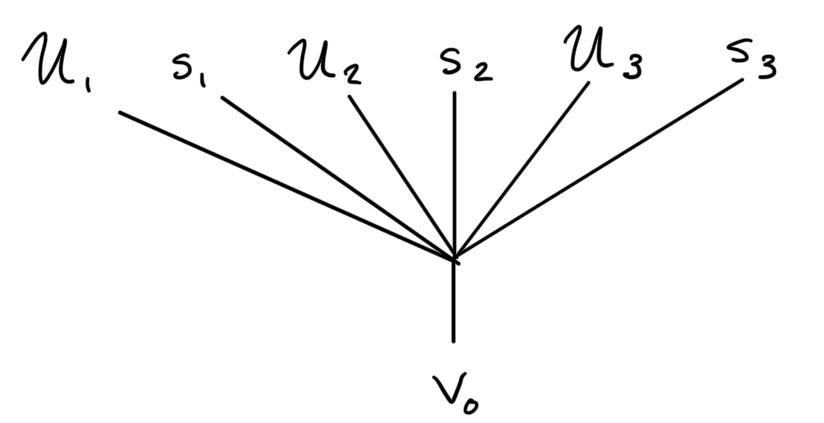}
			
			\caption{}\label{basgraph}
		\end{figure}
		A \emph{centered, unweighted operation} is defined to be an \emph{allowable concatenation} of these basic allowable labelled graphs.
		
		\begin{wrapfigure}{r}{5cm}
			\includegraphics[width = 5cm]{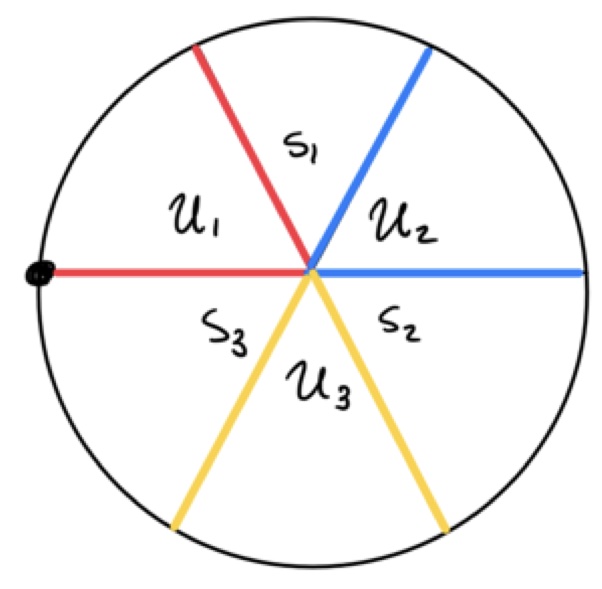}
			
			\caption{}\label{basgraph2}
		\end{wrapfigure}
		
		To define centered unweighted operations, we need to define \emph{allowable concatenations}. To do this, note first that each algebra element used to label the basic graph in Figure~\ref{basgraph} has an initial idempotent and a final idempotent. Reading clockwise around the central vertex, the initial idempotent of a one algebra element is precisely the final idempotent of the previous one. 
		
		Assign to each idempotent a color. The labelling of the basic graph in Figure~\ref{basgraph} induces a coloring of its branches corresponding to the idempotent in which each branch lies, as in Figure~\ref{basgraph2}. Note that there are $N$ colors, and each sector labelled with $U_i$ bordered by edges colored with the $i$-th color, and each sector labelled with an $s_i$ is bordered by edges colored with the $i$-th and $(i + 1)$-st colors, counted clockwise around the vertex.
		
		 In general, a labelled graph is \emph{allowable in the unweighted sense} if it satisfies the following conditions:
		 \begin{itemize}
		 	\item it is a coherently colored $2N$-valent graph with no cycles
			\item For each internal vertex $v$, if such that if we punch out a sufficiently small disk around each $v$, the labelled graph in this local picture looks like the picture above (sans the root vertex), up to cyclic permutation.
		\end{itemize}
		An \emph{allowable concatenation} of basic graphs is any labelled graph which is allowable in the unweighted sense.

		In practice, we usually suppress the coloring of these graphs, and just label the regions as in 
		
		\begin{figure}[h!]
			\includegraphics[width = 10cm]{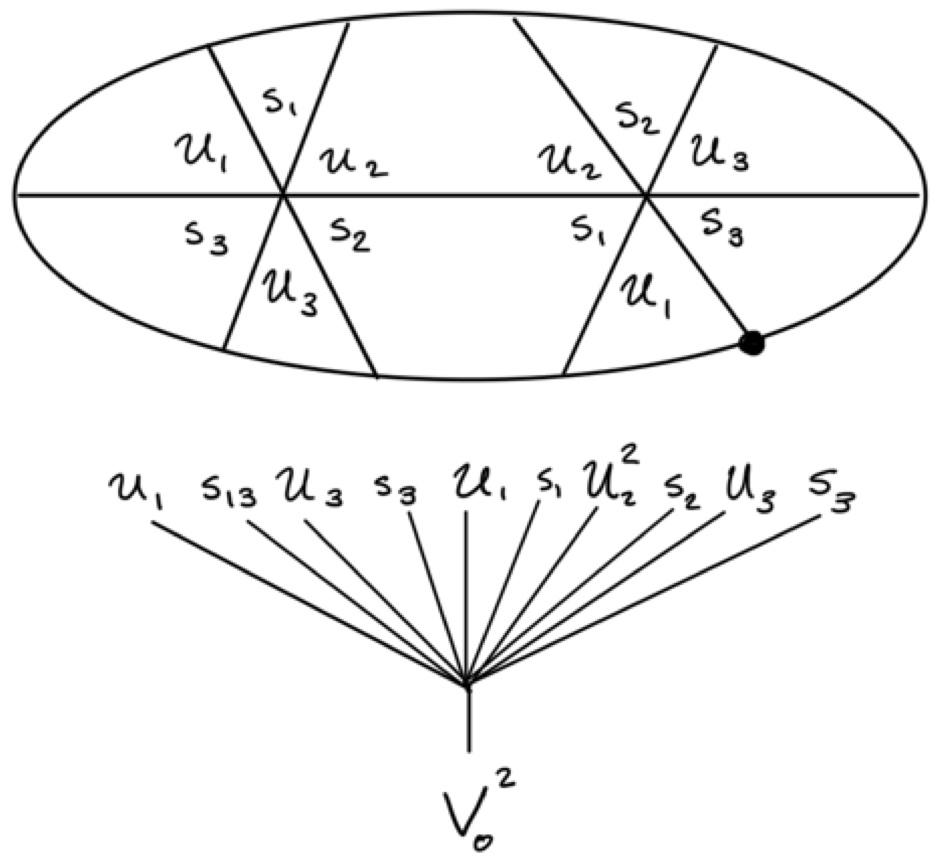}
		\end{figure}\vspace{-5pt}
		
		Extended unweighted operations are defined by adding a sequence of labelled 2-valent vertices to the edge adjacent the root, with labels either all on the left of the root (looking into the disk), in which case the operation is called left-extended, or on the right of the root, in which case the operation is called right extended. For instance, in the case $N = 3$, we could have
		
		\begin{figure}[H]
			\includegraphics[width = 12cm]{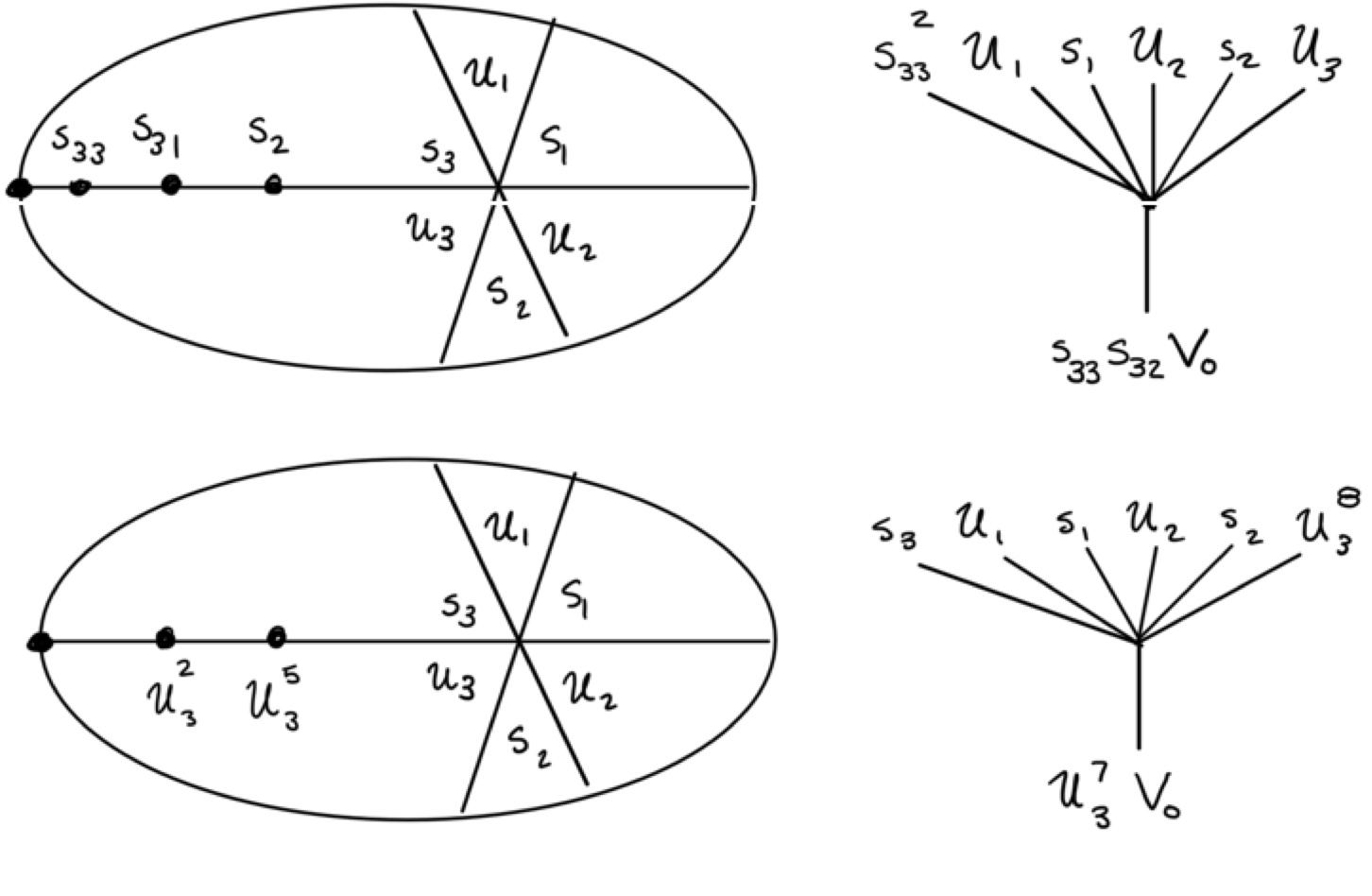}
		\end{figure}
		Note that we are not allowed to insert 2-valent vertices on any other edges. Weighted operations likewise correspond to colored planar graphs, this time allowing certain cycles, namely, for each $1 \leq i \leq N + 1$,
		 \begin{figure}[H]
			\includegraphics[width = 5cm]{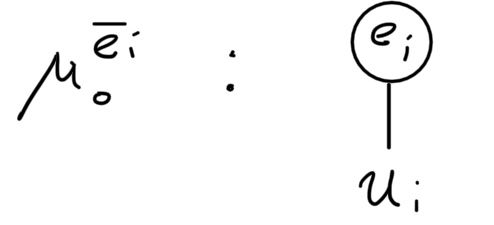}
		\end{figure}
		
		This induces higher multiplications, for instance, in the case $N = 3$,
		
		\begin{figure}[H]
			\includegraphics[width = 9cm]{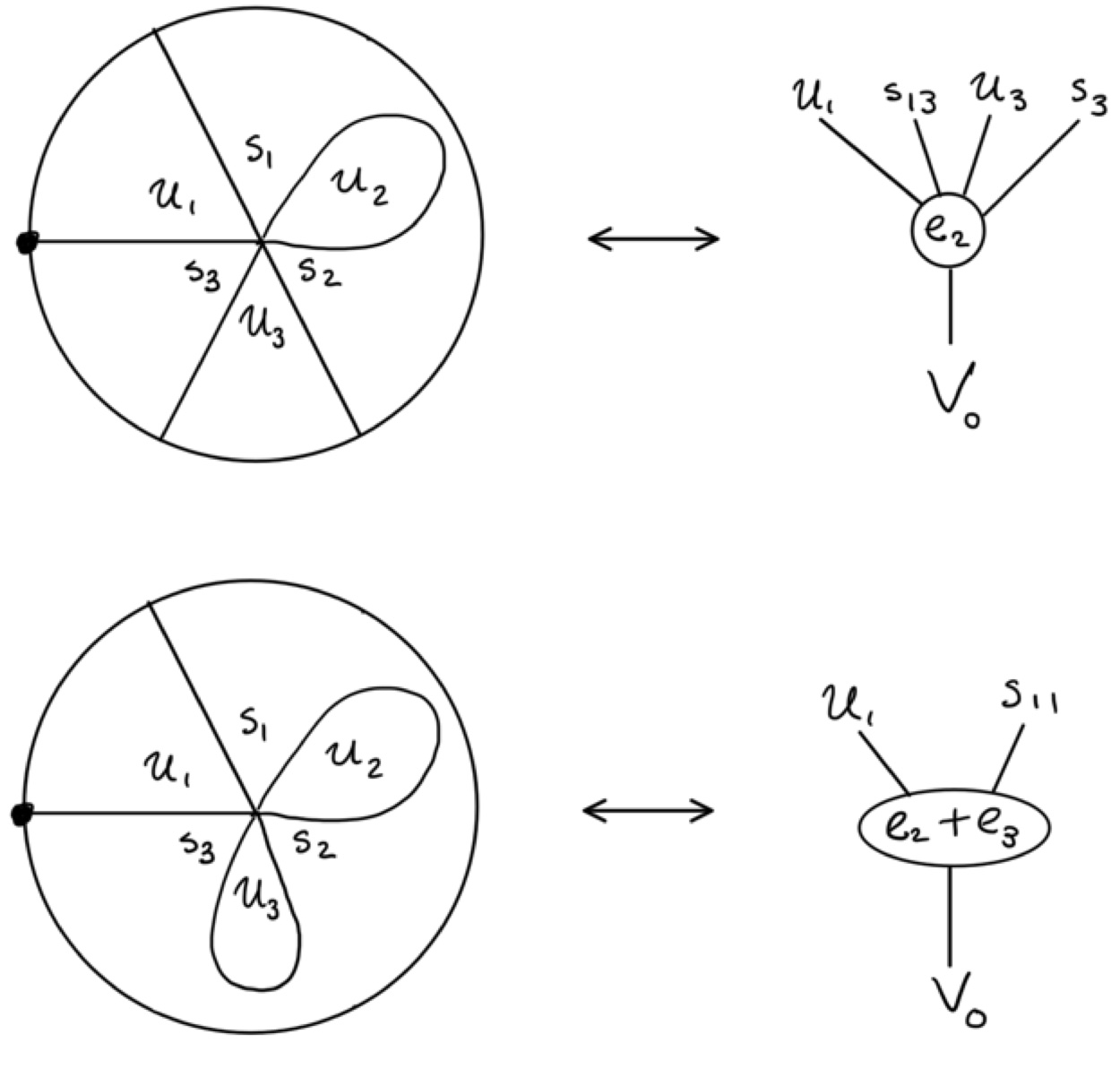}
		\end{figure}
		We call these cycles (for $1 \leq i \leq N$) \emph{petal cycles} which induce \emph{petal weight} in the corresponding trees. In general, these basic petal-weighted operations are of the form
		\[
			\mu_{2N - 2k}^{\w}(a_1, \ldots, a_{2N - 2k}) = V_0
		\]
		where $\w = \sum_{j = 1}^k \e_{i_j}$ for distinct $i_j$'s, and for suitable $a_i$ as detailed above. 
		
		We also have weighted higher multiplications corresponding to \emph{internal cycles}, for instance, in case $N = 3$, 
		
		\begin{figure}[H]
			\includegraphics[width = 8cm]{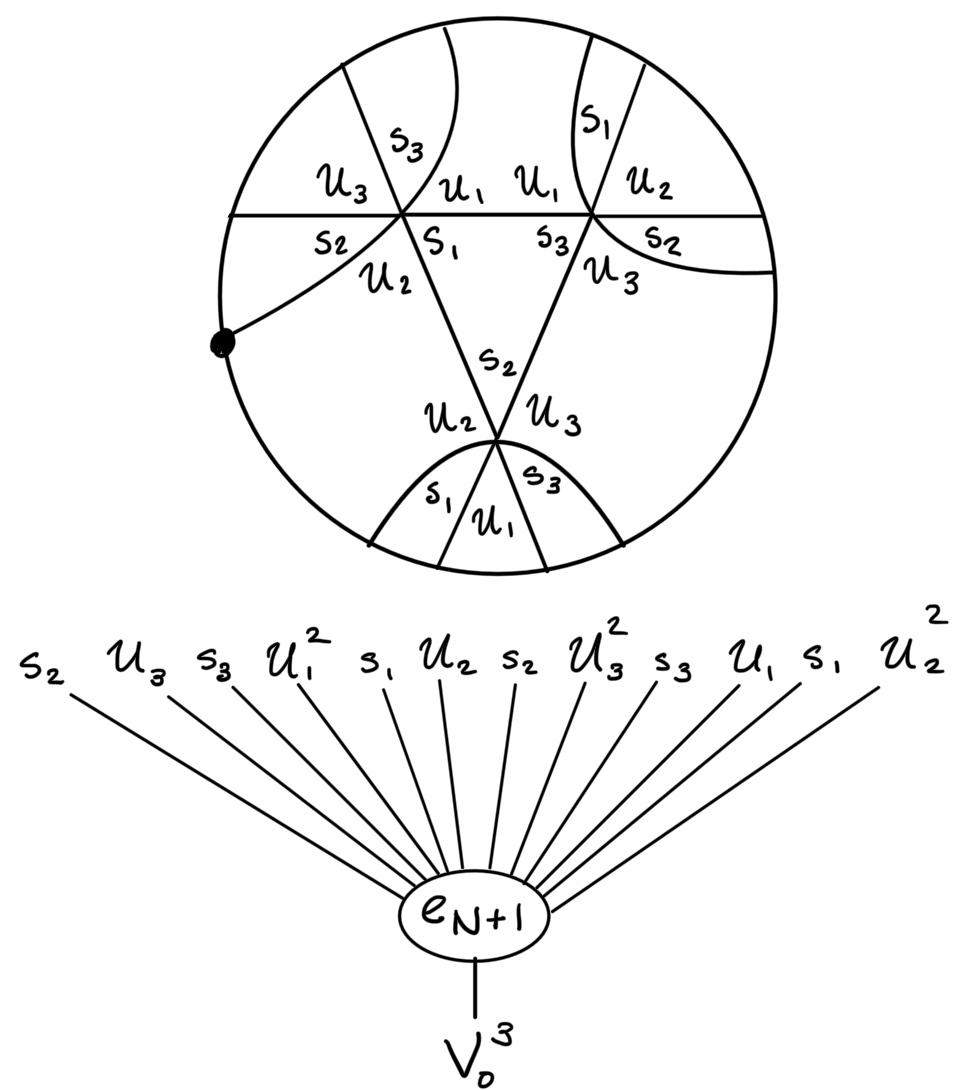}
		\end{figure}
		\hspace{-0.45cm}and the other operations which correspond to moving the basepoint around to any of the other outer vertices. In general, these basic, internally weighted, operations are of the form 
		\[
			\mu_{N (2N - 2)}^{\e_{N + 1}}(a_1, \ldots, a_{N (2N - 2)}) = V_0^N
		\]
		for suitable $a_i$, as detailed above. 
	
	\subsection{Maslov and Alexander gradings for $\A$}\label{aaMas}
	
	Next, since operations are defined in terms of trees, we need to specify order to say exactly which trees describe allowable graphs as defined above. To do this, we need to introduce gradings. The Maslov grading $m$ is a map from the set of labelled trees to $\Z$. It is required to satisfy the following property:
	\begin{equation}\label{aa1}
		m(\mu_n^{\w}(a_1, \ldots, a_n)) = \sum_{i = 1}^n m(a_i) + m(\w) + n - 2
	\end{equation}
	where $\w$ is the weight and $a_1, \ldots, a_n \in \A$. We set
	\begin{align*}
		m(U_i) = m(s_{ij}) &= 0 \text{ for each $i, j$} \\
		m(\e_i) &= 2 \text{ for each } 1 \leq i \leq N + 1 \\
		m(\e_0) &= - (2N - 2) \\
		m(V_0) &= 2N - 2
	\end{align*}
	Notice that in each case, this fits with~\eqref{aa1}, and that the $m(\e_{N + 1}) = 3$ is independent of choice of $N$. 
	
	The Alexander grading is a vector or homological grading. Let $X$ be $S^1$ minus $2N$ points, one corresponding to each endpoints of a $\beta$-arc. The 0-th homology of $X$ is just $\Z^{2N}$. Write the generators $\overline{k}$ for $1 \leq k \leq 2N$. The Alexander grading is a map $A: \A \to H_0(X; \Z) \cong \Z^{2N}$, such that 
	\begin{align*}
		A(U_i) &= \overline{2i - 1} \text{ for each } 1 \leq i \leq N \\
		A(s_{ij}) &= \sum_{k = i}^{j - 1} \overline{2k} \\
		A(\e_i) &= \overline{2i - 1} \text{ for each } 1 \leq i \leq N \\
		A(\e_{N + 1}) &= \sum_{k = 1}^{N} \overline{2k} \\
		A(V_0) &= \sum_{k = 1}^{2N} \overline{k}
	\end{align*}
	This notion is similar to the notion of initial and final idempotents for an algebra element, but distinct. For instance, $U_i$ starts and ends in the $i$-th idempotent and has Alexander grading in the corresponding (i.e. $(2i - 1)$st component). However, idempotents keep track only of which slot the algebra element starts and ends in, whereas the Alexander grading captures all the slots the algebra element covers. In the case of $s_{ij}$, for instance, the initial and final idempotents are $i$ and $j$, respectively, and while $A(s_{ij})$ does have components in the corresponding slots ($2i$ and $2j$, respectively) it also has components in all corresponding even slots, providing a fuller picture.
	
	The Alexander grading of a sequence can also be defined, as:
	\[
		A(a_1, \ldots, a_n) = \sum_{i = 1}^n a_i.
	\]
	The Alexander grading of a weighted sequence is
	\[
		A(\w, a_1, \ldots, a_n) = A(\w) + A(a_1, \ldots, a_n).
	\]
	We require that operations preserve the Alexander grading.
	
	Graphically, all this looks like
	\begin{figure}[H]
		\includegraphics[width = 7cm]{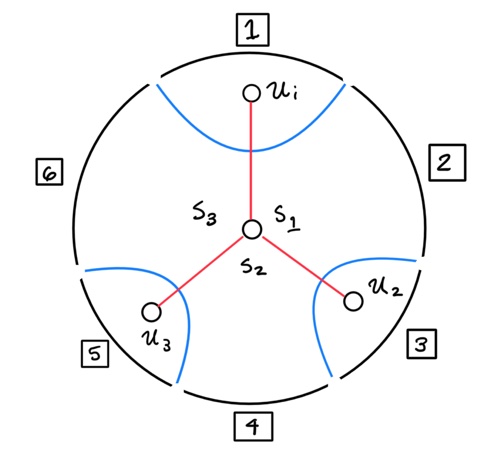}
		
		\caption{}
	\end{figure}\vspace{-10pt}
	\hspace{-0.4cm}where the boxed numbers are the Alexander grading corresponding to each component of the cut-out $S^1$ along the boundary, and we can see why each algebra element fits into the Alexander grading to which it is assigned. As usual, we are just looking at $N = 3$ as an example, but this also holds for other $N$. 
	
	\subsection{$\A_{\infty}$-relations}
	
	We are nearly ready to verify the $\A_{\infty}$ relations for $\A$. The next step is to give a condition for an unweighted tree to correspond to an allowable graph (i.e. to an unweighted operation). Before we can state the condition, we need some further definitions. A sequence $(a_1, \ldots, a_n)$ of elements of $\A$ is \emph{left extended} if the initial element $a_1$ is not basic, but can be written as a product of more than one of the basic elements ($U_i$'s and $s_{i}$'s)\footnote{While we usually view any $s_{ij}$ as basic, even if $j > i + 1$ (in which case $s_{ij}$ can be decomposed as a product of more than one $s_i$), we only want to look at single $s_i$'s for this case.}; it is \emph{right extended} if $a_n$ can be written as a product of more than one of the basic elements of $\A$. The sequence is \emph{centered} if it is neither left nor right extended.
	
	The condition is as follows:
	\begin{theorem}\label{aa2}
		\emph{(Unweighted higher multiplications)} Let $T$ be an unweighted tree whose inputs are labelled with elements of $\A$, say $a_1, \ldots, a_n$, counted from left to right, with $n > 2$. Then $T$ represents an operation if and only if the following conditions hold:
		\begin{enumerate}[label = (\roman*)]
			\item \emph{(The idempotents match up)} The initial idempotent of $a_i$ is the final idempotent of $a_{i - 1}$ for each $1 < i \leq n$; 
			
			\item \emph{(The Maslov grading works out)} $n = j (2N - 2) + 2$ where $j \in \N$;
			
			\item \emph{(Extensions)} The sequence $(a_1, \ldots, a_n)$ is either centered, or left or right extended, but not both at once; 
			
			\item \emph{(The Alexander grading works out)} The condition varies depending on whether the sequence $(a_1, \ldots, a_n)$ is centered, left extended, or right extended. 
			\begin{itemize}
				\item Centered: $A(a_1, \ldots, a_n) = j \cdot \sum_{k = 1}^{2N} \overline{k}$;
				
				\item Left extended: Write $a_1 = a_1' \cdot a_1''$, where $a_1'' = U_i$ or $s_i$ for some $i$ (i.e. it is a truly basic element). Then we require that $A(a_1, \ldots, a_n) = A(a_1') + j \cdot \sum_{k = 1}^{2N} \overline{k}$;
				
				\item Right extended: Write $a_n = a_n'' \cdot a_n'$ where $a_n'' = U_i$ or $s_i$ for some $i$. Then we require that $A(a_1, \ldots, a_n) = A(a_n') + j \cdot \sum_{k = 1}^{2N} \overline{k}$;
			\end{itemize} 
			
			\item \emph{(Correct output)} The output is $V_0^j$ if the $(a_1, \ldots, a_n)$ is centered $a_1' V_0^j$ if it is left extended, and $V_0^j a_n'$ if it is right extended;
		\end{enumerate}
	\end{theorem}
	
	\begin{proof}
		Any suitable graph yields a labelled tree by reading off the labels on the sectors of the graph counting clockwise from the root; for instance
		\begin{center}
			\includegraphics[width = 10cm]{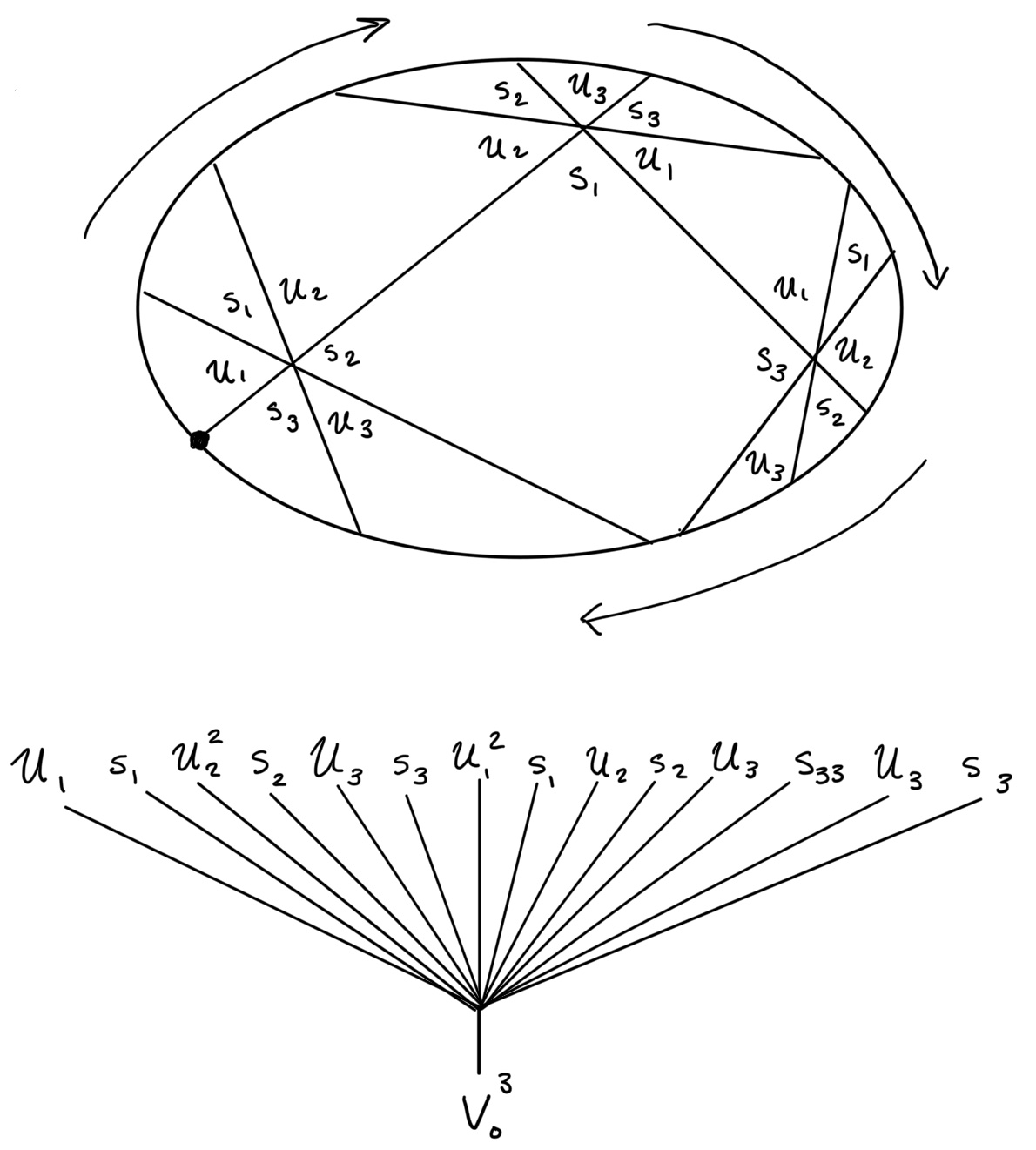}
		\end{center}
		which in this case, is an unweighted $\mu_{14}$. Because of how we constructed our graphs, the resulting tree clearly satisfies (i)-(v).
		
		\begin{figure}[h]
			\includegraphics[width = 5cm]{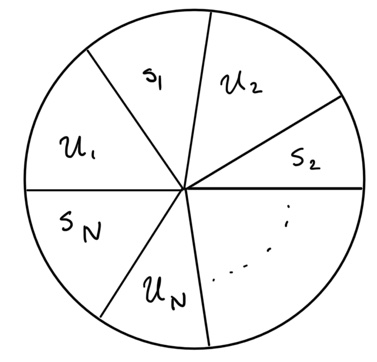}
			
			\caption{}\label{aa23'}
		\end{figure}
		
		For the other direction, suppose a tree $T$ satisfies conditions (i)-(v), and is centered, with the inputs labelled by $(a_1, \ldots, a_n)$. Then we can construct a simply connected, $2N$-valent graph for $T$ inductively, in the following way. First of all, if $j = 1$, so $n = 2N$, then $T$ determines some choice of root vertex in the graph of Figure~\ref{aa23'}, that is, it uniquely determines a suitable graph. Now suppose that for all $j' < j$, if $T$ is a centered tree satisfying (i)-(v), we can define a unique, suitable corresponding graph, and look at a $T$ with $j(2N - 2) + 2$ inputs. Because $T$ has $j(2N - 2) + 2$ inputs and covers exactly $j\cdot 2N$ slots in the Alexander grading, there must be some subsequence of $2N - 2$ consecutive basic elements among the $(a_1, \ldots, a_n)$, say $a_{i + 1}, \ldots, a_{i + 2N - 1}$. Look at $a_i$ and $a_{i + 2N}$, and write
		\begin{align*}
			a_i &= a_i' a_i'' \\
			a_{i + 2N} &= a_{i + 2N}'' a_{i + 2N}',
		\end{align*}
		where $a_i''$ and $a_{i + 2N}''$ are basic elements of $\A$. It is now clear that the initial idempotent of $a_{i}''$ is the same as the final idempotent of $a_{i + 2N}''$; so now we have a sequence $a_i'', a_{i + 1}, \ldots, a_{i + 2N - 1}, a_{i + 2N}''$ of basic elements, which by the rules on idempotents, must be some cyclic permutation of $U_1, s_1, \ldots, U_N, s_N$. Hence, it determines some choice of root vertex on the single-vertex graph in Figure~\ref{aa23'}. Remove the sequence $a_i'', a_{i + 1}, \ldots, a_{i + 2N - 1}, a_{i + 2N}''$ from $a_1, \ldots, a_n$, to get $a_1, \ldots a_{i - 1}, a_i', a_{i + 2N}', a_{i + 2N + 1}, \ldots, a_n$. Look at the tree $T'$ with these elements as inputs, and with output $V_0^{j - 1}$. Then it is clear that this $T'$ still satisfies (i)-(v), but the output power of $V_0$ is less than $j$; hence the inductive hypothesis holds and we can get a unique, suitable rooted graph corresponding to $T'$. Now attach the root of the graph corresponding to $a_i'', a_{i + 1}, \ldots, a_{i + 2N - 1}, a_{i + 2N}''$ to the edge between $a_i', a_{i + 2N}'$. This is clearly a graph that corresponds to $T$. Looking back at the way we got from graphs to trees in the first part, this is clearly unique.
		
		Now, consider a suitable left-extended tree $T$, and write $a_1 = a_1' a_1''$, where $a_1''$ is a basic element. The output of $T$ is $a_1' V_0^j$ for $j$ as in the statement. The tree $T'$ with inputs $(a_1'', a_2, \ldots, a_n)$ and output $V_0^j$ is a centered tree satisfying (i)-(v) and affords a unique graph. Now decompose $a_1'$ into a product of basic elements, say $a_1 = p_1 \cdots p_k$, and mark the initial of the graph for $T'$ with points labelled $p_1, \ldots, p_k$; for instance:
		\begin{center}
			\includegraphics[width = 10cm]{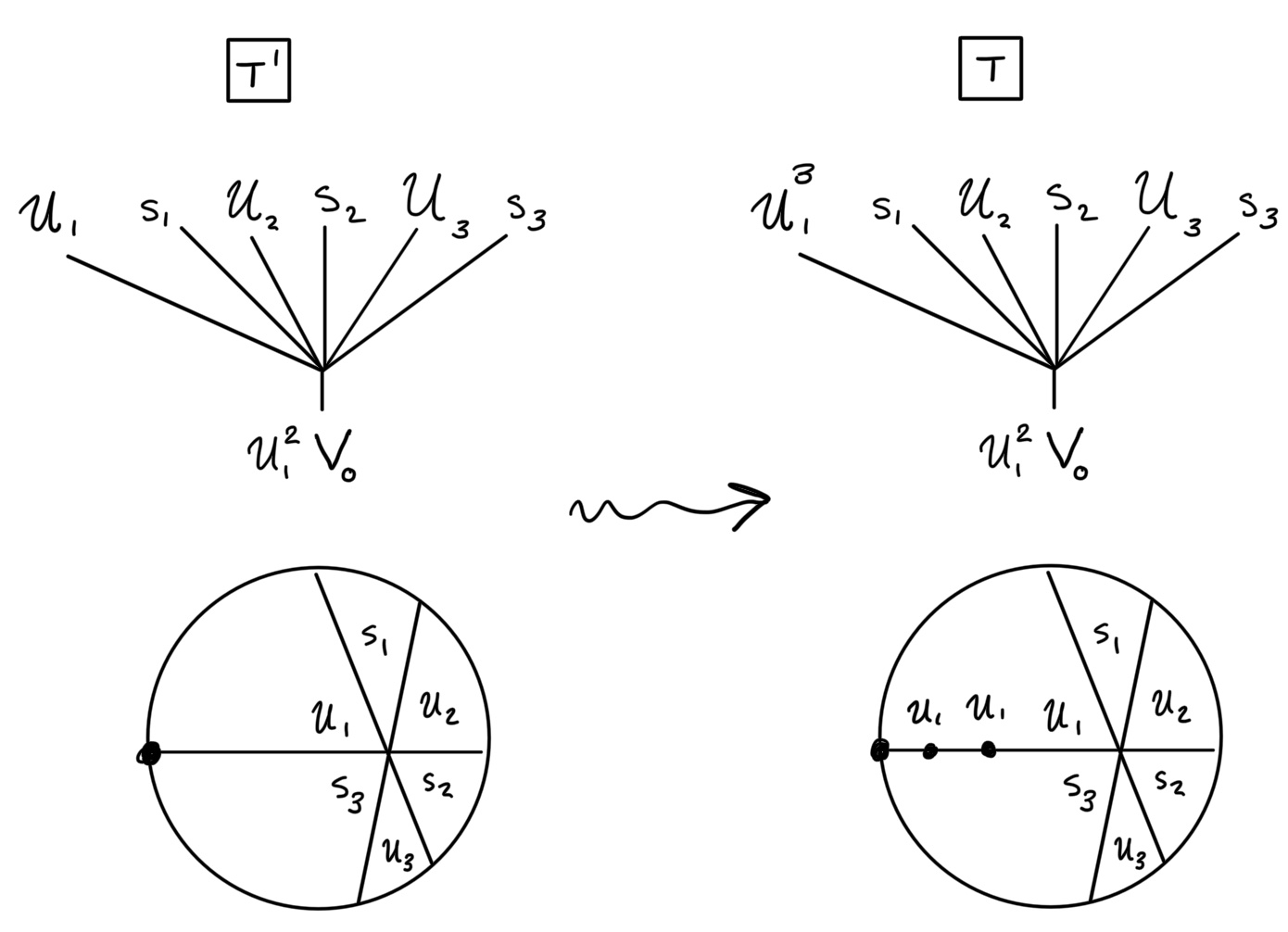}
		\end{center}
		This gives the graph for $T$, which is clearly allowable. The situation for a right-extended tree is analogous.
	\end{proof}
	
	Before we can verify the $\A_{\infty}$-relations for $\A$, we need to define the notion of \emph{augmented graphs}, which we will use in the proof. There are two types. The first is an allowable graph in which we choose one internal edge (that is, an edge not intersecting the boundary) and draw a dotted line from this edge to the boundary. Note that the direction of drawing the edge matters, so
	\begin{center}
		\includegraphics[width = 12cm]{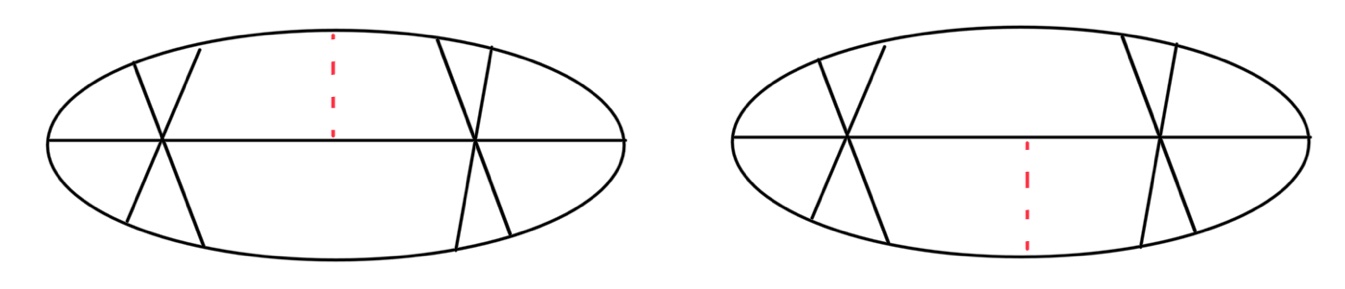}
	\end{center} 
	are distinct. We are also allowed to draw this dotted line on the on an edge adjacent to an ``extended'' vertex, for instance
	\begin{center}
		\includegraphics[width = 10cm]{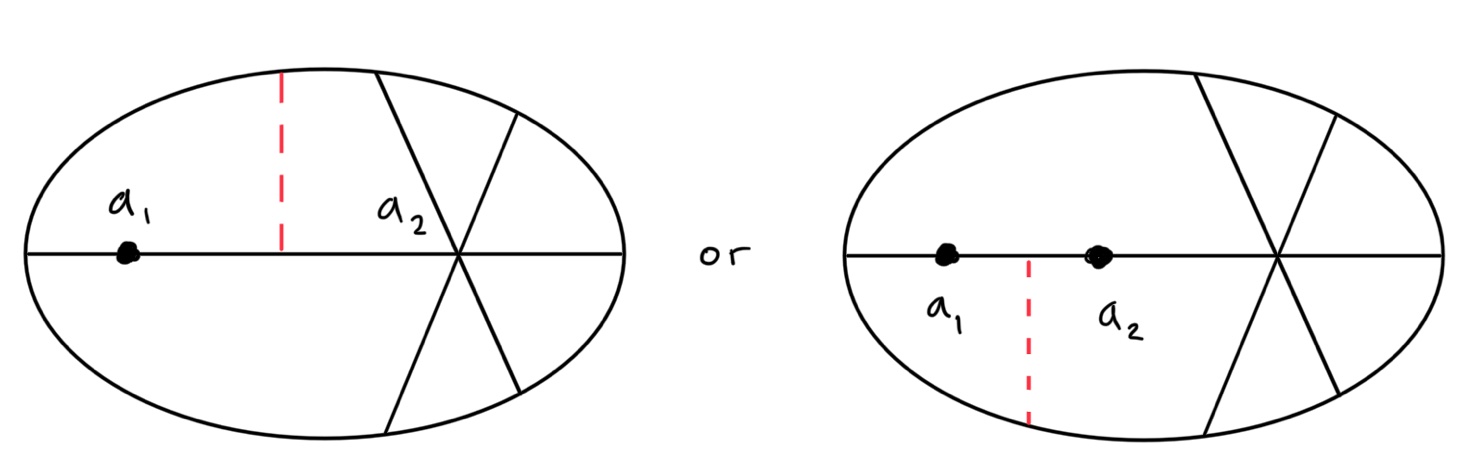}
	\end{center}
	Note that we can only draw this dotted edge on the side of the extension (so in the first case, only on the top, and in the second, only on the bottom). We will prove that that these augmented graphs are in one to one correspondence with pulls in Theorem~\ref{aa7}. 
	
	The other kind of augmented graph are two vertex, two component graphs of the form 
	\begin{center}
		\includegraphics[width = 9cm]{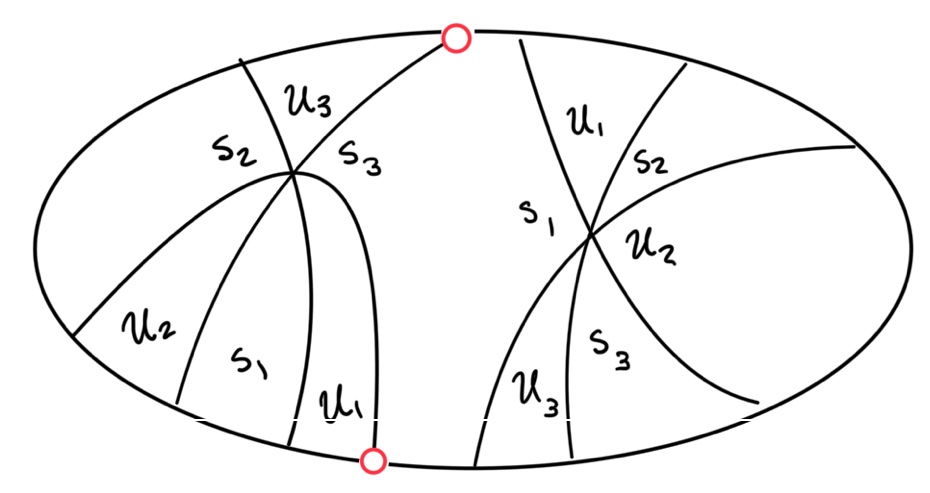}
	\end{center}
	Here, we are again using $N = 3$ as an example. All augmented graphs of the second kind (for $N = 3$) can be obtained from this one by permuting the labels of the first component one step counterclockwise, and the labels of the second component one step clockwise. 
	
	The hollow roots in the graph above indicate the two possible ways to resolve this into allowable compositions; this graph above corresponds to the tree
	\begin{center}
		\includegraphics[width = 9cm]{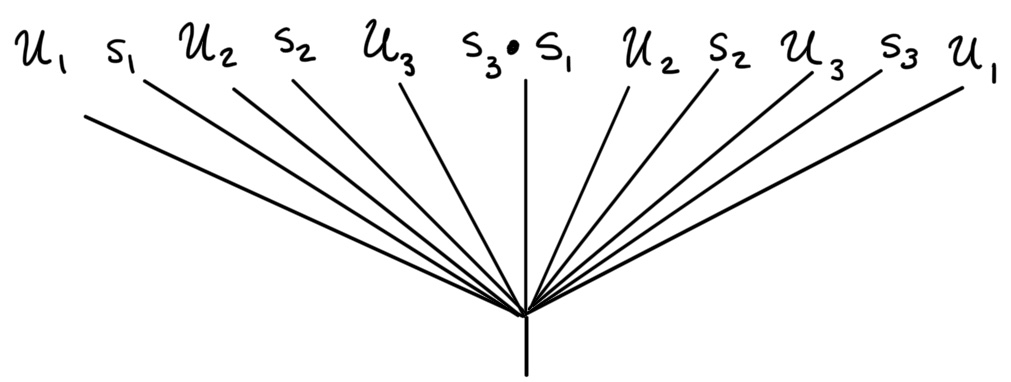}
	\end{center}
	and the two ways to resolve it are
	\begin{center}
		\includegraphics[width = 7cm]{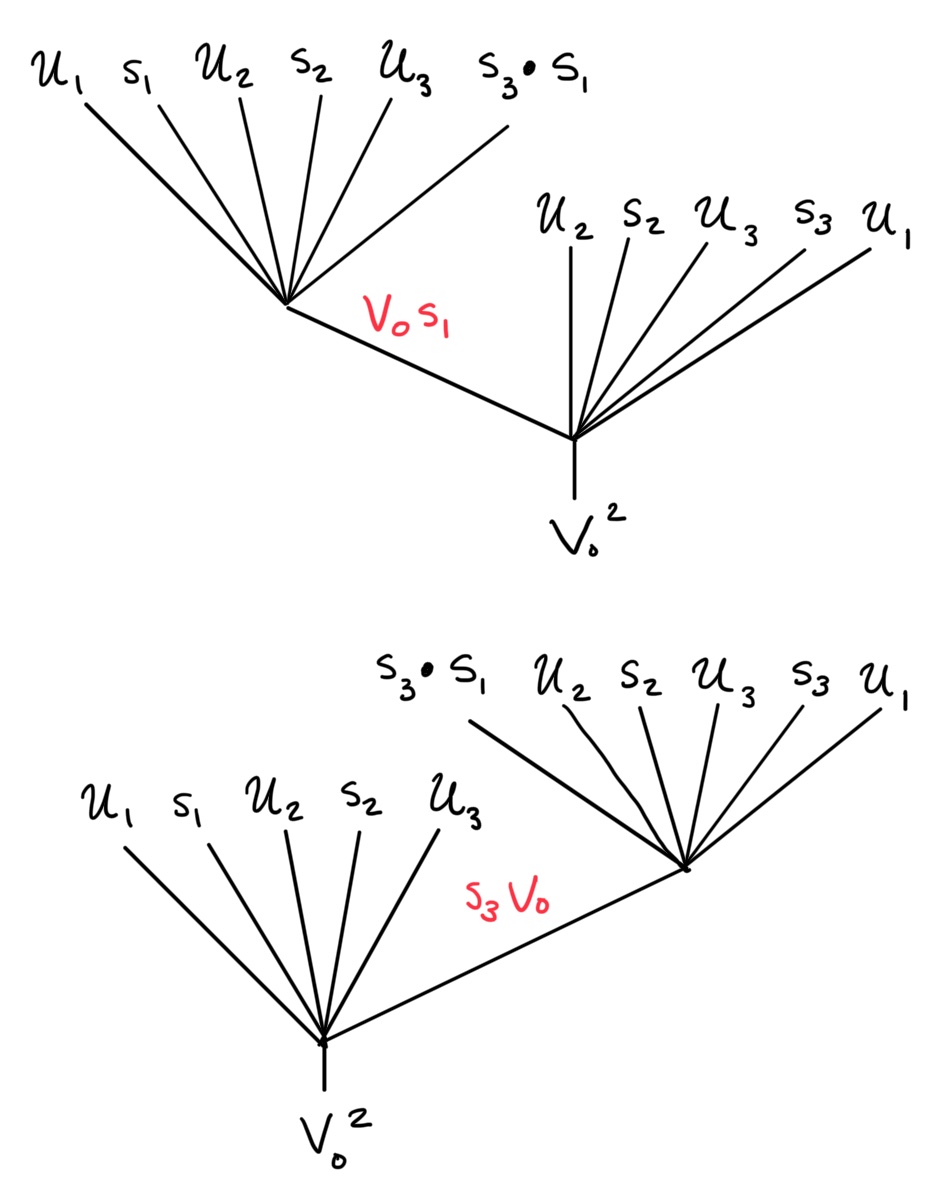}
	\end{center}
	Thus, the second kind of augmented graph are in one to one correspondence with this particular kind of split. A tree admits this particular kind of split if and only if satisfies (i), (iii), and (iv) from Theorem~\ref{aa2}, has $2 (2N) - 1$ inputs -- that is, one extra from the required number from (ii) of Theorem~\ref{aa2} -- and has no multipliable pairs. Essentially, the multipliable pair (which leads to the extra input) is the pair of elements on the two ends. But while those correspond to adjacent elements in the graph, they are not adjacent in the tree.
	
	 We have actually just proven:
	\begin{lemma}\label{aa19}
		Let $T$ be a tree satisfying conditions (i), (iii) and (iv) from Theorem~\ref{aa2}, and with $n = 2 \cdot (2N - 2) + 3$, with no multipliable pairs. Then $T$ corresponds to a unique augmented graph of the first kind, and hence, to a cancelling pair of splits as pictured above.
	\end{lemma}
	
	We are now ready to prove the $\A_{\infty}$-relations for unweighted trees labelled with elements of $\A$. 
	
	\begin{theorem}\label{aa3}
		\emph{(Unweighted $\A_{\infty}$-relations)} The $\A_{\infty}$ relations for unweighted trees labelled with elements of $\A$ are satisfied -- that is, the sum over all ways to add an edge to a given tree $T$ is always zero.
	\end{theorem}
		
	\begin{proof}
	
		First of all, the $\A_{\infty}$ relations for the usual $\mu_2$ -- a product of multipliable elements -- are trivial, so we can restrict our attention to $\A_{\infty}$ relations corresponding to trees with at least three inputs.
		
		 To prove our claim for this case, we are going to show that non-vanishing $\A_{\infty}$ relations in $\A$ (i.e. trees $T$ for which at least one of the trees corresponding to adding an edge to $T$ is a non-zero operation) are in one-to-one correspondence with augmented graphs as described above. 
		
		Note first that any augmented graph of the first kind yields exactly two suitable graphs (i.e. operations, or compositions of operations) in the following ways: 
		\begin{center}
			\includegraphics[width = 12cm]{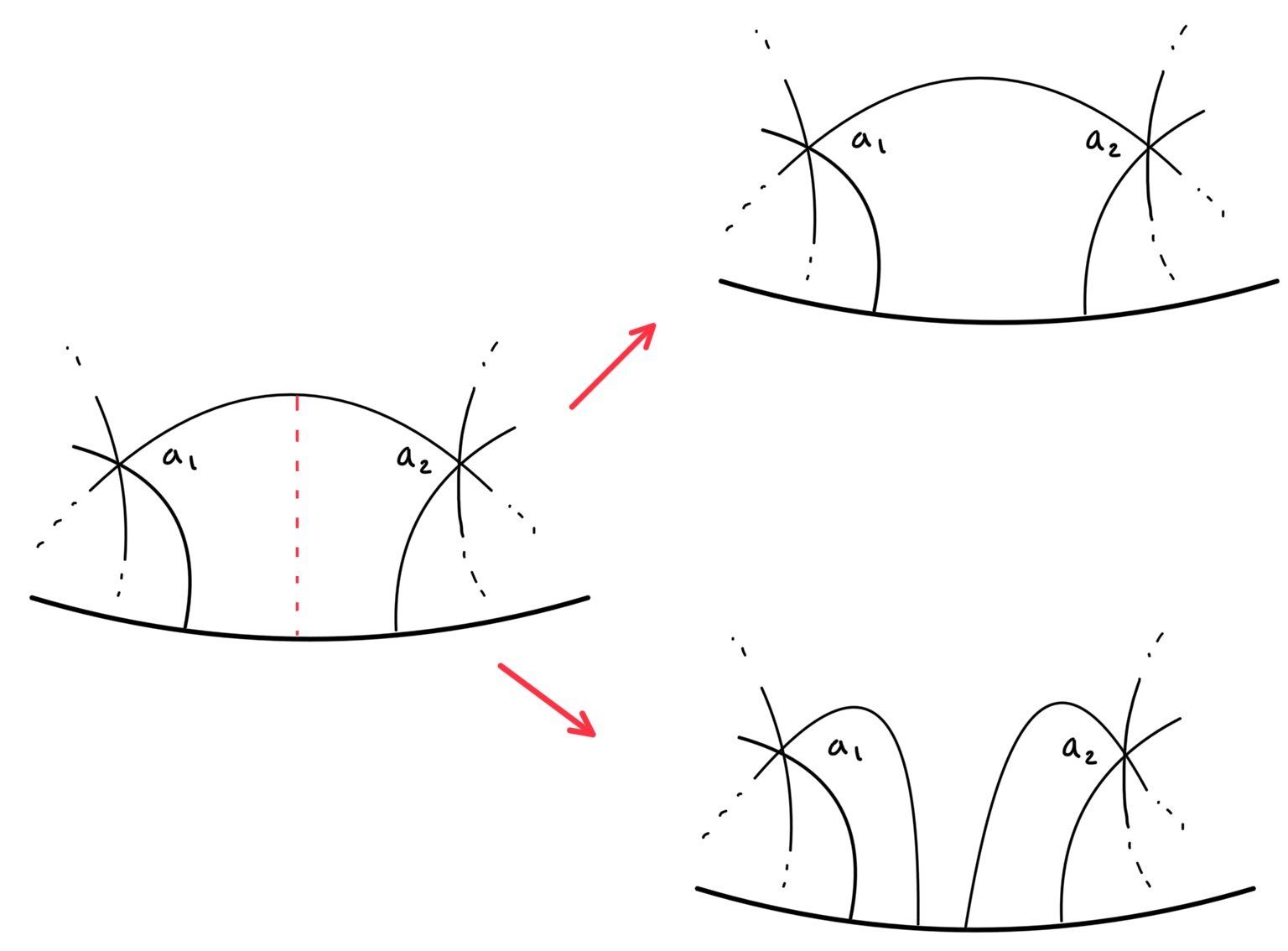}
		\end{center}
		which corresponds to the $\A_{\infty}$ relation
		\begin{center}
			\includegraphics[width = 8cm]{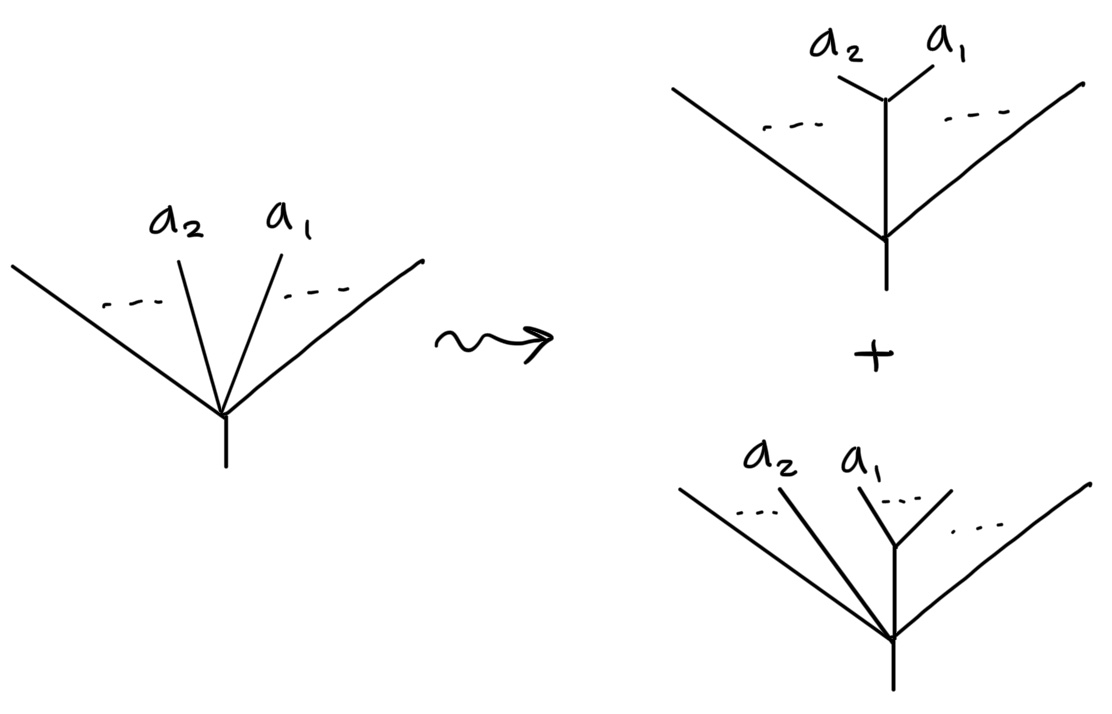}
		\end{center}
		The $a_1, a_2$, are assumed to be multipliable, and as usual, we were using $N = 3$ as a shorthand for general $N$, specifically in depicting the generally $2N$-valent vertices above. That the second graph above has two disjoint components (and hence, corresponds to a composition i.e. a split) follows from the fact that our graphs do not contain cycles, and therefore, cutting an edge gives two components. Also, we are assuming without loss of generality that the root is on the component adjacent to $a_2$. 
		
		That there are non-trivial inputs on either side of the internal operation containing $a_1$ follows from the fact that the root is not on the edge containing $a_2$. This implies, in particular, that the only way to obtain a composition of either of the forms
		\begin{center}
			\includegraphics[width = 7cm]{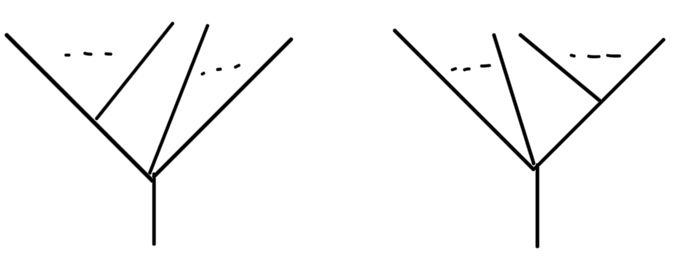}
		\end{center}
		as the result of an $\A_{\infty}$-relation is in the particular kind of split dealt with in Lemma~\ref{aa19}, that is, from the second kind of augmented graph.
		
		Alternately, if the dotted edge is drawn between two vertices on the initial edge (at least one of which must be extended), we get:
		\begin{center}
			\includegraphics[width = 12cm]{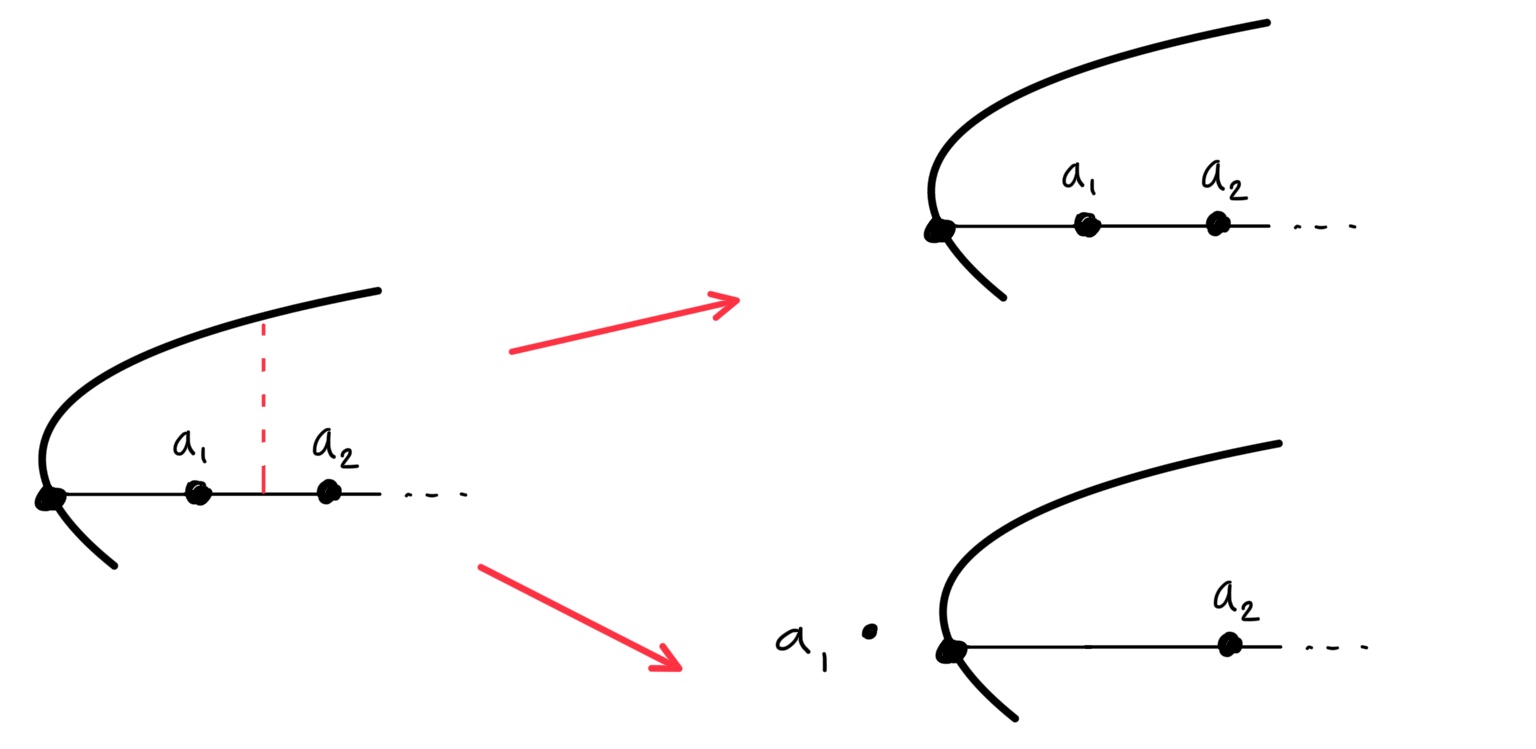}
		\end{center}
		which, on the level of trees, corresponds to 
		\begin{center}
			\includegraphics[width = 10cm]{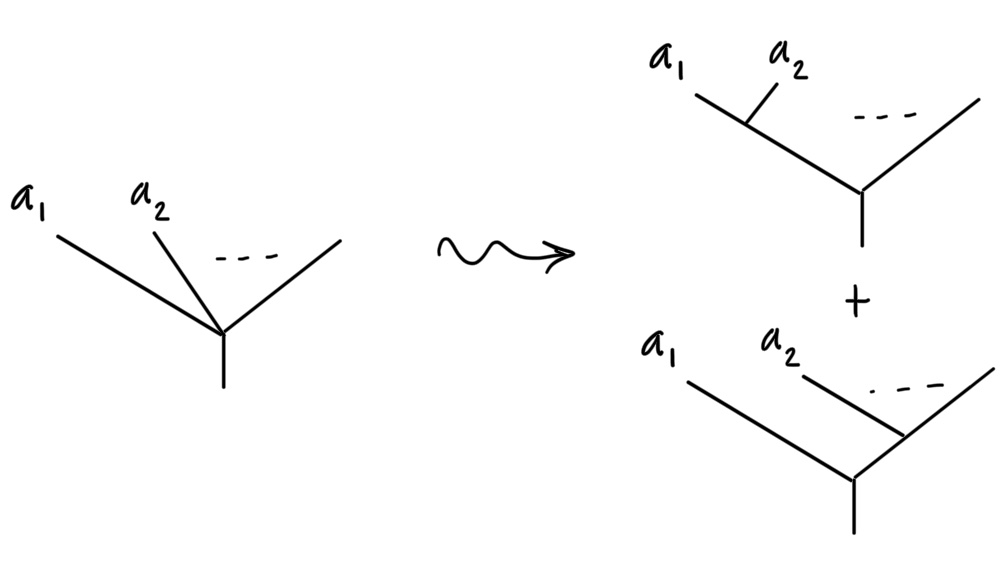}
		\end{center}
		Note that the $a_2$ vertex is not necessarily 2-valent; it can be the first of the $2N$-valent vertices, even though the 2-valent case is pictured above. The case in which we cut on the initial edge for a right-extended tree is also analogous.
		
		By the definition of augmented graphs (of the first kind) the graphs appearing on right hand sides of the pictures above were allowable. Hence, the corresponding trees on the right hand sides of the pictures above are operations. We can therefore summarize and reframe what we have just shown in the following lemma:
				
		\begin{lemma}\label{aa4}
			Each augmented graph of the first kind corresponds to a tree $T$ satisfying conditions (i), (iii), and (iv) from Theorem~\ref{aa2}, and in place of (ii), satisfies
			\begin{enumerate}[label = (\roman*')] \addtocounter{enumi}{1}
				\item The number of inputs is $n = j(2N - 2) + 3$, where $j \in \Z_{\geq 1}$, and $T$ contains a single multipliable pair;
			\end{enumerate}
			In this case, $T$ admits exactly one pull. In addition, any tree $T$ satisfying these conditions corresponds to a unique augmented graph of the first kind (and hence admits exactly one pull).
		\end{lemma}
		
		\begin{remark}\label{aa5}
			\emph{As usual, when we say $T$ ``admits'' a split or pull, we mean that the tree $T'$ obtained from $T$ by this method of adding an internal edge represents a non-zero operation.}
		\end{remark}
		
		The first part of Lemma~\ref{aa4} was proven in the paragraph above. The second part follows because a tree admits a pull if and only if it has a single multipliable pair. To get an augmented graph from an arbitrary tree satisfying the conditions from Lemma~\ref{aa4}, we pull together the two multipliable elements in $T$ to get bona fide operation $T'$, which has a well-defined, allowable graph. Then we draw a dotted line between the two elements which are multiplied in $T'$ but not in $T$, for instance:
		\begin{center}
			\includegraphics[width = 13cm]{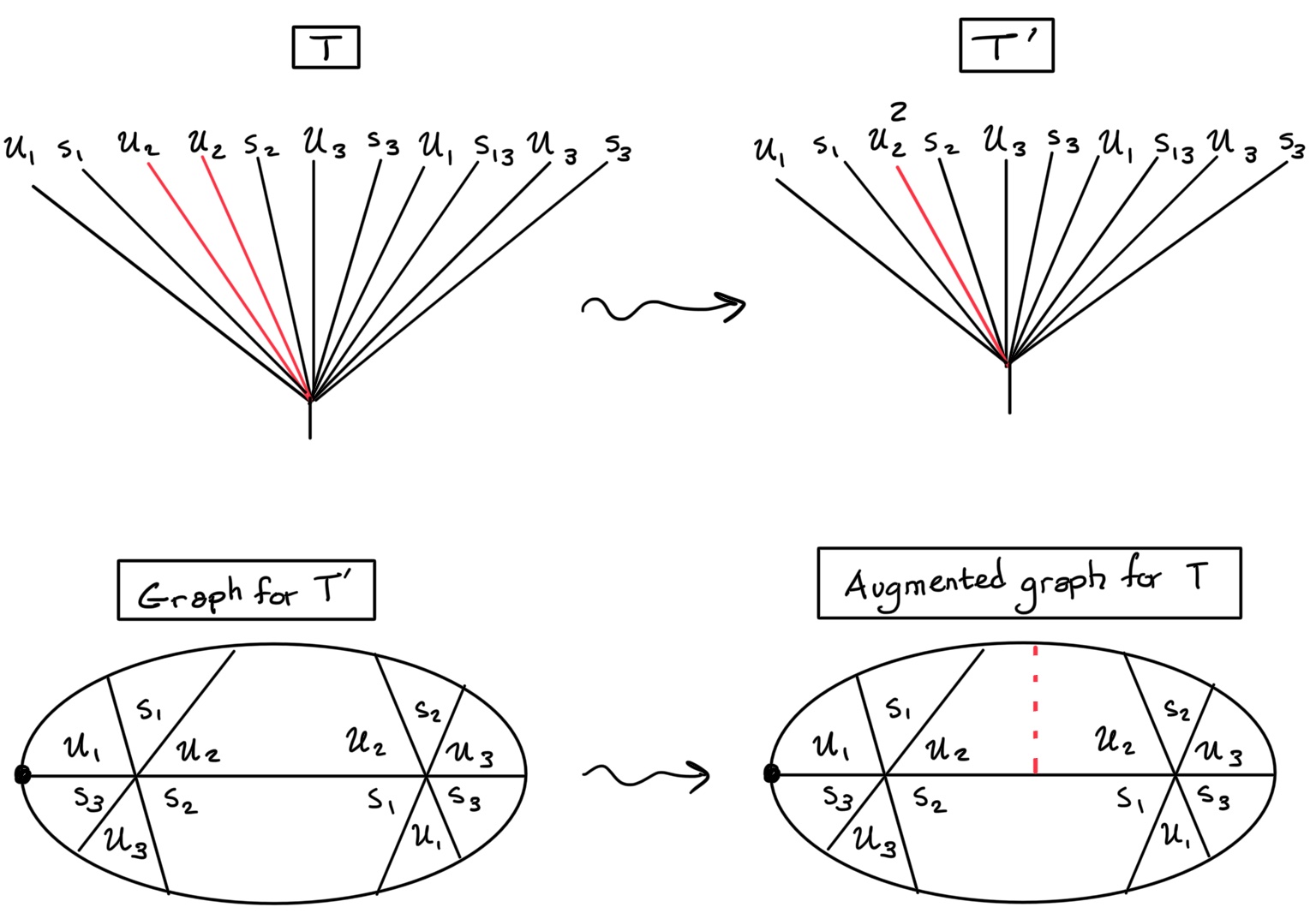}
		\end{center}
		It is also clear that by reversing this process, we can get back from an augmented graph to a tree of the form from Lemma~\ref{aa4}.
		
		Lastly, we need to show that if $T$ admits a split that is not of the particular kind dealt with in Lemma~\ref{aa19}, then $T$ is of the form from Lemma~\ref{aa4}. Suppose $T$ admits such a split. Then drawn as trees, the resulting compositions look like
		\begin{center}
			\includegraphics[width = 10cm]{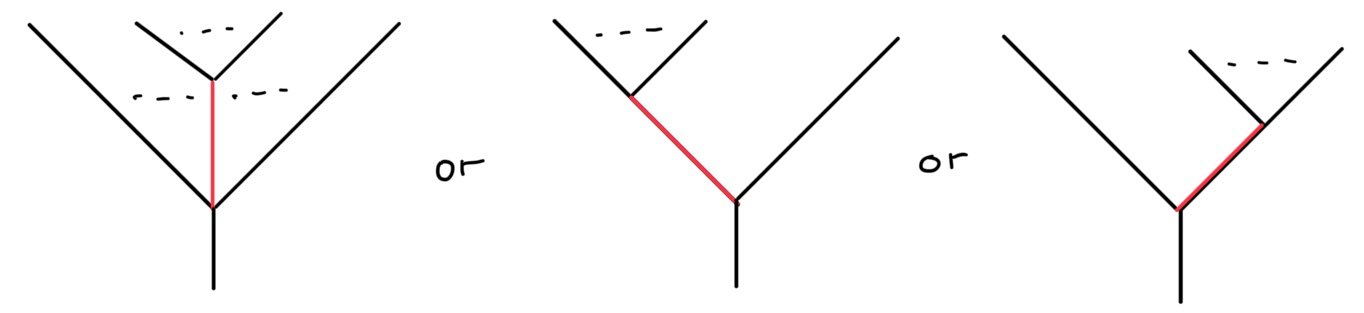}
		\end{center}
		In the first case, inner operation must be either left or right extended. If it were not, then the resulting tree would be of the form
		\begin{center}
			\includegraphics[width = 4cm]{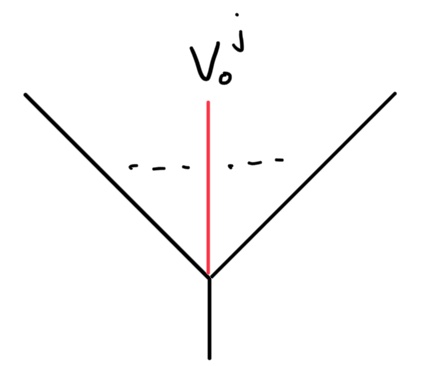}
		\end{center}
		i.e. it would not be an operation because it is not sufficiently rigid, since by the condition on Maslov gradings, operations cannot have multipliable pairs. In the other cases, the inner operation can be centered; there is no problem with a $\mu_2$ involving a raw power of $V_0$.
		
		By an argument analogous to the one in the proof of Theorem~\ref{aa2}, we can show that the compositions pictured above each correspond to one of the following types of graph:
		\begin{center}
			\includegraphics[width = 12cm]{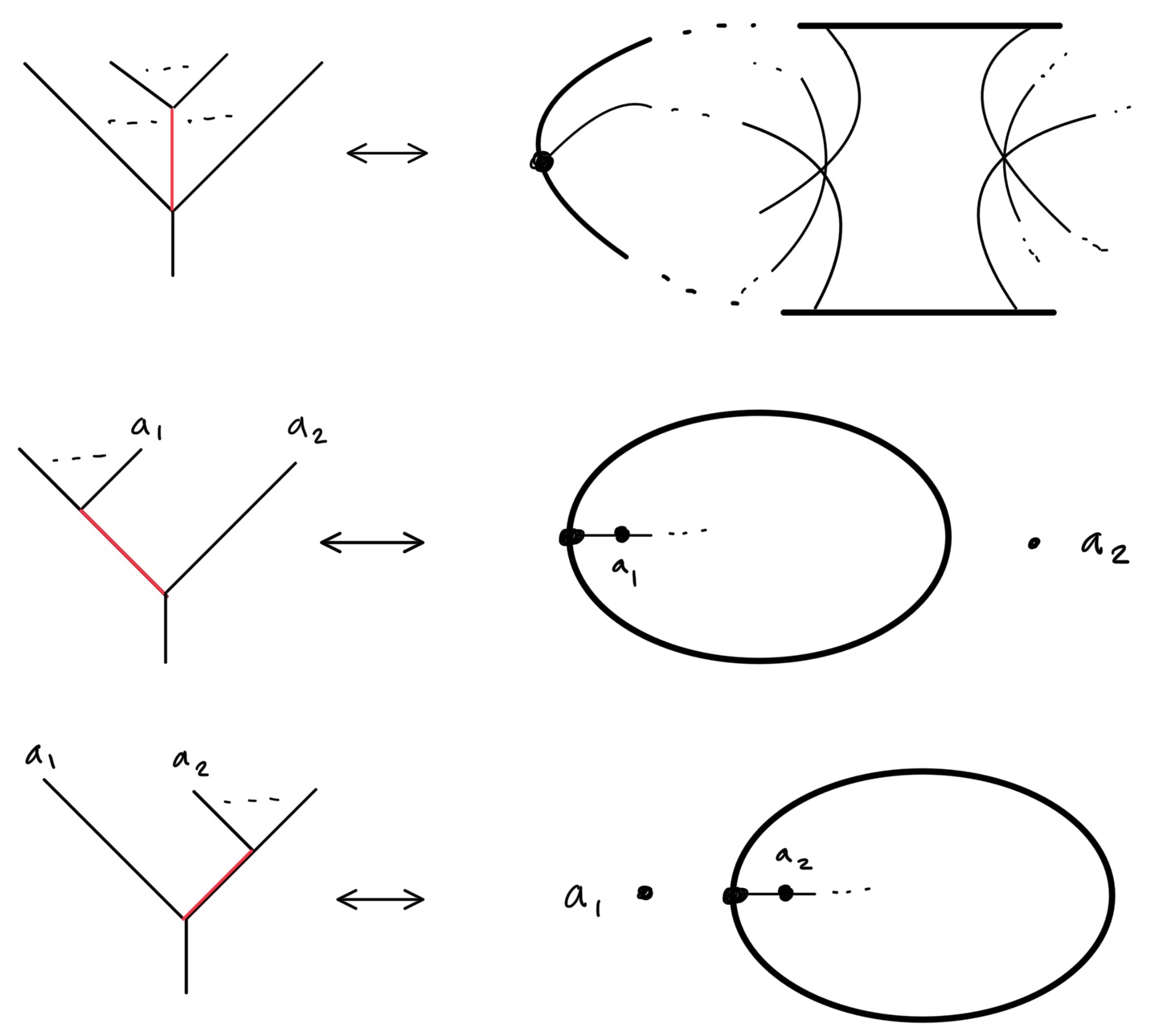}
		\end{center}
		Again, in the second two cases, the pictured vertex in the main graph can be $2N$-valent or $2$-valent; this changes nothing, and will be obvious from the composition. Notice also that each of these graphs has exactly two disjoint components. (We are counting the singleton elements on the left and right of the second and third, respectively, as formal ``connected graphs''.) In the second and third cases, the two-component graphs pictured above clearly came from pushing out the dotted edge of an augmented graph, that is,
		\begin{center}
			\includegraphics[width = 10cm]{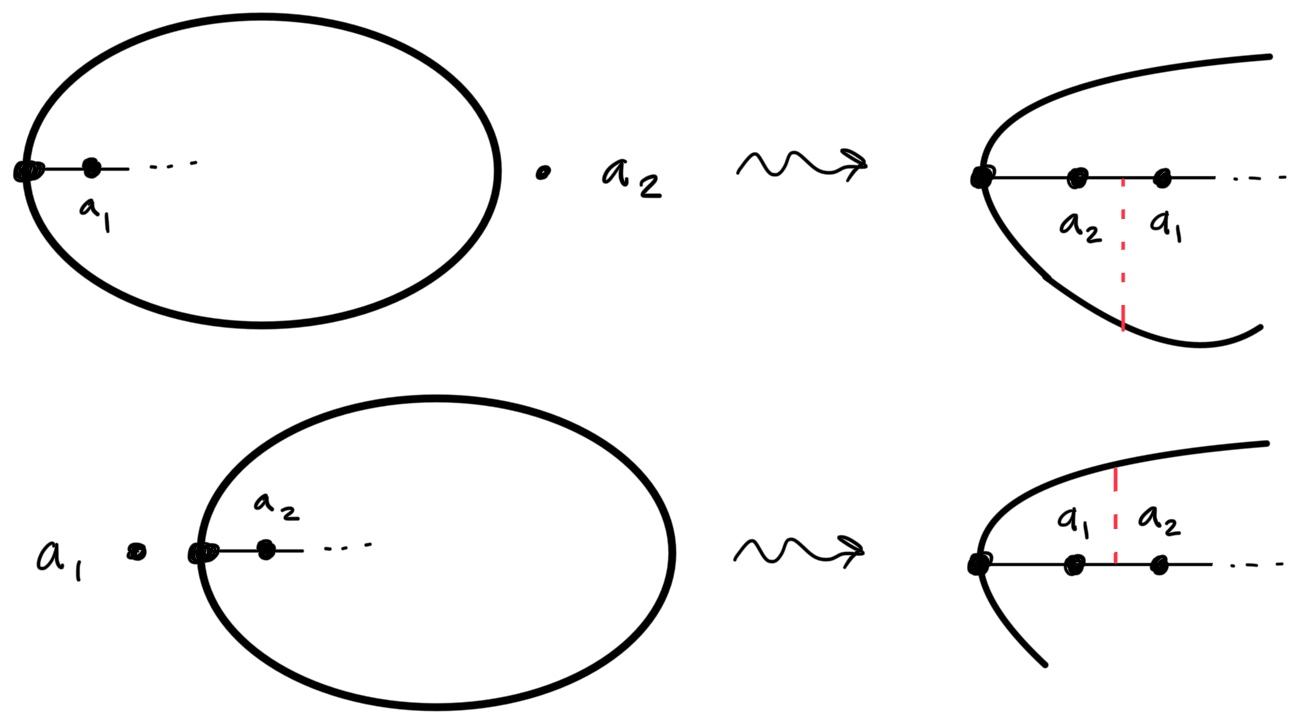}
		\end{center}
		drawn above with corresponding (original, given) trees. Likewise, in the first case, because the internal operation of the composition was either left or right extended, we can use this to unambiguously determine an augmented graph that gives rise to our given two component graph by pushing out an edge. For instance
		\begin{center}
			\includegraphics[width = 12cm]{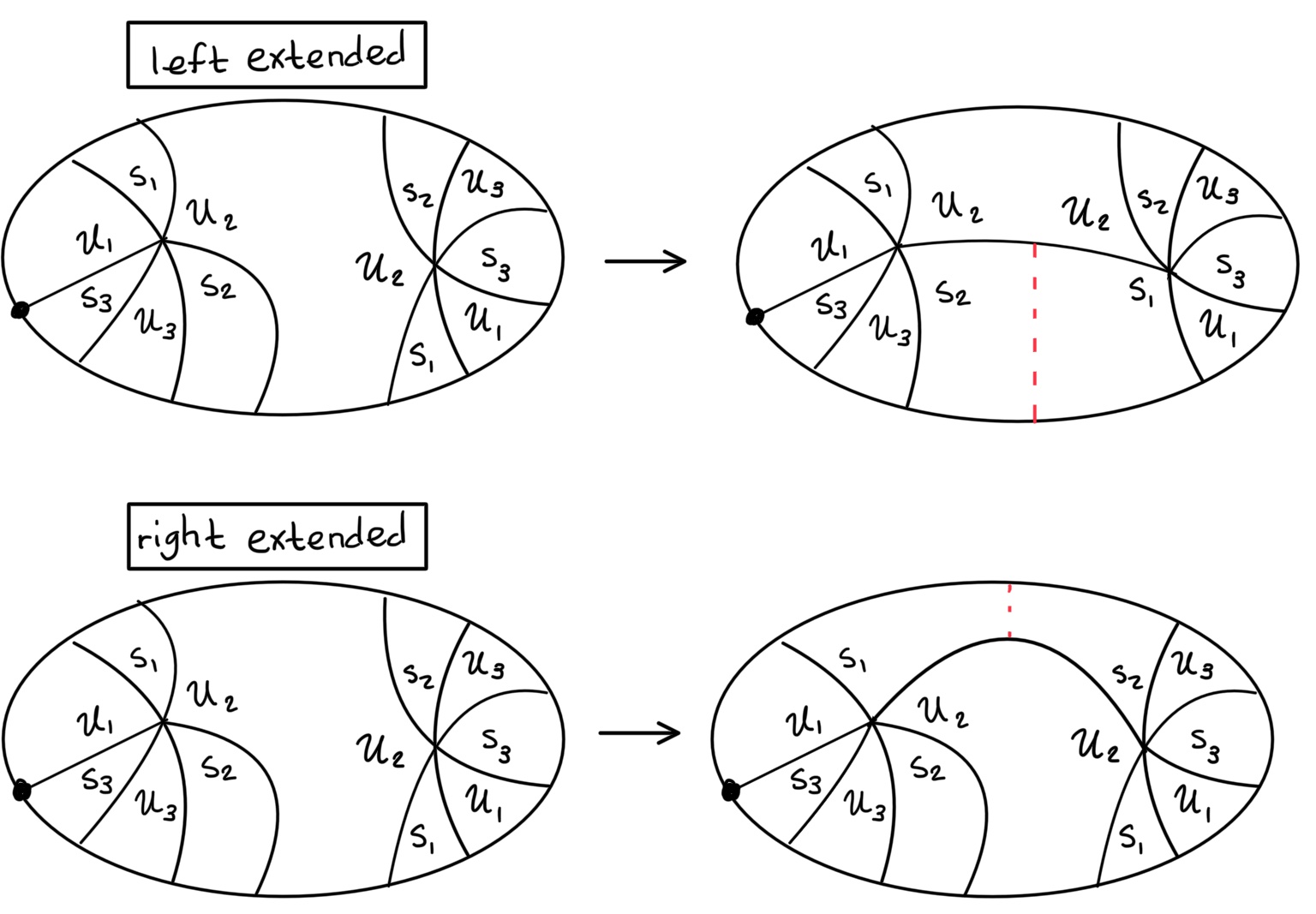}
		\end{center}
		Hence, in each of the three cases, our split corresponds to a unique augmented graph of the first kind. By Lemma~\ref{aa4}, these augmented graphs each give rise to a unique pull. This means that splits (of this kind) and pulls always cancel. The other kind of splits always cancel in pairs, as discussed before. Since these are the only possible ways to get non-zero terms in the $\A_{\infty}$-relation for a given unweighted tree, we have proved our claim.
	\end{proof}
	
	\begin{remark}\label{aa6}
		\emph{We did not really need to show that non-vanishing $\A_{\infty}$ -relations in the unweighted case were in one-to-one correspondence with augmented graphs. We could just have shown that a tree admits a pull or a split if and only if it satisfies the conditions of Lemma~\ref{aa4}, and that in this case it admits exactly one of each -- without mentioning augmented graphs at all. However, the augmented graphs will be important in the weighted case, which is why we include the subsidiary part of the argument above.}
	\end{remark}
	
	The next step is to run through the same process, allowing cycles. There are two kinds of cycles. The first are \emph{petal cycles}, that is $\e_i$-weight for $1 \leq i \leq N$, which, in the trees and graphs, correspond to e.g. 
	\begin{center}
		\includegraphics[width = 9cm]{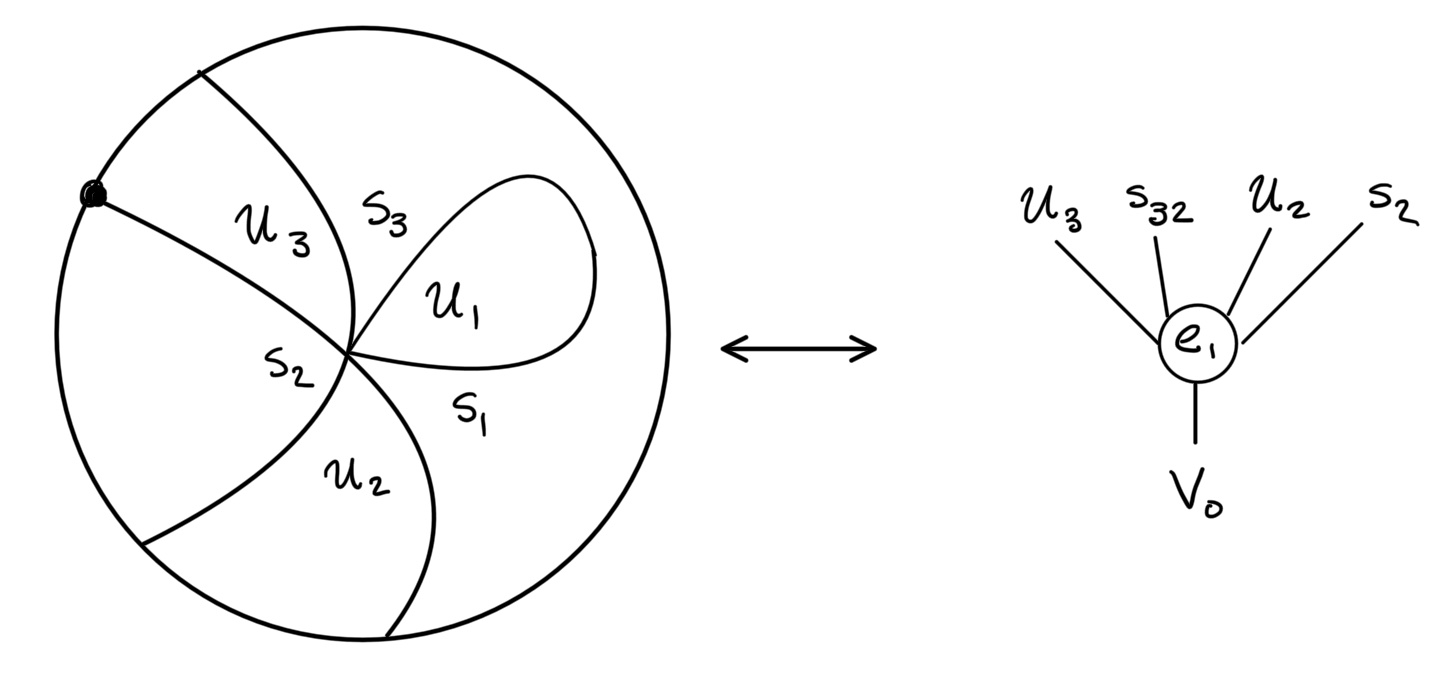}
	\end{center}
	The second are \emph{internal cycles}, i.e. $\e_{N + 1}$-weight, which shows up as e.g. (in the case $N = 3$ and with a particular choice of root)
	\begin{center}
		\includegraphics[width = 12cm]{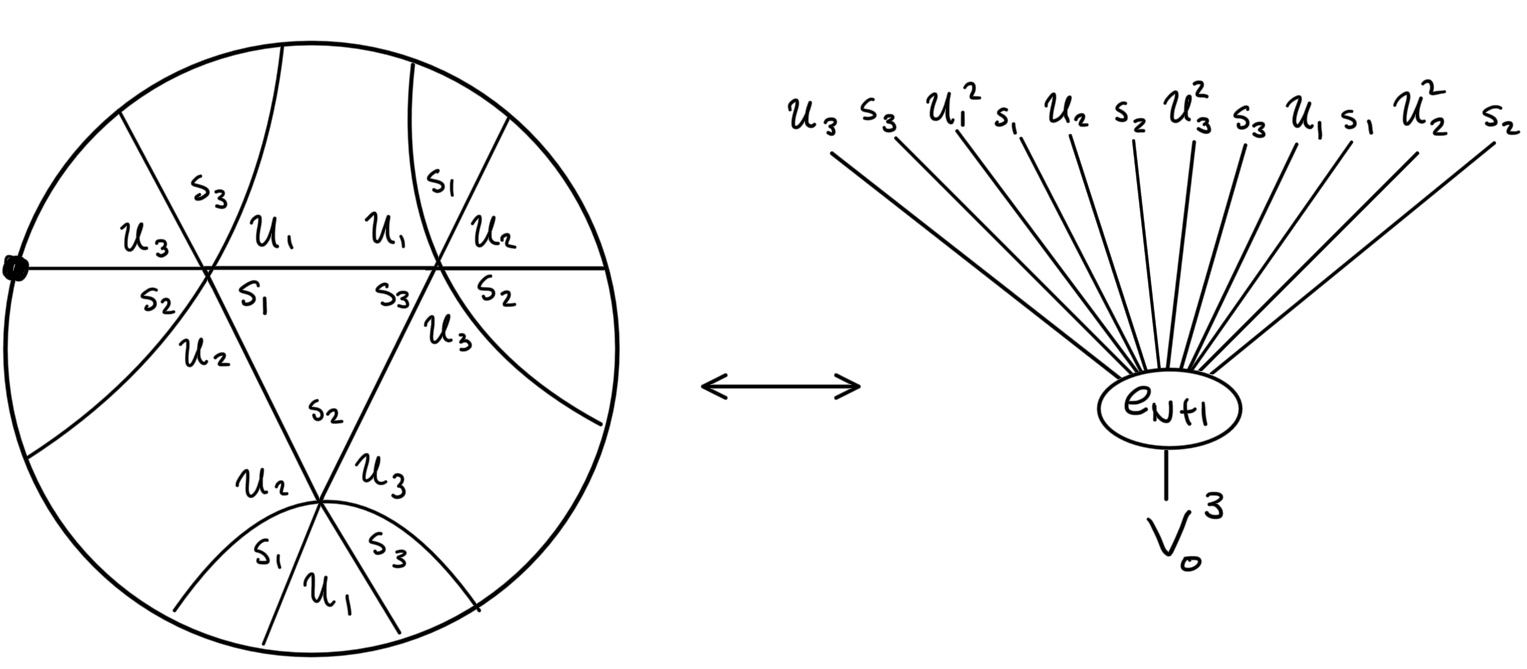}
	\end{center} 
	
	Before we can consider the case with cycles, we need to modify our definition of left / right extensions. In the unweighted case, a tree was extended whenever left- (resp. right-) most element was non-basic. However, in the weighted case, this is not sufficient. For instance the weighted sequence $\w = \e_2+ \e_3,$ with $(a_1, a_2) = (s_{11}, U_1)$ gives tree and graph,
	\begin{center}
		\includegraphics[width = 9cm]{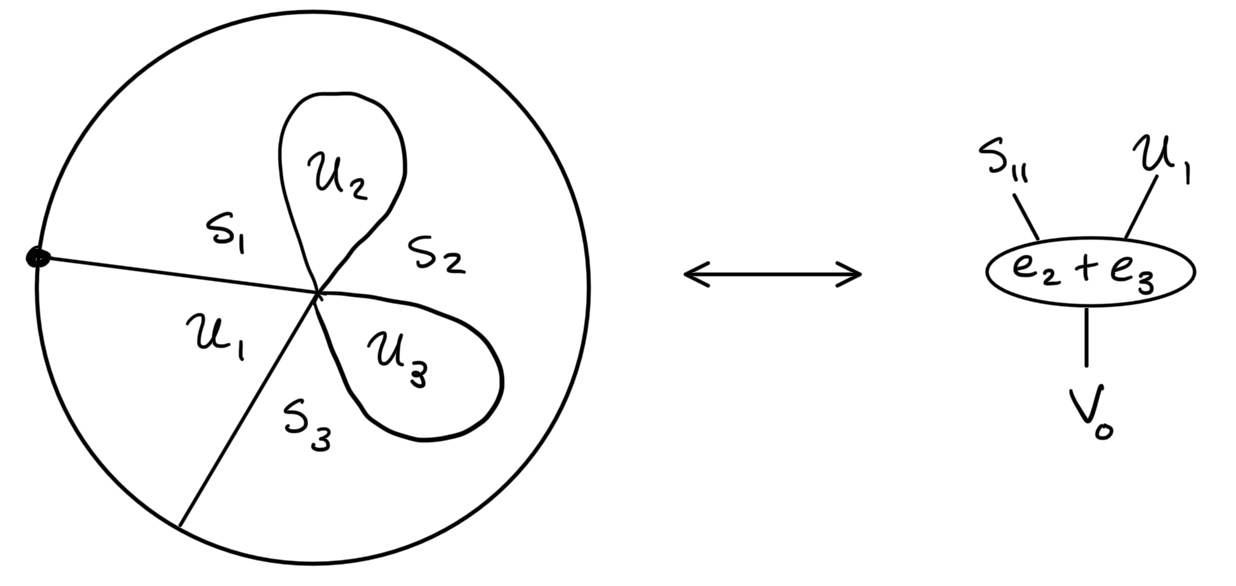}
	\end{center}
	which is a centered graph (there are no 2-valent vertices) so this condition is clearly not enough. Recall that in the weighted case, the other way of identifying non-centered trees was that these had unbalanced Alexander grading, i.e 
	\begin{equation}\label{aa8}
		A(a_1, \ldots, a_n) = j \cdot \sum_{k = 1}^{2N} \overline{k} + \text{ some non-zero term}
	\end{equation}
	where $j$ was assumed maximal in the sense that the non-zero term does not include another factor of $\sum_{k = 1}^{2N} \overline{k}$. Using the fact that $A(\e_i) = \overline{2i - 1}$ for each $1 \leq i \leq N$, the analogous condition to~\eqref{aa8}, for the petal-weighted case, is 
	\begin{equation}\label{aa9}
		A(\w, a_1, \ldots, a_n) = j \cdot \sum_{k = 1}^{2N} \overline{k} + \text{ some non-zero term}
	\end{equation}
	where again, $j$ is maximal in the sense that the non-zero term does not contain an extra factor of $\sum_{k = 1}^{2N} \overline{k}$. We also need to define the notion of an \emph{offset term}. Look at a weighted sequence $(\w, a_1, \ldots, a_n)$ satisfying~\eqref{aa9}. If we can write $a_1 = \alpha a_1'$, with $\alpha, a_1' \in \A$ and
	\[
		A(\w, a_1', a_2, \ldots, a_n) = j \cdot \sum_{k = 1}^{2N} \overline{k}.
	\]
	or if we can write $a_n = a_n' \alpha$, with $\alpha, a_n' \in \A$ and 
	\[
		A(\w, a_1, \ldots, a_{n - 1}, a_n') = j \cdot \sum_{k = 1}^{2N} \overline{k},
	\]
	then we say that $(\w, a_1, \ldots, a_n)$ has a \emph{well-defined offset term}, namely $\alpha$. If a weighted sequence satisfies~\eqref{aa9} and has a well-defined offset term, then it is \emph{left extended} if the offset is on the left, and right-extended if the offset is on the right. If, on the other hand, 
	\begin{equation}\label{aa10}
		A(\w, a_1, \ldots, a_n) = j \cdot \sum_{k = 1}^{2N} \overline{k},
	\end{equation}
	then this is a \emph{centered} sequence. There are sequences which are not centered, left, or right extended, but these do not correspond to suitable graphs, so we do not consider them. 
	
	We define the \emph{total magnitude} of a weight 
	\[
		\w = \sum_{j = 0}^{N + 1} k_j \e_{j}
	\] as $|\w| = k = \sum_{j = 0}^{N + 1} k_j$.
	
	Finally, recall that operations are still \emph{defined} to correspond to suitable (and now possibly non-simply-connected) graphs. For petal cycles, the condition for a tree to be an operation is as follows.
	
	\begin{theorem}\label{aa7}
		\emph{(Weighted higher operations)} Let $T$ be a tree with inputs $a_1, \ldots, a_n \in \A$, with $n \geq 0$, and some non-zero weight $\w$ which is some sum of the $\{\e_i\}_{i = 1}^{N + 1}$ with total magnitude $k$. Then $T$ represents an operation if and only if it satisfies the following conditions:
		\begin{enumerate}[label = (\roman*)]
			\item \emph{(The idempotents match up)} The initial idempotent of $a_i$ is the final idempotent of $a_{i - 1}$ for each $1 < i \leq n$; 
			
			\item \emph{(The Maslov grading works out)} $n = j (2N - 2) + 2 - 2k $ where $j \in \N$ and again, $|\w| = k$;
			
			\item \emph{(Extensions)} The sequence $(\w, a_1, \ldots, a_n)$ is either centered, left or right extended -- in particular, this means that there cannot be an offset on both left and right ends of the sequence, or somewhere in the middle; 
			
			\item \emph{(The Alexander grading works out)} If $(\w, a_1, \ldots, a_n)$ is centered, then we require it satisfy~\eqref{aa10}; if it is left or right extended, we require that it satisfy~\eqref{aa9};
						
			\item \emph{(Correct output)} The output is $V_0^j$ if the $(a_1, \ldots, a_n)$ is centered, and $\alpha V_0^j$ if it is extended, where $\alpha$ is the offset term, as usual;
		\end{enumerate}
	\end{theorem}
	
	\begin{proof}
		Again, it is clear how to read off a tree from an allowable graph (this time allowing cycles) -- just count inputs traveling counterclockwise along the boundary, as before. That the resulting tree satisfies (i) and (iii) is obvious. For (v), we just \emph{define} the output to be $V_0^j$, where e$j$ is the number of vertices. For (ii), we count the number, $j$, of vertices, and note that each petal cycle drops the number of sectors (and hence, $n$) by 2 without changing $j$ (from the unweighted case), and each internal cycle likewise drops the number of sectors by two, without changing $j$. Since each $\e_i$ adds 2 to the Maslov grading of a sequence, and the output is still $V_0^j$, the Maslov gradings work out. For (iv), note that each $2N$-valent vertex contributes a $j \cdot \sum_{k = 1}^{2N} \overline{k}$ to the Alexander grading of the tree, and the 2-valent vertices (if there are any) contribute the off-set term, so the Alexander grading is fine.
		
		For the reverse, start with a tree $T$ satisfying (i)-(v), with some non-zero weight. (If there is no weight, we are back to the case of Theorem~\ref{aa2}, in which case we are done.) The issue is how to ``distribute'' the weight in the graph so as to make it an operation which gives back $T$ when we read off the elements from each sector (as in the previous paragraph). For each petal weight, we can (by the condition on Alexander gradings) find a unique spot in the tree where a $U_i$ is missing between a (multiplied) pair $s_{i - 1}\cdot s_{i}$. Define $T_1$ by replacing this $s_{(i - 1) (i + 1)}$ with the sequence $(s_{i - 1}, U_i, s_i)$, and removing the $\e_i$ from the weight. This removes all petal weights from the tree. Now we only have (possibly) internal weights. Notice also that $T_1$ also satisfies all the conditions (i)-(v) from the statement. 
		
		We now  need to show that any tree satisfying conditions (i)-(v), with only $\e_{N + 1}$-weight, has a (unique) corresponding graph. This is obvious in the base case, which for $N = 3$, looks like
		\begin{center}
			\includegraphics[width = 12cm]{aa38}
		\end{center}
		up to choice of root. Now assume we have proven this for all $n' < n$, and $T_1$ is a tree satisfying (i)-(v), with only $\e_{N + 1}$ weight and $n$ inputs. Write $(k \e_{N + 1}, a_1, \ldots, a_n)$ for the weighted input sequence of $T_1$. The condition (ii) on $n$ forces there to be at least one ``extremal subsequence'', that is, a maximal sequence of basic elements. The meaning of maximal will depend on $k$. Notice that in the unweighted case, $n$ depended only on the number of vertices ($j$). In this case, every additional $\e_{N + 1}$ drops the number of inputs by 2, because essentially two of the sectors get folded together to make the canonical graph (pictured above in the case $N = 3$). This means that, depending on $k$, the maximal length sequence of basic elements  will take one of the following forms
		 \begin{enumerate}
		 	\item (The chunk corresponds to an unweighted piece of graph) A sequence of basic elements, specifically, a length $2N - 2$ subsequence $a_{i + 1}, \ldots, a_{i + 2N - 1}$ which was the middle $2N - 2$ elements of a cyclic permutation of $U_1, s_1, \ldots, U_N, s_N$ (e.g, in the case $N = 3$, one such sequence would be $s_2, U_3, s_3, U_1$). 
			
			 If we are dealing with this case again, the subsequence determines a labeled graph with a single $2N$-valent vertex. We can excise the subsequence from $a_1, \ldots, a_n$ as we did in the proof of Theorem~\ref{aa2} to get a shorter sequence $(k \e_{N + 1}, b_1, \ldots, b_{n - (2N - 2)})$ which still satisfies the conditions (i)-(v), and being of length $< n$, corresponds to a graph (with internal cycles). Attaching the smaller subgraph corresponding to the excised subsequence to the correct edge as we did there, we have a graph for $T_1$, which is clearly unique.
			
			\item (The chunk corresponds to a piece of the graph containing internal cycle) The other option is that the maximal sequence of basic elements (the extremal portion of the tree) is a subsequence, of length either $N (2N - 2) -2$ or $N (2N - 2) - 3$, which consists of the middle terms of a cyclic permutation of the sequence from the outside of the basic graph for internal weight. In the case $N = 3$, this basic graph is pictured above.
			
			  We can use this subsequence to label the outer sectors of the basic graph for internal weight, except for two or three (adjacent) ones, so for instance
		 \begin{center}
		 	\includegraphics[width = 7cm]{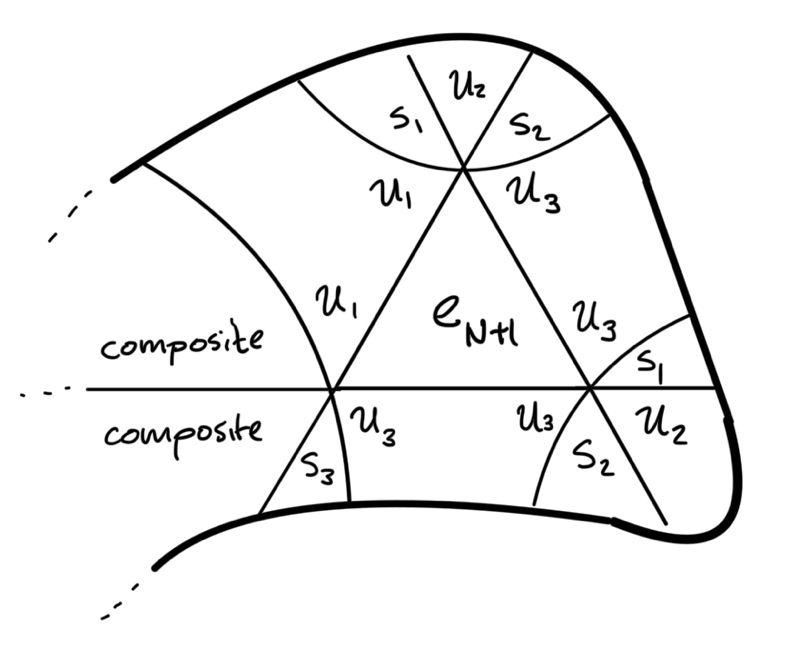}
		 \end{center}
		 which could appear if this internal cycle was attached to an unweighted basic chunk, e.g. 
		 \begin{center}
		 	\includegraphics[width = 8cm]{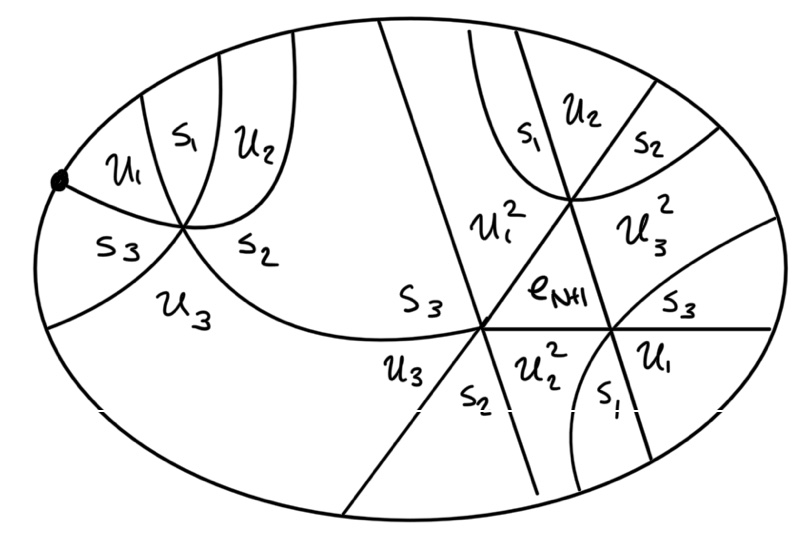}
		 \end{center}
		 On the other hand, this maximal sequence could be much shorter, as in the case
		 \begin{center}
		 	\includegraphics[width = 12cm]{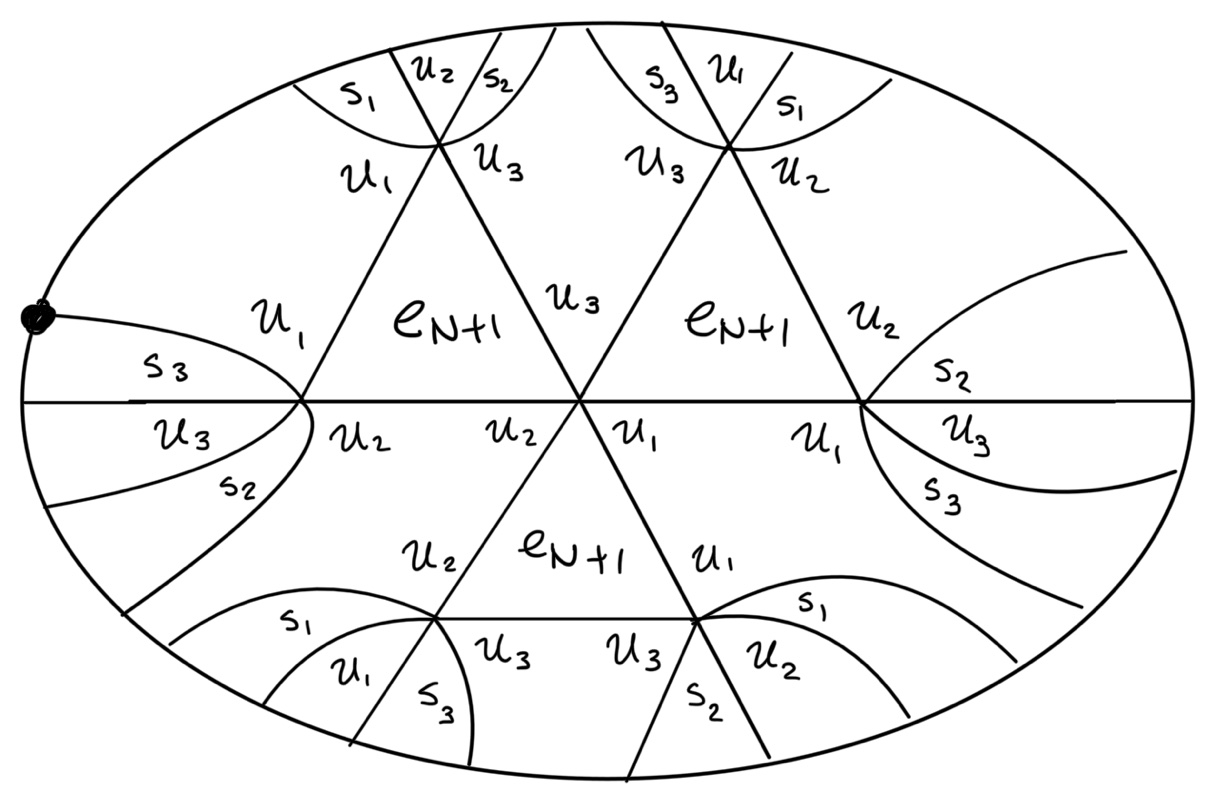}
		 \end{center}
		 Here, the maximal sequence is of length $7 = N (2N - 2) - (2N - 1)$ (so we lost all of the possible labels around the vertex attached to the rest of the graph). 
		 
		 This is what was meant by the ambiguity in the length of the subsequence. Viewed on the level of the graph we are going to construct, if this ``extremal'' chunk of the graph is attached to an unweighted bit, then we are looking at the first picture, and the subsequence is of length $N (2N - 2) - 2$ -- we only lose the two labels attached to the adjacent unweighted bit. On the other hand, if it is directly attached to one or more internally weighted chunks, the maximal sequence of unaffected elements (which look like a sequence from the basic internally weighted pictures) becomes shorter, down to $N(2N - 2) - (2N - 1)$. This is not an issue, however, since in each case, we can still retrieve enough information about this extremal bit to excise it and reduce to a shorter graph, which is what we want.
		 
		  Now write the subsequence as $a_{i + 1}, \ldots, a_{i + m -1}$. The the condition on the idempotents forces $a_i$ and $a_{i + m}$ to be the product of the two missing elements from the canonical internal weight chunk, multiplied, on the left and the right respectively, by some other elements $a_i'$ and $a_{i +m}'$. There are two subcases to consider, corresponding to the two pictures above. 
		  
		  In the first subcase, the initial idempotent of $a_i'$ and the final idempotent of $a_{i + m}'$ match up (and $m = N(2N - 2) + 1$) and we can consider the tree $T_2$ with weighted input sequence 
		  \[
		  ((k - 1) \e_{N + 1}, a_1, \ldots, a_{i - 1}, a_i', a_{i + m}', a_{i + m + 1}, \ldots, a_n),
		  \]
		  and output $V_0^{j - N}$. This will still satisfy (i)-(v), and since it satisfies the inductive hypothesis, it gives back a well-defined graph. Attaching the chunk from the first picture above to the correct edge of this graph, we get back a graph for $T_1$. 
		  
		  	In the second subcase, the idempotents of $a_i'$ and $a_{i + m}'$ do not match up (because, from the perspective of the graph we are trying to construct, there is an $s_j$ between them), but they will be of the form
		  	\begin{align*}
		  		a_i' &= U_{j - 1}^p \\
				a_{i + m}' &= U_j^q,
			\end{align*}
			where $p, q \in \N$, and $j, (j - 1)$ are counted $\emm N$. In this case, we consider the tree $T_2$ with weighted input sequence 
			\[
				((k - 1) \e_{N + 1}, a_1, \ldots, a_{i - 1}, a_i', s_j, a_{i + m}', a_{i + m + 1}, \ldots, a_n).
			\]
			and output $V_0^{j - (N - 1)}$. This still satisfies (i)-(v) and is of length strictly less than $n$, so we have a graph for $T_2$. Attaching the chunk from the second picture above to the correct edge of this graph, we get back a well-defined graph for $T_1$.
		\end{enumerate}
		
		This now gives us a way to associate a graph to any internally weighted tree $T_1$ satisfying (i)-(v). To conclude, we recall how $T_1$ was initially obtained from our original $T$ (which included petal weights in addition to internal weights). In place of each petal weight $\e_i$ (and corresponding $\alpha s_{(i - 1) (i + 1)} \beta$) we inserted an $\alpha s_{i - 1}, U_i, s_i \beta$ into the sequence in the proper place. (There could be other terms multiplying the $s_{(i - 1)(i + 1)}$, but this is the only relevant part.) This corresponds to a portion of the graph for $T_1$ that looks like
		\begin{center}
			\includegraphics[width = 5cm]{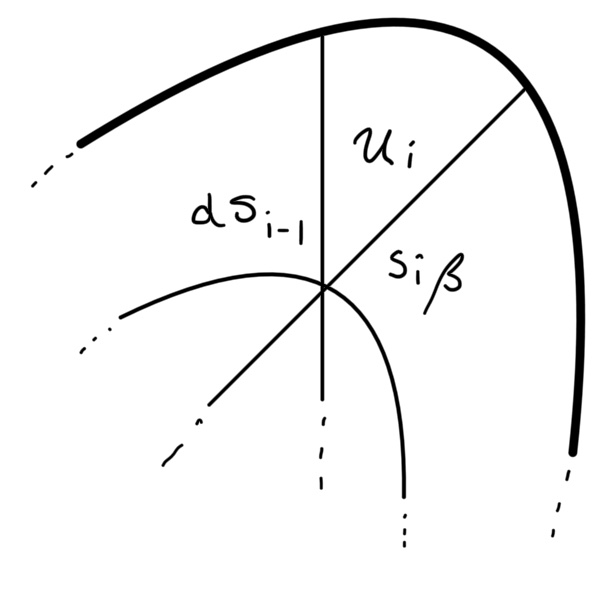}
		\end{center}
		Replacing each such section (which originally corresponded to a petal weight) with 
		\begin{center}
			\includegraphics[width = 5cm]{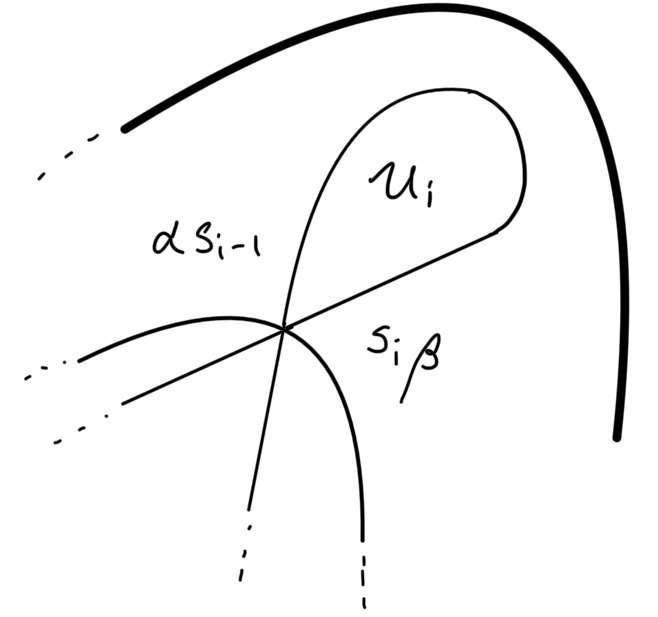}
		\end{center}
		The input sequence corresponding to this new graph is the same as the one for $T_1$, except with each subsequence $\alpha s_{i - 1}, U_i, s_i \beta$ (which originally corresponded to a petal weight) replaced with $\alpha s_{(i - 1)(i + 1)} \beta$, with an extra $\e_i$ in the weight; in other words, it is precisely the original $T$.
	\end{proof}
	
	Next come the $\A_{\infty}$-relations. There are, in this case, three ways to add an edge to a tree $T$: a push, a pull, and a split. First, we have the following analogue of Lemma~\ref{aa4}:
	
	\begin{lemma}\label{aa11}
		A weighted tree $T$ admits a pull precisely when it satisfies conditions (i), (iii), and (iv) from Theorem~\ref{aa7}, and in place of (ii), satisfies
		\begin{enumerate}[label = (\roman*')] \addtocounter{enumi}{1}
			\item $n = j(2N - 2) + 3 - 2k$, where $j \in \Z_{\geq 1}$, $n$ is the number of inputs, and $k$ is the total magnitude of the weight of $T$;
		\end{enumerate}
		In this case, $T$ admits exactly one pull, and $T$ can be expressed uniquely as an augmented graph. In addition, each augmented graph corresponds to  a tree satisfying these conditions.
	\end{lemma} 
	
	This condition is exactly the same as Lemma~\ref{aa4}, except that (ii)' is modified to fit the weighted case. The proof of the part about pulls is obvious, and is omitted. The only thing that needs verification is the part about augmented graphs, but this follows exactly as in the proof of Lemma~\ref{aa4}, and can also be omitted.
	
	In order to verify the $\A_{\infty}$ relations, we just need to show that each split or push corresponds to an augmented graph, and that each augmented graph corresponds to a split or push, but not both. The second part is easier:
	
	\begin{lemma}\label{aa12}
		Every augmented graph gives rise to exactly one split or one push, but not both. 
	\end{lemma}
	
	\begin{proof}
		Choose an augmented graph $G$ of the first kind. The dotted line in $G$ be draw in one of five ways. First we have the three options from the unweighted case -- the dotted line going from a true internal edge to the boundary, or from the initial edge to the boundary in the left extended case, and from the final edge to the boundary in the right extended case. We also have the following two new options:
		\begin{center}
			\includegraphics[width = 8cm]{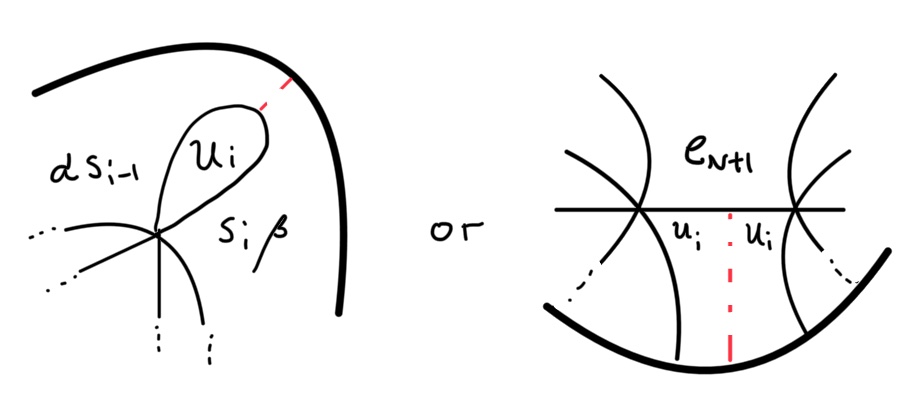}
		\end{center}
		To get the pull from these diagrams (as we did in Lemma~\ref{aa11}) we just erase the dotted line. For each augmented graph, there is a corresponding tree given by reading off the elements clockwise around the boundary from the root in the usual way, and counting the dotted line as an unmultiplied multipliable pair. For the graphs above, the multipliable pairs are $(\alpha s_{i - 1}, s_i \beta)$ and $(U_i, U_i)$, respectively. As per Lemma~\ref{aa11}, these trees satisfy the conditions (i), (iii), and (iv), and have one too many inputs, corresponding to the multipliable pairs pictured above.
		
		For our present situation, we look back at the original augmented graph, and ``open out'' the edge abutting the dotted line, in the following way
		\begin{center}
			\includegraphics[width = 11cm]{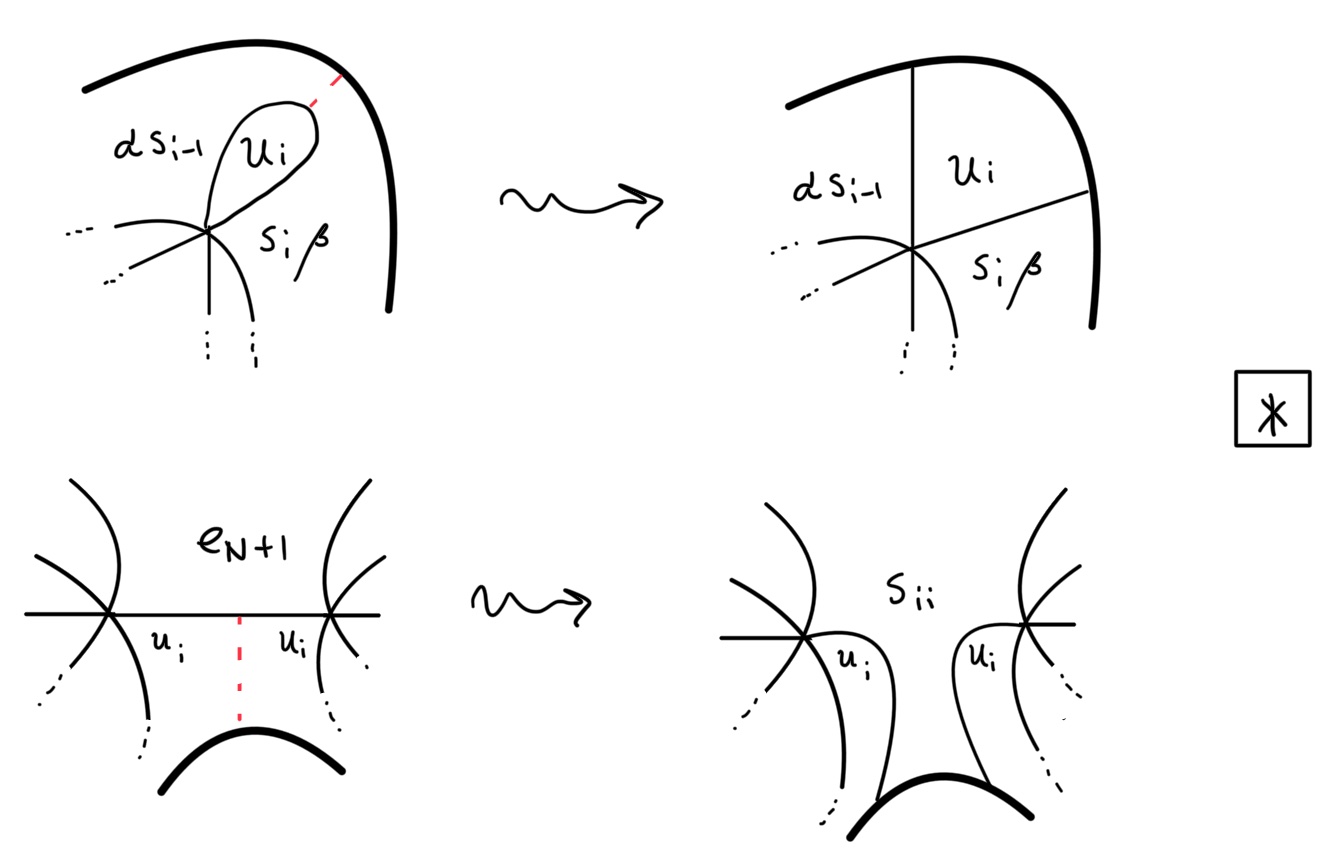}
		\end{center}
		In terms of trees, these are pushes, that is 
		\begin{center}
			\includegraphics[width = 10cm]{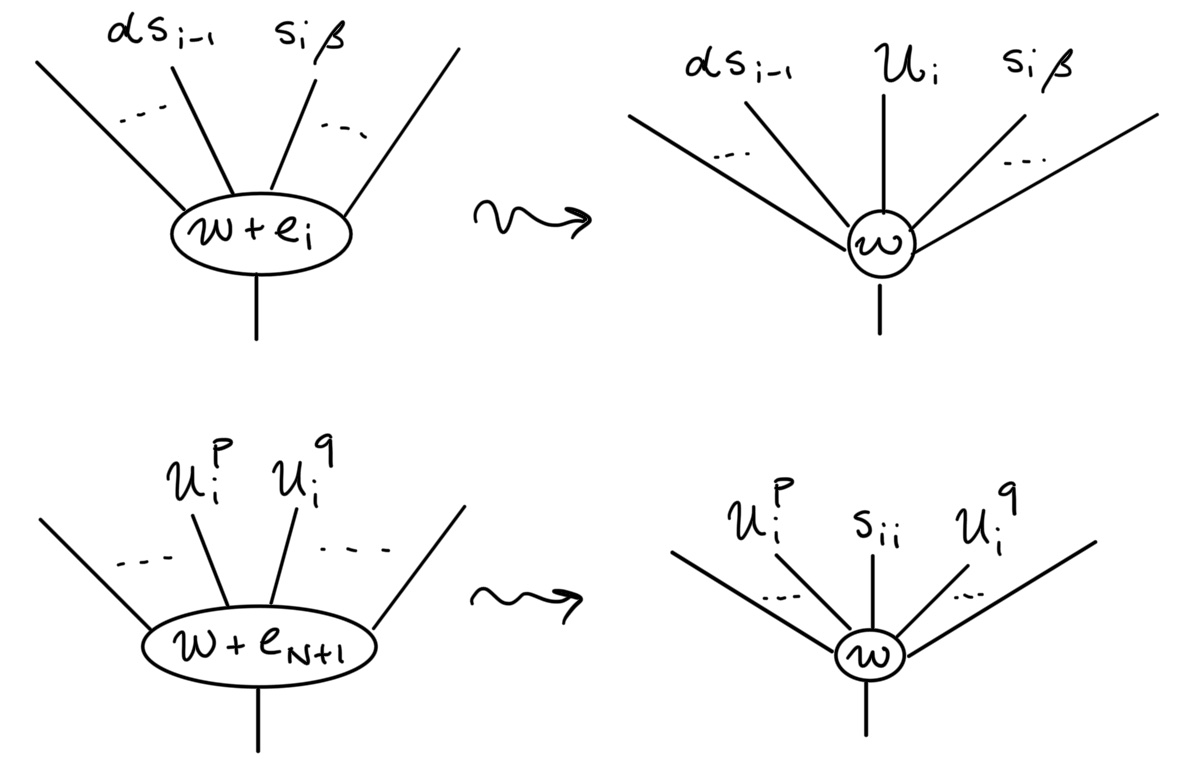}
		\end{center}
		By the remark at the end of the previous paragraph, these are bona fide operations. In the first case, this follows by a counting argument as in the proof of the $\A_{\infty}$-relations for the unweighted case. For the other two, this is just because we have balanced out the Maslov grading by deleting a weight term and adding an input, without messing up any of the other conditions. 
		
		That the second kind of augmented graph always gives rise to the same specific kind of split follows exactly as in the unweighted case.
	\end{proof}
	
	Now we need to verify the other direction, namely:
	
	\begin{lemma}\label{aa13}
		Every split or push gives rise to a unique augmented graph.
	\end{lemma}
	
	\begin{remark}\label{aa14}
		\emph{Graphically, pushing a cycle out of an augmented graph looks like the two starred pictures from the proof of the previous lemma. The goal of the ``push'' part of the proof of Lemma~\ref{aa13} is to retrieve this picture from the tree for a push.}
	\end{remark}
	
	\begin{proof}
		That a split gives rise gives rise to a unique augmented graph follows as in the proof of Theorem~\ref{aa3}; the weight has no effect on the proof. 
		
		For the case of pushing out a basic element of weight (i.e. $\e_i$ for $1 \leq i \leq N + 1$) we need to remember what exactly it means for a tree $T$ to admit a push. Writing the weighted input sequence for $T$ as $(\w, a_1, \ldots, a_n)$, $T$ admits a push (of weight $\e_i$) if and only if we can find $1 \leq j \leq n$ so that the tree $T'$ with the same output as $T$, and input sequence 
		\begin{equation}\label{aa15}
			(\w - \e_i, a_1, \ldots, a_{j - 1}, U_i, a_j, \ldots, a_n)
		\end{equation}
		satisfies (i)-(v) from Theorem~\ref{aa7} and is therefore an operation. Suppose first that $1 \leq i \leq N$ -- that is, we are pushing out petal weight. Then in order for the sequence from~\eqref{aa15} to determine an operation, we need $a_{j - 1} = \alpha s_{i - 1}$ and $a_j = s_i \beta$, for $ \alpha, \beta \in \A$, left-multipliable with $s_{i - 1}$ and right-multipliable with $s_i$, respectively. This means that the original weighted input sequence for $T$ is:
		\begin{equation}\label{aa16}
			(\w, a_1, \ldots, s_{i - 1}, s_{i}, \ldots, a_n).
		\end{equation}
		It must still satisfy (i), (iii), (iv), and (v) from Theorem~\ref{aa7}, and comparing the sequences from~\eqref{aa15} and~\eqref{aa16}, we see that it must also satisfy (ii)' from Lemma~\ref{aa11}. This means that $T$ corresponds to a unique augmented graph. 
		
		Now, consider the case where $i = N + 1$. Then since the idempotents in the sequence from~\eqref{aa15} work out, and it contains no multipliable pairs, there must be some $1 \leq k \leq N$, as well as $p, q \in \N$, so that $a_{j - 1}= U_{\ell}^p, \: a_{j} = U_{\ell}^q$. Remembering that $U_{N + 1}$ is actually a sum of one term in each idempotent, and in any given situation all the terms except the one in the correct idempotent will drop out, the weighted input sequence for $T'$ looks like
		\begin{equation}\label{aa17}
			(\w - \e_{N + 1}, a_1, \ldots, U_{\ell}^p, s_{\ell \ell}, U_{\ell}^q, \ldots, a_n)
		\end{equation}
		This determines a well defined graph, and the portion of this graph corresponding to the $U_{\ell}^p, s_{\ell \ell}, U_{\ell}^q$ segment looks like
		\begin{center}
			\includegraphics[width = 4cm]{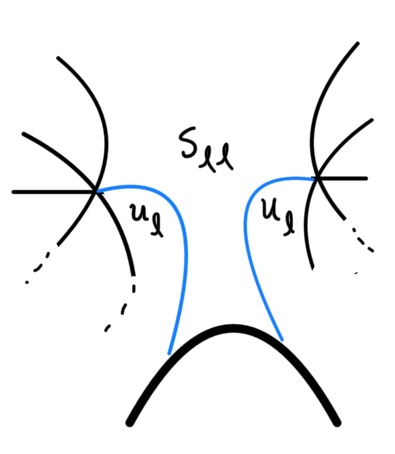}
		\end{center}
		(Here, the coloring on the two bottom edges is to facilitate the next part of the proof, and does not indicate anything else.) To get an augmented graph for $T$, zip the two colored edges together, and add a dotted edge out to the boundary, as
		\begin{center}
			\includegraphics[width = 4cm]{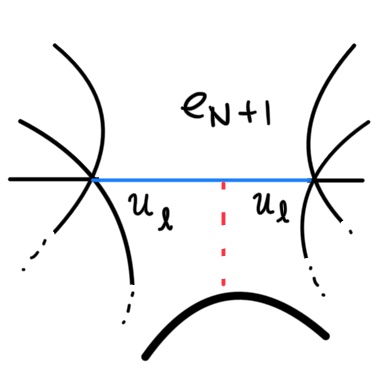}
		\end{center}
		This clearly gives rise to a tree, say $T''$. Zipping these edges together (and adding the dotted edge, to keep the two powers of $U_{\ell}$ separate) added an $\e_{N + 1}$-weight and deleting the $s_{\ell\ell}$, but did not change anything of the other inputs or weight from the bona fide suitable graph for $T'$. This means that these are the only ways $T''$ differs from $T'$, that is, $T'' = T$, and we have found an augmented graph for $T$, as desired. This is clearly unique, because the place where we are allowed to push out the $\e_{N + 1}$ is well-defined, which completely determines $T'$, and hence the augmented graph. Thus, the proof is complete. 
	\end{proof}
	
	Now, recall that by definition, all terms of $\A_{\infty}$-relation corresponding to a given tree $T$ vanish identically unless $T$ admits a push, a split, or a pull. Combining Theorem~\ref{aa3} with Lemmas~\ref{aa11},~\ref{aa12}, and~\ref{aa13}, we have now shown that in all cases when one term of the $\A_{\infty}$-relation for some $T$ is non-vanishing, there is a unique cancelling term, that is:
	
	\begin{theorem}\label{aa18}
		$\A$ satisfies the $\A_{\infty}$-relations, and is a bona fide $\A_{\infty}$-algebra. 
	\end{theorem}
	
	\section{The $\beta$-bordered algebra $\B$}\label{bb1}
	
	\subsection{Algebra elements and basic operations}
	
		The algebra $\B$ from Theorem~\ref{duality} is constructed as follows. $\B$ is an $R = \F[V_0, \ldots, V_{N + 1}]$-module with generators $\{\rho_i\}_{i = 1}^{N}$ and $\{\sigma_i\}_{i = 1}^{N}$. Simple multiplication (i.e. $\mu_2$) is defined by
		\begin{align*}
			\rho_i \rho_j &= 0 \text{ for each } i, j \\
			\sigma_i \sigma_j &= 0 \text{ for each } i, j \\
			U_0 &:= \sum_{i = 1}^N \rho_i \sigma_i \rho_{i + 1} \sigma_{i + 1} \cdots \rho_{i + N - 1} \sigma_{i + N - 1} \\
			&\qquad+ \sum_{i = 1}^N \sigma_i \rho_{i + 1} \sigma_{i + 1} \rho_{i + 2} \cdots \sigma_{i + N - 1} \rho_i
		\end{align*}
		where, in the last definition, all indices are counted $\emm N$. We also assume the $V_i$ to be central elements which can be multiplied by anything. (Again, as in the previous section, two algebra elements are said to be \emph{multipliable} or \emph{can be multiplied} if and only if they give non-zero product.) The $\rho_i$ each live in the $i$-th idempotent (i.e. the initial and final idempotents are both $i$) and each $\sigma_i$ has initial idempotent $i$ and final idempotent $i + 1$. We have $\de \rho_i = V_i$ for each $1 \leq i \leq N$, and one new weighted $\mu_0$, namely
		\begin{center} 
			\includegraphics[width = 5cm]{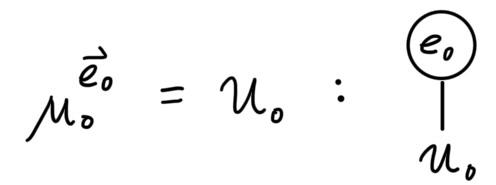}
		\end{center}
		There are no other weighted operations. This weighted $\mu_0$ does not affect any of the $\A_{\infty}$-relations for higher multiplication defined below -- it only adds the requirement that $\de U_0 = 0$, which is not a surprise. Therefore, we disregard it in the subsection that follows, in which we prove the $\A_{\infty}$ relations for unweighted operations. 
		
		\subsection{Higher operations and $\A_{\infty}$-relations}
		
		We also define additional nonzero $\mu_k$ for $1 \leq k \leq N$, according to the following rules.
		
		First, recall that a sequence of algebra elements $(\tau_k, \ldots, \tau_1)$ in $\B$ is \emph{in an allowable sequence of idempotents} if the initial idempotent of $\tau_{i + 1}$ equals the final idempotent of $\tau_i$, for each $i$. We then define the \emph{length of an algebra element $\tau \in \B$} by decomposing it into a product of $\sigma$s and $\rho$s: this decomposition is unique, and the length of $\tau$ is defined to be the number of $\sigma$s that appears int his product. (We really want to say the length is the difference between the final and initial idempotents of $\tau$, but if there are more than $N$ $\sigma$s in the product that makes up $\tau$, then we will be short some number of $N$s.) We then define the \emph{total length} of $(\tau_k, \ldots, \tau_1)$ as the sum of the lengths of each of the individual elements,
		\[
			\ell (\tau_k, \ldots, \tau_1) = \sum_{j = 1}^{k} \ell (\tau_k).
		\]
		We say that the sequence $(\tau_k, \ldots, \tau_1)$ \emph{stretches too far} if 
		\begin{enumerate}[label = (S\arabic*)]
			\item for some $j > 1$, the sequence $(\tau_k, \ldots, \tau_j)$ is of total length $\geq N$, \emph{and}
					
			\item for some $j < k$, the sequence $(\tau_j, \ldots, \tau_1)$ is of total length $\geq N$;
		\end{enumerate}
		If one of these two fails -- say (S1) -- then counting backwards from the final idempotent of $\tau_k$, the point at which we have backed up exactly $N$ steps will be somewhere in the decomposition (into $\rho$s and $\sigma$s) of $\tau_1$. We define this to be the \emph{cut point}. For instance:
		\begin{center}
			\includegraphics[width = 8cm]{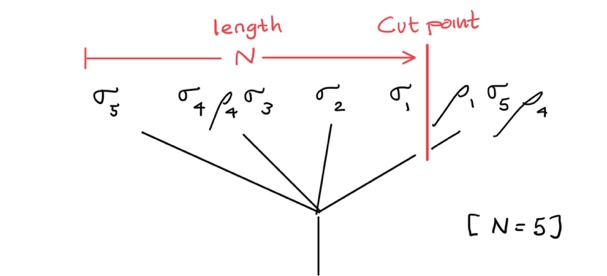}
		\end{center}
		is the cut point in one cases where (S1) fails. If (S2) fails, then counting forwards from the initial idempotent of $\tau_1$, the point at which we have gone exactly $N$ steps is somewhere in the decomposition of $\tau_k$. In this case we define this to be the cut point. For instance
		\begin{center}
			\includegraphics[width = 8cm]{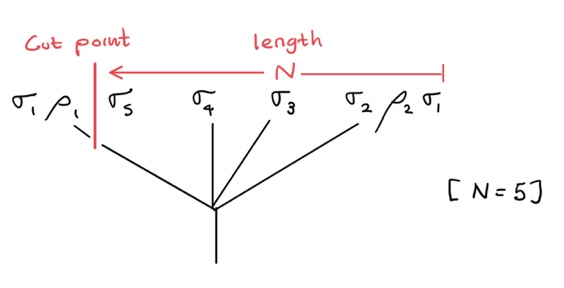}
		\end{center}
		is the cut point for one cases where (S2) fails. If both (S1) and (S2) fail -- that is, $(\tau_k, \ldots, \tau_1)$ does not stretch too far no matter which direction we count from -- we have ambiguity in the definition of the cut point. In this case, we just define it counting back from the final idempotent of $\tau_k$. After we have defined the operations, we will show that this ambiguity does not matter.  
		
		We say that $(\tau_k, \ldots, \tau_1)$ is \emph{left flanked} if $\tau_k = \rho_{\lambda'} \tau_k'$, where $\lambda'$ is the final idempotent of $\tau_k$. The sequence is \emph{right flanked} if $\tau_1 = \tau_1' \rho_{\lambda}$, where $\lambda$ is the initial idempotent of $\tau_1$. It is \emph{doubly flanked} if it is both left and right flanked. (We do not like things that are doubly flanked.) 
		
		A sequence is \emph{impermissibly flanked} if
		\begin{enumerate}[label = (F\arabic*)]
			\item $(\tau_k, \ldots, \tau_1)$ is left-flanked, and satisfies (S2) but not (S1);
			
			\item $(\tau_k, \ldots, \tau_1)$ is right-flanked, and satisfies (S1) but not (S2);
			
			\item $(\tau_k, \ldots, \tau_1)$ is doubly flanked;
		\end{enumerate}
		The difficulty in each of these cases, is that the cut point would have to be defined counting from the same end on which the sequence is flanked. As we will see later, impermissible flanking always forces $\mu_k(\tau_k, \ldots, \tau_1) = 0$. 
		
		More generally, a sequence $(\tau_k, \ldots, \tau_1)$ with $\ell(\tau_k \ldots, \tau_1) \geq N$ is said to be \emph{allowable} only when
			\begin{enumerate}[label = (\roman*)]
				\item $k \leq N$; 
				
				\item $(\tau_k, \ldots, \: \tau_1)$ is in an allowable sequence of idempotents, meaning in particular, no jumps;
				
				\item The sequence does not stretch too far;
				
				\item The sequence is not impermissibly flanked;
				
				\item Each $\tau_j$ ($1 < j < k$), is of the form 
				\[
					\sigma_{\lambda_j} \rho_{\lambda_{j}} \sigma_{\lambda_j - 1} \cdots \rho_{\lambda_{j - 1} + 1} \sigma_{\lambda_{j - 1} + 1}
				\]
				where we may possibly have $\lambda_j = \lambda_{j - 1} + 1$, i.e. $\tau_j = \sigma_{\lambda_j}$, and we also require
				\begin{align*}
					\tau_1 &= \begin{cases}\sigma_{\lambda_1} \rho_{\lambda_{1}} \sigma_{\lambda_1 - 1} \cdots \rho_{\lambda_{0} + 1} \sigma_{\lambda_{0}} \\ \text{ or} \\ \sigma_{\lambda_1} \rho_{\lambda_{1}} \sigma_{\lambda_1 - 1} \cdots \rho_{\lambda_{0} + 1} \sigma_{\lambda_{0}} \rho_{\lambda_0} \end{cases} \\ \\
					\tau_k &= \begin{cases} \sigma_{\lambda_k} \rho_{\lambda_{1}} \sigma_{\lambda_1 - 1} \cdots \rho_{\lambda_{k} + 1} \sigma_{\lambda_{k}} \\ \text{ or} \\ \rho_{\lambda_k + 1}\sigma_{\lambda_k} \rho_{\lambda_{1}} \sigma_{\lambda_k - 1} \cdots \rho_{\lambda_{k - 1} + 2} \sigma_{\lambda_{k - 1} + 1} \end{cases}
				\end{align*}
				For instance, with $N = 3$, allowable sequences include
				\begin{align*}
					(\sigma_3, \sigma_2, \sigma_1), \: &(\rho_2 \sigma_1, \sigma_3, \sigma_2) \\
					(\sigma_2, \sigma_1 \rho_1 \sigma_3 \rho_3 \sigma_2), &\text{ and }(\rho_1 \sigma_3 \rho_3 \sigma_2 \rho_2 \sigma_1),
				\end{align*}
				but not
				\begin{align*}
					(\sigma_2 \rho_2 \sigma_1), \:&(\rho_1\sigma_3, \sigma_2, \sigma_1\rho_1) \\
					(\sigma_1\rho_3, \sigma_3, \sigma_2), &\text{ or } (\rho_1, \rho_1, \rho_1),
				\end{align*}
				etc. 
			\end{enumerate}
		
		When defining $\mu_k$s, there are two general cases to consider.
		\begin{enumerate}[label = \textbf{Case \arabic*:}]
			\item $\ell(\tau_k, \ldots, \tau_1) < N$. Then we have only $\mu_2$s and $\mu_1$s. In this case, associativity still holds (since there are no higher multiplications for such sequences as these), and we have the following lemma:
	
	\begin{lemma}\label{rad20}
		\begin{enumerate}[label = (\alph*)]
			\item $\de( \sigma_{j} \rho_{j} \sigma_{j -1} \cdots \rho_{i + 1} \sigma_{i}) = 0$ for all $i \leq j$ and where the product alternates between $\rho_{\lambda}$'s and $\sigma_{\lambda}$'s in decreasing order;
			
			\item $\de (\sigma_j \rho_j \cdots \sigma_{i} \rho_i) = \sigma_j \rho_j \cdots \sigma_{i} V_i$ for each $i \leq j$, with alternating product as in (a);
			
			\item $\de (\rho_i \sigma_{i - 1} \cdots \rho_{j + 1} \sigma_j) = V_i \sigma_{i - 1} \cdots \rho_{j + 1} \sigma_j$, for each $j < i$, product alternating as above;
			
			\item $\de ( \rho_j \sigma_{j - 1} \cdots \rho_{i + 1} \sigma_i \rho_i) = V_j \sigma_{j - 1} \cdots \rho_{i + 1} \sigma_i \rho_i + \rho_j \sigma_{j - 1} \cdots \rho_{i + 1} \sigma_i V_i$.
			
			\item $\de (\tau V) = (\de \tau) V$ where $\tau$ is a product of $\rho$'s and $\sigma$'s and $V$ is some product of $\{V_i\}_{i = 1}^N$. 
		\end{enumerate}
		\end{lemma}
		
		\begin{proof}[Proof of Lemma~\ref{rad20}]
			This is by induction on the number of terms in the product (which could be more than $j - i$ since we are counting $\emm N$). 
		\end{proof}
		
		Lemma~\ref{rad20} suffices to verify the $\A_{\infty}$ relation for any tree with inputs $\tau_k, \ldots, \tau_1$ such that $\ell(\tau_k, \ldots, \tau_1) < N$. In the remaining cases, we will verify $\A_{\infty}$-relations for trees labelled with a sequences of length $\geq N$.

			\item $\ell(\tau_k, \ldots, \tau_1) \geq N$ and $k = 1$ -- that is, we are dealing with a differential. This differential, $\de \tau$, is defined to be non-zero precisely in the following cases:
			\begin{enumerate}[label = \textbf{Case 2.\arabic*}]
				\item $\tau$ is unflanked;	
				
					Write $\tau = \sigma_{j} \rho_j \sigma_{j - 1} \cdots \rho_{i + 1} \sigma_i$. Then (cutting from the left and the right, respectively) we can write
					\begin{equation}\label{algFull20}
						\tau = \begin{cases} \underbrace{(\sigma_j \rho_j  \cdots \rho_{j + 2} \sigma_{j + 1})}_{\text{length } N} \cdot \underbrace{\rho_{j + 1} \sigma_{j + 1} \cdots \rho_{i + 1} \sigma_i}_{\text{length } \ell(\tau) - N} & \text{ from the left} \\ \underbrace{\sigma_j \rho_j \cdots \sigma_{i} \rho_i}_{\text{length } \ell(\tau) - N} \cdot \underbrace{(\sigma_{i + N - 1} \rho_{i + N - 1} \cdots \rho_{i + 1} \sigma_i)}_{\text{length } N} & \text{from the right}\end{cases}
					\end{equation}
					Then we define
					\begin{equation}\label{algFull19}
					\de \tau = \underbrace{\prod_{\lambda \neq (j + 1)} V_{\lambda}}_{N \text{ total } V \text{s}} \cdot \underbrace{(\rho_{j + 1} \sigma_{j + 1} \cdots \rho_{i + 1} \sigma_i)}_{\text{length } \ell(\tau) - N} + \underbrace{\prod_{\lambda \neq i} V_{\lambda}}_{N \text{ total } V \text{s}} \cdot \underbrace{(\sigma_j \rho_j \cdots \sigma_{i} \rho_i)}_{\text{length } \ell(\tau) - N}
					\end{equation}

				\item $\tau$ is left-flanked;
				
					Write $\tau = \rho_{j + 1} \sigma_{j} \rho_j \sigma_{j - 1} \cdots \rho_{i + 1} \sigma_i$. Then we cut from the right, and write
					\[
						\tau = \underbrace{\rho_{j + 1} \sigma_j \rho_j \cdots \sigma_{i} \rho_i}_{\text{length } \ell(\tau) - N} \cdot \underbrace{(\sigma_{i + N - 1} \rho_{i + N - 1} \cdots \rho_{i + 1} \sigma_i)}_{\text{length } N}
					\]
					and define 
					\begin{equation}\label{algFull22}
						\de \tau = V_{j + 1} \cdot \underbrace{\sigma_{j} \rho_j \sigma_{j - 1} \cdots \rho_{i + 1} \sigma_i}_{\text{length } \ell(\tau)} + \underbrace{\prod_{\lambda \neq i} V_{\lambda}}_{N \text{ total } V \text{s}} \cdot \underbrace{(\rho_{j + 1} \sigma_j \rho_j \cdots \sigma_{i} \rho_i)}_{\text{length } \ell(\tau) - N}
					\end{equation}
				
				\item $\tau$ is right-flanked;
				
					Write $\tau = \sigma_{j} \rho_j \sigma_{j - 1} \cdots \rho_{i + 1} \sigma_i \rho_i$. Then we cut from the left, and write
					\[
						\tau = \underbrace{(\sigma_j \rho_j  \cdots \rho_{j + 2} \sigma_{j + 1})}_{\text{length } N} \cdot \underbrace{\rho_{j + 1} \sigma_{j + 1} \cdots \rho_{i + 1} \sigma_i \rho_i}_{\text{length } \ell(\tau) - N}
					\]
					and define 
					\begin{equation}\label{algFull23}
						\de \tau = V_i \cdot \underbrace{\sigma_{j} \rho_j \sigma_{j - 1} \cdots \rho_{i + 1} \sigma_i}_{\text{length } \ell(\tau)} + \underbrace{\prod_{\lambda \neq (j + 1)} V_{\lambda}}_{N \text{ total } V \text{s}} \cdot \underbrace{(\rho_{j + 1} \sigma_{j + 1} \cdots \rho_{i + 1} \sigma_i \rho_i)}_{\text{length } \ell(\tau) - N}
					\end{equation}
				
				\item $\tau$ is doubly flanked;
				
					Write $\tau = \rho_{j + 1}\sigma_{j} \rho_j \sigma_{j - 1} \cdots \rho_{i + 1} \sigma_i \rho_i$. In this case, we just define
					\begin{equation}\label{algFull21}
						\de \tau = V_{j + 1} \cdot \sigma_{j} \rho_j \sigma_{j - 1} \cdots \rho_{i + 1} \sigma_i \rho_i + V_i \cdot \rho_{j + 1}\sigma_{j} \rho_j \sigma_{j - 1} \cdots \rho_{i + 1} \sigma_i.
					\end{equation}
			\end{enumerate}
			
			\begin{remark}\label{algFull24}
				\emph{Essentially, in both cases 2.2 and 2.3, the first term of the differential is calculated as in the length $< N$ case, and the second term is calculated analogous to the more general length $\geq N$ case, by at a point length $N$ from the unflanked side differentiating that length-$N$ piece, and leaving the rest fixed.}
						
				\emph{In Case 2.1, there is no analogous nonzero differential from the length $< N$ case, so both terms are by analogy to the length $\geq N$, $k > 1$ cases from below.}
						
				\emph{In Case 2.3, there is no analogous non-zero differential from the length $\geq N$, $k > 1$ cases, so both terms are calculated as they would be if $\ell(\tau) < N$.}
			\end{remark}
			
			\item $\ell(\tau_k, \ldots, \tau_1) \geq N$, $k > 1$. If $(\tau_k, \ldots, \tau_1)$ is not allowable, then $\mu_k(\tau_k, \ldots, \tau_1) = 0$. If it is allowable, then we define $\mu_k(\tau_k, \ldots, \tau_1)$ in the following way.
			\begin{enumerate}[label = \textbf{Case 3.\arabic*:}]
				\item $(\tau_k, \ldots, \tau_1)$ is unflanked and stretches too far from the left, but not the right;
				
				Then we can compute the cut-point counting up from the right, and write
				\begin{equation}\label{algFull10}
					(\tau_k, \ldots, \tau_1) = (\tau_k' | \sigma_{i + N - 1} B_{i + N - 1} \cdots B_{i + 1} \sigma_i)
				\end{equation}
				where $|$ is the cut point, $i$ is the initial idempotent of $\tau_1$, each $B_j = \rho_j$ or denotes a comma, but not both, and $\tau_k'$ is the leftover part of the decomposition of $\tau_k$ to the left of the cut-point. In this case, we define
				\begin{equation}\label{algFull6}
					\mu_k(\tau_k, \ldots, \tau_1) = \tau_k'\cdot V_{N + 1} \prod \{ V_j : B_j = \rho_j\}
				\end{equation}
				For instance, when $N = 5$, we could have
				\begin{center}
					\includegraphics[width = 6cm]{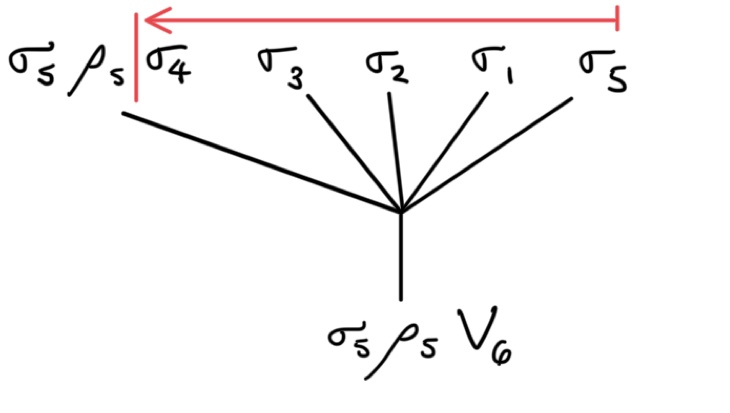}
				\end{center}
				
				\item $(\tau_k, \ldots, \tau_1)$ is unflanked and stretches too far from the right, but not the left;
				 
				 Then the cut-point is defined counting down from the final idempotent of $\tau_k$. This means we can write
				\begin{equation}\label{algFull11}
					(\tau_k, \ldots, \tau_1) = (\sigma_{i} B_i \sigma_{i - 1} B_{i - 1} \cdots B_{i - N + 2} \sigma_{i - N + 1} | \tau_1'),
				\end{equation}
				where $i + 1$ is the final idempotent of $\tau_k$, the $|$ denotes the cut point, each $B_{j}$ is either a comma or $\rho_j$ (but not both), and $\tau_k$' is the leftover part of the decomposition of $\tau_1$ which lies to the right of the cut point. We then define
				\begin{equation}\label{algFull4}
					\mu_k(\tau_k, \ldots, \tau_1) = V_{N + 1} \prod \{ V_j : B_j = \rho_j\} \cdot \tau_1'.
				\end{equation}
				For instance, with $N = 5$, we have
				\begin{center}
					\includegraphics[width = 5cm]{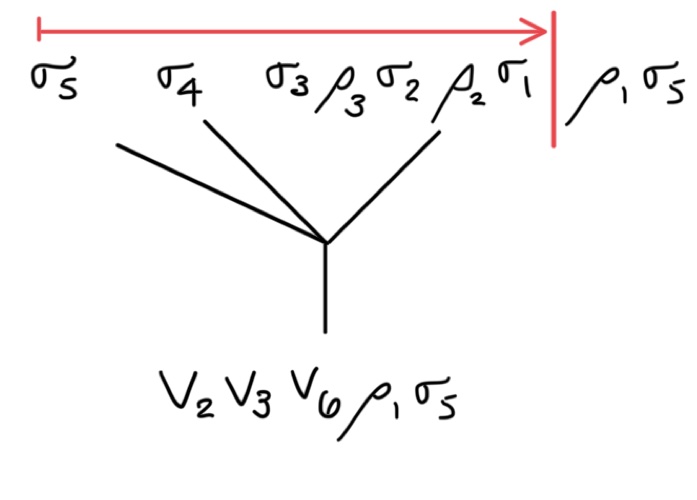}
				\end{center}
				
				\item $(\tau_k, \ldots, \tau_1)$ is unflanked and does not stretch too far in either direction;
				
				Then we can write $(\tau_k, \ldots, \tau_1)$ as in~\eqref{algFull10} or as in~\eqref{algFull11}. In fact, in this case the decompositions given by counting the cut from the left versus from the right differ in a very predictable way:
				
				\begin{lemma}\label{algFull12}
					Suppose $k \geq 1$, $\ell(\tau_k, \ldots, \tau_1) \geq N$, and the sequence is unflanked, as above. Then we can count the cut both from the left and from the right, to get
					\begin{equation}\label{algFull8}
					(\tau_k, \ldots, \tau_1) = \begin{cases} (\sigma_{j - 1} B_{j - 1} \sigma_{j - 2} \cdots B_{j + 1} \sigma_{j} | \tau_1' ) & {\small \text{ (from the left) }} \\ (\tau_k' | \sigma_{i + N - 1} B_{i + N - 1}' \sigma_{i + N - 2} \cdots B_{i + 1}' \sigma_i) & {\small \text{(from the right)}} \end{cases}
				\end{equation}
				where $i < j \emm N$ are the initial and final idempotents of the sequence, respectively. 
				
				Writing the sequence in this way, $B_{\lambda} = B_{\lambda}'$ for each $\lambda$ which appears more than once, and in fact, writing $\Lambda = \sigma_{j - 1} \rho_{j -1} \cdots \rho_{i + 1} \sigma_i$, we have
				\begin{equation}\label{algFull14}
					\tau_1' = \rho_j \Lambda, \quad \tau_k' = \Lambda \rho_i.
				\end{equation}
			\end{lemma}
				
			\begin{remark}\label{algFull25}
				\emph{The statement of Lemma~\ref{algFull12} is true even when $k = 1$ (i.e. we are dealing with a $\mu_1$) not just in the current case. However, in that case, it is patently obvious, so we exclude that from the proof.}
			\end{remark}
					
			\begin{proof}
				First, note that on the overlap of the $\sigma / B$ parts of the two decompositions, the $B_{\lambda}$ and $B_{\lambda}'$ have to match, because these are just decompositions of the same sequence. Hence, if $(\tau_k, \ldots, \tau_1)$ is of length $N$ (i.e the ``overlap'' is the whole sequence) then $\tau_1'$ and $\tau_k'$ are trivial, and the two ways of writing the sequence in~\eqref{algFull8} are actually the same.
				
				 Next, suppose $\ell(\tau_k, \ldots, \tau_1) > N$ (i.e. suppose we have $B_{\lambda}$ other than those in the overlap of the two decompositions from~\eqref{algFull8}). Because the total sequence does not stretch too far in either direction, it is either a $\mu_1$ -- in which case \emph{all} $B_{\lambda}$s and $B_{\lambda}'$s are $\rho_{\lambda}$s -- or it is of length $\leq 2N$. Hence, every idempotent appears at most twice in this sequence. 
				
				Suppose now that $\lambda$ is \emph{not} in the overlap, but that there is both a $B_{\lambda}$ and a $B_{\lambda}'$. This means $\lambda < j \emm N$ and $\lambda > i \emm N$; we put in the exclusive bound because if e.g. $\lambda = j$, then there is a $B_j'$ in the 2nd decomposition of~\eqref{algFull8}, but no $B_j$ the first, because that is one step past the final idempotent of $\tau_k$, and likewise when $\lambda = i$. Because the sequence does not stretch too far in either direction, when $\lambda$ appears on the right, it is part of the right-most term of the sequence (otherwise $(\tau_k, \ldots, \tau_1)$ would satisfy (S1)), and when it appears on the left, it is part of the left-most term (otherwise $(\tau_k, \ldots, \tau_1)$ would satisfy (S2)). In terms of idempotents, this means that:
				\begin{align*}
					i < \lambda &\leq i' \emm N; \\
					j' \leq \lambda &< j \emm N,
				\end{align*}
				where $i'$ denotes the final idempotent of $\tau_1$ and $j'$ denotes the initial idempotent of $\tau_k$, respectively. Moreover, if $\lambda = i'$, then the sequence $(\tau_{k}, \ldots, \tau_2)$ would be of length $\geq N$, i.e. $(\tau_k, \ldots, \tau_1)$ would stretch too far from the left, and likewise, if $\lambda = j'$, then $(\tau_k, \ldots, \tau_1)$ stretches too far from the right. Hence, we actually have
				\begin{equation}\label{algFull9}
					i < \lambda < i' \emm N \text{ and } j'< \lambda < j \emm N.
				\end{equation}
				Now decompose $\tau_1, \tau_k$ into products of $\rho$s and $\sigma$s. Because~\eqref{algFull9} says that $\lambda$ appears in the interior of the decompositions of both $\tau_1$ or $\tau_k$, it must correspond to a $\rho_{\lambda}$ both in $\tau_1$ and in $\tau_k$. Hence, $B_{\lambda} = \rho_{\lambda}$ in the top line of~\eqref{algFull8}, where the $\lambda$ is to the far left, and $B_{\lambda}' = \rho_{\lambda}$ in the bottom line, where the $\lambda$ is to the far right. 
	
				We can now apply the decomposition from~\eqref{algFull8} to get 
				\begin{align*}
					\tau_1' &= B_j \sigma_{j - 1} B_{j - 1} \cdots B_{i + 1} \sigma_i \\
					&= \rho_j \sigma_{j - 1} \rho_{j -1} \cdots \rho_{i + 1} \sigma_i\\
					\tau_k' &= \sigma_{j - 1} B_{j - 1} \cdots B_{i + 1} \sigma_i B_i \\
					&= \sigma_{j - 1} \rho_{j - 1} \cdots \rho_{j + 1} \sigma_i \rho_i
				\end{align*}
				which is precisely~\eqref{algFull14}.
				\end{proof}
			
				With this done, we define
				\begin{equation}\label{algFull13}
				\begin{split}
					\mu_k(\tau_k, \ldots, \tau_1) &= \tau_k' V_{N + 1} \prod\{ V_{\lambda} : j < \lambda < i \emm N, B_{\lambda} = \rho_{\lambda}\} \\ &\qquad +\tau_1' V_{N + 1} \prod \{V_{\lambda} : j < \lambda < i \emm N, B_{\lambda} = \rho_{\lambda} \}
				\end{split}
				\end{equation}
				Using Lemma~\ref{algFull12}, this becomes
				\begin{equation}\label{algFull15}
				\begin{split}
					\mu_k(\tau_k, \ldots, \tau_1) &= \Lambda \rho_i \cdot V_{N + 1} \prod\{ V_{\lambda} : j < \lambda < i \emm N, B_{\lambda} = \rho_{\lambda}\} \\
					&\qquad+ \rho_j \Lambda \cdot V_{N + 1} \prod \{V_{\lambda} : j < \lambda < i \emm N, B_{\lambda} = \rho_{\lambda} \}
				\end{split}
				\end{equation}
				So for instance, with $N = 5$, we could have
				\begin{center}
					\includegraphics[width = 12cm]{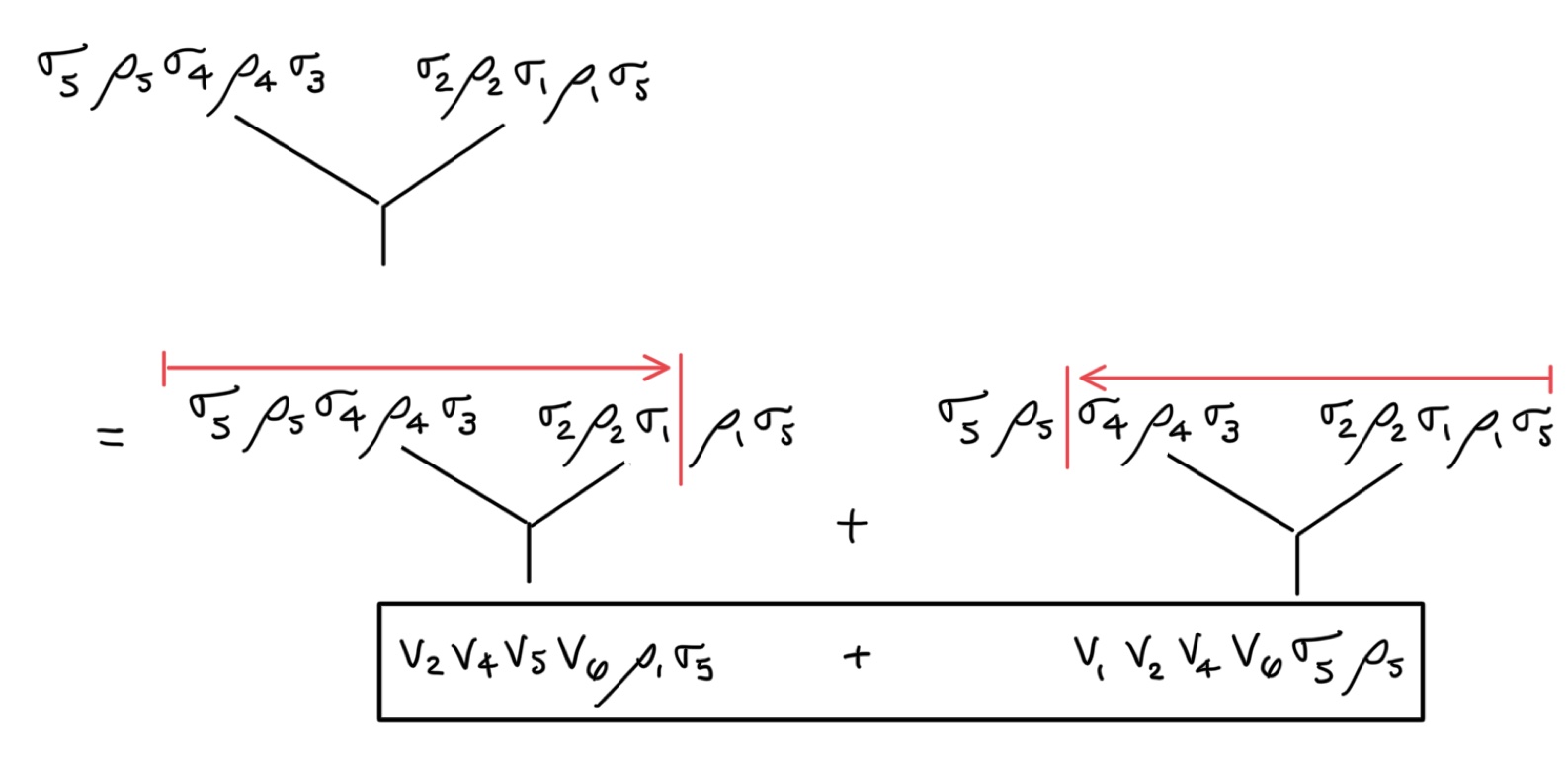}
				\end{center}
				
				\item $(\tau_k, \ldots, \tau_1)$ is left-flanked but of length $N$;
				
				Here, we can technically cut from either direction, but the cut-point is trivial since the sequence is of length $N$. We can just write
				\[
					(\tau_k, \ldots, \tau_1) = (\rho_{i} \sigma_{i + N - 1} B_{i + N - 1} \sigma_{i + N - 2} B_{i + N - 2} \cdots B_{i + 1} \sigma_{i} ),
				\]
				where again, each $B_j$ is either a comma or a $\rho_j$ but not both. Define
				\begin{equation}\label{algFull5}
					\mu_k(\tau_k, \ldots, \tau_1) =\rho_{i} V_{N + 1} \prod \{ V_j : B_j = \rho_j\}.
				\end{equation}
				so for instance, in the case $N = 5$, we could have
				\begin{center}
					\includegraphics[width = 5cm]{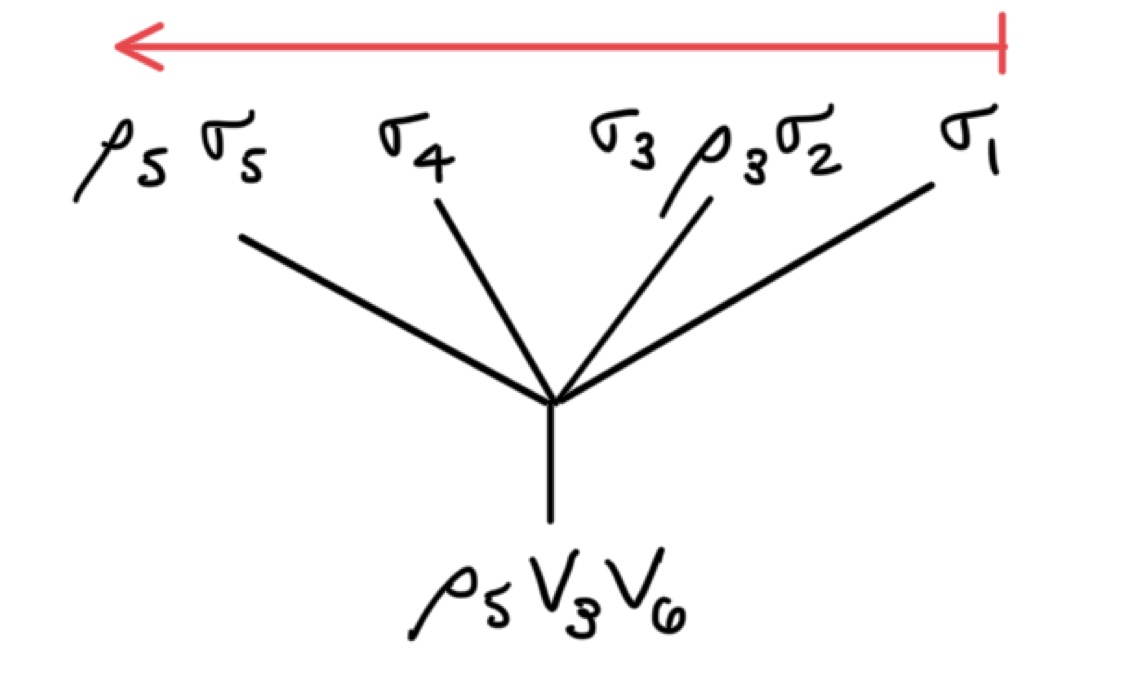}
				\end{center}
				
				\item $(\tau_k, \ldots, \tau_1)$ is right flanked but of length $N$;
				
				As in the previous case, we write
				\[
					(\tau_k, \ldots, \tau_1) =  (\sigma_{i + N - 1} B_{i + N - 1} \sigma_{i + N - 2} B_{i + N - 2} \cdots B_{i + 1} \sigma_{i} \rho_i),
				\]
				and define
				\begin{equation}\label{algFull16}
					\mu_k(\tau_k, \ldots, \tau_1) = \rho_{i} V_{N + 1} \prod \{ V_j : B_j = \rho_j\}
				\end{equation}
				
				\item $(\tau_k, \ldots, \tau_1)$ is left-flanked and of length $> N$;
				
				 Because $(\tau_k, \ldots, \tau_1)$ is not inadmissibly flanked, is not right flanked and it must not stretch too far counting from the right; in this case we write 
				\[
					(\tau_k, \ldots, \tau_1) = (\tau_k' | \sigma_{i + N - 1} B_{i + N - 1} \cdots B_{i + 1} \sigma_i)
				\]
				as in~\eqref{algFull10}, but this time with the understanding that the $\tau_k'$ ends with a $\rho_j$ (where $j$ is the final idempotent of $\tau_k$). We then define $\mu_k(\tau_k, \ldots, \tau_1)$ using~\eqref{algFull6}.
				
				\item $(\tau_k, \ldots, \tau_1)$ is right-flanked and of length $> N$;
				
				Again, because the sequence is not impermissibly flanked, it must not be flanked on the left, and does not stretch too far from that direction, so we can write $(\tau_k, \ldots, \tau_1)$ as in~\eqref{algFull11}, but with the understanding that in this case, $\tau_1'$ begins with a $\rho$, and define
				\[
					\mu_k(\tau_k, \ldots, \tau_1) = V_{N + 1} \prod \{ V_j : B_j = \rho_j\} \cdot \tau_1'.
				\]
				as in~\eqref{algFull4}.
			\end{enumerate}
		\end{enumerate}
		
		The point in Cases 3.4 through 3.7 is that except in the trivial case (when the sequence is of length $N$, so there is no real cut) we never want to define the cut counting down from the same end of the sequence where the sequence is flanked. Suppose the sequence $(\tau_k, \ldots, \tau_1)$ is of length $> N$, and we did define the cut counting from the same direction in which the sequence was flanked. Then (assuming the sequence did not stretch too far in that direction) we would get a cut that looked like
			\[
				(\tau_k, \ldots, \tau_1) = \begin{cases} (\rho_{i + 1} \sigma_i B_i \cdots B_{i - N + 2} \sigma_{i - N + 1} | \tau_1') & {\small \text{ if we had left flanking and counted from the left}} \\ (\tau_k' | \sigma_{i + N - 1} B_{i + N - 1} \cdots B_{i + 1} \sigma_i \rho_i) & {\small \text{ if we had right flanking and counted from the right}} \end{cases}
			\]
			where $\tau_1' = \rho_{i + 1} \omega_1$ and $\tau_k' = \omega_k \rho_i$, for some $\omega_1, \omega_k \in \B$. By analogy with the definitions above, in left flanked case we should have
			\begin{align*}
				\mu_k(\tau_k, \ldots, \tau_1) &= \rho_{i + 1} \cdot V_{N + 1} \prod \{ V_j : B_j = \rho_j\} \cdot \tau_1' \\
				&= \rho_{i + 1} \tau_1' \cdot V_{N + 1} \prod \{ V_j : B_j = \rho_j\} \\
				&{\small \text{(Because all the $V_{\lambda}$ are central)}} \\
				&= \rho_{i + 1} \cdot \rho_{i + 1} \omega_1 \cdot V_{N + 1} \prod \{ V_j : B_j = \rho_j\} \\
				&= 0
			\end{align*}
			because $\rho_{\lambda}$ cannot be multiplied with one another. Likewise, in the right flanked case,
			\[
				\mu_k(\tau_k, \ldots, \tau_1) = \tau_k' \cdot V_{N + 1} \prod \{ V_j : B_j = \rho_j\} \cdot \rho_{i} = 0.
			\]
			We avoid this by never counting the cut from the same direction a sequence is flanked. This is also why we rule out doubly flanked sequences -- because in this case, there is no way to define the cut. 
			
			The next step is to verify the $\A_{\infty}$ relations. Again, the $\A_{\infty}$ relations, as stated in~\eqref{defs1}, also say that for any tree $T$ labelled with elements of $\B$, the sum of all ways to add an edge to $T$ is zero. Adding an edge corresponds either to a differential, or to pulling together $j$ adjacent elements of the sequence $\tau_k, \ldots, \tau_1$ with which the input edges of $T$ are labeled, that is:
			\begin{center}
				\includegraphics[width = 7cm]{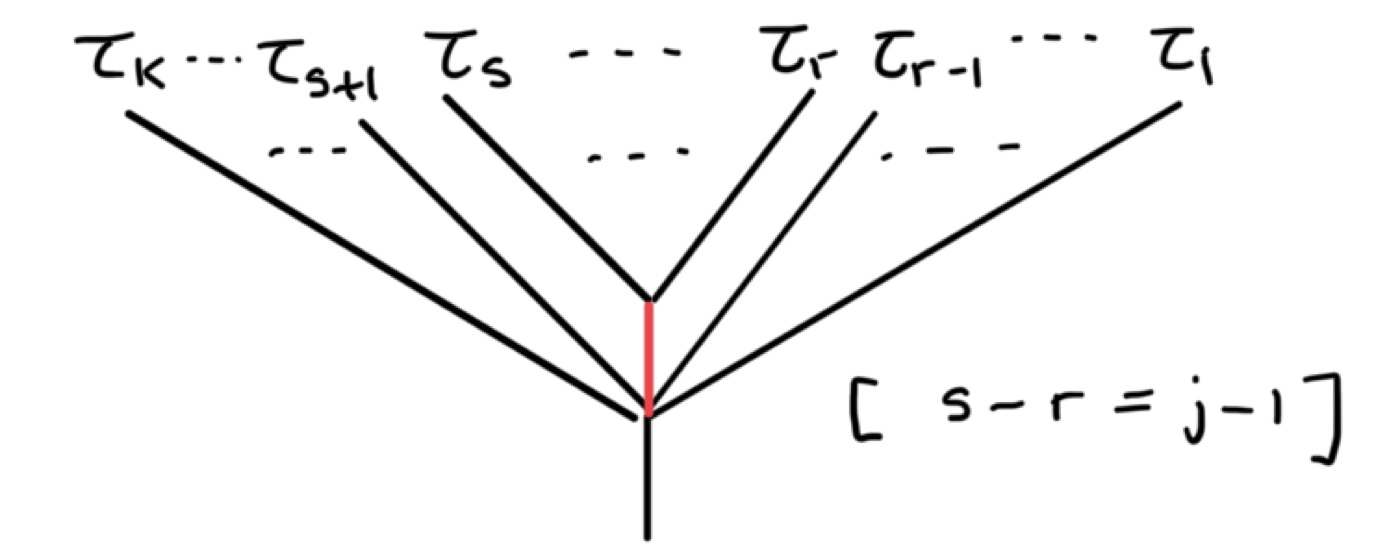}
			\end{center}
			So the $\A_{\infty}$ relation for $T$ counts the number of ways to pull together $j$ elements and get a non-zero composition of algebra operations.
						
			\begin{proposition}\label{algFull17}
				$\B$ satisfies the $\A_{\infty}$ relations, and is therefore a bona fide $\A_{\infty}$-algebra.
			\end{proposition}
			
			\begin{proof}
				Look at a sequence $\taa := (\tau_k, \ldots, \tau_1)$. We have two notions of length. The first is $k$, the actual number of elements in the sequence. The second is $\ell := \ell(\taa)$. If $\ell < N$, then the only operations we have are the usual $\mu_2$ and the $\mu_1$, so the $\A_{\infty}$ relations for $\taa$ are trivial. For the remainder of the proof, we will deal with only $\ell \geq N$. 
				
				If $k = 1$, then we go through Cases 2.1-2.4, and it is clear that all double differentials vanish. 
				
				Suppose now that $k = 2$, $A = \ell (\tau_2)$ and $B = \ell (\tau_1)$. We are assuming with $\ell = A + B \geq N$, so there are the following cases to consider
				\begin{enumerate}[label = \textbf{Case II.\arabic*}]
					\item $A, B < N$: This is just the Leibniz rule, as in the low-length case. 
					
					\item $A \geq N, \: B < N, \: \ell < 2N$; There are 16 cases, corresponding to whether $\tau_1, \tau_2$ are flanked, and whether each is left, right, or doubly flanked. The verification in each case follows directly from the rules for high-length products, above.
					
					\item $A \geq N, \: B < N, \: \ell \geq 2N$; Same situation as Case II.2. The only difference is that since $\ell \geq 2N$, the relations may involve terms of the form
					\begin{center}
						\includegraphics[width = 3cm]{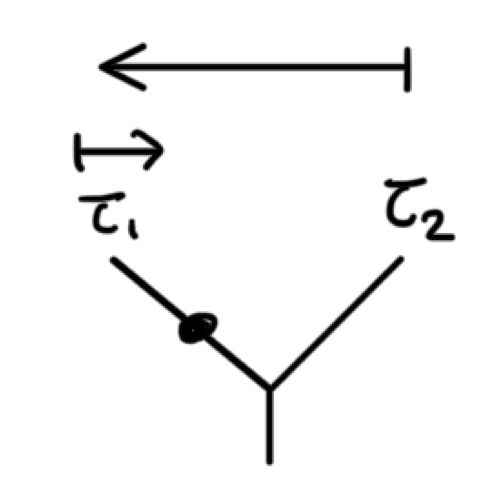}
					\end{center}
					in the case where $\tau_1$ is not right flanked and $\tau_2$ is not left-flanked. Here, again, the arrow denotes the direction in which the cut point is computed if the operation is of length 
					
					\item $A < N, \: B \geq N, \: \ell < 2N$; Analogous to Case II.2;
					
					\item $A < N, \: B \geq N, \: \ell \geq 2N$; Analogous to Case II.3
					
					\item $A, B \geq N$ (so $\ell \geq 2N$ is implied); This is the straight Leibniz rule again, since $\taa$ stretches too far in both directions so high-length products are all trivial;
				\end{enumerate}
			Now suppose $k > 2$. Then there are three cases to consider. The first option is that $\taa$ contains a multipliable pair (in the traditional sense). No high-length multiplication can involve this pair; hence, the $\A_{\infty}$ relation is trivial, unless this is the only multipliable pair, in which case it looks like
			\begin{center}
				\includegraphics[width = 10cm]{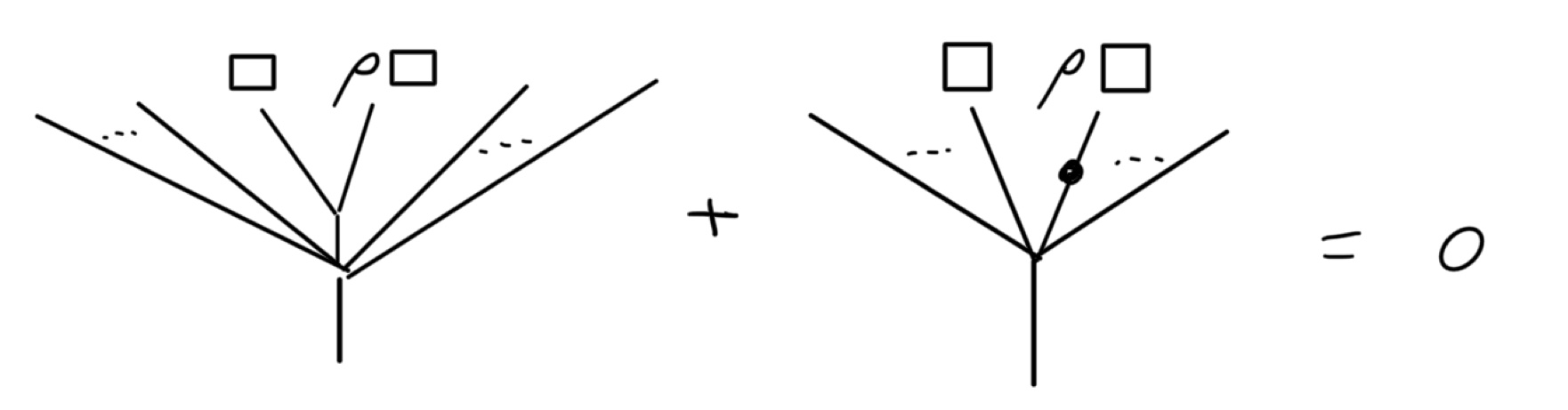}
			\end{center}
			Here, we assumed without loss of generality that the right term of the multipliable pair was the one with an exposed $\rho$; the other case is analogous. 
			
			The second option is that one of the terms in the $\A_{\infty}$ relation has the second operation a $\mu_2$, i.e. one of the terms looks like
			\begin{center}
				\includegraphics[width = 6cm]{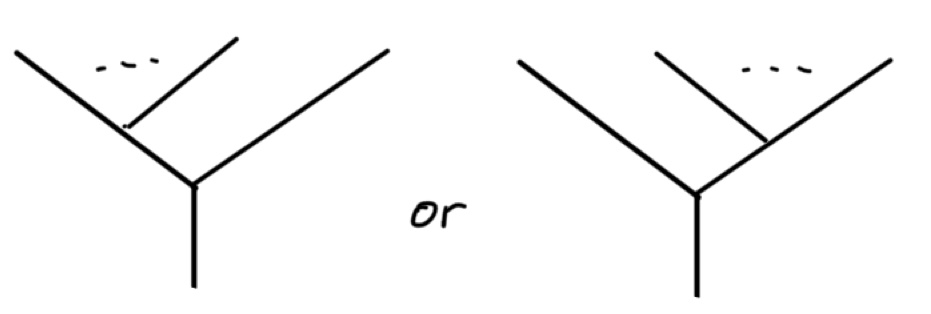}
			\end{center}
			We are going to verify the $\A_{\infty}$ relations for cases in which a term like the left picture appears, since the other situation is completely analogous. In what follows, we will the internal lengths in this diagram as
			\begin{align*}
				A &:= \ell(\tau_k) \\
				B &:= \ell(\tau_{k - 1}, \ldots, \tau_3) \\
				C &:= \ell(\tau_2) \\
				D &:= \ell(\tau_1),
			\end{align*}
			that is,
			\begin{center}
				\includegraphics[width = 8cm]{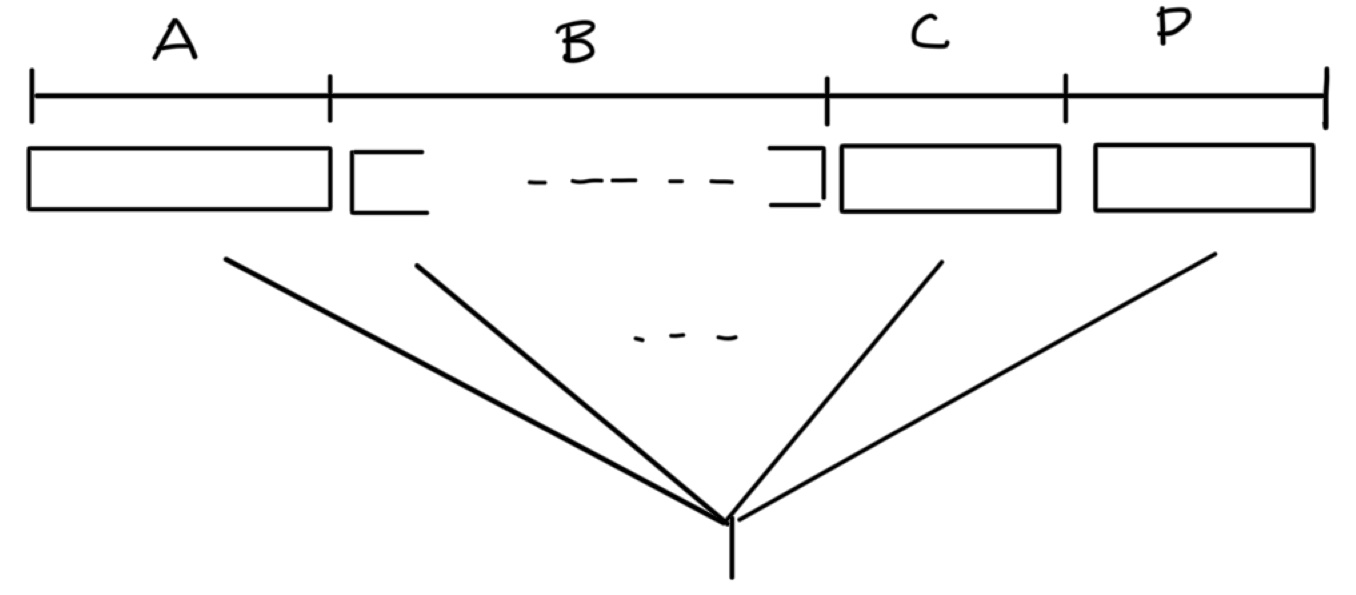}
			\end{center}
			We subdivide the cases according to the cases 3.1-3.7 from the definition of higher multiplications:
			\begin{enumerate}[label = \textbf{Case 3.\arabic*}]
				\item For this case the base assumptions are $A + B \geq N$, $B + C < N$ (so $A + B + C \geq N$ is implied). Within this, we need to subdivide both by the various possibilities for the internal lengths, and by whether $\tau_1$ is unflanked, left / right flanked, or doubly flanked. 
				
				\begin{enumerate}[label = (3.1.\arabic*)]
					 
					\item $\tau_1$ unflanked;
					\begin{enumerate}[label = (3.1.1.\arabic*)]
						\item $A, D < N, \: B + C + D < N$ ($\ell < 2N$ implied);
						\begin{center}
							\includegraphics[width = 6cm]{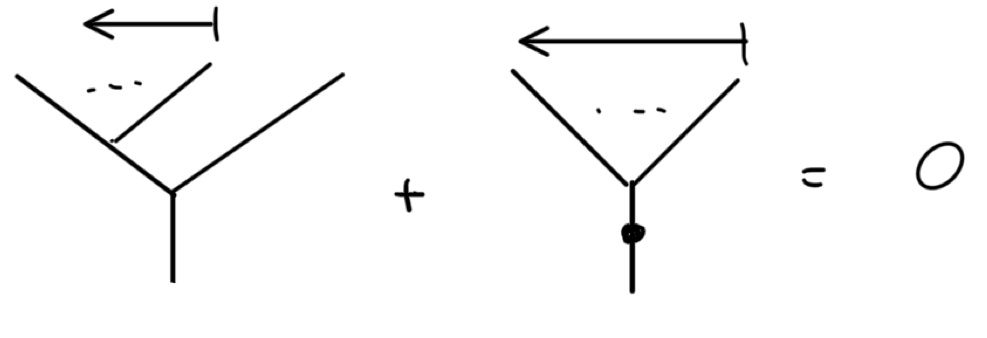}
						\end{center}
						
						\item $A, D < N, \: B + C + D \geq N, \: \ell < 2N$;
						\begin{center}
							\includegraphics[width = 6cm]{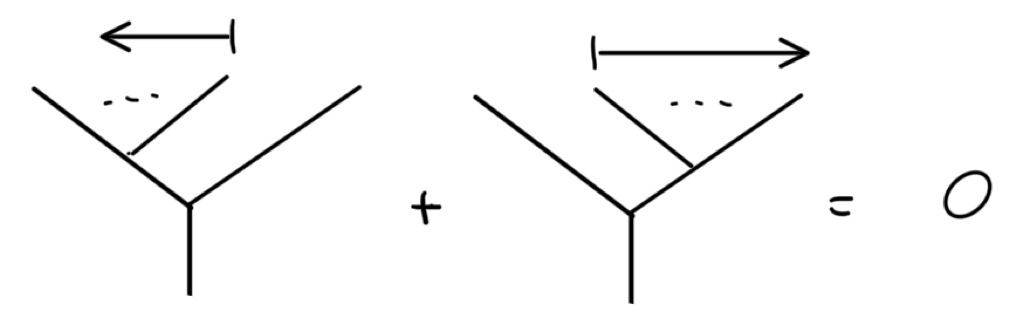}
						\end{center}
						
						\item $A, D < N, \: B + C + D \geq N, \: \ell \geq 2N$; 
						\begin{center}
							\includegraphics[width = 8cm]{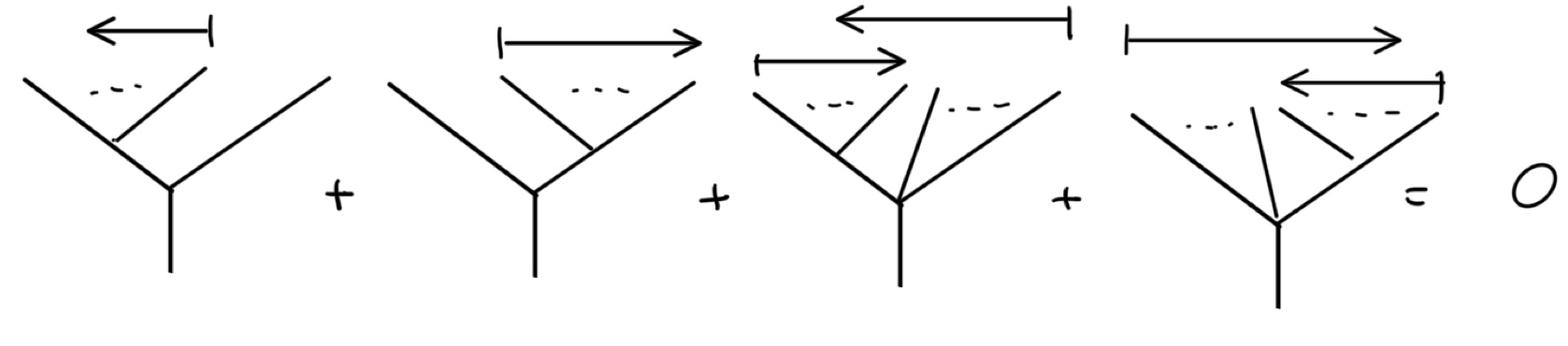}
						\end{center}
						
						\item $A < N, \: D \geq N$ ($B + C + D \geq N$ and $\ell \geq 2N$ implied);
						\begin{center}
							\includegraphics[width = 6cm]{algFull30}
						\end{center}
						
						\item $A \geq N, \: D < N, \: B + C + D < N, \: \ell < 2N$;
						\begin{center}
							\includegraphics[width = 6cm]{algFull29}
						\end{center}
						
						\item $A \geq N, \: D < N, \: B + C + D < N, \: \ell \geq 2N$;
						\begin{center}
							\includegraphics[width = 8cm]{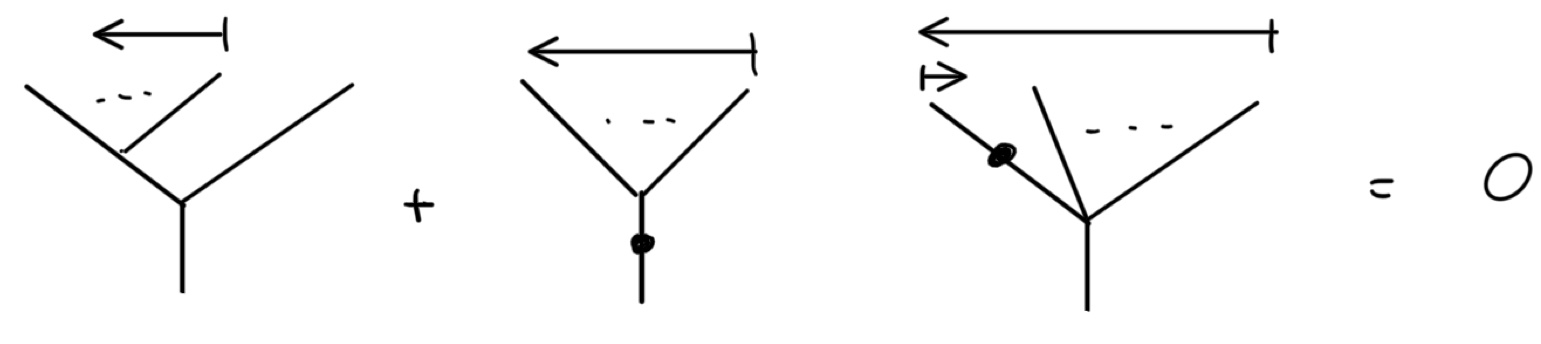}
						\end{center}
						
						\item $A \geq N, \: D < N, \: B + C + D \geq N$ ($\ell \geq 2N$ implied);
						\begin{center}
							\includegraphics[width = 6cm]{algFull30}
						\end{center}
						
						\item $A, D \geq N$ ($B + C + D \geq N$ and $\ell \geq 2N$ implied);
						\begin{center}
							\includegraphics[width = 6cm]{algFull30}
						\end{center}
					\end{enumerate}
					
					\item $\tau_1$ right flanked, i.e. the tree looks like
					\begin{center}
						\includegraphics[width = 4cm]{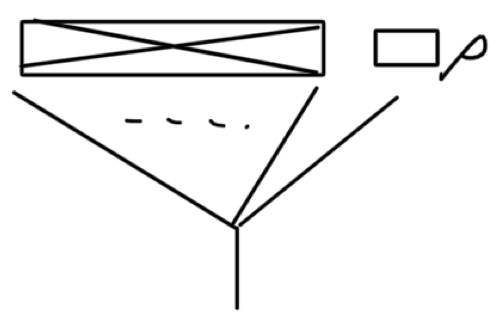}
					\end{center}
					\begin{enumerate}[label = (3.1.2.\arabic*)]
						\item $A, D < N, \: B + C + D < N$ ($\ell < 2N$ implied);
						\begin{center}
							\includegraphics[width = 7cm]{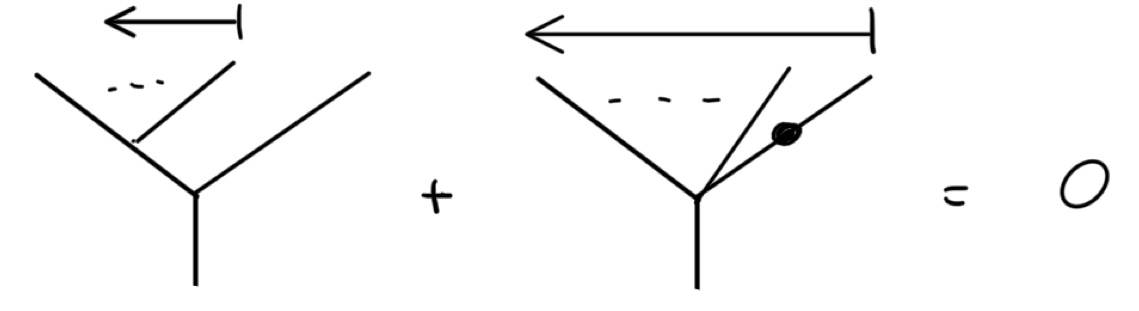}
						\end{center}
						
						\item $A, D < N, \: B + C + D \geq N, \: \ell < 2N$; Same relation as (3.1.1.2);
						
						\item $A, D < N, \: B + C + D \geq N, \: \ell \geq 2N$; Same relation as (3.1.1.2) -- the other two terms from (3.1.1.3) drop out because of flanking;
						
						\item $A < N, \: D \geq N$ ($B + C + D \geq N$ and $\ell \geq 2N$ implied); Same relation as (3.1.1.2)
						
						\item $A \geq N, \: D < N, \: B + C + D < N, \: \ell < 2N$; Same as (3.1.2.1);
						
						\item $A \geq N, \: D < N, \: B + C + D < N, \: \ell \geq 2N$; Same as (3.1.2.1);
						
						\item $A \geq N, \: D < N, \: B + C + D \geq N$ ($\ell \geq 2N$ implied); Same as (3.1.1.2);
						
						\item $A, D \geq N$ ($B + C + D \geq N$ and $\ell \geq 2N$ implied); Same as (3.1.1.2);
					\end{enumerate}
					
					\item $\tau_1$ left flanked, i.e. the tree looks like
					\begin{center}
						\includegraphics[width = 4cm]{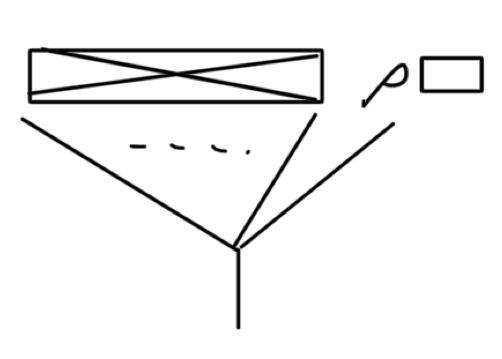}
					\end{center}
					There are non-zero terms in the $\A_{\infty}$ relation if and only if $B + C + D < N$, in which case the two terms are
					\begin{center}
						\includegraphics[width = 8cm]{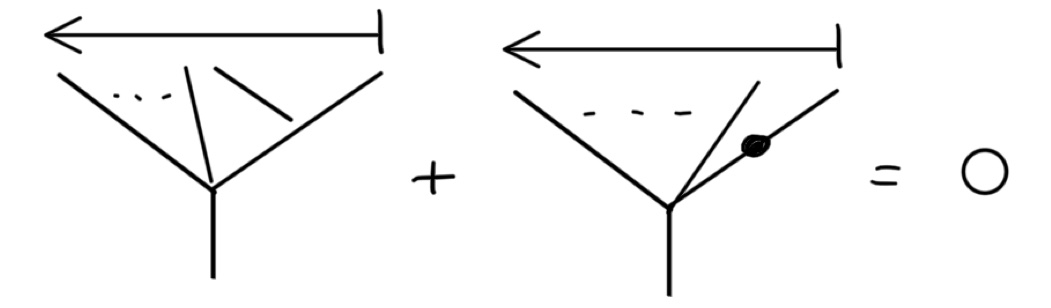}
					\end{center}
					
					\item $\tau_1$ is doubly flanked, i.e. the tree looks like
					\begin{center}
						\includegraphics[width = 4cm]{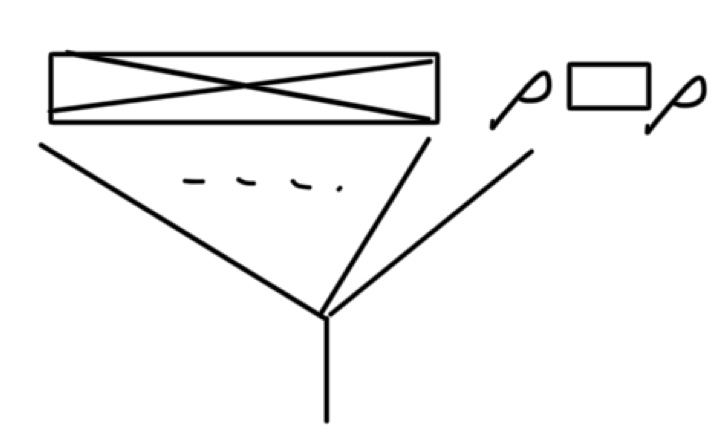}
					\end{center}
					In this case, if we evaluated the inner high-length product from the right, i.e. 
					\begin{center}
						\includegraphics[width = 3cm]{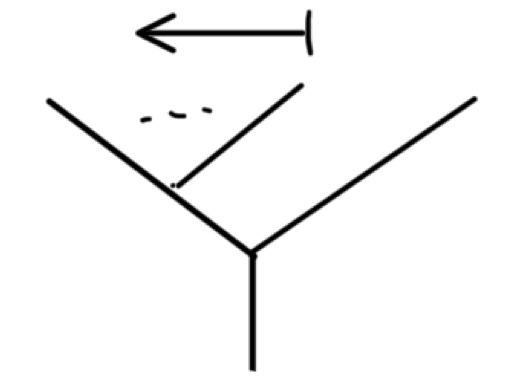}
					\end{center}
					Then we would get adjacent $\rho$s; hence, we would have to evaluate from the left. But this is impossible since $A+ B \geq N$. Hence, the $\A_{\infty}$ relation for a string of this form does not have a term with $\mu_2$ as the external operation, and will be dealt with in the last part of this proof.
				\end{enumerate}

				\item For this case the base assumptions are $A + B < N, \: B + C \geq N$ (so again, $A + B + C \geq N$ is implied). Also, the inner multiplication will be evaluated from the left, the two elements multiplied by the outer $\mu_2$ will not be multipliable (in the usual sense). Hence, the string $\mu(\tau_k, \ldots, \tau_2), \tau_1$ must be of length $\geq N$, which gives $\ell \geq 2N$. It follows that $C + D \geq N$, and also $D < N$ since the outer (high length) $\mu_2$ will be evaluated from the right. Subdividing as above:
				\begin{enumerate}[label = (3.2.\arabic*)]
					\item $\tau_1$ unflanked; The only thing we can control is whether $C < N$ or $\geq N$. But this does not actually affect the diagrams that appear in the $\A_{\infty}$ relation because even when $C \geq N$, the (one or two) term(s) of the differential do not give non-zero compositions. Hence in both cases, the $\A_{\infty}$ relation is
					\begin{center}
						\includegraphics[width = 7cm]{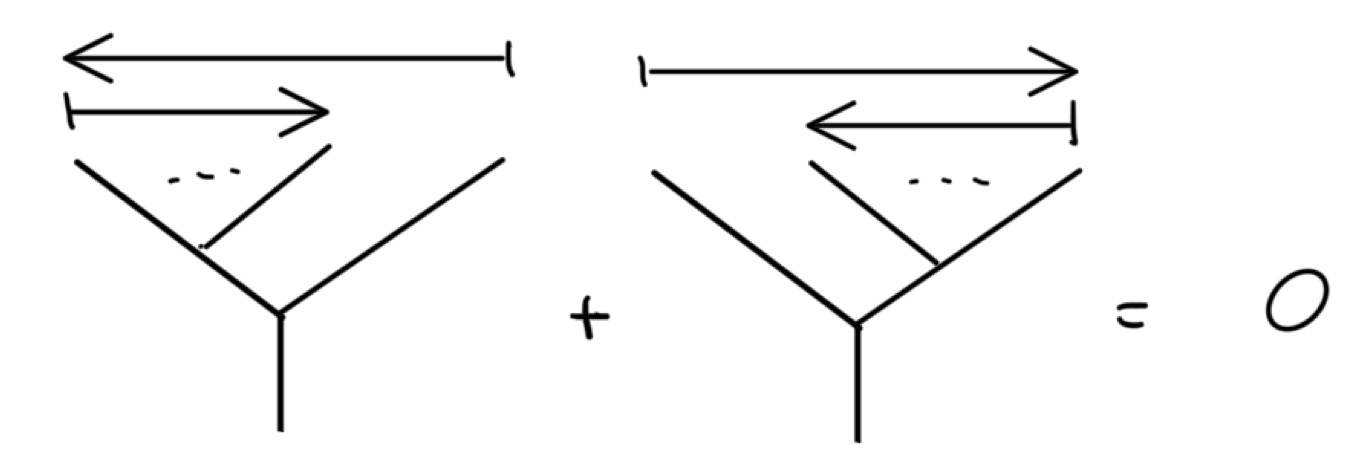}
					\end{center}
					
					\item $\tau_1$ right flanked; In this case, since the inner multiplication is evaluated from the left, the resulting sequence will be doubly flanked. Since it is also of length $\geq N$, it does not contribute, and this case is trivial.
					
					\item $\tau_1$ left flanked; Here, there are the same two subcases as in Case (3.2.1), but here both of those trees vanish for flanking reasons (because the outer multiplication is still high length).
					
					\item $\tau_1$ doubly flanked; Again, both contributing terms vanish in all subcases.
				\end{enumerate}
				
				\item For this case, the base assumptions are $A + B, B + C < N$ and $A + B + C \geq N$.
				\begin{enumerate}[label = (3.3.\arabic*)]
					\item $\tau_1$ unflanked;
					\begin{enumerate}[label = (3.3.1.\arabic*)]
						\item $D < N, \: B + C + D < N$ ($\ell < 2N$ implied);
						\begin{center}
							\includegraphics[width = 6cm]{algFull29}
						\end{center}
						
						\item $D < N, \: B + C + D \geq N, \: \ell < 2N$;
						\begin{center}
							\includegraphics[width = 6cm]{algFull30}
						\end{center}
						
						\item $D < N, \: B + C + D \geq N, \: \ell \geq 2N$; 
						\begin{center}
							\includegraphics[width = 10cm]{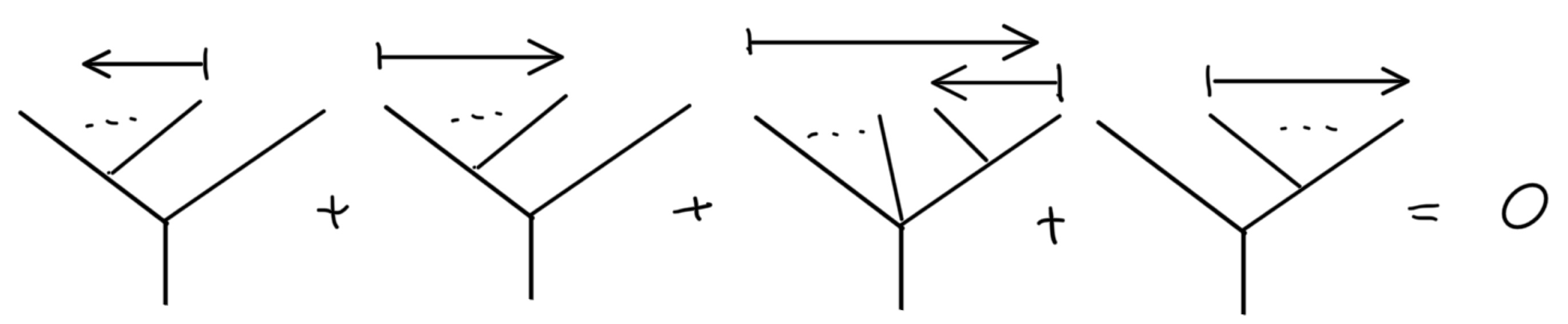}
						\end{center}
						
						\item $D \geq N$ ($B + C + D \geq N$ and $\ell \geq 2N$ implied);
						\begin{center}
							\includegraphics[width = 6cm]{algFull30}
						\end{center}
					\end{enumerate}
					
					\item $\tau_1$ right flanked;
					\begin{enumerate}[label = (3.3.2.\arabic*)]
						\item $D < N, \: B + C + D < N$ ($\ell < 2N$ implied);
						\begin{center}
							\includegraphics[width = 6cm]{algFull38}
						\end{center}
						
						\item $D < N, \: B + C + D \geq N, \: \ell < 2N$; Same as (3.1.1.2);
						
						\item $D < N, \: B + C + D \geq N, \: \ell \geq 2N$; Same as (3.1.1.2);
						
						\item $D \geq N$ ($B + C + D \geq N$ and $\ell \geq 2N$ implied); Same as (3.1.1.2);
					\end{enumerate}
					
					\item $\tau_1$ left flanked; 
					\begin{enumerate}[label = (3.3.3.\arabic*)]
						\item $B + C + D < N$;
						\begin{center}
							\includegraphics[width = 8cm]{algFull38}
						\end{center}
						
						\item $B + C + D \geq N$; Same as (3.1.1.2)
					\end{enumerate}	
					
					\item $\tau_1$ doubly flanked; Because the overall string is right flanked, we can never evaluate the outer multiplication from the right. This means that in both the $B+ C < N$ and $\geq N$ cases, the picture is precisely the one from Case (3.3.3.2), above.
				\end{enumerate}
				 
				\item Here, since $(\tau_k, \ldots, \: \tau_2)$ is left-flanked and length $N$, the situation is much simpler. The base assumption is $A + B + C = N$. 
				\begin{enumerate}[label = (3.4.\arabic*)]
					\item $\tau_1$ unflanked; There are three cases:
					\begin{enumerate}[label = (3.4.1.\arabic*)]
						\item $D < N, \: B + C + D < N$; Same as (3.1.1.1);
						
						\item $D < N, \: B + C + D \geq N$; Same as (3.1.1.2);
						
						\item $D \geq N$; Same as (3.1.1.2);
					\end{enumerate}
					
					\item $\tau_1$ right flanked -- same single case as in (3.4.1);
					
					\item $\tau_1$ left flanked -- vanishes;
					
					\item $\tau_1$ doubly flanked -- vanishes;
				\end{enumerate}
				
				\item Here, $(\tau_k, \ldots, \: \tau_2)$ is right-flanked and the base assumption is again $A + B + C = N$. 
				\begin{enumerate}[label = (3.5.\arabic*)]
					\item $\tau_1$ unflanked; There are no cases, just
					\begin{center}
						\includegraphics[width = 8cm]{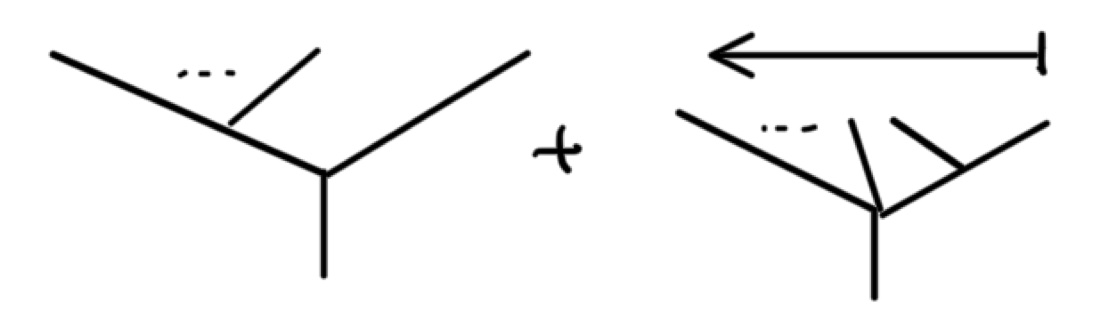}
					\end{center}
					
					\item $\tau_1$ right flanked -- same single case as in (3.5.1);
					
					\item $\tau_1$ left flanked -- all terms vanish because too much flanking;
					
					\item $\tau_1$ doubly flanked -- vanishes; 
				\end{enumerate}
				
				\item Here, $(\tau_k, \ldots, \: \tau_2)$ is left flanked and the base assumptions are $A + B + C \geq N, \: B+ C < N$, or else there is no contributing term of the form we want. 
				\begin{enumerate}[label = (3.6.\arabic*)]
					\item $\tau_1$ unflanked; We have the same 8 subcases as in Case (3.1.1):
					\begin{enumerate}[label = (3.6.1.\arabic*)]
						\item $A, D < N, \: B + C + D < N$ ($\ell < 2N$ implied);
						\begin{center}
							\includegraphics[width = 8cm]{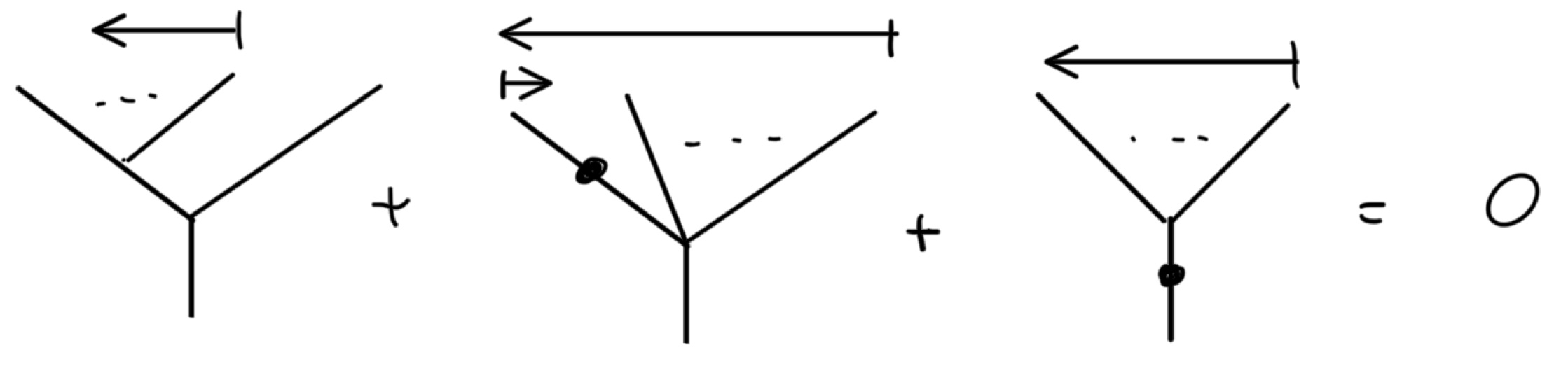}
						\end{center}
						
						\item $A, D < N, \: B + C + D \geq N, \: \ell < 2N$; Same as (3.1.1.2);
						
						\item $A, D < N, \: B + C + D \geq N, \: \ell \geq 2N$; Same as (3.1.1.2);
						
						\item $A < N, \: D \geq N$ ($B + C + D \geq N$ and $\ell \geq 2N$ implied); Same as (3.1.1.2);
						
						\item $A \geq N, \: D < N, \: B + C + D < N, \: \ell < 2N$; Same as (3.6.1.1);
						
						\item $A \geq N, \: D < N, \: B + C + D < N, \: \ell \geq 2N$; Same as (3.6.1.1);
						
						\item $A \geq N, \: D < N, \: B + C + D \geq N$ ($\ell \geq 2N$ implied); Same as (3.1.1.2);
						
						\item $A, D \geq N$ ($B + C + D \geq N$ and $\ell \geq 2N$ implied); Same as (3.1.1.2);
					\end{enumerate}
					
					\item $\tau_1$ right flanked: then there are only two cases:
					\begin{enumerate}[label = (3.6.2.\arabic*)]
						\item $B + C + D < N$;
						\begin{center}
							\includegraphics[width = 8cm]{algFull37}
						\end{center}
						
						\item $B + C + D \geq N$; all terms vanish and teh relation is trivial;
					\end{enumerate}
					
					\item $\tau_1$ left flanked -- No contributing term with $\mu_2$ as the external operation;
					
					\item $\tau_1$ doubly flanked -- same;
				\end{enumerate}
				
				\item Here, $(\tau_k, \ldots, \: \tau_2)$ is right flanked and the base assumptions are again $A + B + C \geq N, \: A + B < N$, or else there is no contributing term of the form we want. 
				\begin{enumerate}[label = (3.7.\arabic*)]
					\item $\tau_1$ unflanked; The only possible cases are
					\begin{enumerate}[label = (3.7.1.\arabic*)]
						\item $B + C + D < N$; Same as (3.3.3.1);
						
						\item $B + C + D \geq N$; Same as (3.1.1.2);
					\end{enumerate}
					
					\item $\tau_1$ right flanked -- no cases, just the same  relation as in Case (3.1.1.2);
					
					\item $\tau_1$ left flanked -- all terms vanish because too much flanking;
					
					\item $\tau_1$ doubly flanked -- same as (3.7.3); 
				\end{enumerate}
			\end{enumerate}
			
			Now what is left is the case in which $k \geq 3$, $\taa$ contains no multipliable pairs, and no term of the $\A_{\infty}$ relation for $\taa$ has a $\mu_2$ as the second operation of the composition. This means that for each non-vanishing term of the $\A_{\infty}$ relation, both operations in the composition (internal and external) must be high length. Therefore, if the first operation was $\mu_j(\tau_s, \ldots, \tau_r)$, with $1 < r \leq s < k$ (i.e. was not based at either the far left or the far right of the tree) then we would end up with a multipliable pair in the resulting sequence $(\tau_k, \ldots, \tau_{s + 1}, \tau', \tau_{r - 1}, \ldots, \tau_1)$. This cannot happen, since the outer operation is at least a $\mu_3$. 
			
			Hence, we can now restrict ourselves to the case when $r = 1$ or $s = k$ or both. If $r = 1, \: s = k$, then this is just the statement that $\de \circ \mu_j = 0$ for any high length operation $\mu_j$. This follows easily from the definitions of the high length $\mu_j$. Now suppose $r > 1, \: s = k$. 
			
			Because we are also that the composition we are looking at does not  have a $\mu_2$ as the outer operation, the picture is
			\begin{center}
				\includegraphics[width = 5cm]{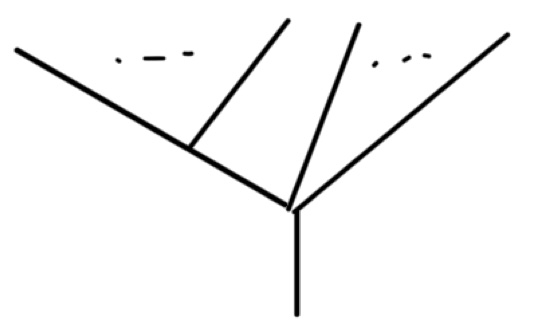}
			\end{center}
			that is, $s = k, \: r \geq 3$. In particular, this means that the given operation (if it is non-vanishing) must be evaluated in the directions
			\begin{center}
				\includegraphics[width = 5cm]{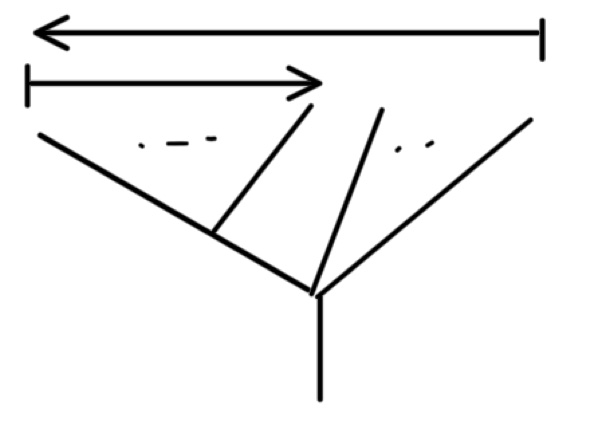}
			\end{center}
			Writing the internal lengths as
			\begin{align*}
				A' &:= \ell(\tau_k, \ldots, \tau_{r + 1}) \\
				B' &:= \ell(\tau_r) \\
				C & := \ell(\tau_{r - 1}, \ldots, \tau_1),
			\end{align*}
			that is,
			\begin{center}
				\includegraphics[width = 5cm]{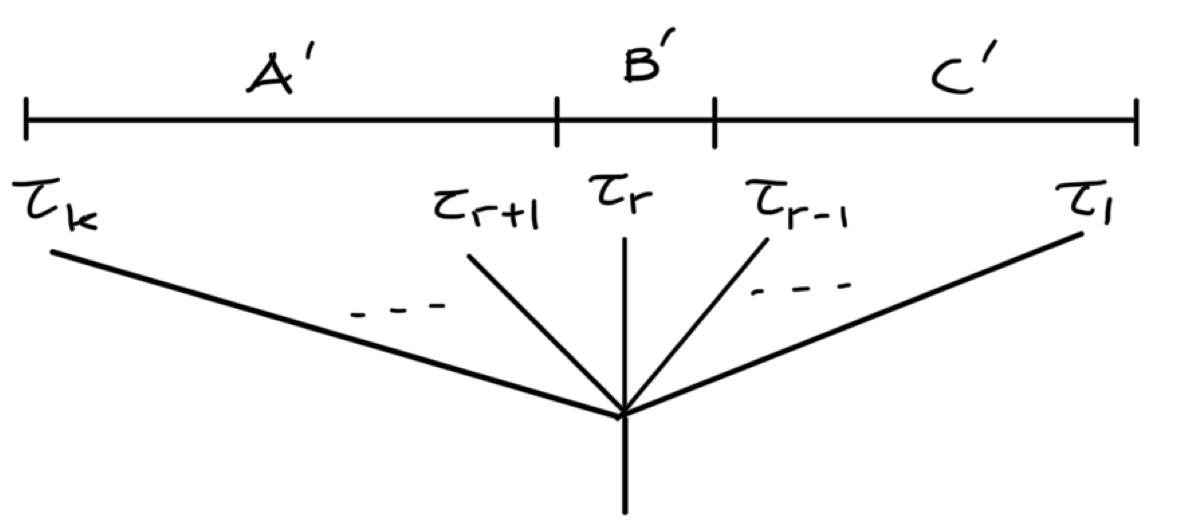}
			\end{center}
			the composition we are looking at is non-vanishing if and only if $A' < N, \: A' + B' \geq N,$ and $C' < N$. This means that the $\A_{\infty}$ relation is precisely
			\begin{center}
				\includegraphics[width = 8cm]{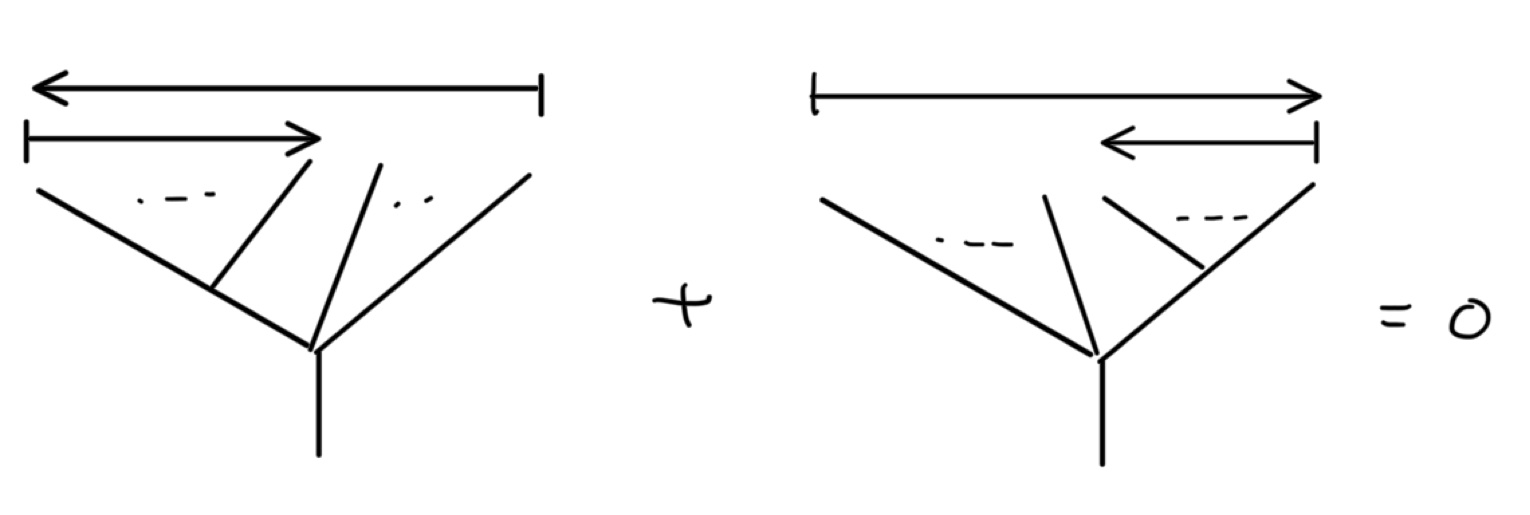}
			\end{center}
			The other case ($r = 1$, $s < k$) is entirely analogous, and will be omitted. Thus the proof is complete, and $\B$ is a bona fide $\A_{\infty}$ algebra.
		\end{proof}
		
		\subsection{Maslov and Alexander gradings for $\B$}
		We define
		\begin{align*}
			m(\rho_i) = m(\sigma_i) &= -1 \text{ for each } i \\
			m(V_i) &= - 2 \text{ for each } 1 \leq i \leq N + 1 \\
			m(\e_{N + 1}) &= 1 \\
			m(U_0) &= - 2N
		\end{align*}
		(This last one is technically not a definition, but a computation from the definition of $U_0$ and the rules for the $\rho_i$ and $\sigma_i$.) We know from Section~\ref{aaMas} that also
		\begin{align*}
			m(V_0) &= 2N - 2 \\
			m(\e_i) &= 2 \text{ for each } 1 \leq i \leq N + 1 \\
			m(\e_0) &= - (2N - 2)
		\end{align*}
		Note that all of the new definitions still satisfy~\eqref{aa1}. That all the higher operations satisfy~\eqref{aa1} follows from the conditions (on length, etc.) which we placed on the operations when we defined them.
		
		For the Alexander gradings, we have
		\begin{align*}
			A(\rho_i) &= A(U_i) = \overline{2i - 1} \text{ for each } 1 \leq i \leq N \\
			A(\sigma_i) &= A(s_i) = \overline{2i} \text{ for each } 1 \leq i \leq N \\
			A(\e_0) &= \sum_{k = 1}^{2N} \overline{k} \\
			A(V_i) &= A(U_i) = \overline{2i - 1} \text{ for each } 1 \leq i \leq N \\
			A(V_{N + 1}) &= A(\e_{N + 1}) = \sum_{k = 1}^N \overline{2k}
		\end{align*}
		Again, it is clear that all operations on $\B$ preserve the Alexander grading. This therefore defines a suitable $\A_{\infty}$-algebra with operationally bounded higher multiplications, satisfying all the properties we need it to.

	\section{The AA bimodule}\label{aabim}
	
		The goal of this section is to define the AA bimodule $\:_{\B} Y_{\A}$ which, with the $DD$-bimodule $\:^{\A} X^{\B}$ defined in Section~\ref{DDbim}, will fit into the duality relation of Theorem~\ref{duality}.
		
			$Y$ is defined to be a $\F[V_0, \ldots, V_{N + 1}]$-module, generated by the intersections of $\alpha$ and $\beta$ arcs in the diagram
			\begin{figure}[H]
				\includegraphics[width = 5cm]{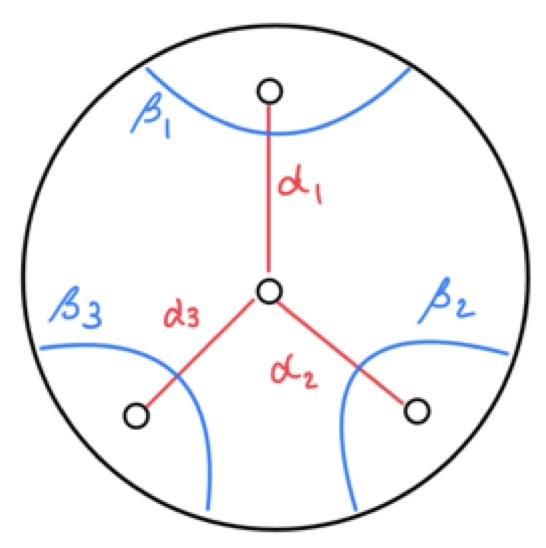}
				
				\caption{}
			\end{figure}
			
			\begin{wrapfigure}{r}{2cm}
				\includegraphics[width = 2cm]{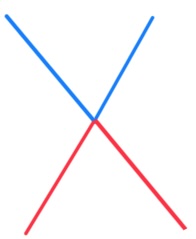}
				
				\caption{}\label{bim1'}
			\end{wrapfigure}
			
			We notate these as $\x = \{i\}$ for $1 \leq i \leq N$, where $\x = \{i\}$ is in the $i$-th idempotent. 
			
			Recall from Section~\ref{aabimdef} that $\A_{\infty}$ operations on an AA bimodule $\:_{\B} Y_{\A}$ are in one to one correspondence with corollas with a designated central leaf, and labelled to the right of this leaf with $\A$-elements and on the left with $\B$-elements. See Figure~\ref{bimdefpic}. This is the basic definition. In this section we are going to manually define the non-zero operations, and show that they can be written in terms of this definition.
			
			We do that as follows. In order to define a valid AA bimodule that satisfies the $\A_{\infty}$-relations, we will use the following steps:
			\begin{enumerate}
				\item Define the set of non-zero operations manually in a way specific to the case we are dealing with -- \emph{not} in terms of corollas, for this step;
				
				\item Show that each operation, as manually defined in (1), can be expressed as a labelled corolla;
				
				\item Exhibit a condition (Proposition~\ref{bim10}) which a corolla $T$ must satisfy in order to correspond with a non-zero operation, as defined in (1);
				
				\item Show that the AA bimodule with generators as above, nad non-zero operations -- now written, properly, in terms of corollas -- as in Proposition~\ref{bim10} satisfies the $\A_{\infty}$ relations.
			\end{enumerate}
			
			For step (1), we define non-zero bimodule operations to be certain allowable types of planar graphs in the disk. A non-zero operation, often called just an \emph{operation} is defined as follows. Start with a 4-valent planar graph embedded in $\D$ which admits a coherent coloring in two colors -- in our case, red and blue, with red bounding $\A$-sectors and blue bounding $\B$-sectors -- such that each vertex looks like Figure~\ref{bim1'}. We require that each red sector is labelled with an $\A$-algebra element, each blue sector is labelleed with a $\B$-algebra element, and each mixed sector is labelled with a number $1 \leq i \leq N$, chosen so that the idempotent along the incoming (i.e. left) red / blue strand matches the initial idempotent of the $\A$- or $\B$-label in the all-red or all-blue sector, and the idempotent along the outgoing (i.e. right) red / blue strand matches the final idempotent of the $\A$- or $\B$-label in this sector. With these conditions, the three allowable types of labelling are
			\begin{figure}[H]
				\includegraphics[width = 12cm]{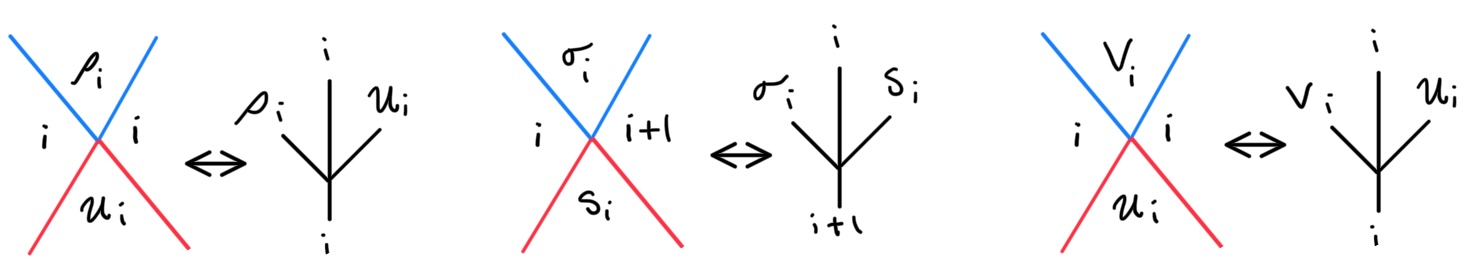}
				
				\caption{}\label{basic}
			\end{figure} 
			We define each basic graph in Figure~\ref{basic} to correspond to a corolla as notated above. These are the only three types of allowable basic graphs. A 4-valent graph, colored and labeled as described in the previous paragraph, is an operation if and only if it satisfies the conditions below. Where possible, we specify the tree (or partial tree) corresponding to each conditions. We will say 
			\begin{enumerate}[label = Op. \arabic*]
				\item:\label{op4} The graph does not contain \emph{composite} regions -- i.e. regions that are bordered by a boundary edge, red edge, and blue edge, and more than one of at least one of these kinds; for instance:
				\begin{figure}[H]
					\includegraphics[width = 6cm]{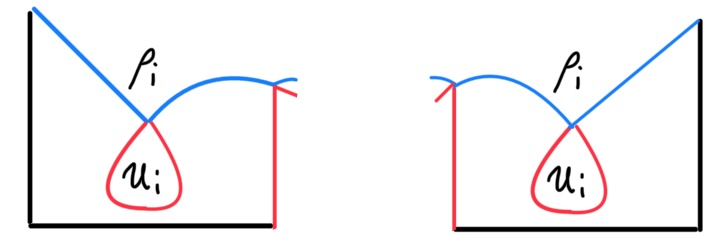}
					
					\caption{}
				\end{figure}\vspace{-10pt}
				where the composite regions are the unlabeled ones. The difficulty with these regions is that the extra edges can be pushed out the boundary -- that is, graphs containing composite regions do not correspond to operations.
				
				This has two notable consequences:
				\begin{enumerate}[label = (\alph*)]
					\item All weights must be isolated away from the boundary, so, for petal weights, as in
				\begin{center}
					\includegraphics[width = 6cm]{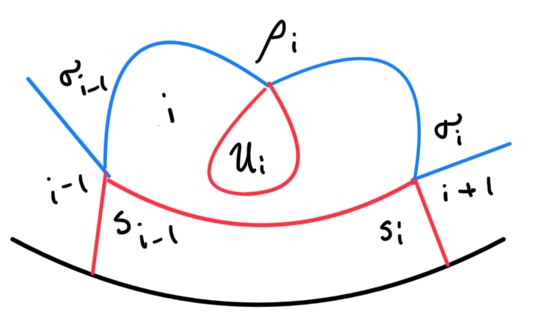}
				\end{center}
				and not like
				\begin{center}
					\includegraphics[width = 6cm]{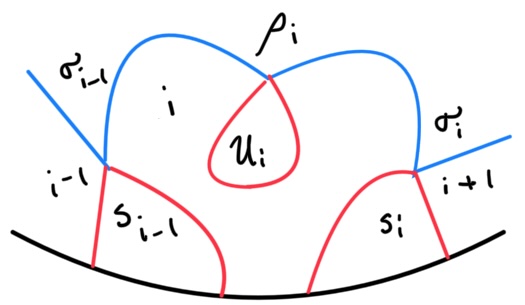}
				\end{center}
				Likewise for central weights, we allow graphs such as
				\begin{center}
					\includegraphics[width = 6cm]{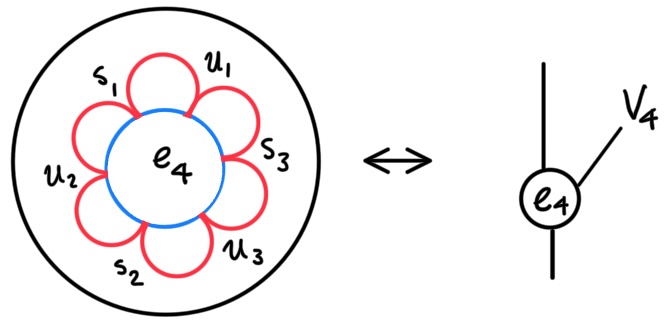}
				\end{center} 
				or
				\begin{center}
					\includegraphics[width = 6cm]{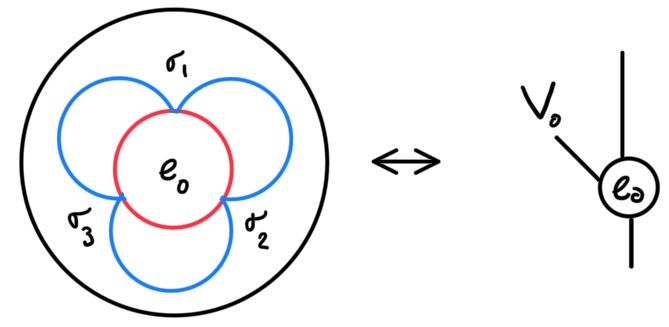}
				\end{center}
				but not
				\begin{center}
					\includegraphics[width = 7cm]{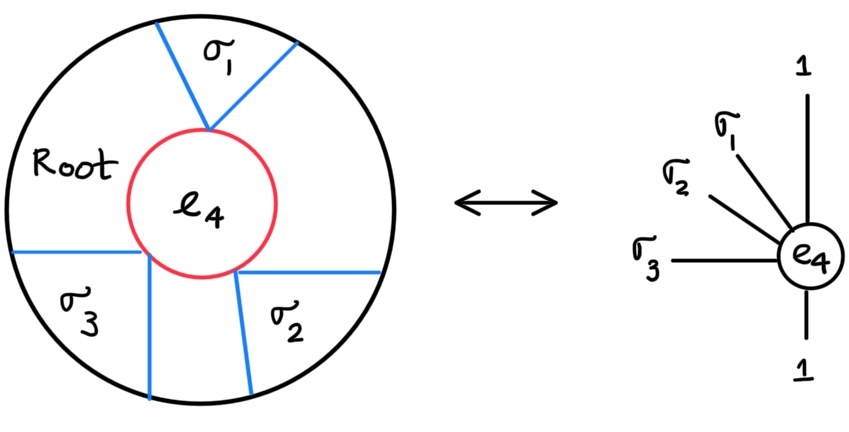}
				\end{center}
				
				\item here are no unmultiplied multipliable pairs; more specifically, if we have $a, a' \in \A$ labelling adjacent regions 
				\begin{center}
					\includegraphics[width = 4cm]{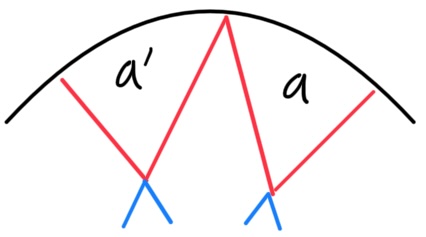}
				\end{center}
				then $a\cdot a' = a' \cdot a = 0$, and if $b, b' \in \B$ labelling adjacent regions
				\begin{center}
					\includegraphics[width = 4cm]{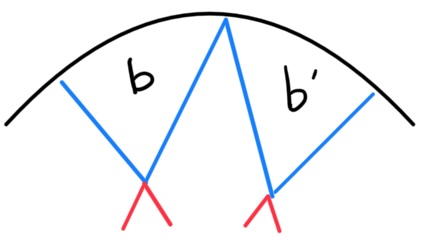}
				\end{center}
				then $b \cdot b' = b' \cdot b = 0$.
				\end{enumerate}
				
				\item:\label{firstop} One vertex of each edge must be internal, so we cannot have e.g.
				\begin{center}
					\includegraphics[width = 3cm]{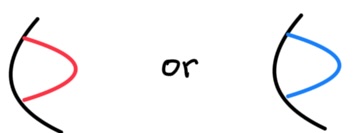}
				\end{center}
				
				\item: The graph must be rooted in the following sense: at least one of the edges intersects the boundary, which by the preceding conditions implies that there must be at least one region of the form
					\begin{center}
						\includegraphics[width = 1.5cm]{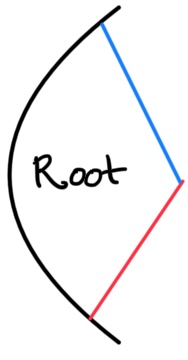}
					\end{center} 
					We choose one such region as the designated region or \emph{root} idempotent. 
				
				\item: The sequence of elements in a given region must be multipliable, read either left to right or clockwise according to the orientation of the main disk; as in
				\begin{center}
					\includegraphics[width = 8cm]{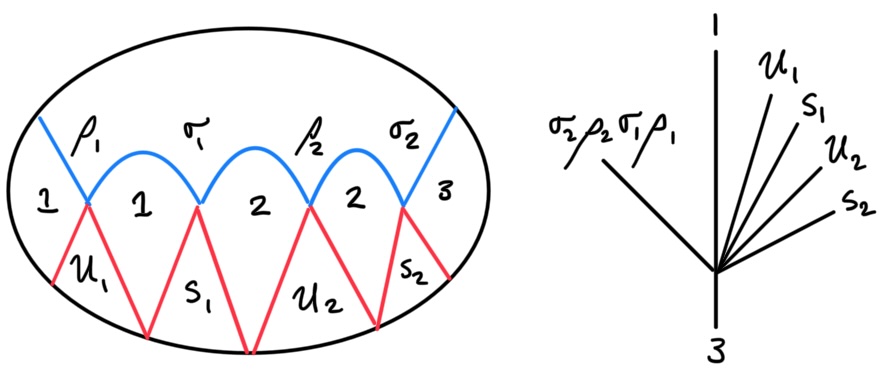}
				\end{center}
				
				\item: Cycles are of one of the following three forms:
				\begin{enumerate}
					\item $\e_i$ for $1 \leq i \leq N$ i.e. petal weight:
					\begin{center}
						\includegraphics[width = 5cm]{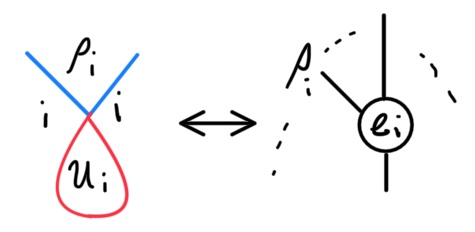}
					\end{center}
					
					\item $\e_{N + 1}$ i.e. central $\A$-cycle:
					\begin{center}
						\includegraphics[width = 8cm]{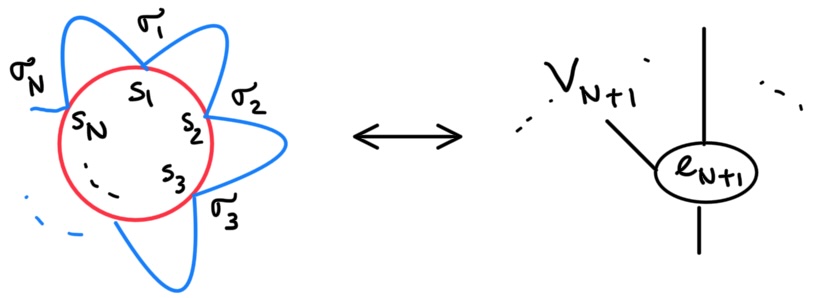}
					\end{center}
					
					\item $\e_0$ i.e. central $\B$-cycle:
					\begin{center}
						\includegraphics[width = 8cm]{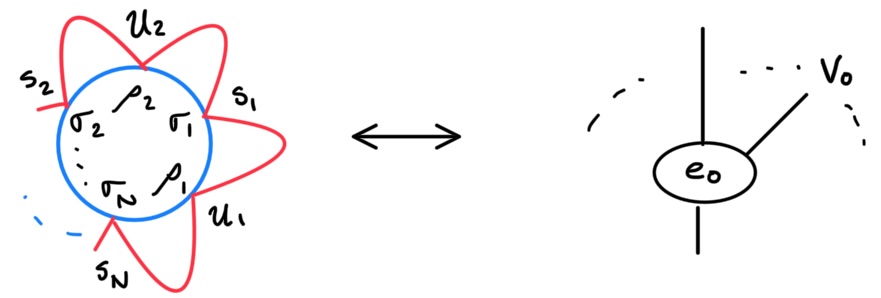}
					\end{center}
				\end{enumerate} 
				
				\item:\label{lastop} (Balancing) Any vertex of the form
				\begin{center}
					\includegraphics[width = 1cm]{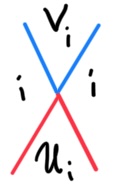}
				\end{center}
				is a petal cycle on the far side, i.e.
				\begin{center}
					\includegraphics[width = 10cm]{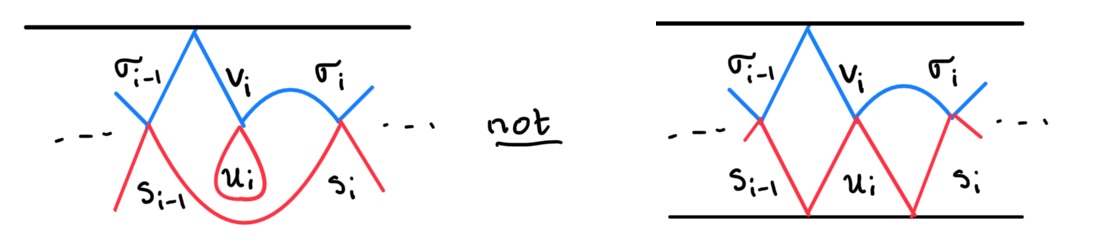}
				\end{center}
			\end{enumerate}
			
			Before we talk about how to associate a tree to any allowable graph, we need to talk about Maslov and Alexander gradings. We make the following definition for weighted input sequences: for $a_1, \ldots, a_n \in \A, \x \in Y, b_1, \ldots, b_k \in \B$ and $\w$ the weight,
			\begin{equation}\label{bim1}
				m(\w, b_k, \ldots, b_1, \x, a_1, \ldots, a_n) = m(\w) + m(\x) + \sum_{i = 1}^n m(a_i) + \sum_{i =1}^k m(b_i)
			\end{equation}
			We also stipulate that $m(\x) = 0$ for each $\x \in Y$. We define the Alexander grading of an input sequence as
			\begin{equation}\label{bim6}
				A(\w, b_k, \ldots, b_1, \x, a_1, \ldots, a_n) = A(\w) + \sum_{i = 1}^n A(a_i) + \sum_{i =1}^k A(b_i)
			\end{equation}
			
			We also need to talk about the kind of modified matching between elements on the $\A$-side and the $\B$ side that will be crucial to deciding which trees correspond to operations. Let $\eta_1, \ldots, \eta_r, \tau_1, \ldots \tau_{\ell},$ be the basic inputs on the $\A$ and $\B$-sides, respectively, so in the case of the tree,
				\begin{center}
					\includegraphics[width = 3cm]{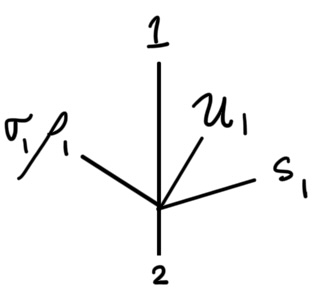}
				\end{center}
				we have 
				\begin{align*}
					\eta_1 &= U_1, \: \eta_2 = s_1\\
					\tau_1 &= \rho_1, \: \tau_2 = \sigma_1
				\end{align*}
				For weighted trees, we also keep track of $\w_1, \ldots, \w_t$, the basic components of the weight $\w$ (that is, with $\w_j = \e_{i_j}$ for each $1\leq j \leq t$, and some $0 \leq i_j \leq N + 1$). 
				
				A weighted input sequence $(\w, b_k, \ldots, b_1, \x, a_1, \ldots, a_n)$ \emph{admits a matching} if the corresponding sequence of basic elements, which we now write as
				\[
					(\e_{i_1}, \ldots, \e_{i_p}, q \cdot \e_0, \eta_1, \ldots, \eta_r, \tau_1, \ldots, \tau_{\ell})
				\]
				where 
				\begin{itemize}
					\item $1 \leq i_j \leq N + 1$ for each $1 \leq j \leq p$,
					\item $q \in \N$,
					\item We always count $\e_i$ on the $\A$-side for $1\leq i \leq N + 1$, and count $\e_0$ on the $\B$-side
				\end{itemize}
				satisfy the following conditions:
				\begin{enumerate}[label = Mat. \arabic*:]
					\item $p + r = q + \ell$ -- write the common sum as $\alpha$;
					
					\item There exists a one-to-one indexing map $f: \{\e_{i_1}, \ldots, \e_{i_p}, \eta_1, \ldots, \eta_r\}\to \{1, \ldots, \alpha\} $, such that if $r' < r''$, then $f(\eta_{r'}) < f(\eta_{r''})$;
					
					\item There exists a one-to-one indexing map $g: \{\e_0, \ldots, \e_0, \tau_1, \ldots, \tau_r\}\to \{1, \ldots, \alpha\}$ (where we include $q$ copies of $\e_0$, all of which are counted as distinct), such that if $\ell' < \ell''$, then $g(\tau_{\ell'}) < g(\tau_{\ell''})$;
					
					\item The resulting matching 
					\begin{equation}\label{bim13}
						h := (g^{-1} f):  \{\e_{i_1}, \ldots, \e_{i_p}, \eta_1, \ldots, \eta_r\} \to \{\e_0, \ldots, \e_0, \tau_1, \ldots, \tau_r\}
					\end{equation}
					preserves Alexander grading and never matches one weight with another weight. (We need to include this second part of the condition because $\e_0$ and $\e_{N + 1}$ have components in every Alexander grading.)
					
					\item $h$ matches any $V_i$ with the corresponding $\e_i$ for each $0 \leq i \leq N$
					
					\item $h$ matches $e_i$ with either $V_i$ or $\rho_i$ for $1 \leq i \leq N$, and only with $V_i$ for $i = 0$ or $N + 1$.
				\end{enumerate} 
				An input sequence \emph{admits a matching} if and only if we can construct $h$, as in~\eqref{bim13}, with the stipulated properties. 
				
				Here is an example of a tree whose input sequence admits a matching, and how the matching is constructed:
				\begin{center}
					\includegraphics[width = 8cm]{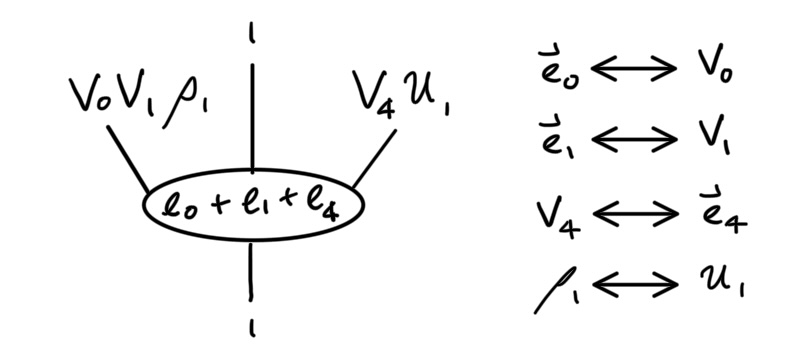}
				\end{center}
				Here is an example which does not admit a matching
				\begin{center}
					\includegraphics[width = 5cm]{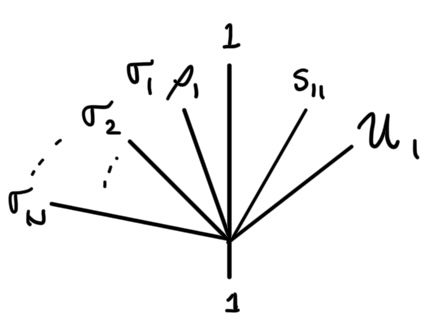}
				\end{center}
				because we are not allowed to rearrange the $\tau_i$ and $\eta_i$. See also Remark~\ref{bim12} for comparison between the condition ``admits a matching'' and the idempotent and Alexander grading conditions for operations given below. 
				
				We can also define the input sequence corresponding to an allowable graph, by reading the $\A$-inputs counterclockwise along the boundary from the root sector, and reading the $\B$-inputs clockwise along the boundary from the root sector, and adding in weight wherever it is found. We will see in Lemma~\ref{bim3} that every allowable graph admits a matching.
				
				Here is an example of an allowable weighted graph, and how to read the inputs off a graph, and compute the matching:
				\begin{center}
					\includegraphics[width = 13cm]{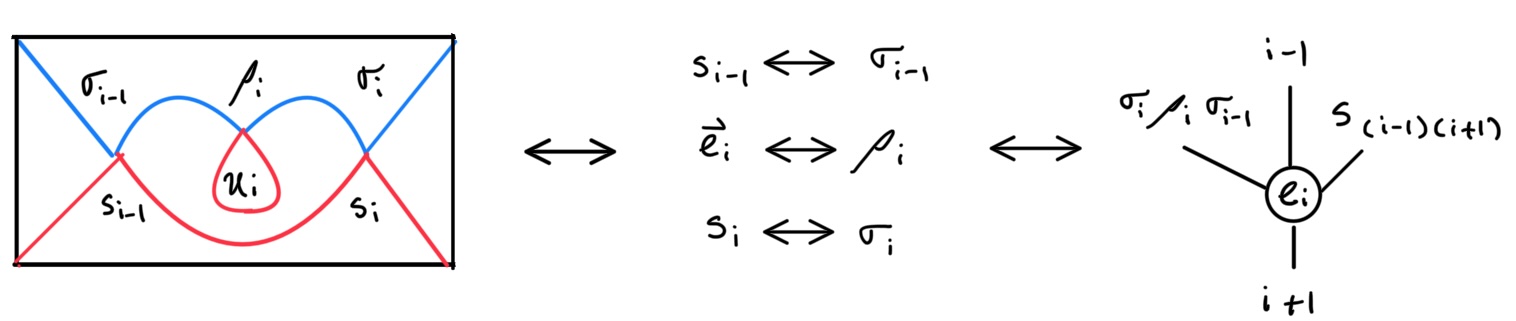}
				\end{center}
				with the requirement that the matching go in order -- i.e. if $\eta_1$ is matched to $\tau_7$, then $\eta_2$ cannot be matched to $\tau_3$. 
			
			With this in mind, we now show that:
			 
			\begin{lemma}\label{bim3}
				Any allowable graph $G$ with $\A$-inputs $a_1, \ldots, a_n$, initial element $\x \in Y$, $\B$-inputs $b_1, \ldots, b_k$ and weight $\w$ satisfies
			\begin{equation}\label{bim2}
				m(\w) + m(\x) + \sum_{i = 1}^n m(a_i) + \sum_{i =1}^k m(b_i) + k + n - 1 = 0
			\end{equation} 
			and 
			\begin{equation}\label{bim7}
				A(\w, b_k, \ldots, b_1, \x, a_1, \ldots, a_n) = 0 \emm2
			\end{equation}
			and weighted input sequence admits a matching in the sense described above.
			\end{lemma}
			
			\begin{remark}\label{bim4}
				\emph{Note that \emph{if we could} associate an operation $m(\: m_{k|n}^{\w} (b_k, \ldots, b_1, \x, a_1, \ldots, a_n)) = \y$ to our graph $G$, then since $m(\x) = m(\y) = 0$, Lemma~\ref{bim3} implies that this operation satisfies the bimodule equivalent of~\eqref{aa1}, namely}
				\begin{equation}\label{bim5}
					m(\: m_{k|n}^{\w} (b_k, \ldots, b_1, \x, a_1, \ldots, a_n)) = m(\w) + m(\x) + \sum_{i = 1}^n m(a_i) + \sum_{i =1}^k m(b_i) + k + n - 1
				\end{equation}
				\emph{Notice that the ``$-1$'' instead of ``$-2$'' is because we have one extra input in addition to the $\A$- and $\B$-inputs, namely the $\x$.}
			\end{remark}
			
			\begin{proof}[Proof of Lemma~\ref{bim3}]
				To get a matching from $G$, just go vertex by vertex clockwise from the root and match inputs across vertices. The only things to be careful of are:
				\begin{itemize}
					\item $U_i$ in a given sector counts as $\e_i$ in the matching if it is part of a petal cycle, and as $U_i$ otherwise. 
					
					\item $(U_1, s_1, \ldots, U_N, s_N)$ or a cyclic permutation of this sequence counts as $V_0$ if it is isolated about an $\e_0$ cycle, as in
					\begin{center}
						\includegraphics[width = 8cm]{bim9}
					\end{center}
					
					\item $(\sigma_1, \ldots, \sigma_N)$ or a cyclic permutation of this sequence counts as a $V_{N + 1}$ if it is isolated about a $\e_{N + 1}$ cycle, as in
					\begin{center}
						\includegraphics[width = 8cm]{bim8}
					\end{center}
					
					\item $(\rho_1, \sigma_1, \ldots, \rho_N, \sigma_N)$ or a cyclic permutation of this sequence counts as an $\e_0$ if it is in a central cycle;
					
					\item $(s_1, \ldots, s_N)$ or a cyclic permutation of this sequence counts as an $\e_{N + 1}$ if it is in a central cycle;
				\end{itemize}
 			
				For~\eqref{bim2}, we are going to use induction on the number of vertices in a given allowable graph. This is obvious for the basic graphs
				\begin{center}
					\includegraphics[width = 6cm]{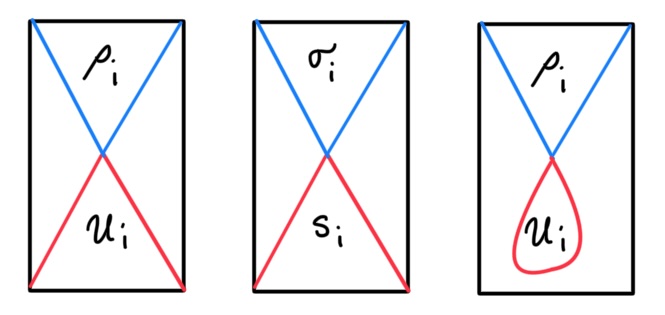}
				\end{center}
				by a straight calculation. 
				
				Next, note that if every edge of an \emph{allowable} graph is attached to a central cycle, the graph is a disjoint union of the basic central cycles
				\begin{center}
					\includegraphics[width = 6cm]{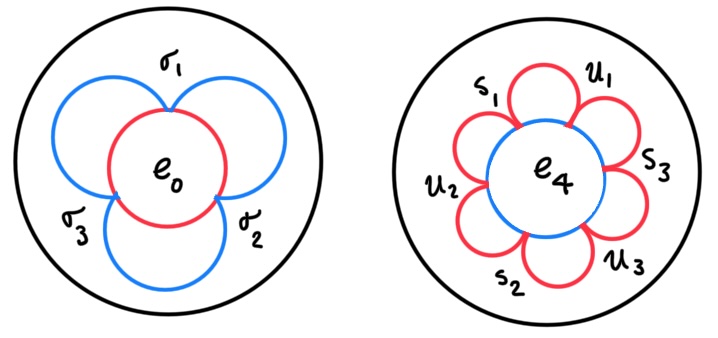}
				\end{center}
				(here drawn in the case $N = 3$). So in this case, the graph can look like
				\begin{center}
					\includegraphics[width = 12cm]{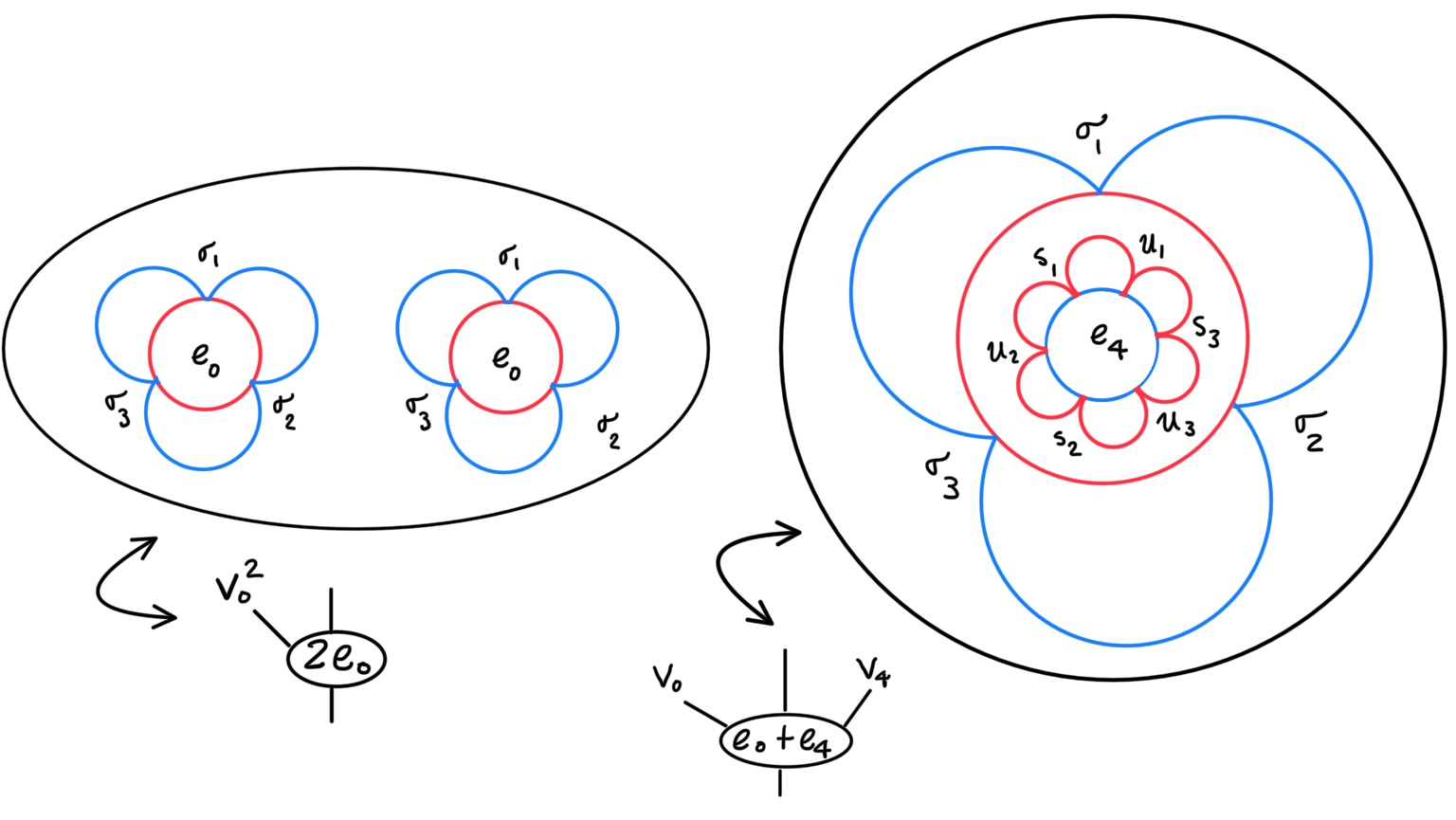}
				\end{center}
				
				If $G$ is one of the two basic central cycles, then it satisfies~\eqref{bim2}, since
				\begin{align*}
					m(e_0) + m(V_0) + 2 - 2 &= - (2N - 2) + (2N - 2) + 2 - 2 = 0 \\
					m(e_{N + 1}) + m(V_{N + 1}) + 2 - 2 &= 2 - 2 + 2 - 2 = 0
				\end{align*}
				for the first and second cases, respectively. For higher weight (where every edge is attached to a central cycle) the left and right sides of~\eqref{bim2} are just finite linear combinations of the left and right sides of the two lines above -- with matching coefficients on left and right -- so that~\eqref{bim2} is still satisfied.
				
				In fact, this allows us to remove all central cycles from an allowable graph, without affecting either side of~\eqref{bim2} or~\eqref{bim7}.
				
				We can therefore assume without loss of generality that our allowable graph $G$ does not contain central cycles. Suppose now that every allowable $k$-vertex graph satisfies~\eqref{bim2} and~\eqref{bim7}. Let $G$ be an allowable graph with $(k + 1)$ vertices. Pick one of the vertices to excise, and assume first that this is a vertex of the form below:
				\begin{figure}[H]
					\includegraphics[width = 1cm]{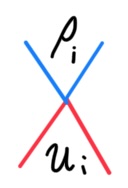}
					
					\caption{}\label{bim16'}
				\end{figure} 
				If the $U_i$ in Figure~\ref{bim16'} is a petal-weight, then the surrounding part of the graph must be of the form
				\begin{center}
					\includegraphics[width = 6cm]{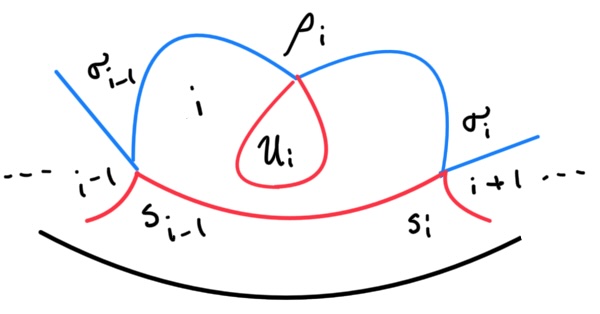}
				\end{center}
				This portion contributes
				\begin{equation}\label{bim8}
					2 + m(\sigma_{i -1}) + m(\rho_i) + m(\sigma_{i}) + m(\e_i) + m(s_{i - 1}) + m(s_i) = 2 - 3 + 2 = 1
				\end{equation}
				to the left hand side of~\eqref{bim2}, where the first ``2'' comes because this part has two inputs, and where we are assuming the two $s$-terms are not part of a cycle. Replace this portion of the graph with
				\begin{center}
					\includegraphics[width = 6cm]{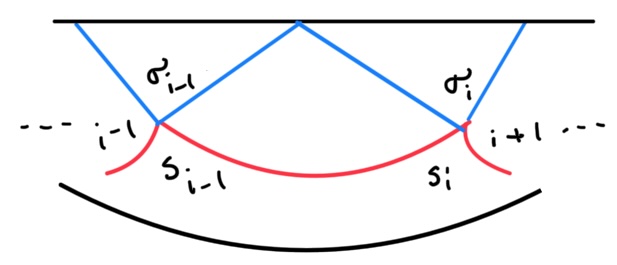}
				\end{center}
				This now contributes
				\begin{equation}\label{bim9}
					3 + m(\sigma_{i - 1}) + m(\sigma_i) + m(s_{i - 1}) + m(s_i) = 3 - 2 = 1.
				\end{equation}
				to the left hand side of~\eqref{bim2}.
				
				Likewise, if the $U_i$ does not come from a petal-cycle, then the part of the graph surrounding this vertex looks like one of
				\begin{center}
					\includegraphics[width = 12cm]{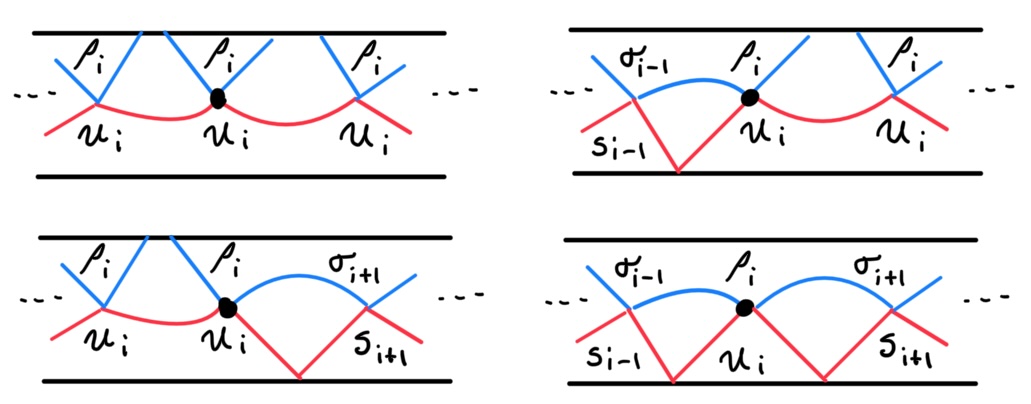}
				\end{center}
				As above, we can calculate the contribution from each one of these to the left hand side of~\eqref{bim2}, then remove the vertex in question, and recalculate to see that the total contribution of this part of the graph does not change. The resulting graph is also an allowable graph, and hence, we can apply the induction hypothesis.
				
				If the vertex looked like a 
				\begin{figure}[H]
					\includegraphics[width = 1.5cm]{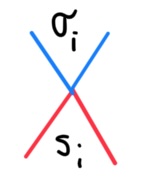}
				\end{figure}
				\hspace{-0.4cm}then we can go through the options for adjacent vertices, and in each case, excise this vertex in such a way as to not change the sum on the left hand side of~\eqref{bim2}. This means that we can apply the induction hypothesis in this case, as well.
			\end{proof}
			
			Notice that in the limited cases, such as some of the ones treated above, it is reasonably obvious how to associate a tree to a graph. In general, we do that as follows. 
			
			\begin{proposition}\label{bim10}
				\emph{(Identifying operations in $Y$)} Any allowable graph $G$ corresponds to a tree with input sequence 
				\begin{equation}\label{bim12}
					(\w, b_k, \ldots, b_1, \x, a_1, \ldots, a_n)
				\end{equation}
				and output $\y \in Y$ satisfying the following properties
				\begin{enumerate}[label = (\roman*)]
					\item \emph{(Idempotents)} 
					\begin{itemize}
						\item The initial idempotent of $a_i$ (resp. $b_i$) is the final idempotent of $a_{i - 1}$ (resp. $b_{i - 1}$) for each $1 < i \leq n$ (resp. $1 < i \leq k$);
						
						\item The initial idempotent of $a_1$ (resp. $b_1$) is the idempotent of $\x$;
						
						\item The final idempotent of $a_n$ (resp. $b_k$) is the idempotent of $\y$;
					\end{itemize}
					
					\item \emph{(Maslov grading)} The input sequence satisfies~\eqref{bim2};
					
					\item \emph{(Alexander grading)} The input sequence satisfies~\eqref{bim7} and admits a matching;
				\end{enumerate}
				Under these conditions, the associated tree is unique with this property. Moreover, every weighted tree the input sequence of which satisfies (i)-(iii) determines a unique allowable graph.
			\end{proposition}
			
			\begin{remark}\label{bim12}
				\emph{Notice that the condition ``admits a matching'' is stronger than the condition we place on Alexander grading (namely~\eqref{bim7}, above). The matching matches elements with the same Alexander grading, so trees that admit a matching automatically satisfy~\eqref{bim7}. However, the second tree discussed at the end of the definition of matchings,}
				\begin{center}
					\includegraphics[width = 5cm]{bim33}
				\end{center}
				\emph{satisfies conditions (i) and (ii), as well as~\eqref{bim7}, but does not admit a matching. We include the explicit Alexander grading condition~\eqref{bim7} in the above for clarity rather than strictly for completeness.}
			\end{remark}
			
			\begin{proof}[Proof of Proposition~\ref{bim10}]
				Look at a given allowable graph $G$. We can read off the input sequence by looking at the labels of the all-red and all-blue sectors of each vertex. If a sector is part of a cycle, then it goes to weight, otherwise it goes to the appropriate side. The $\x \in Y$ comes from the idempotent of the root. That the resulting input sequence $(\w, b_k, \ldots, b_1, \x, a_1, \ldots, a_n)$ satisfies (ii) and (iii) follows directly from Lemma~\ref{bim3}, so we only need to verify that 
				\begin{enumerate}[label = (\alph*)]
					\item We can choose an arrangement of $(\w, b_k, \ldots, b_1, \x, a_1, \ldots, a_n)$ satisfies (i), and 
					
					\item There is a single well-defined idempotent corresponding to an element $\y \in Y$, which is the final idempotent of both $a_n$ and $b_k$, and hence the output element for the tree $T$ corresponding to $G$.
				\end{enumerate}
				The trick is that (a) comes from (iii) from the statement. This works by induction on the number of vertices in $G$. The base case is the graphs
				\begin{center}
					\includegraphics[width = 12cm]{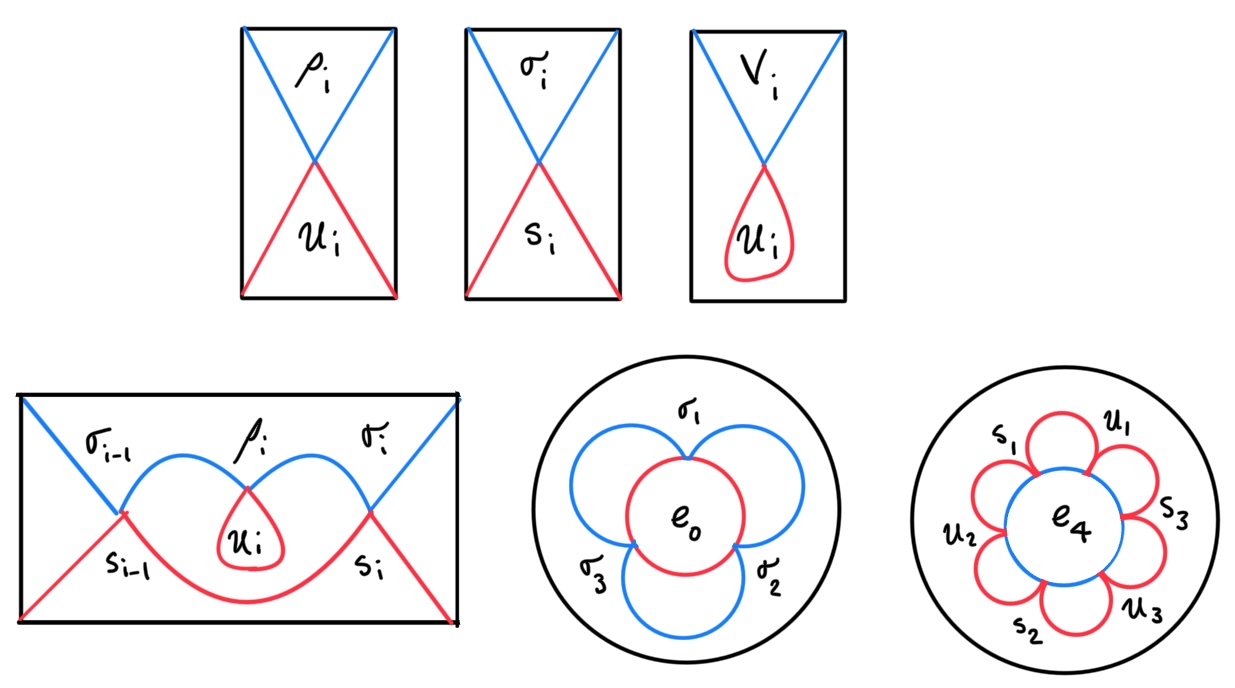}
				\end{center}
				 These clearly satisfy (a) and (b), corresponding to the trees
				 \begin{center}
				 	\includegraphics[width = 9cm]{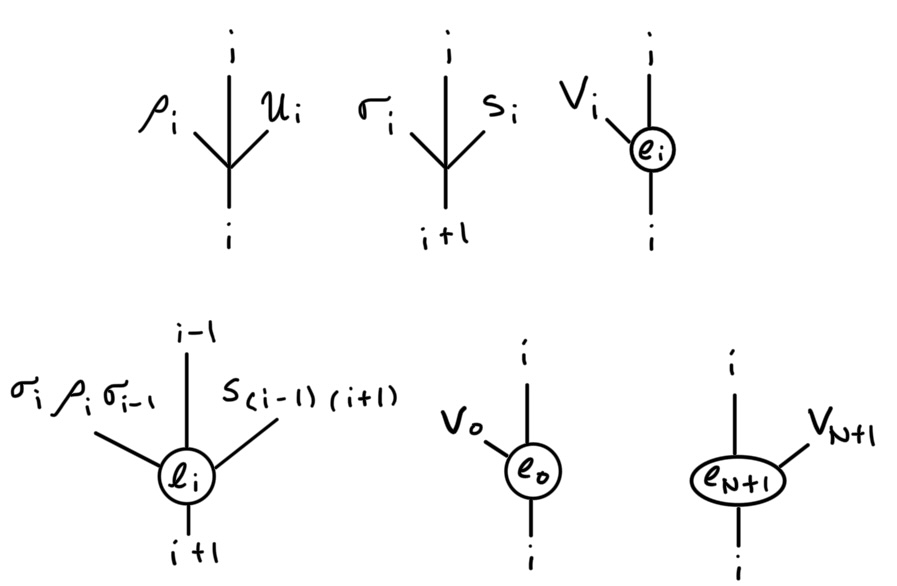}
				 \end{center}
				
				Now suppose we know that there is an arrangement of the inputs satisfying (i), and a well-defined final idempotent for all allowable graphs with $\ell$ or fewer vertices $(\w, b_k, \ldots, b_1, \x, a_1, \ldots, a_n)$ with $k \leq k_0, n \leq n_0$. Look at a graph with $\ell + 1$ vertices. Then as in the proof of Lemma~\ref{bim3}, we can excise a single vertex (which we can assume without loss of generality is not part of a cycle, and is not the final vertex if there is one) without changing that this is an allowable graph, $G'$. By the induction hypothesis, we can get an (allowable) tree $T'$ corresponding to $G'$. Then looking at the point in $G'$ from which we excised our original vertex, we can find the corresponding point in $T'$, and add back the appropriate algebra elements (properly multiplied) to get an allowable tree $T$ for $G$.
				
				The second half of the proof deals with getting an allowable graph back from a given tree $T$, the inputs of which satisfy (i)-(iii). We use the given matching to get a graph $G$ for $T$ by induction on $\alpha$, the number of inputs on either side, counting weights (with notation as in 	the definition of the matching). 	
				
				We now need to verify that 
				\begin{enumerate}
					\item $G$ satisfies Op. 1-6 (i.e. is a bona fide operation);
					
					\item $G$ gives back the original tree $T$ when we read off the elements as we did in the first half of the proof;
				\end{enumerate}
				We do (1) first. This graph is clearly rooted -- we just choose the initial idempotent $\x$ as the root, which is well-defined because the graph was constructed inductively. The sequence of elements in a given region is multipliable by construction of $G$ -- as we walked down the pairs of $\A$- and $\B$-elements, we connected vertices as
				\begin{center}
					\includegraphics[width = 3cm]{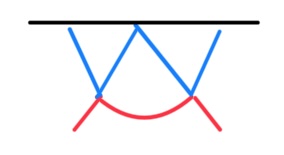}
				\end{center}
				when the adjacent $\A$ elements are multipliable and as
				\begin{center}
					\includegraphics[width = 3cm]{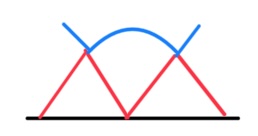}
				\end{center}
				when the $\B$ elements are multipliable with the exception of petal weights, which show up in the obvious way:
				\begin{center}
					\includegraphics[width = 6cm]{bim11}
				\end{center}
				This also gives us back the fact that all multipliable pairs are multiplied. That petal weights and central weights are isolated from the boundary also follows from the construction of the graph. The only thing that is left is balancing, and this follows from stipulation Mat. 5 from the definition of matchings.
				
				Now we can use the first part of the proof to get back an allowable tree $T'$ from $G$, and $T'$ will have the same sequence of basic elements and weights (in the same order) as $T$. Also, there are no unmultiplied multipliable pairs of elements in $T'$. The only thing that could go wrong is that there are such pairs in $T'$ (because everything else is the same by construction). But $T'$ satisfies~\eqref{bim2}, and if $T$ had one or more multipliable pairs (with everything else the same as $T'$, as it has) then it could not also satisfy~\eqref{bim2}. Since both $T$ and $T'$ satisfy~\eqref{bim2}, it follows that since everything else matches, they must be identical. 
			\end{proof}
			
			The last step of the section is to verify the $\A_{\infty}$ relations for $Y$ -- that is, show that given any tree $T$, the sum over all ways to add an edge to $T$ is zero. Again, the four ways to add an edge to a tree are a pull, a split, a push, or a differential, as noted in Section~\ref{aabimdef}. It is also important to remember that in our case, there are only non-zero differentials on the $\B$-side. 
			
			First, we note that:
			\begin{lemma}\label{ver5}
				Start with a tree $T$ and look at the left hand side of~\eqref{bim2}, namely
				\begin{equation}\label{ver6}
					m(\w) + \sum_{i = 1}^n m(a_i) + \sum_{i =1}^k m(b_i) + k + n - 1
				\end{equation}
				corresponding to the inputs and weight of $T$. If $T$ gives rise to a non-zero operation (or a composition of non-zero operations, as for a split) by one of the four actions above, then this quantity is $+1$. In particular, a pull, push, or differential shifts the quantity of~\eqref{ver6} by $-1$. 
			\end{lemma}
			
			\begin{proof}
				For a pull, we recall that 
				\begin{equation}\label{ver7}
					m(\mu_j^{\w}(a_1, \ldots, a_j)) = \sum m(a_j) + m(\w) + j - 2
				\end{equation}
				Let $T$ be the initial tree and $T'$ the one that results after pulling together $a_1, \ldots a_j$ (which is an adjacent string of elements on either the $\A$-side or the $\B$-side). Write $C$ for then the expression of~\eqref{ver6}, for $T$, and $D$ for the expression on the right hand side of~\eqref{ver7}. Then the expression of~\eqref{ver6} is, for $T'$,
				\[
					C + D - j + 1 - \sum m(a_j) - m(\w) = C - 1.
				\]
				If $T'$ corresponds to a non-zero operation, then $C = +1$.
				
				Now look at a split, where by splitting $T$, we obtain a composition of trees each corresponding to a \emph{valid} bimodule operation. Suppose the operation on the top pulls together $a_1, \ldots, a_j$ and $b_1, \ldots b_{j'}$, in addition to some weight $\w'$, so:
				\begin{center}
					\includegraphics[width = 9cm]{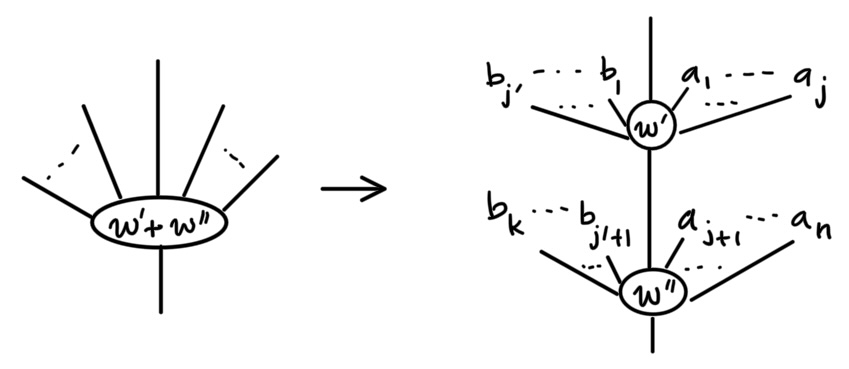}
				\end{center}
				Because $T'$ corresponds to a valid operation, we have
				\begin{equation}\label{ver8}
					m(\w') + \sum_{i = 1}^j m(a_i) + \sum_{i =1}^{j'} m(b_i) + j + j' - 1	= 0 
				\end{equation}
				Likewise,
				\begin{equation}\label{ver9}
					m(\w'') + \sum_{i = j + 1}^n m(a_i) + \sum_{i = j' + 1}^{k} m(b_i) + (n - j) + (k - j') - 1 = 0
				\end{equation}
				Summing the quantities from~\eqref{ver8} and~\eqref{ver9} tells us that the left hand side of~\eqref{ver6} is $+1$. 
				
				For a push and a differential, the argument is analogous to the one for a pull, and is omitted. 
			\end{proof}
			
			As in the case of $\A$, we introduce the notion of \emph{augmented graphs}, which will correspond to non-trivial relations. An augmented graph is a graph satisfying~\ref{firstop} through~\ref{lastop}, and contains exactly one composite region. We first show that every such graph corresponds to a pair of operations. 
			
			There are four basic types of composite regions, each of which corresponds to a particular type of relation.
			\begin{enumerate}
				\item Pull vs. split:
				\begin{center}
					\includegraphics[width = 6cm]{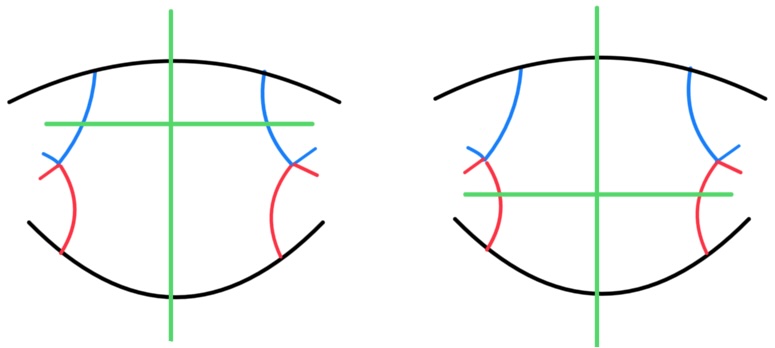}
				\end{center}
				In the first case, we are assuming the pair of elements in the adjacent blue sectors are multipliable, and the pair in the adjacent red sectors is not; in the second, we are assuming the reverse. The two green lines represent the two actions possible for this graph. We also assume for a moment that both pairs are not multipliable; we will deal with the other case in a moment. 
				
				For now, look at the first picture, and replace it with the graph:
				\begin{center}
					\includegraphics[width = 4cm]{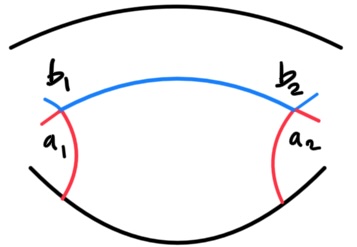}
				\end{center}
				Since everything else about the original augmented graph was allowable, and we have now eliminated the objectionable composite region, this is an allowable graph, and corresponds to a nonzero operation $T'$. 
				
				To get back a tree $T$ corresponding to the original graph $G$, find the pair of elements that were multiplied to get from $G$ to $G'$, and unmultiply them on the level of trees. Since nothing else has changed (and nothing else changed to go from $G$ to $G'$), this tree $T$ corresponds to $G$. 
				
				But notice that because the elements $b_1$ and $b_2$ (which will show up on the $\B$-side of $T$) are unmultiplied, we can also modify $T$ by splitting:
				\begin{center}
					\includegraphics[width = 9cm]{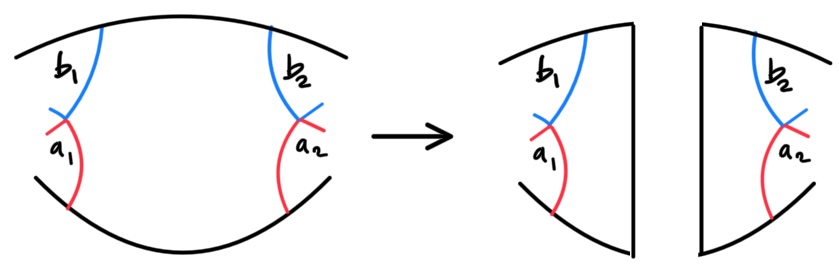}
				\end{center}
				which corresponds to:
				\begin{center}
					\includegraphics[width = 6cm]{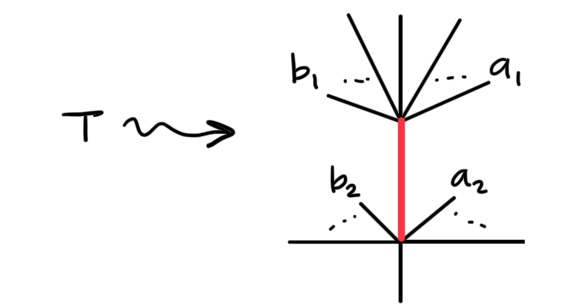}
				\end{center}
				To conclude that both trees (top and bottom) correspond to non-zero operations, just note that the two graphs (left and right) of the split are each allowable, and hence correspond to non-zero operations. Since the trees corresponding to these sections of the graph are just $T'$ and $T''$, the conclusion follows. 
				
				The case where the two $\B$-elements are multipliable, but not the pair of $\A$-elements, is analogous.
				
				Finally, we have to deal with the case where both pairs are multipliable only way that both pairs would be multipliable is
				\begin{center}
					\includegraphics[width= 10cm]{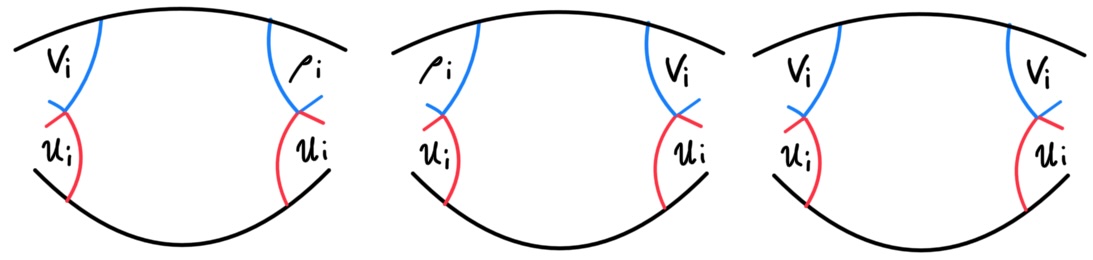}
				\end{center} 
				We do not consider this case, because the augmented graph $G$ would correspond to a tree $T$ with 
				\[
					m(\w) + m(\x) + \sum_{i = 1}^n m(a_i) + \sum_{i =1}^k m(b_i) + k + n - 1 = 2
				\]
				which, by Lemma~\ref{ver5}, cannot yield a non-zero operation (or composition of non-zero operations) by any of the four actions.
				
				\item Pull vs. push, with central orbit:
				\begin{center}
					\includegraphics[width = 6cm]{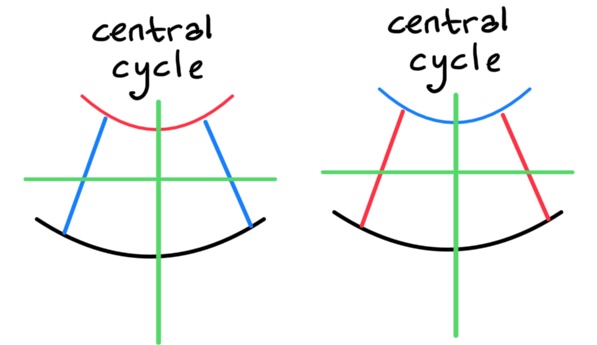}
				\end{center}
				We are assuming, of course, that all other regions are correctly zipped up, e.g., in the case $N = 3$ and for the $\alpha$-bordered orbit:
				\begin{center}
					\includegraphics[width = 8cm]{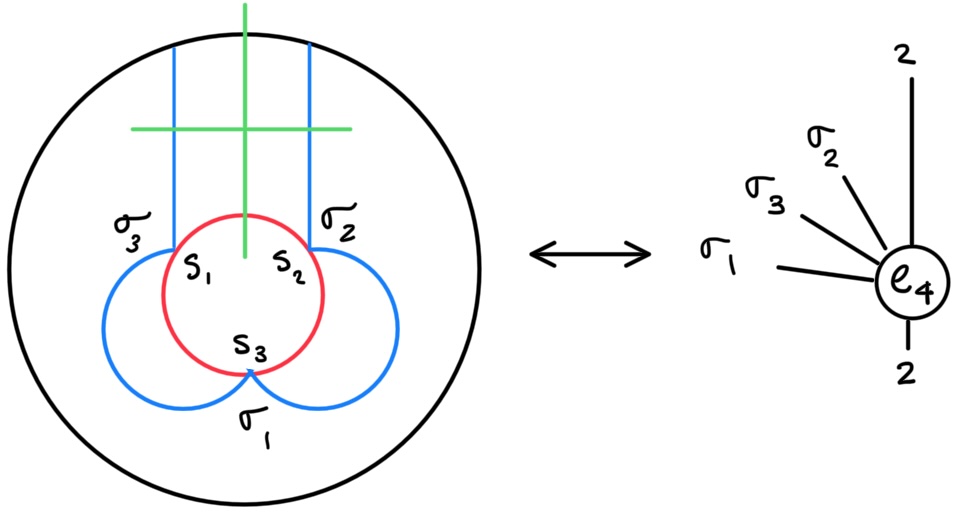}
				\end{center}
				This means that if we zip the two edges together, we get
				\begin{center}
					\includegraphics[width = 8cm]{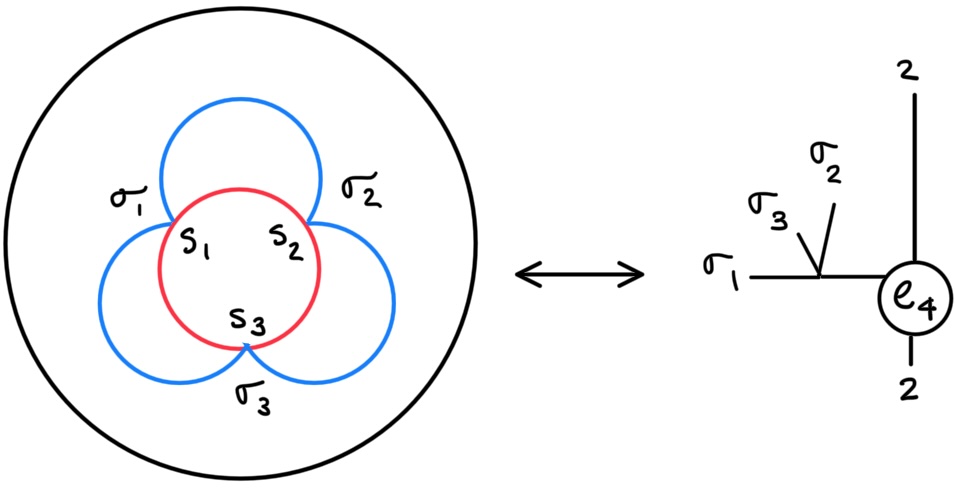}
				\end{center}
				which is just the standard $m_{1|0}^{\e_0}(V_0, \x) = \x$, where $\x$ our root of choice. 
				
				Likewise, if we unzip the red edge (with a little imagination) we get
				\begin{center}
					\includegraphics[width = 10cm]{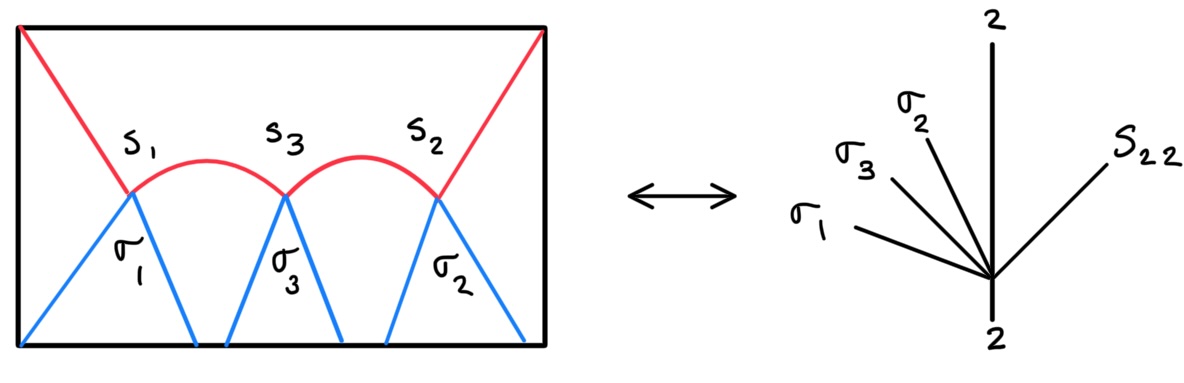}
				\end{center}
				It is clear that the original augmented graph $G$ corresponds to the tree
				\begin{center}
					\includegraphics[width= 2cm]{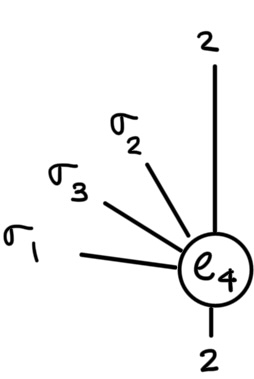}
				\end{center}
				or the analogous graph, for general $N$. The case for the $\beta$-bordered orbit is analogous.
				
				\item Pull vs. push, with radial orbit: 
				\begin{center}
					\includegraphics[width = 4cm]{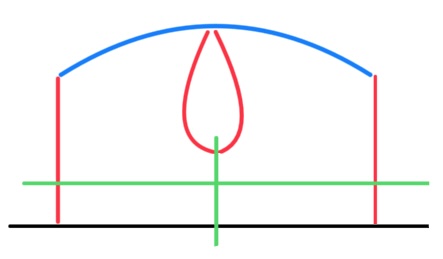}
				\end{center}
				The two allowable graphs obtained from $G$ are:
				\begin{center}
					\includegraphics[width = 9cm]{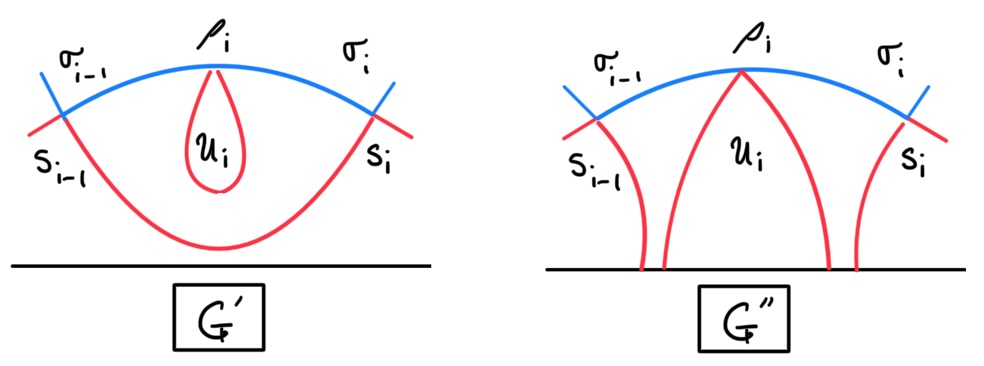}
				\end{center}
				Using the first one, we can get a tree $T'$ corresponding to a non-zero operation, then unzipping the two $s_i$ at the correct location (in the tree) we can get back a tree $T$ that clearly corresponds to the original graph $G$. Pushing out the radial orbit as in the second picture gives another allowable graph, $G''$, and drawing the three corresponding trees together: 
				\begin{center}
					\includegraphics[width = 14cm]{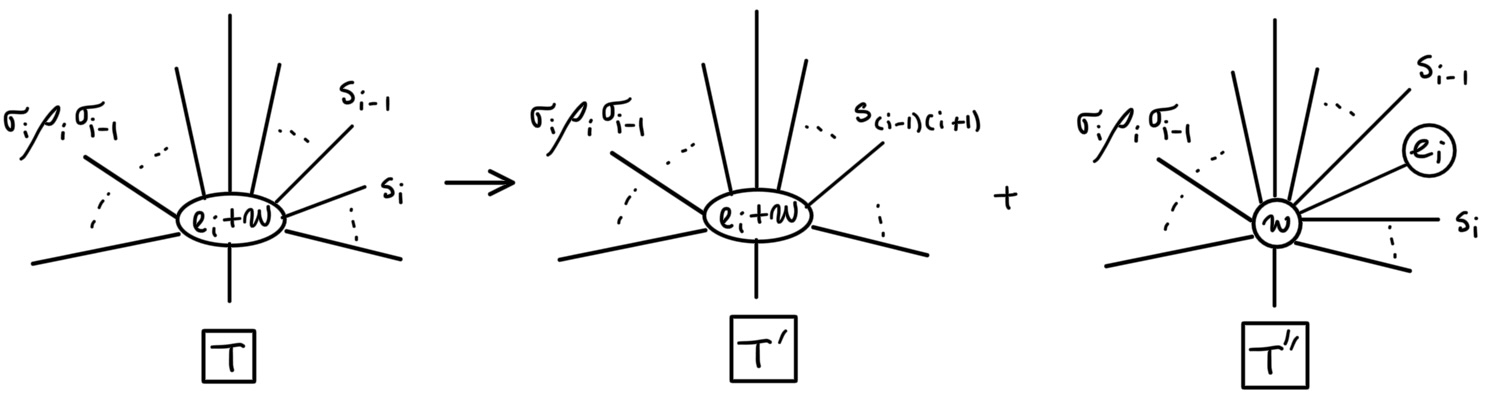}
				\end{center}  
				These are all the trees for $G', G''$, and $G$, above. Hence, again, $G$ gives rise to a cancelling pair of operations.
				
				\item Differential vs. push:
				\begin{center}
					\includegraphics[width = 6cm]{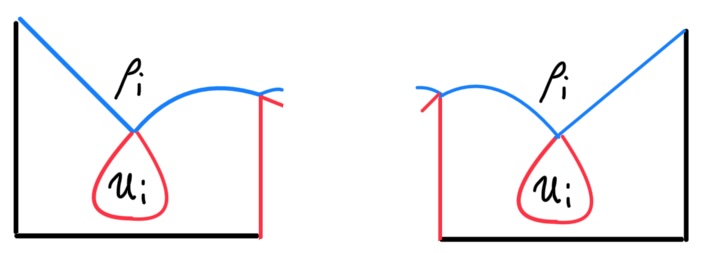}
				\end{center}
				The other way this can appear is
				\begin{center}
					\includegraphics[width = 1.5cm]{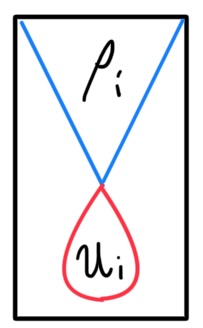}
				\end{center}
				but the trees are managed the same way, so we will restrict our attention to the leftmost picture above.
				
				For this, the two resulting graphs (from resolving the augmented region) are
				\begin{center}
					\includegraphics[width = 7cm]{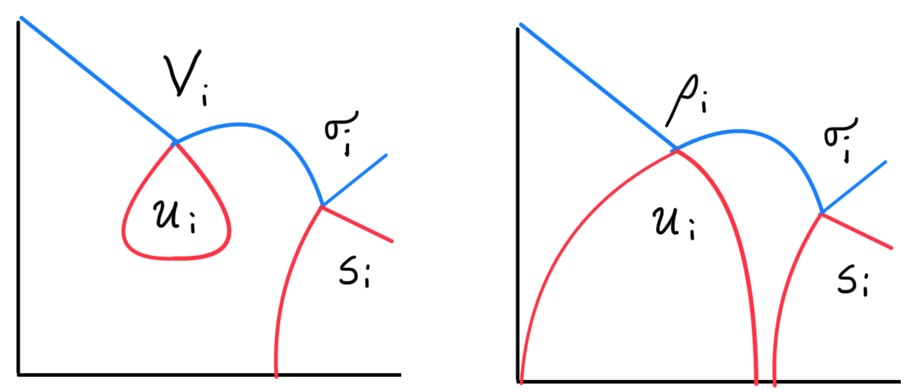}
				\end{center}
				Again, since these are allowable graphs, they correspond to a pair of non-zero operations, which have corresponding (allowable) trees $T'$ and $T''$. From either of these, we can get back a tree $T$ for the original graph $G$. That the original method of getting $T \rightsquigarrow T'$ and $T \rightsquigarrow T''$ correspond to a differential and a push, respectively, is clear on inspection.
			\end{enumerate}
			
			\begin{remark}\label{ver10}
				\emph{Our graph $G$ is never allowed to have ``composite'' regions. The issue is not that we cannot define a corresponding tree: we can, by the exact same methods used above when we associated a tree to an allowable graph in Lemma~\ref{bim10}. The issue is that when we look at the the resulting tree $T$, and evaluate the expression from~\eqref{ver6}, we will get $+2$ (assuming that the rest of the graph $G$ is allowable). By Lemma~\ref{ver5}, this means that no way of ``resolving'' this region (by applying one of our four actions to the corresponding tree and considering the graph corresponding to the result) will yield an allowable graph. Hence, no graph $G$ with such a ``composite'' region can correspond to a non-trivial relation.}
			\end{remark}
			
			The main lemma in the proof of the $\A_{\infty}$-relations for $Y$ is now as follows:
			\begin{lemma}\label{ver3}
				The non-trivial relations $T$ (i.e. trees where adding some single edge gives a non-trivial operation) are in one-to-one correspondence with augmented graphs as described above. 
			\end{lemma} 
			
			\begin{proof}
				We have already shown that augmented graphs each correspond to a non-trivial relation. Now start with a non-trivial relation, some tree $T$. Any resulting tree $T'$ (from adding an edge to $T$) is a non-trivial operation, and therefore satisfies the conditions of Proposition~\ref{bim10}, and corresponds to an allowable graph $G'$. To get from $T$ to $T'$, we performed one of the four allowable actions (pull, split, push, or differential). We can mark, on $T$ and $T'$, where this operation happened, and the corresponding location on $G'$. By the discussion above, going from $T$ to $T'$ corresponds to zipping up a composite region. We can walk back through the process above to obtain an augmented graph $G$ with one composite region $R$, such that one of the ways of zipping $R$ up to get an allowable graph gives back exactly $G'$. Note that $G$ has the same number of vertices (with the same labels) as $G'$, and $T$ has the same sequence of basic inputs (recalling that $V_0, \: V_{N + 1}$ are not regarded as basic, but composite) up to differentiation (possibly $\rho_i$ in $G$ but $V_i$ in $G'$) and push-outs of $\e_i$'s. This means that $T$ corresponds to the augmented graph $T'$. 
			\end{proof}
			
			\begin{proposition}\label{ver4}
				$Y$ satisfies the $\A_{\infty}$ relations and is a valid $AA$-bimodule. 
			\end{proposition}
			
			\begin{proof}
				The trick is that each augmented graph can be resolved into an allowable graph in precisely two ways. In the one-to-one correspondence between relations and augmented graphs from Lemma~\ref{ver3}, this shows us that any tree $T$ which has a non-trivial relation can be resolved into a non-zero operation in exactly two ways; that is, the $\A_{\infty}$ relations hold. 
			\end{proof}
	
	\section{The DD bimodule}\label{DDbim}
	
	The goal of this section is to construct the DD bimodule $\:^{\A} X^{\B}$ and verify the $\A_{\infty}$ relations for it. The generators of the DD bimodule $X$ are the complements of generators of $Y$, that is, writing 
	\[
		\overline{\x} = \overline{\{i\}} = \{1 \ldots, N\} \smallsetminus \{i\},
	\]
	$X$ is generated by $\{\overline{\x}: \x \text{ a generator of $Y$}\}$. For the operations, recall that a DD bimodule over $\A$ and $\B$ (e.g. $X$) is just a type-D module over the $\A_{\infty}$ algebra $\A \otimes \B$, with operations defined as in Section~\ref{ddorigin}. The operation on $X$ is
	\begin{equation}\label{dd1}
		\delta^1: X \to \A \otimes \B \otimes X
	\end{equation}
	given by
	\[
		\delta^1(\overline{\x}) = U_i \otimes \rho_i \otimes \x + s_i \otimes \sigma_i \otimes \y
	\]
	where $\overline{\x} = \overline{\{i\}}, \: \overline{\y} = \overline{\{i + 1\}}$, and $1 \leq i \leq N$. In terms of trees, this looks like
	\begin{center}
		\includegraphics[width = 7cm]{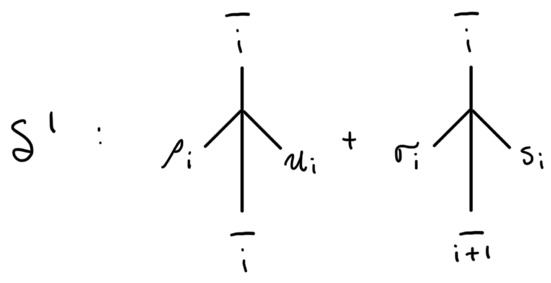}
	\end{center}
	
	In order to explicitly discuss these operations, we need to state what diagonal we are using to define operations on $\A \otimes \B$. There are two points in the construction of a diagonal where we are required to make an explicit choice. First, when choosing the seeds for the weighted portion of the algebra, and second, during the inductive step of the acyclic models construction. We already specified in Section~\ref{diagsfirst} that in each acyclic models step, we would choose the right-moving tree. Refer back to the last portion of Section~\ref{diagsfirst} for further details. 
	
	It remains to specify which seeds we choose for the construction of the weighted part of the digagonal. We do so now. The key point in motivating our choice is that each $\e_i$ popsicle only appears in one of the two algebras, determined by the color of arc which intersects the boundary circle corresponding to $\e_i$ -- so the $\e_0$ popsicle $\Psi_0^{\e_0}$ appears in $\B$, and the $\Psi_0^{\e_i}$, $1 \leq i \leq N + 1$, only appear in $\A$. With this in mind, we choose as seeds:
	\begin{equation}\label{diag2}
		\Gamma^{0, \e_0} ( \Psi_0^{\e_0}) = \top \otimes \Psi_0^{\e_0}, \qquad \Gamma^{0, \e_i} (\Psi_0^{\e_i}) = \Psi_0^{\e_i} \otimes \top \text{ for each } 1 \leq i \leq N + 1,
	\end{equation}
	or in terms of trees,
	\begin{center}
		\includegraphics[width = 10cm]{diag3}
	\end{center}
	In our case, we also stipulate that the $Y_1, Y_2 \in R$ are set equal to 1.
	
	The rest of this section is devoted to verifying~\eqref{dd2} for this particular choice of $\delta^1$. The first step will be to identify which summands $(\mu_n^{\w} \otimes \I) \circ \delta^n$ do not vanish, for a given $\overline{\x} = \overline{\{i\}}$. 
	
	Before we do this, however, we need to recall the definition of \emph{total length} of a sequence of algebra elements $\{a_1, \ldots, a_n\}$, initially given defined in Section~\ref{bb1}. Again, length is only defined for a sequence $\{a_1, \ldots, a_n\}$ where the initial idempotent of $a_i$ is the final idempotent of $a_{i - 1}$ for each $i$. To get the length of such a sequence, start $a_1$, and go along the list adding $1$ each time we pass into to a new idempotent. 
	
	For instance, $U_1, s_1, U_2$ has length 2, and $\sigma_1, \sigma_2, \sigma_3, \sigma_1$ has length 4 (where we are assuming $N = 3$ to make this allowable from an idempotent standpoint). Notice again that the notion of length is distinct from the number of basic inputs for a given sequence; so for instance $(\sigma_2 \rho_2 \sigma_1, \sigma_3)$ has length 3, but has 4 basic inputs.

	We are now ready to state the main lemma.
	
	\begin{lemma}\label{dd3}
		Start with $\overline{\x} = \overline{\{i\}}$. The only non-vanishing $(\mu_n^{\w} \otimes \I) \circ \delta^n$ are those which give as output
		\begin{enumerate}[label = (\alph*)]
			\item $U_{i} \otimes V_i \otimes \overline{\x}$ for $1 \leq i \leq N$;
			
			\item $V_0 \otimes U_0\otimes \overline{\x}$;
		\end{enumerate}
		Moreover, each of these can be obtained in exactly two ways.
	\end{lemma}
	
	\begin{proof}
		In this proof we are not going to try to determine all the possible operations $\mu_n^{\w}$; since these come from the weighted diagonal, they get very complicated very quickly. Instead, since we are looking at $(\mu_n^{\w} \otimes \I) \circ \delta^n$, the trick here is that the input sequence of $\mu_n^{\w}$ (the operation on $\A \otimes \B$) must come from repeated applications of $\delta_1$, and hence, 
		\begin{enumerate}
			\item Each input on the $\A$-side and the $\B$-side is basic (no composite / pre-multiplied inputs);
			
			\item There are the same number of basic inputs on the $\A-$ and $\B$-sides (this one is more or less obvious);
			
			\item\label{matchme}The Alexander gradings of the $k$-th input on the $\A$-side and the $k$-th input on the $\B$-side match, for each $1 \leq k \leq n$ (this is the important condition);
		\end{enumerate} 
		The key takeaway from these conditions -- especially~\eqref{matchme} -- is that once we determine the $\A$-side of the inputs, we have now determined the $\B$-side of the inputs, as well; each $\A$-input has a unique partner on the $\B$-side, and vice-versa. For instance, if the sequence of $\A$-inputs for $\mu_n^{\w} \otimes \I$ is $s_1, U_2, s_2, s_3$, then we know the sequence of  $\B$-inputs \emph{must} be $\sigma_1, \rho_2, \sigma_2, \sigma_3$, and the full sequence is
		\[
			(\overline{\x}, s_1 \otimes \sigma_1, U_2\otimes \rho_2, s_2 \otimes \sigma_2, s_3 \otimes \sigma_3)
		\]
		Now, notice that if $a_1, a_2$ are basic elements in $\A$ with $a_1 a_2 \neq 0$, then their partner elements $b_1, b_2$ on the $\B$-side will have $b_1 b_2 = 0$, and vice versa. Indeed, 
		\begin{itemize}
			\item $U_i \cdot U_i \neq 0$ and $\rho_i \cdot \rho_i = 0$;
			
			\item $s_i \cdot s_{i + 1} \neq 0$ and $ \sigma_{i + 1} \sigma_i = 0$;\footnote{Recall that multiplication in $\B$ goes right to left rather than left to right, which is why we reverse the order. (In this particular case, it is of course also true that $\sigma_i \cdot \sigma_{i + 1} = 0$.)}
			
			\item $\sigma_i \rho_i \neq 0$ and $U_i s_i = 0$;
			
			\item $\rho_{i + 1} \sigma_i \neq 0$ and $s_i U_{i + 1} = 0$;
		\end{itemize}
		We do not mention the $V_i$, because these are not generators of $\A$ or $\B$, but rather elements of the ground ring $\F[V_0, \ldots, V_{N + 1}]$.
		
		We claim that this implies that if an input sequence admits simple multiplication on one side, then it must be of the form
		\begin{center}
			\includegraphics[width = 10cm]{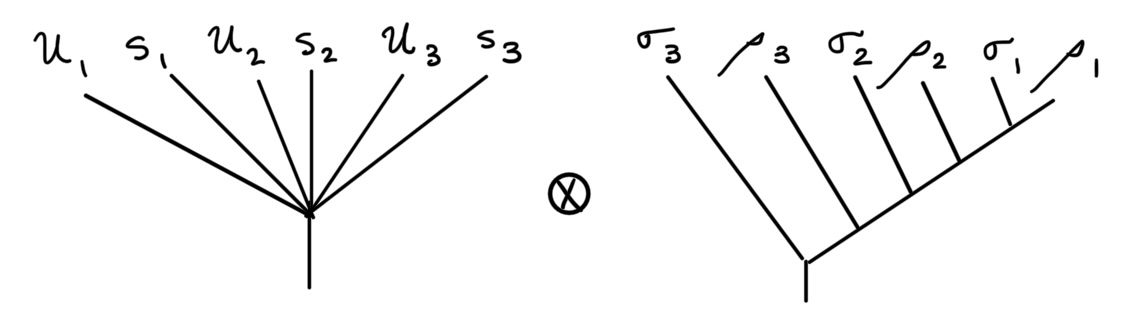}
		\end{center}
		or
		\begin{center}
			\includegraphics[width = 6cm]{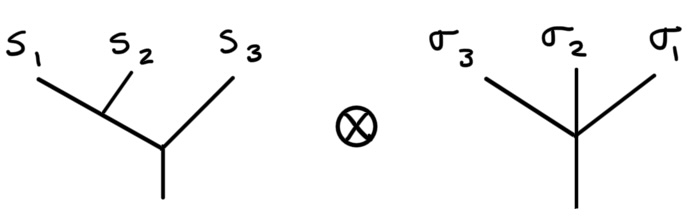}
		\end{center}
		where of course we are using $N = 3$ as an example, and we allow for cyclic permutation of the inputs. We now need to prove this claim, which we do by induction on the number of (basic) inputs. Notice that we are not currently dealing with the cases $n = 0$ and $n = 1$. We will do that at the very end. 
		
		If the number of basic inputs is $< N$, then all operations on the $\B$-side must be compositions of $\mu_2$s. By the above discussion, if the composition of $\mu_2$'s does not vanish on the $\B$-side, then \emph{any} composition of $\mu_2$'s we choose will vanish on the $\A$-side. A $\mu_2$ can only ever multiply adjacent elements, and all adjacent elements must eventually be multiplied in our composite. Any tree corresponding to a non-vanishing operation on the $\A$-side, with fewer than $N$ basic inputs, must involve some $\mu_2$. Since all adjacent elements on the $\A$-side are un-multipliable, we will eventually end up multiplying unmultipliable elements, and the entire expression will vanish. For instance:
		\begin{center}
			\includegraphics[width = 6cm]{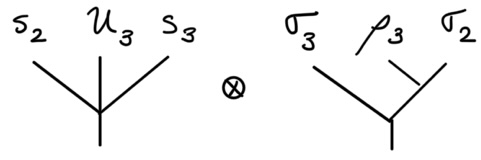}
		\end{center}
		and there is no way (weighted or unweighted) to pull together the three elements on the $\A$-side that does not vanish. Now suppose the number of basic inputs is $N$. If the $\A$-input sequence is $s_1, s_2, \ldots, s_N$, then the $\B$-input sequence has to be $\sigma_1, \ldots, \sigma_N$, and the only non-vanishing pair of trees from any $\mu_N^{\w}$ is
		\begin{center}
			\includegraphics[width = 7cm]{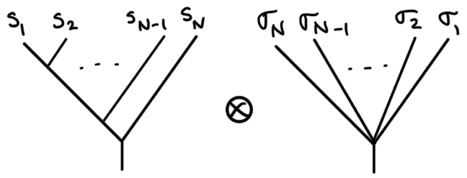}
		\end{center}
		Any other pair of trees, weighted or not, would have to include some lower multiplication on the $\B$-side, which will vanish. That this pair of trees \emph{does} appear in $\Gamma_*^{N, 0} \Psi_0^{N}$, follows by induction on $n$, and since we have stipulated that $\Gamma_*^{**}$ is a right-moving diagonal. In fact, by this argument, the pairs
		\begin{center}
			\includegraphics[width = 10cm]{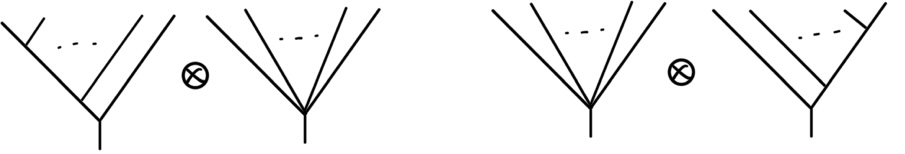}
		\end{center} 
		both appear in the sum that makes up $\Gamma_*^{n, 0} \Psi_0^{n}$, for any $n \in \N$. 
		
		So we have now proved the claim for sequences with up to (and including) $N$ inputs. Suppose the sequence of basic inputs has between $N$ and $2N$ elements. Then again, the tree on the $\A$-side must involve some internal $\mu_2$. Even if the $\A$-side has some (non-zero) weighted higher multiplication, as in
		\begin{center}
			\includegraphics[width = 7cm]{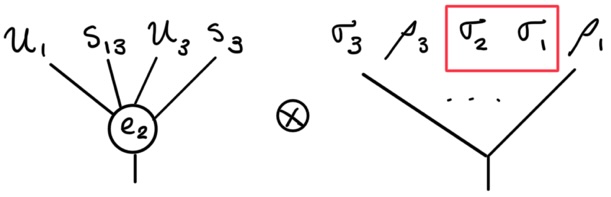}
		\end{center}
		the $\A$-tree must include a $\mu_2$, as pictured above; that is, the sequence of $\A$-inputs must have at least one multipliable pair. (We are not specifying the $\B$-tree yet, but we indicate what the problem will be.) If this multipliable pair on the $\A$-side is a pair of adjacent $U_i$'s, we are done, because there is no non-zero multiplication on the $\B$-side that includes adjacent $\rho_i$'s. If these are $s_i, s_{i + 1}$, then we have a corresponding $\sigma_i, \sigma_{i +1}$ on the $\B$-side. The only way the tree on the $\B$-side does not immediately vanish (because the $\sigma_i, \sigma_{i + 1}$ are not multipliable) is if they form part of a \emph{length} $N$ subsequence as in
		\[
			(\sigma_1, \sigma_2, \ldots, \sigma_N), (\sigma_1, \rho_2, \sigma_2, \ldots, \sigma_N), \text{ or } (\sigma_1, \rho_1, \sigma_2, \rho_2, \sigma_3, \ldots, \sigma_N).
		\]
		So far, if this is the case, there is no problem; since the sequence of inputs contains at least $N$ elements, this could happen. Even the corresponding partner elements on the $\A$-side might be multipliable (for example if we were looking at a weighted higher multiplication, which could involve fewer than $2N$ basic elements). But once this subsequence of $\B$-inputs has been pulled together, there are fewer than $N$ basic elements left over, as in
		\begin{center}
			\includegraphics[width = 6cm]{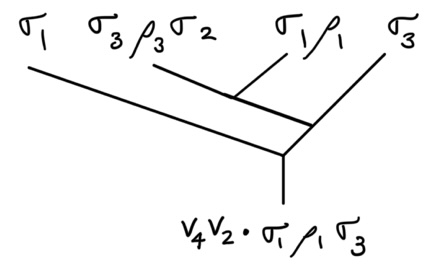}
		\end{center}
		Because there are no $\mu_n^{\w}$ on the $\B$-side involving fewer than $N$ basic input elements, the multiplication on the $\B$-side must be a composition of $\mu_2$'s, that is, all adjacent elements on the $\B$-side must now be multipliable. This means that outside the partner subsequence on the $\A$-side for our chosen subsequence, on which we cannot comment, all adjacent elements are non-multipliable, that is
		\begin{center}
			\includegraphics[width = 10cm]{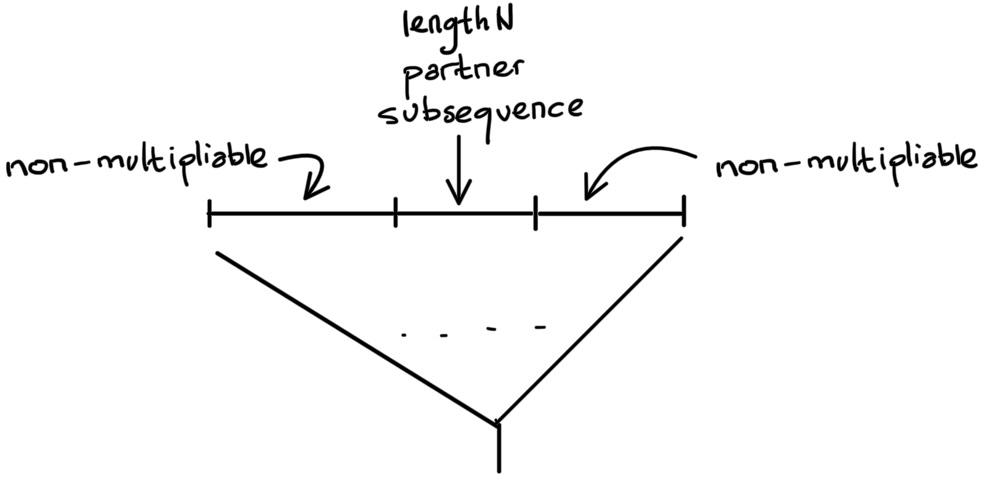}
		\end{center}
		We are going to look at the $\A$-side and argue that given this set-up, there is no tree (composite or tree) that can pull these basic elements together in a way that does not vanish. The issue here is not the \emph{number} of basic inputs on the $\A$-side, it is the \emph{length} of the sequence. While it is possible to have composite trees with nonmultipliable elements on the sides and fewer than $2N$ basic input elements, e.g.
		\begin{center}
			\includegraphics[width = 2cm]{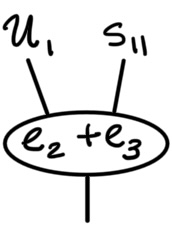}
		\end{center}
		the input sequences of any such tree must have length $\geq 2N$. By the condition on matching idempotents (that basic elements on the $\A$-side and $\B$-side are always paired off) from above, it follows that the length of the input sequences on the $\A$-side and $\B$-side are the same. But the sequence on the $\B$-side has length $< 2N$, and hence, the one on the $\A$-side has length less than $2N$, too. This means that there can be no tree with the given inputs on the $\A$-side that has non-zero output. Hence, when $N < n < 2N$, $(\mu_n^{\w} \otimes \I) \circ \delta^n$ always vanishes.
		
		Next we consider terms of the form $(\mu_{2N}^{\w} \otimes \I) \circ \delta^{2N}$ -- that is, situations where the input sequence for our $\mu$ has $2N$ basic elements. By an argument analogous to the one used in the case where there were $N$ basic inputs, the only pair of trees that does not vanish is
		\begin{center}
			\includegraphics[width = 10cm]{dd9}
		\end{center}
		or equivalent, with cyclically permuted sequence of inputs, i.e.
		\begin{center}
			\includegraphics[width = 10cm]{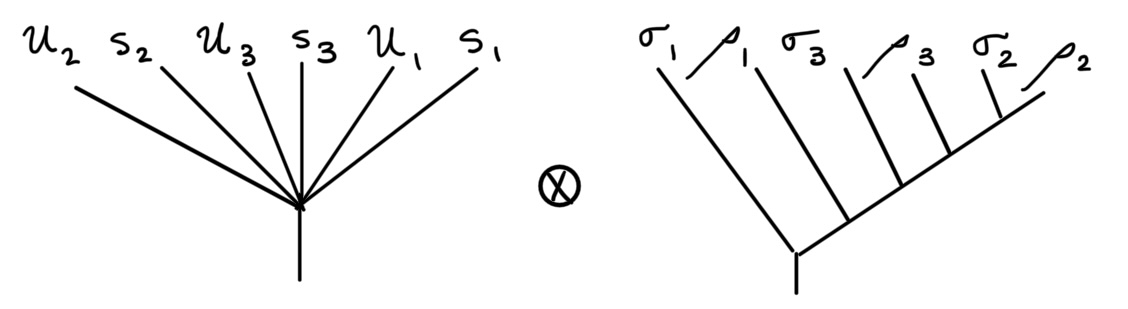}
		\end{center}
		
		We can use an analogous argument and induction to show that $\mu_n^{\w}$ vanishes on all input sequences with $n \not\equiv 0 \emm N$ basic elements. All that remains is sequences with $k N$ basic elements, $k > 1$. We claim that $\mu_{kN}^{\w}$ is also always zero, for $k > 1$ and any $\w$. If $|\w| > 0$, then first of all, because there are no weighted higher multiplications on the $\B$-side, all weight must be on the $\A$-side, and must be $\e_i$'s with $1 \leq i \leq N + 1$. If we pull together a subsequence of (adjacent) elements in a non-zero weighted multiplication, this subsequence will necessarily contain a number of basic inputs which is not a perfect multiple of $N$. By induction down, we can then conclude, as in the case where there were $n \not\equiv 0 \emm N$ basic elements, that if such a multiplication is included (and the total number of basic inputs is a multiple of $N$) then we will have to multiply a pair of nonmultipliable elements. 
		
		Now assume $|\w| = 0$. Then the conclusion holds for dimension reasons. Indeed, $\Gamma^{**}_*$ has to preserve dimension, as defined in~\eqref{diag9}, and for each $\Psi_n^{0}$, we have $\dim \Psi_n^{0} = n - 2$. When we calculate dimension for pairs of trees, we stipulate that $\dim (T \otimes T') = \dim T + \dim T'$. Thus, in general, the two allowable pairs of trees from above:
		\begin{center}
			\includegraphics[width = 10cm]{dd7}
		\end{center}
		clearly have dimension $N - 2$ and $2N - 2$, respectively, the same as the initial trees $\Psi_N^0$ and $\Psi_{2N}^0$ from which they came. 
		
		Now look at a sequence with $n \equiv 0 \emm N$ elements. By an argument analogous to the ones used above, we can conclude that the sequence of inputs can be subdivided into chunks of the form
		\begin{enumerate}
			\item $s_1 \otimes \sigma_1, s_2 \otimes \sigma_2, \ldots, s_N \otimes \sigma_N$, up to cyclic permutation;
			
			\item $U_1 \otimes \rho_1, s_1 \otimes \sigma_1, \ldots, U_N \otimes \rho_N$, again up to cyclic permutation
		\end{enumerate}
		and in order for the full pair of trees not to vanish, these chunks are connected as:
		\begin{center}
			\includegraphics[width = 11cm]{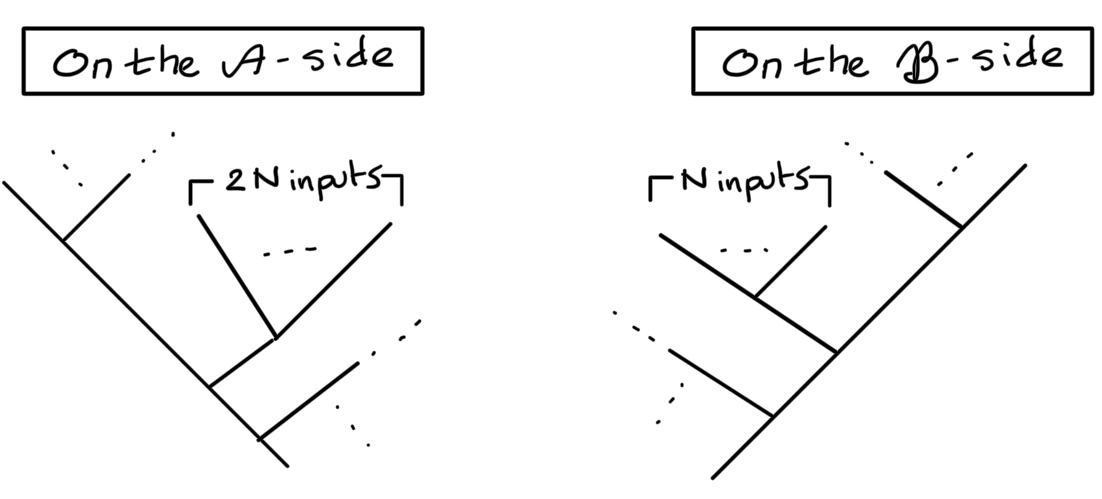}
		\end{center}
		Notice that any element that is connected in via a $\mu_2$ does not contribute to the dimension, since it contributes one input and one vertex. This means that the dimension of a pair that contains $k$ chunks of type (1) and $j$ chunks of type (2) will have dimension
		\begin{equation}\label{dd19}
			\dim (T \otimes T') = k \cdot 2N - k + j \cdot N - j - 1
		\end{equation}
		whereas the dimension of the initial tree will be
		\begin{equation}\label{dd20}
			\dim \Psi_{(2k + j)N}^0 = (2k + j) N - 2,
		\end{equation}
		The right hand sides of~\eqref{dd19} and~\eqref{dd20} will agree if and only if $j + k = 1$. Since $j, k$ are both integers, this means that the only options where the tree is obtainable for dimension reasons are the ones from above, i.e.
		\begin{center}
			\includegraphics[width = 10cm]{dd7}
		\end{center}
		This handles all cases where $n > 1$.
		
		There is one small technical point; namely, that the outputs of e.g.
		\begin{center}
			\includegraphics[width = 10cm]{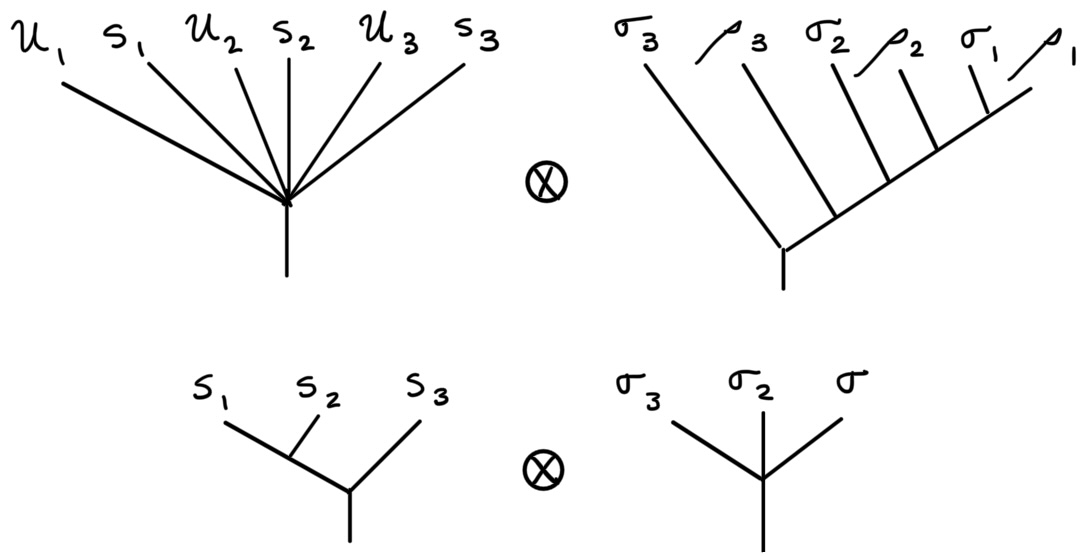}
		\end{center}
		are actually $s_{11} \otimes V_{N + 1}$ and $V_{0} \otimes \sigma_N \rho_N \cdots \sigma_1 \rho_1$, respectively, not $U_{N + 1} \otimes V_{N + 1}, \: V_0 \otimes U_0$, which is what we want. This is, however, only the output from one choice of initial idempotent. When we allow our initial idempotent over all $\{\overline{\x}\}$, and then take the sum of each of the corresponding operations, we do indeed get $U_{N + 1} \otimes V_{N + 1}, \: V_0 \otimes U_0$, as desired.
		
		Let us now consider $n = 0$ and $n = 1$. The first, i.e. $n = 0$ is fairly simple, as it is just the weighted base-case of the weighted diagonal; we do not apply $\delta^1$ at all, and look at what we can get out via $\Gamma_*^{0, \w} \Psi_0^{\w}$. Because the only weighted $\mu_0$s (on either the $\A$-side or the $\B$-side) are $\mu_0^{\e_i}$, with $1 \leq i \leq N + 1$ on the $\A$-side and $i = 0$ on the $\B$-side, we have that
		\begin{equation}\label{dd4}
			((\mu \otimes \nu) \circ \Gamma_*^{0, \w})(\Psi_0^{\w}) = \begin{cases} \mu_0^{\e_i} \otimes \nu(\top) & \w = \e_i \text{ for } 1 \leq i \leq N + 1 \\ \mu(\top) \otimes \mu_0^{\e_0} & \w = \e_0 \\ 0 & \text{ otherwise}\end{cases}
		\end{equation}
		where $\mu$ and $\nu$ are as defined in Section~\ref{tens} above.
		Since idempotents (of the two outputs) have to match up, we have
		\begin{align*} 
			\mu_0^{\e_i} \otimes \nu(\top) &= U_i \otimes V_i \text{ for } 1\leq i \leq N + 1 \\
			 \mu(\top) \otimes \mu_0^{\e_0} &= V_0 \otimes U_{0}
		\end{align*}
		This means that the output of the $n = 0$ term in the $\A_{\infty}$ relation are $V_0 \otimes U_0$ and $\{U_i \otimes V_i\}_{i = 1}^{N + 1}$. 
		
		Now we look at $n = 1$. Because the weighted diagonal is defined using stably weighted trees (i.e. no 2-valent unweighted vertices) and there are no weighted $\mu_1$s on either side, the only non-identically-vanishing operation on $\A \otimes \B$ will be the tensor differential. Since there is no non-vanishing differential on $\A$, this is just
		\[
			\del_{\otimes} = \id_{\A} \otimes \de_{\B}.
		\]
		This means that the non-vanishing terms for $n = 1$ are going to be
		\[
			((\del_{\otimes} \otimes \I)\circ \delta^1)(\overline{\{i\}}) = \del_{\otimes} (U_i \otimes \rho_i) \otimes \overline{\{i\}} = U_i \otimes V_i \otimes \overline{\{i\}}.
		\]
		
		Notice that this also shows that each $U_i \otimes V_i$ with $1 \leq i \leq N + 1$, as well as $V_0 \otimes U_0$, can be obtained in two ways. 
	\end{proof}
	
	We have now additionally verified:
	\begin{proposition}\label{dd5}
		The $\A_{\infty}$-relations,~\eqref{dd2}, hold for $X$, so $X$ is a valid weighted $DD$-imodule. 
	\end{proposition}
	
	Because the only non-vanishing outputs of the left hand side of~\eqref{dd2} are catalogued above, and each appears exactly twice. 

    \bibliography{DualRefs}
    \bibliographystyle{plain}

\end{document}